\documentclass[10pt]{article}
\usepackage{amssymb,amsmath,amsfonts,amsthm,tikz,multirow,wasysym,color,siunitx}
\usetikzlibrary{arrows,calc}

\title{Tilings of the sphere by congruent pentagons IV: Edge combination $a^4b$ with general angles}
\author{Yixi Liao, Junjie Shu, Erxiao Wang\thanks{Corresponding author (wang.eric@zjnu.edu.cn). Research was supported by National Natural Science Foundation of China NSFC-RGC 12361161603 and Key projects of Zhejiang Natural Science Foundation No. LZ22A010003.}\\
	Zhejiang Normal University}

\usepackage[hidelinks, hyperindex]{hyperref}

\usepackage{subfig,overpic}

\DeclareMathOperator{\opt}{\; |\;}

\newcommand{\sub}{\subset}
\newcommand{\pa}{\partial}
\newcommand{\mb}{\mathbf}
\newcommand{\mc}{\mathcal}
\newcommand{\bb}{\mathbb}
\newcommand{\ssum}{\textstyle \sum}
\newcommand\aaa{\alpha}
\newcommand\bbb{\beta}
\newcommand\ccc{\gamma}
\newcommand\ddd{\delta}
\newcommand\eee{\epsilon}
\newcommand\bn{\boldsymbol{n}}

\newcommand{\arcThroughThreePoints}[4][]{
	\coordinate (middle1) at ($(#2)!.5!(#3)$);
	\coordinate (middle2) at ($(#3)!.5!(#4)$);
	\coordinate (aux1) at ($(middle1)!1!90:(#3)$);
	\coordinate (aux2) at ($(middle2)!1!90:(#4)$);
	\coordinate (center) at ($(intersection of middle1--aux1 and middle2--aux2)$);
	\draw[#1]
	let \p1=($(#2)-(center)$),
	\p2=($(#4)-(center)$),
	\n0={veclen(\p1)},       
	\n1={atan2(\y1,\x1)}, 
	\n2={atan2(\y2,\x2)},
	\n3={\n2>\n1?\n2:\n2+360}
	in (#2) arc(\n1:\n3:\n0);
}

\newcommand{\dash}{\hspace{0.1em}\dashrule{0.7}{2.4 1 2.4 1 2.4}\hspace{0.1em}}
\newcommand{\thin}{\hspace{0.1em}\rule{0.7pt}{0.8em}\hspace{0.1em}}
\newcommand{\thick}{\hspace{0.1em}\rule{1.5pt}{0.8em}\hspace{0.1em}}

\newtheorem{theorem}{Theorem}
\newtheorem{lemma}[theorem]{Lemma}

\newtheorem{proposition}{Proposition}
\newtheorem*{theorem*}{Theorem}

\theoremstyle{definition}
\newtheorem*{definition*}{Definition}

\newtheorem*{case*}{Case}

\newtheorem*{subcase*}{Subcase}

\newtheorem*{subsubcase*}{Subsubcase}

\theoremstyle{remark}

\numberwithin{equation}{section}

\begin{document}
    \date{}
	\maketitle

	\begin{abstract}
		We classify edge-to-edge tilings of the sphere by congruent pentagons with the edge combination $a^4b$ and with any irrational angle in degrees: they are three $1$-parameter families of pentagonal subdivisions of the Platonic solids, with $12, 24$ and $60$ tiles; and a sequence of $1$-parameter families of pentagons admitting non-symmetric $3$-layer earth map tilings together with their various rearrangements under extra conditions. Their parameter moduli and geometric data are all computed in both exact and numerical form. The total numbers of different tilings for any fixed such pentagon are counted explicitly. As a byproduct, the degenerate pentagons produce many new non-edge-to-edge quadrilateral tilings. 

\medskip
 
		
	\end{abstract}

	\section{Introduction}
	A lot of progress has been made on spherical tilings in recent years.
	Non-edge-to-edge multihedral tilings of the sphere by regular polygons have been classified in \cite{aehj}. Edge-to-edge monohedral quadrilateral tilings of the sphere have been fully classified in \cite{clm1,lw1,lw2,lw3}.
	In an edge-to-edge tiling of the sphere by congruent pentagons, the pentagon can have five possible edge combinations \cite[Lemma 9]{wy1}: $a^2b^2c,a^3bc,a^3b^2,a^4b,a^5$. Tilings for $a^2b^2c,a^3bc,a^3b^2,a^5$ were classified in \cite{wy3,lwy,wy1,wy2}, whose methods were further extended by Cheung, Luk and Yan to classify $a^4b$-tilings in a preprint \cite{clm2} of \textbf{174 pages}. Our strategy was developed independently and does not rely on their case reductions, reconstruction arguments or classification proof. Agreement between the two classifications provides an independent verification.
	
	The present paper treats pentagons with at least one irrational angle in degrees; the sequel \cite{hllw} treats the rational-angle case using number-theoretic methods. This totally new strategy was designed for this last and the most difficult $a^4b$-pentagon case originally, though our papers applying this strategy to quadrilateral case \cite{lw2,lw3} were completed earlier since the trigonometric Diophantine equations for pentagons are much harder to solve. 
	
	Together, our two papers give a shorter proof ($45+23=$ \textbf{68 pages}) and provide four new results: 
	\begin{enumerate}
		\item full geometric data for all prototiles; 
		\item 3D illustrations of every tiling type; (https://github.com/spherical1tiling2)
		\item explicit counts of the tilings admitted by each prototile; 
		\item new non-edge-to-edge quadrilateral tilings obtained from degenerate pentagons.
	\end{enumerate}
	
	An $a^4b$-pentagon is given by Figure \ref{pentagon}, with normal edge $a$, thick edge $b$ and angles $\aaa$, $\bbb$, $\ccc$, $\ddd$, $\eee$ as indicated. The second picture is the mirror image or flip of the first. The angles determine the orientation. Conversely, the edge lengths and the orientation also determine the angles. So we may present the tiling by shading instead of indicating all angles. Throughout this paper, an $a^4b$-tiling is always an edge-to-edge tiling of the sphere by congruent simple pentagons in Figure  \ref{pentagon}, such that all vertices have degree $\geq$ $3$. We call such a pentagon $\textbf{almost equilateral}$. We assume $a\ne b$ and allow straight angles.

	\begin{figure}[htp]
		\centering
		\begin{tikzpicture}[>=latex,scale=1]

			\begin{scope}[xshift=3cm]

				\draw
				(-54:1) -- (18:1) -- (90:1) -- (162:1) -- (234:1);

				\draw[line width=1.5]
				(234:1) -- (-54:1);

				\node at (90:0.75) {\small $\alpha$};
				\node at (162:0.75) {\small $\beta$};
				\node at (15:0.75) {\small $\gamma$};
				\node at (234:0.75) {\small $\delta$};
				\node at (-54:0.75) {\small $\epsilon$};

			\end{scope}

			\begin{scope}[xshift=6cm]

				\fill[gray!50]  (-54:1) -- (18:1) -- (90:1) -- (162:1) -- (234:1);

				\draw
				(-54:1) -- (18:1) -- (90:1) -- (162:1) -- (234:1);

				\draw[line width=1.5]
				(234:1) -- (-54:1);

				\node at (90:0.75) {\small $\alpha$};
				\node at (162:0.75) {\small $\ccc$};
				\node at (15:0.75) {\small $\bbb$};
				\node at (234:0.75) {\small $\eee$};
				\node at (-54:0.75) {\small $\ddd$};

			\end{scope}

			\begin{scope}[xshift=8cm]

				\draw
				(0,0.4) -- node[above=-2] {\small $a$} ++(1,0);

				\draw[line width=1.5]
				(0,-0.1) -- node[above=-2] {\small $b$} ++(1,0);

			\end{scope}
		\end{tikzpicture}
		\caption{Pentagons with the edge combinations $a^4b$.}
		\label{pentagon}
	\end{figure}

	We use $\alpha^{n_1}\beta^{n_2}\gamma^{n_3}\delta^{n_4}\eee^{n_5}$ to mean a vertex having $n_1$ copies of $\alpha$, $n_2$ copies of $\beta$, etc. The anglewise vertex combination(s), abbreviated as \textbf{AVC}, is the collection of all vertices in a tiling. Then the notation $T(16 \aaa\ddd\eee, 8\bbb^2\ccc, 2\ccc^4)$ means the tiling has exactly $16$ vertices $\aaa\ddd\eee$, $8$ vertices $\bbb^2\ccc$, and $2$ vertices $\ccc^{4}$, and is \textit{uniquely} determined by them. In general there may exist several different tilings with the same set of vertices. Then we use $\{4\aaa\bbb\ccc,14\aaa\ddd\eee,6\bbb^2\ccc,4\bbb\ddd\eee,2\aaa\ccc^{3},2\ccc^{2}\ddd\eee: 3\}$ to mean that there are $3$ different tilings.

	When an almost equilateral pentagon is symmetric about an axis through $\aaa$, i.e. $\bbb=\ccc$, $\ddd=\eee$, it can be divided into two congruent quadrilaterals along this axis, and the corresponding classification has already been completed by \cite{lw1,lw2,lw3}. We translate it as the following theorem. 

\begin{figure}[htp]
		\centering
		\begin{overpic}[scale=0.2]{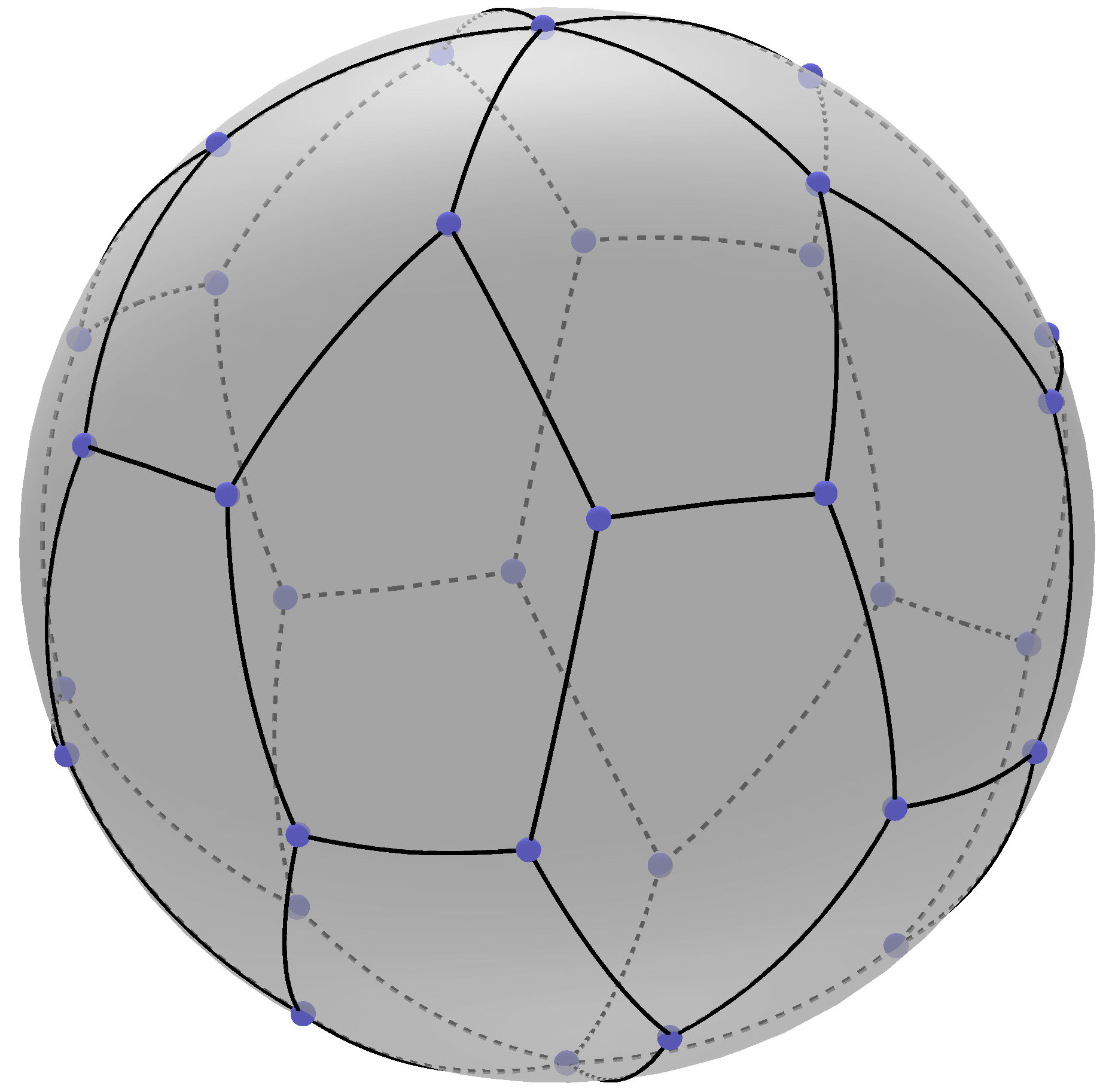}
			\put(90,5){\small $f=20$}
		\end{overpic}   \hspace{30pt}
		\includegraphics[scale=0.202]{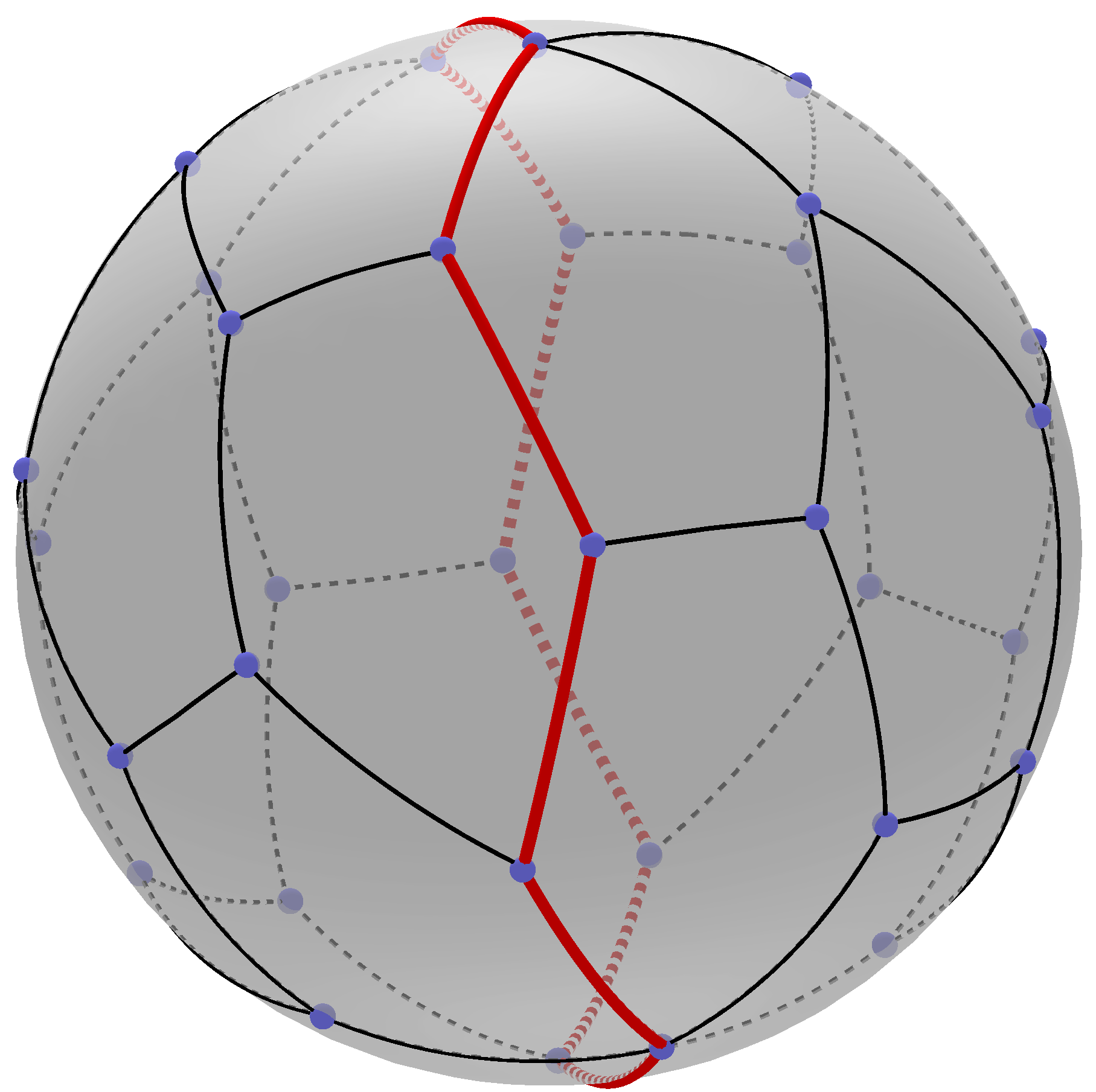}\hspace{30pt}
		\includegraphics[scale=0.22]{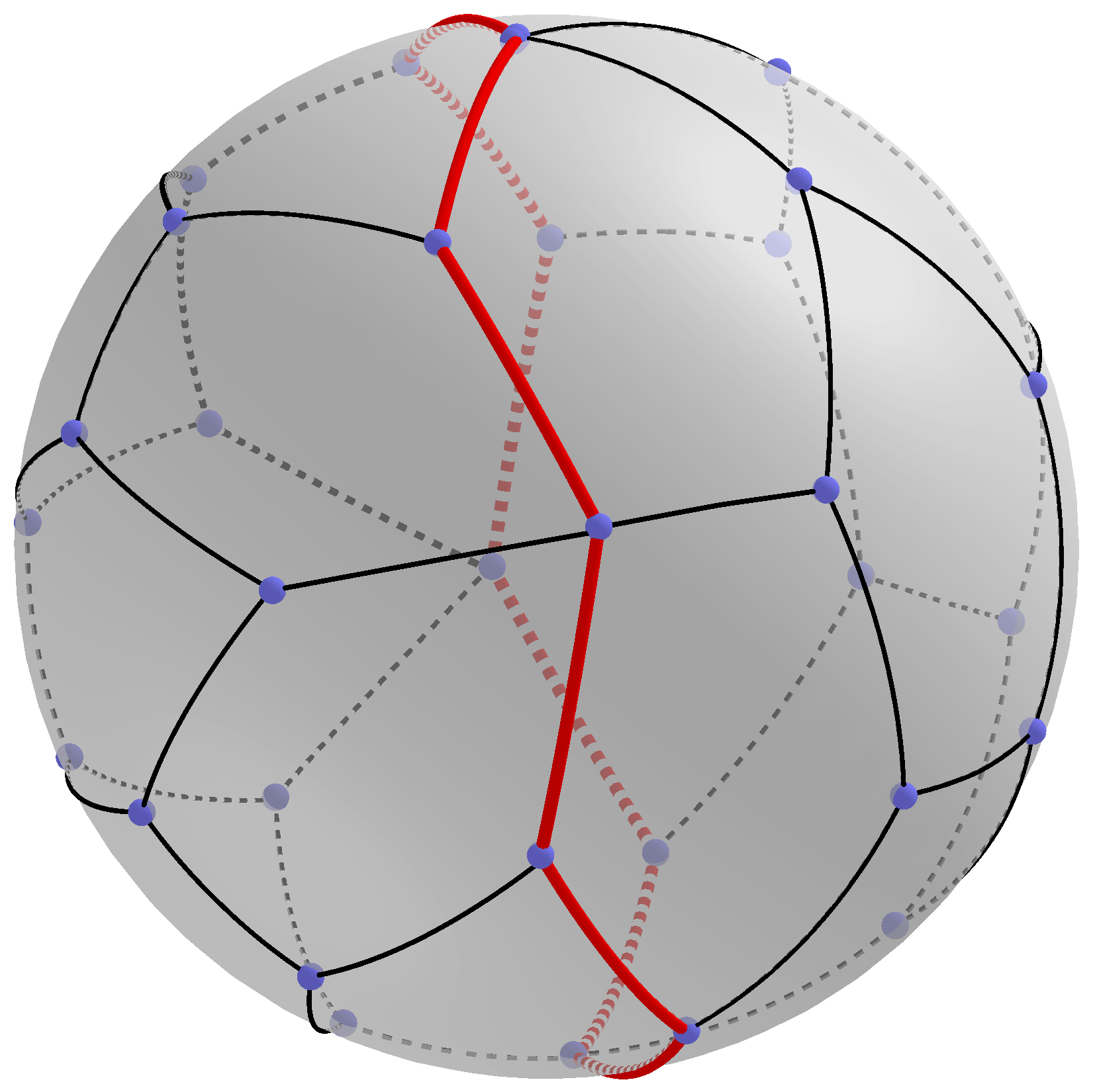}
		\caption{A symmetric $3$-layer earth map tiling and its two standard flip modifications.}
		\label{sym3}
	\end{figure}

	\begin{theorem} \label{theorem1}
		All symmetric $a^4b$-tilings are:
		\begin{enumerate}
			\item A $1$-parameter family of symmetric $a^4b$-pentagonal subdivisions of the tetrahedron with $12$ tiles $T(12\aaa\ddd\eee,4\bbb^3,4\ccc^3)$ (see the first picture of Figure \ref{poufen} and Table \ref{tab-1}).

			\item A sequence of unique symmetric $a^4b$-pentagons with five angles $(\tfrac8f,\,1-\tfrac4f,\,1-\tfrac4f,\,\tfrac12+\tfrac2f,\,\tfrac12+\tfrac2f)\pi$ and
			\begin{equation*}
				\cos a=1-2\left(\frac{\sqrt{5}-1}{4\cos{\frac{4\pi}f}}\right)^2,\,\,\,\,
				\cos b=2 \left(\frac{(3-\sqrt5)\cos^2{\frac{4\pi}f}+\sqrt5-2}{\cos{\frac{4\pi}f}}\right)^2 -1,
			\end{equation*}
			each admitting a symmetric $3$-layer earth map tiling $T(2m\,\aaa\bbb\ccc,2m\,\bbb\ddd\eee,2m\,\ccc\ddd\eee,2\aaa^{m})$ by $f=4m$ tiles for any $m\ge4$, among which each odd $m=2k+1$ case admits two standard flip modifications (see Figure \ref{sym3}):

			$T(4k\,\aaa\bbb\ccc,(4k+2)\bbb\ddd\eee,(4k+2)\ccc\ddd\eee,2\aaa^{k+1}\bbb,2\aaa^{k+1}\ccc)$, and

			$T((4k+2)\aaa\bbb\ccc,(4k+1)\bbb\ddd\eee,(4k+1)\ccc\ddd\eee,\aaa^{k+1}\bbb,\aaa^{k+1}\ccc,2\aaa^k\ddd\eee)$ .
		\end{enumerate}
	\end{theorem}

	The pentagonal subdivision introduced in \cite{wy1} replaces each triangular face of the regular tetrahedron, octahedron or icosahedron by three congruent pentagons determined by a freely chosen point $V$ (see Figure \ref{div_triangle}). The exact range of $V$ was studied in great detail in \cite{lwy}.
	Only when $V$ is on the perpendicular bisector of $AB$ will the pentagon become almost equilateral (see Figure \ref{poufen}).  The parameter $t$ in Table \ref{tab-1} represents the signed spherical distance of $V$ to $AB$.

	\begin{figure}[htp]
		\centering
		\begin{tikzpicture}[scale=1.5]
			\foreach \a in {0,1,2}
			{
				\begin{scope}[rotate=120*\a]

					\draw
					(0:1) arc (30:90:{sqrt(3)});

					\coordinate (M) at (60:{sqrt(3)-1});
					\coordinate (X) at (0:1);

					\pgfmathsetmacro{\th}{18}

					\pgfmathsetmacro{\ph}{-atan((1/3)*(cot(\th+30)))}
					\coordinate (A1) at ({\ph}:{sqrt(2)});
					\coordinate (B1) at ({\ph+180}:{sqrt(2)});

					\pgfmathsetmacro{\ps}{-atan((1/3)*(cot(\th+150)))}
					\coordinate (A2) at ({\ps}:{sqrt(2)});
					\coordinate (B2) at ({\ps+180}:{sqrt(2)});

					\coordinate (V) at (36.8:0.595);
					\coordinate (W) at (-42:0.95);
					\coordinate (W2) at (120-42:0.95);

					\arcThroughThreePoints[line width=1.5]{X}{A1}{V}
					\arcThroughThreePoints[line width=1.5]{W}{A2}{X}
					\arcThroughThreePoints[dashed]{V}{M}{W2}

					\draw
					(0,0) -- (V);

					\draw[gray!50]
					({-sqrt(3)+1},0) -- (1,0);

				\end{scope}
			}

			\fill
			(1,0) circle (0.05);
			\filldraw[fill=white]
			(0,0) circle (0.05);
			\filldraw[fill=gray]
			(36.8:0.595) circle (0.05);

			\node at (0:1.15) {$B$};
			\node at (220:0.18) {$A$};
			\node at (20:0.5) {$V$};

		\end{tikzpicture}
		\caption{Pentagonal subdivision of a regular triangle face.}
		\label{div_triangle}
	\end{figure}

	\begin{figure}[htp]
		\centering
		\includegraphics[scale=0.045]{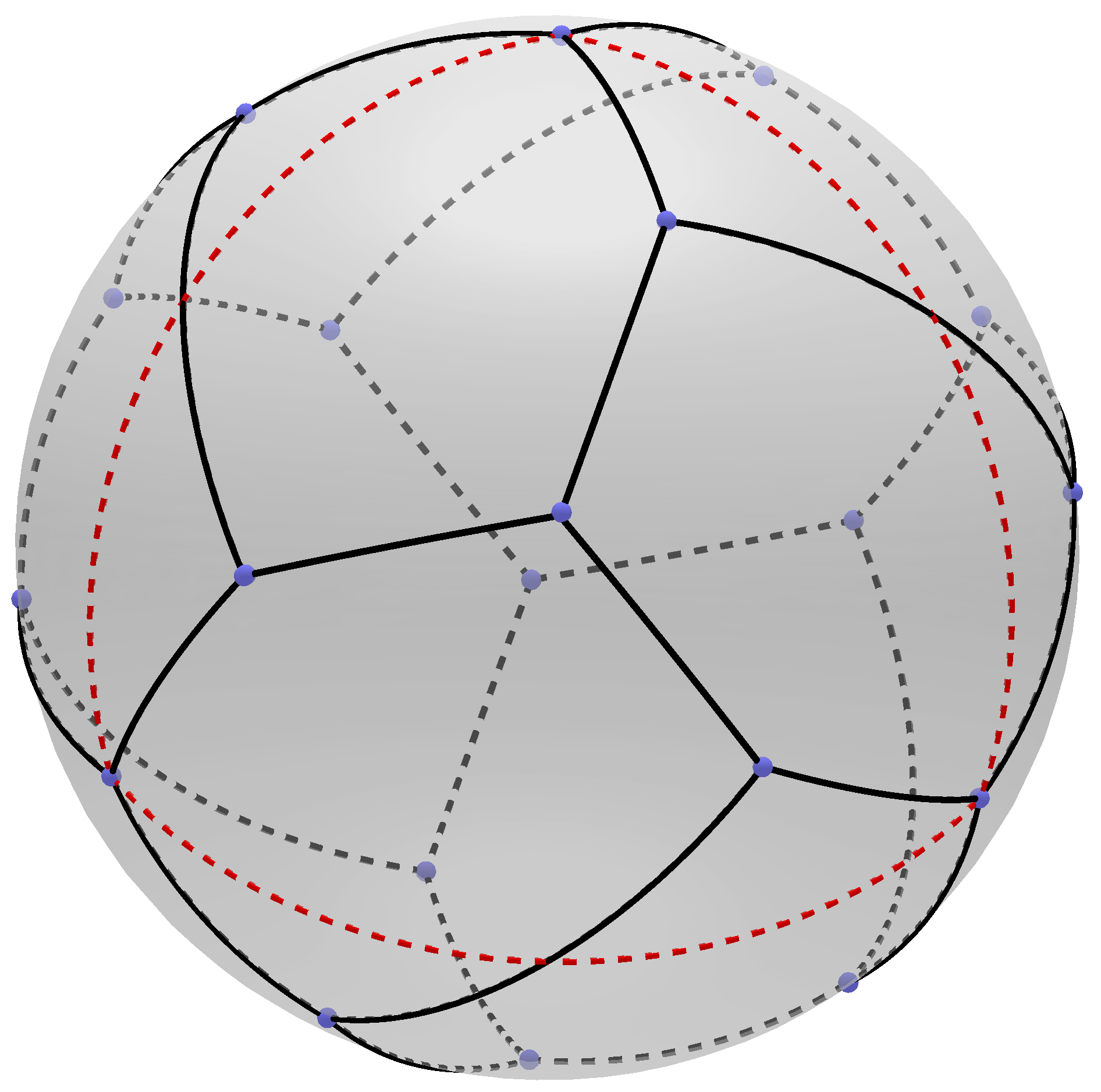}\hspace{30pt}
		\includegraphics[scale=0.047]{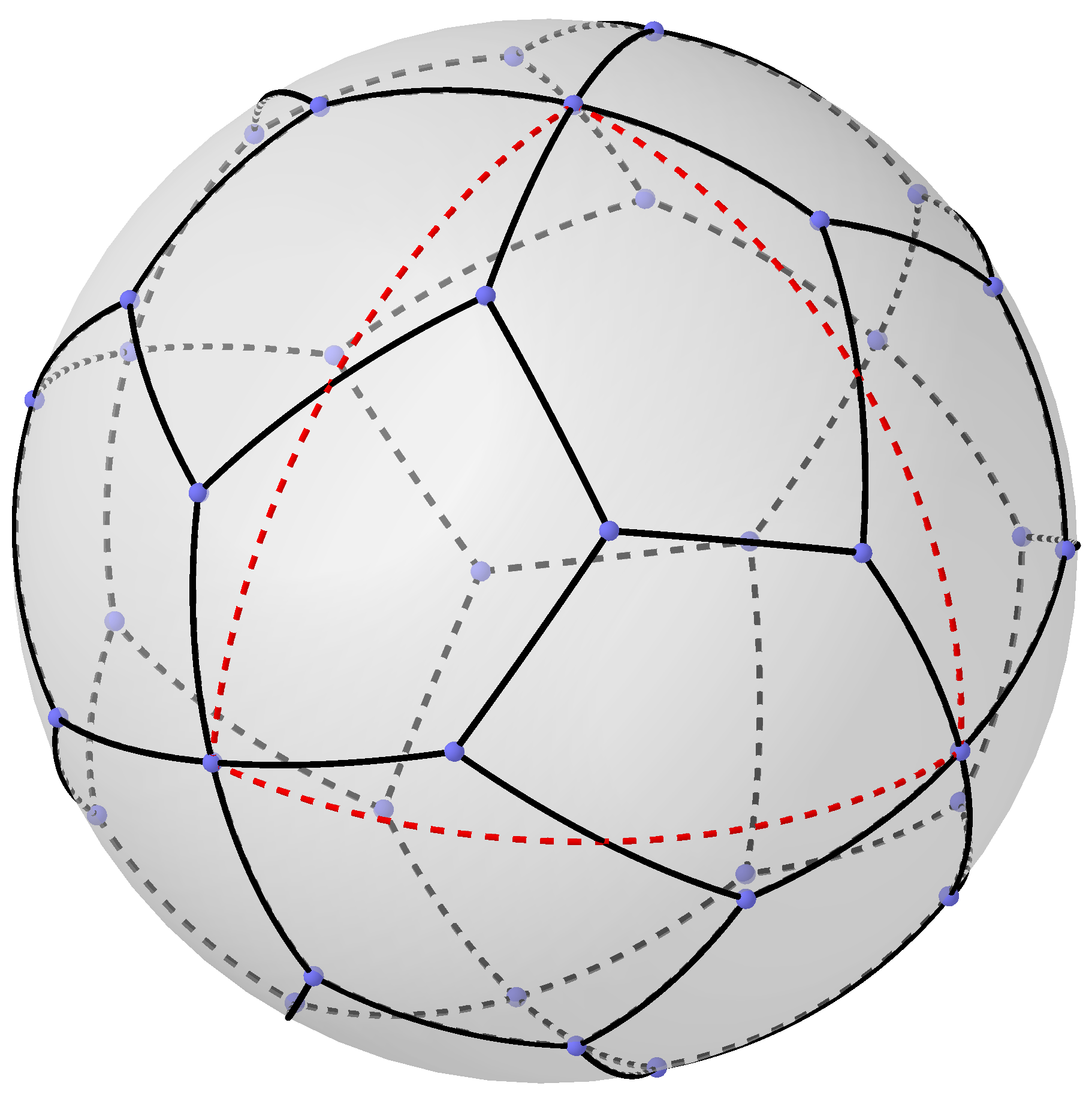}\hspace{30pt}
		\includegraphics[scale=0.048]{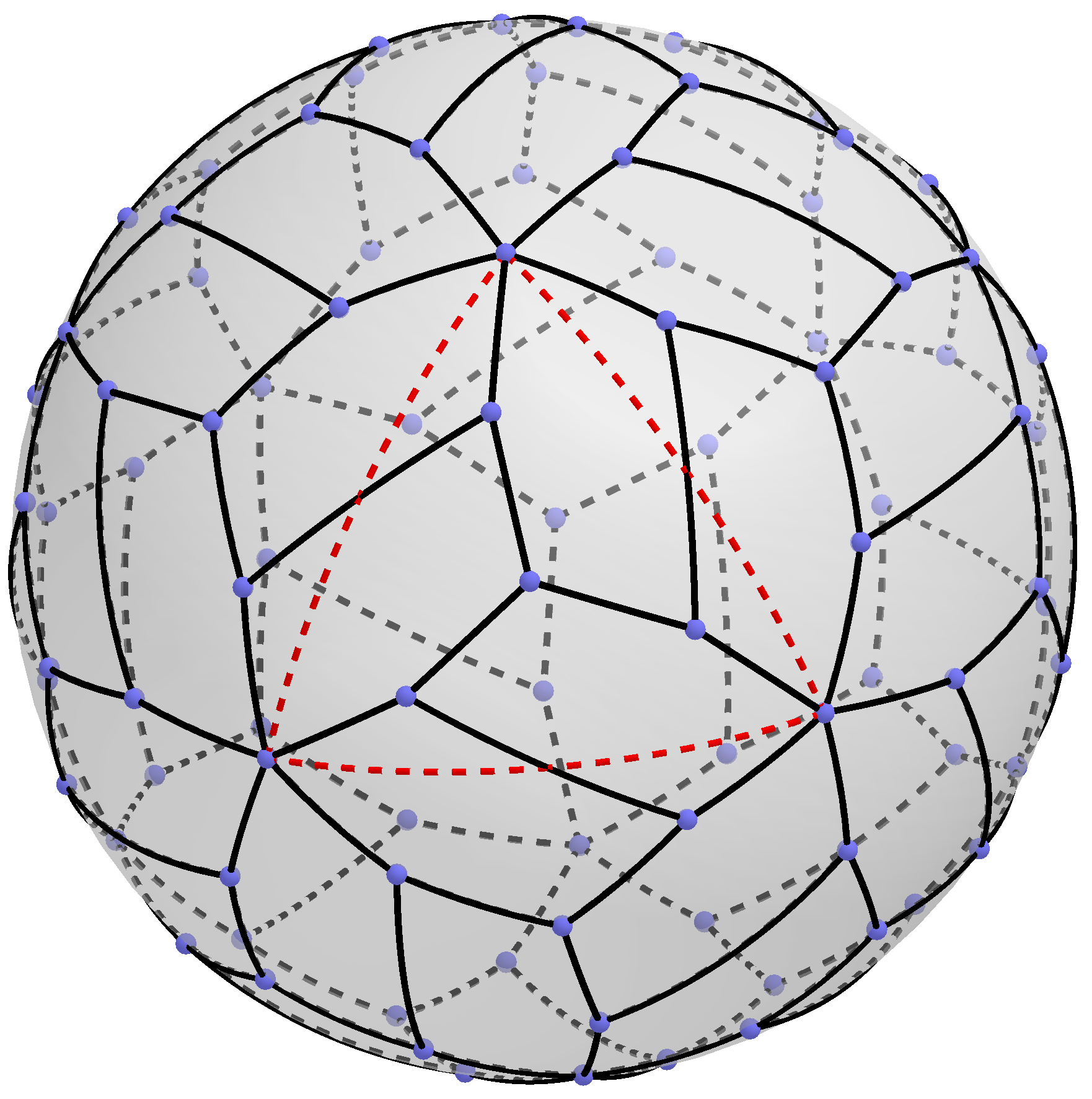}
		\caption{Pentagonal subdivision of the tetra-, octa- and icosa-hedron.}
		\label{poufen}
	\end{figure}

	\begin{table}[htp]
		\begin{center}
			\caption{The geometric data of three $1$-parameter families of $a^4b$-pentagonal subdivisions.}\label{tab-1}
			\bgroup
			\resizebox{\textwidth}{30mm}{
				\begin{tabular}{ll}
					\hline
					\multicolumn{1}{c}{Vertices and tilings} & \multicolumn{1}{c}{$\alpha$, $\bbb$, $\ccc$, $\ddd$, $\eee$, $a$, $b$} \\
					\hline
					\multirow{1}{*}{$T(12\aaa\ddd\eee$, $4\bbb^3$, $4\ccc^3)$} &
					$\bbb=\ccc=\frac{2\pi}{3}$, $\cos a=\frac{\sqrt{6} \cos t }{3}$, $b=\frac{\pi}{2}-2t$, $\cos \frac{\aaa}{2}= \frac{\sqrt{2}\sin t}{\sqrt{1+2\sin^2 t}}$,\\[3pt]
					& $\cos \eee=\frac{1-\sqrt{2}\cos t\cos\frac{b}{2}}{\sqrt{1+2 \sin^{2} t}\, \sin\frac{b}{2}}$, $\ddd=2\pi-\aaa-\eee$, $t\in(-\frac\pi4,\frac\pi4)\backslash\{t_0\}$, $t_0=\arctan\frac{\sqrt{5}}{5}$. The pentagon is equilateral at $t_0$. \\[5pt]
					\hline
					\multirow{1}{*}{$T(24\aaa\ddd\eee$, $8\ccc^3$, $6\bbb^4)$}&$\bbb=\frac{\pi}{2}$, $\ccc=\frac{2\pi}{3}$, $\cos a=\frac{\sqrt{18+6 \sqrt{3}} \cos t}{6}$, $\cos\frac{b}{2}=\frac{1}{12}(\sqrt{3+\sqrt{3}}( \sqrt{(4-2\sqrt{3}) \cos^2 t}+4 \cos t)+6 \sin t)$,\\[3pt]
					&$\cos\frac{\aaa}{2}=\frac{ \sqrt{3+\sqrt{3}}\sin t}{\sqrt{6-(3+\sqrt{3})\cos^{2} t}}$, $\cos \eee= \frac{ 2-\sqrt{3+\sqrt{3}}\cos t  \cos\frac{b}{2}}{\sqrt{6-(3+\sqrt{3}) \cos^{2}t} \sin\frac{b}{2}}$, $\ddd=2\pi-\aaa-\eee$,\\[3pt]
					& $t\in(t_1,t_2)\backslash\{t_0\}$, $\tan t_1=-\frac{3^{\frac{3}{4}} \sqrt{2}}{\sqrt{18+6 \sqrt{3}}}$, $\tan t_2=\frac{\sqrt{6-2\sqrt{3}}}2$, $\cos t_0$ is the largest real root of\\[3pt]
					& $(150+75\sqrt{3}) x^{4}-\sqrt{18+6 \sqrt{3}} (6+10 \sqrt{3}) x^{3}-(90+78\sqrt{3}) x^{2}+12 \sqrt{18+6 \sqrt{3}} x+18=0 $, $t_0\approx0.0892\pi$.\\[3pt]
					\hline
					\multirow{1}{*}{$T(60\aaa\ddd\eee$, $20\ccc^3$, $12\bbb^5)$} &$\bbb=\frac{2\pi}{5}$, $\ccc=\frac{2\pi}{3}$,  $\cos a=\frac{\sqrt{450+30 \sqrt{75+30 \sqrt{5}}}\, \cos t}{30}$,\\[3pt]
					&$\cos \frac{b}{2}=\frac{2^{\frac{5}{4}} (\sqrt{5}+1) \sqrt{3\sqrt{10-2\sqrt{5}}+\sqrt{15}+\sqrt{3}}\cos t+(\sqrt{5}-1)(2^{\frac{1}{4}} \sqrt{3 \sqrt{10-2\sqrt{5}}-\sqrt{15}-\sqrt{3} }\cos t+ 6(5-\sqrt{5})^{\frac{1}{4}}\sin t )}{24 (5-\sqrt{5})^{\frac{1}{4}}}$,\\[3pt]
					&$\cos\frac{\aaa}{2}=\frac{ \sqrt{6 \sqrt{5-\sqrt{5}}+\sqrt{30}+\sqrt{6}}\sin t}{\sqrt{ 6\sqrt{5-\sqrt{5}}(1+\sin^{2} t)-(\sqrt{30}+\sqrt{6})\cos^{2}t}}$, $\cos\eee=\frac{(\sqrt{5}+1) (5-\sqrt{5})^{\frac{1}{4}}- \sqrt{6 \sqrt{5-\sqrt{5}}+ \sqrt{30}+\sqrt{6}} \cos t\cos\frac{b}{2}}{\sqrt{(6+6\sin^{2}t) \sqrt{5-\sqrt{5}}-(\sqrt{30}+\sqrt{6})\cos^{2}t} \sin\frac{b}{2}}$, $\ddd=2\pi-\aaa-\eee$,\\[3pt]
					& $t\in(t_1,t_2)\backslash\{t_0\}$, $\tan t_0= \frac{\sqrt{3}}{1513800}((\sqrt{75+30 \sqrt{5}}(42050-16820 \sqrt{5})+25230\sqrt{5}) \sqrt{90+6 \sqrt{75+30 \sqrt{5}}}$\\[3pt]
					&$+3 ((11 \sqrt{5}-27) \sqrt{75+30 \sqrt{5}}-5 \sqrt{5}-35)\sqrt{ \sqrt{75+30 \sqrt{5}}(7259325+3566635 \sqrt{5})-36362925\sqrt{5}-70401375})$,\\[3pt]
					&$\tan t_1=\frac{(\sqrt{5}-3)}{4} (\sqrt{15+6 \sqrt{5}}+\sqrt{75+30 \sqrt{5}}-3 \sqrt{5}-7)^{\frac12}$, $\tan t_2=\frac{\sqrt{600-10 \sqrt{6}(5+\sqrt{5})^{\frac{3}{2}}}}{20}$.\\[3pt]
					\hline
			\end{tabular}}
			\egroup
		\end{center}
	\end{table}

\begin{figure}[htp]
    	\centering
    	\begin{overpic}[scale=0.048]{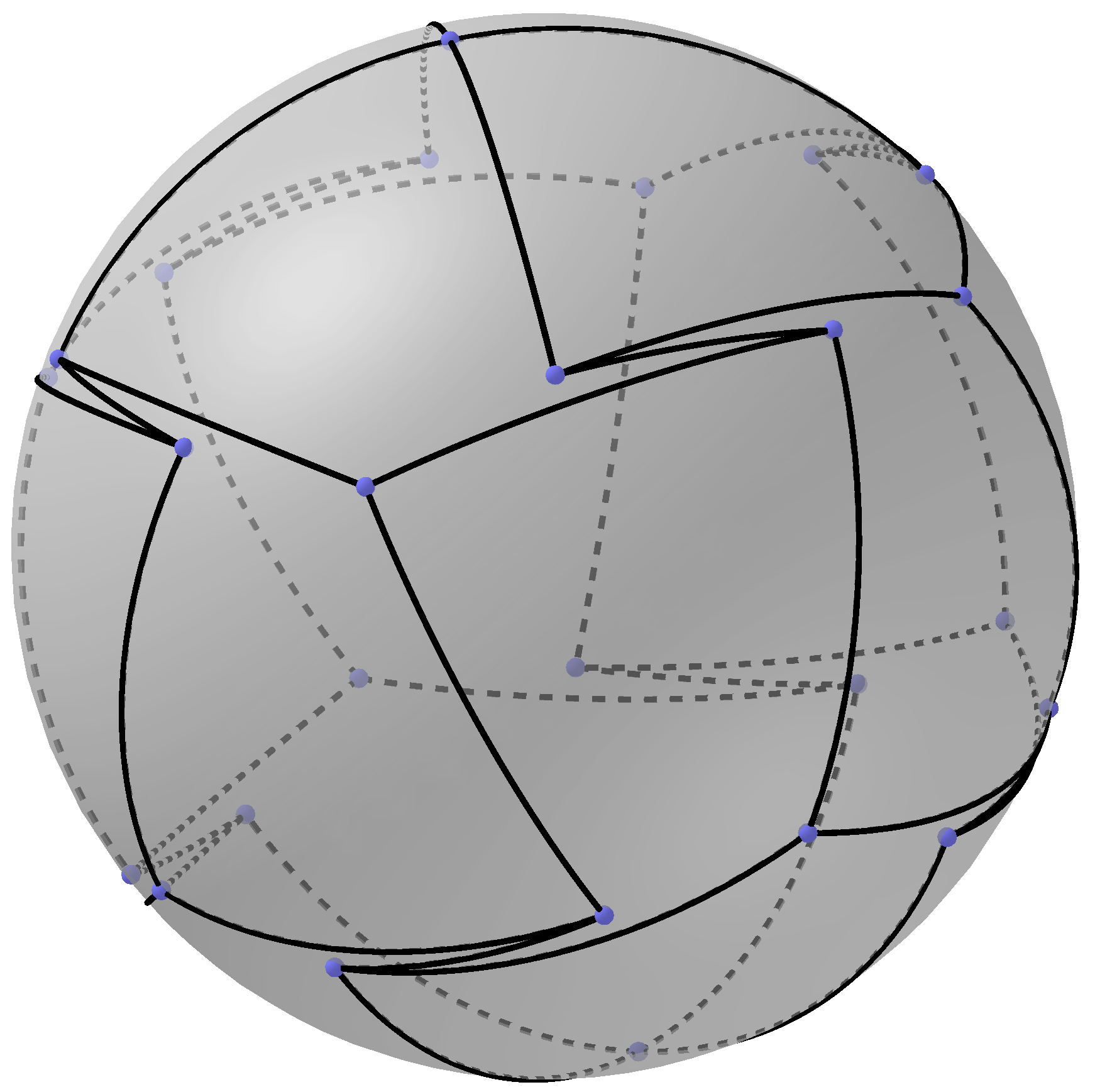} \hspace{20pt}
    		\put(30,-7){\small $\aaa=0.45\pi$}
    	\end{overpic}   \hspace{20pt}
    	\begin{overpic}[scale=0.0480]{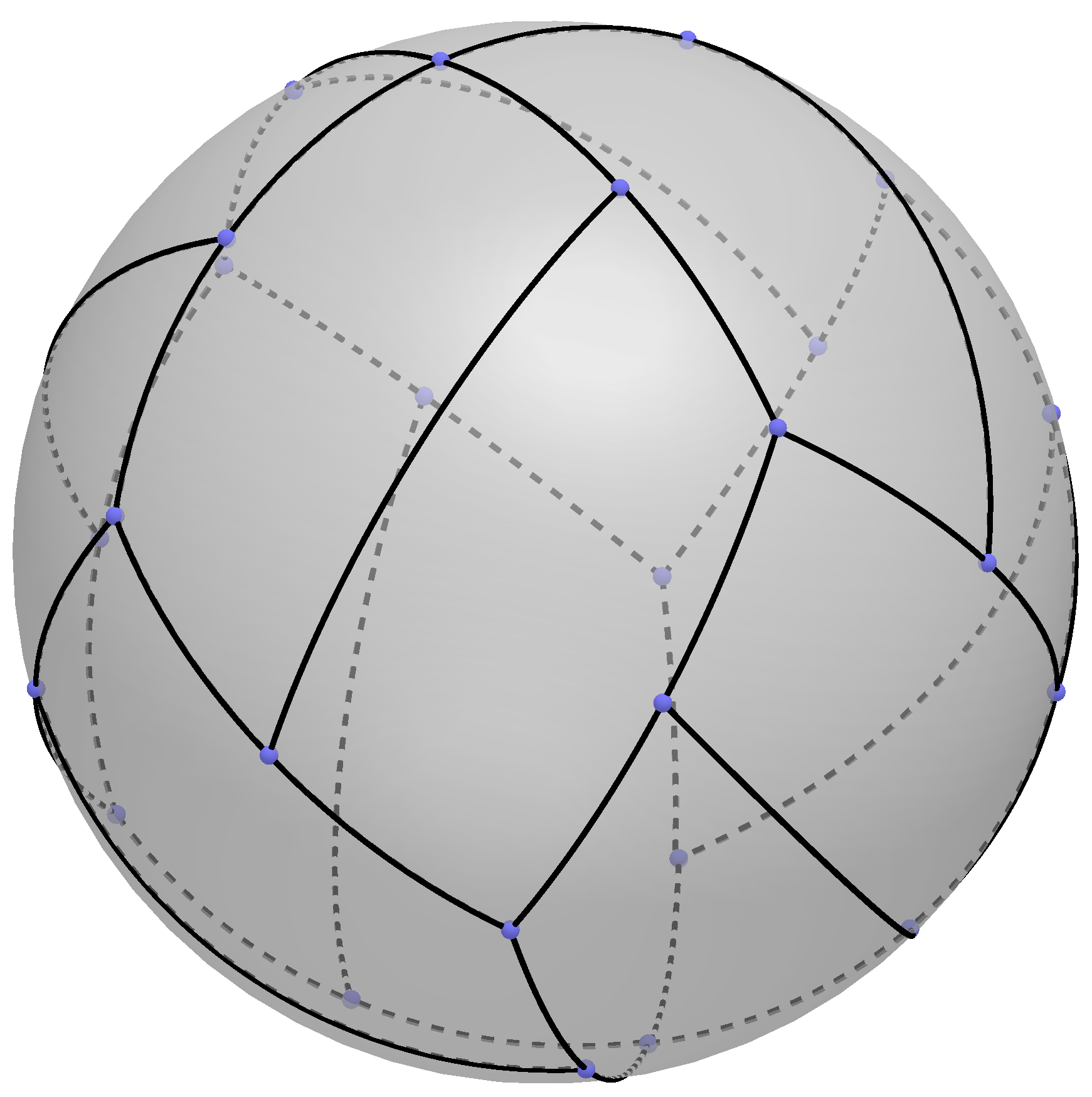} \hspace{20pt}
    		\put(30,-7){\small $\aaa=\pi$}
    	\end{overpic}   \hspace{20pt}
    	\begin{overpic}[scale=0.0498]{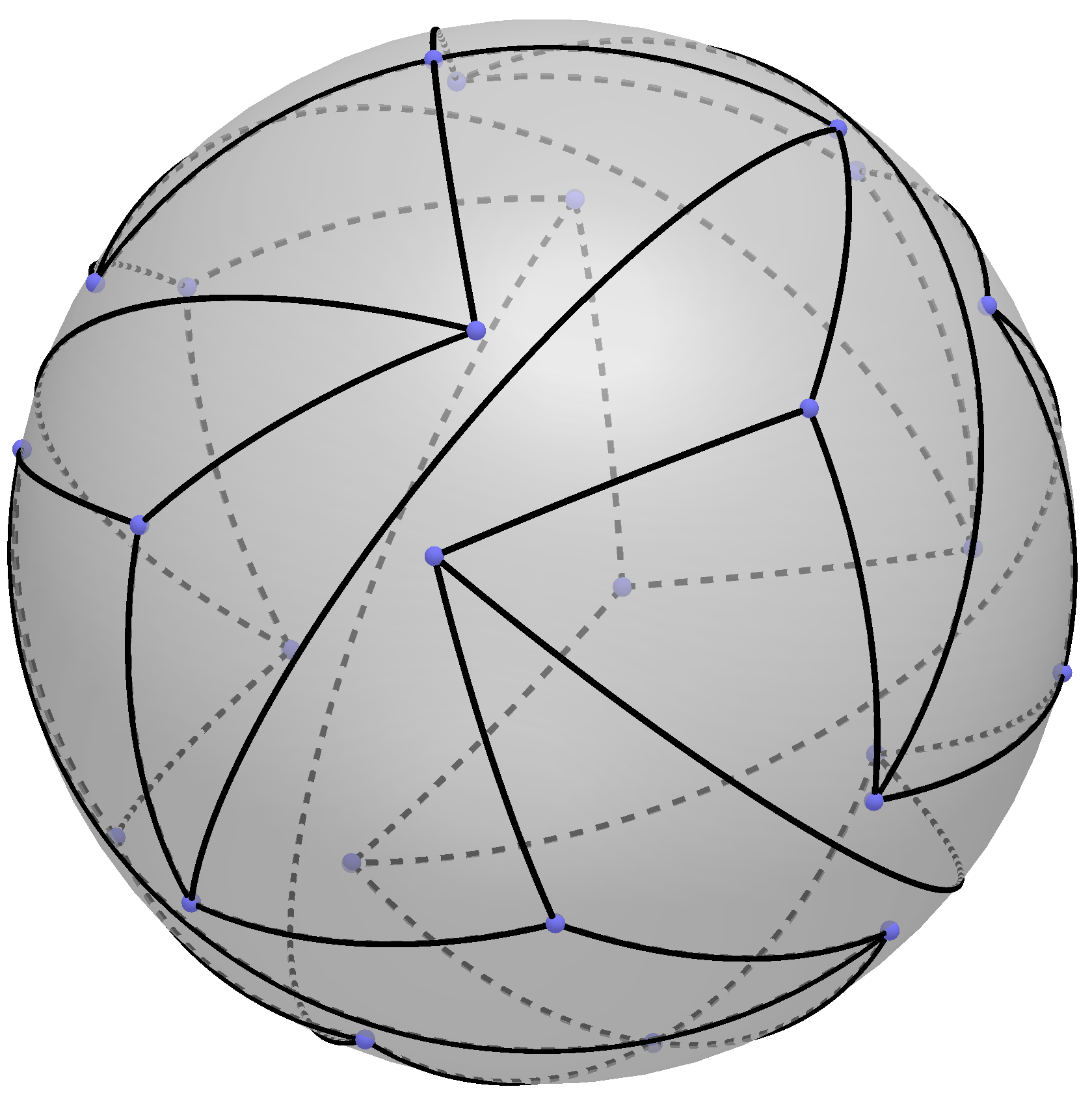}
    		\put(30,-7){\small $\aaa=1.45\pi$}
    	\end{overpic}

    	\caption{The $1$-parameter family of non-symmetric $3$-layer earth map tilings with $f=16$ for three particular parameter values.}
    	\label{f=16 3}
    \end{figure}

\begin{figure}[htp]
		\centering
		\begin{overpic}[scale=0.175]{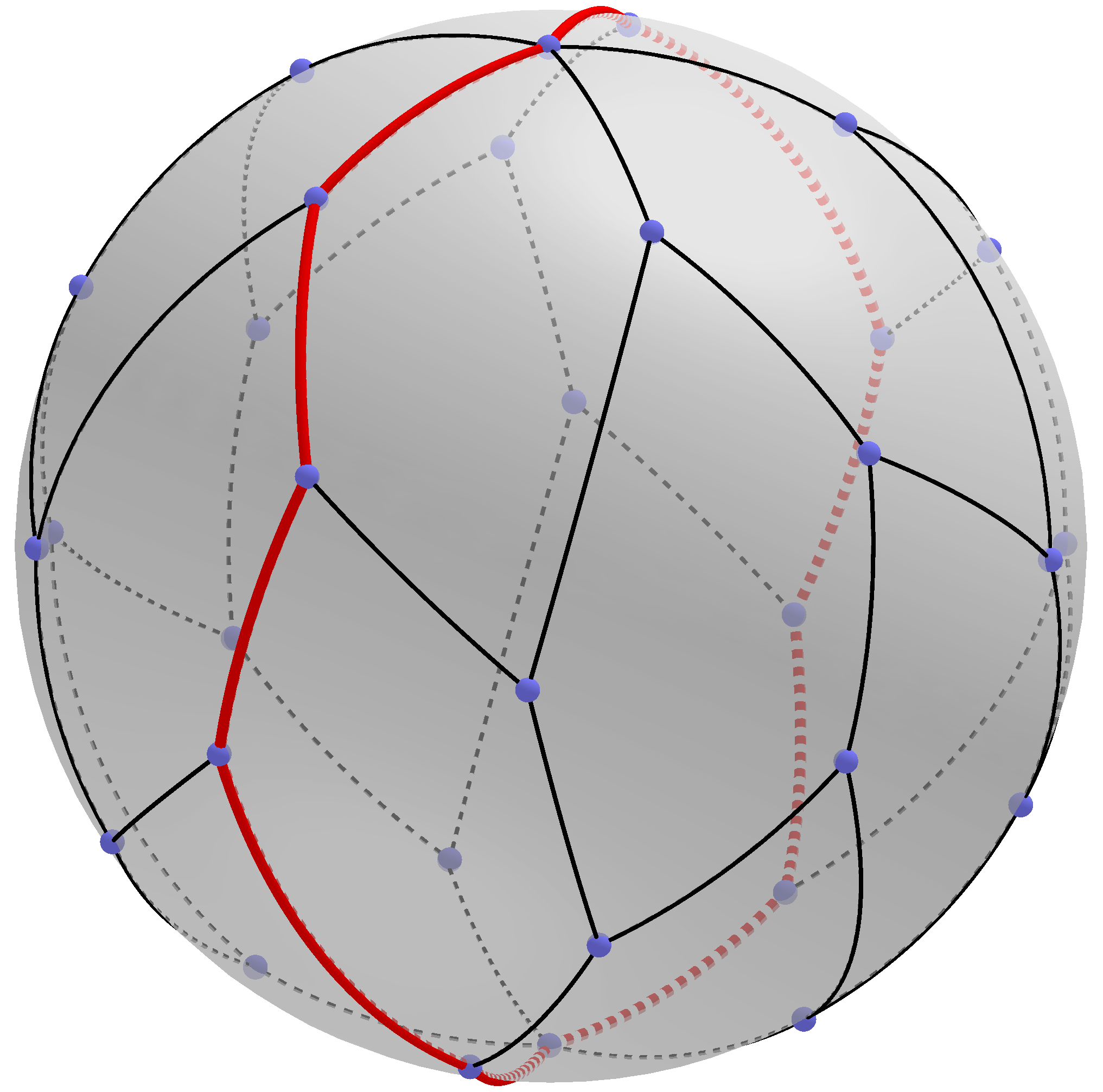}
			\put(120,5){\small $f=20$}
		\end{overpic}   \hspace{60pt}
		\includegraphics[scale=0.183]{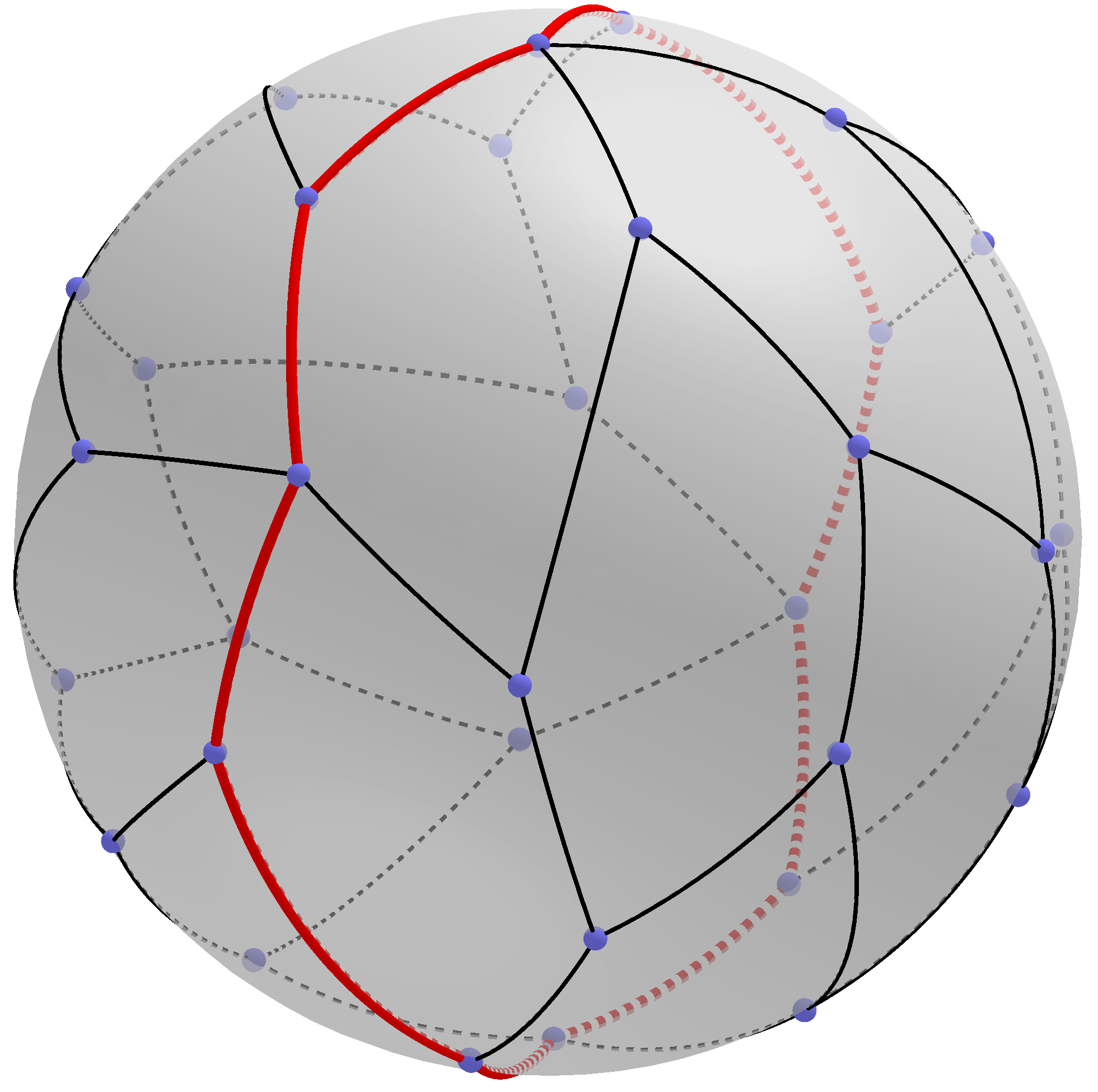}
		\caption{A generic non-symmetric $3$-layer earth map tiling and its standard flip modification.}
		\label{The earth map}
	\end{figure}

\begin{figure}[htp]
		\centering
		\includegraphics[scale=0.125]{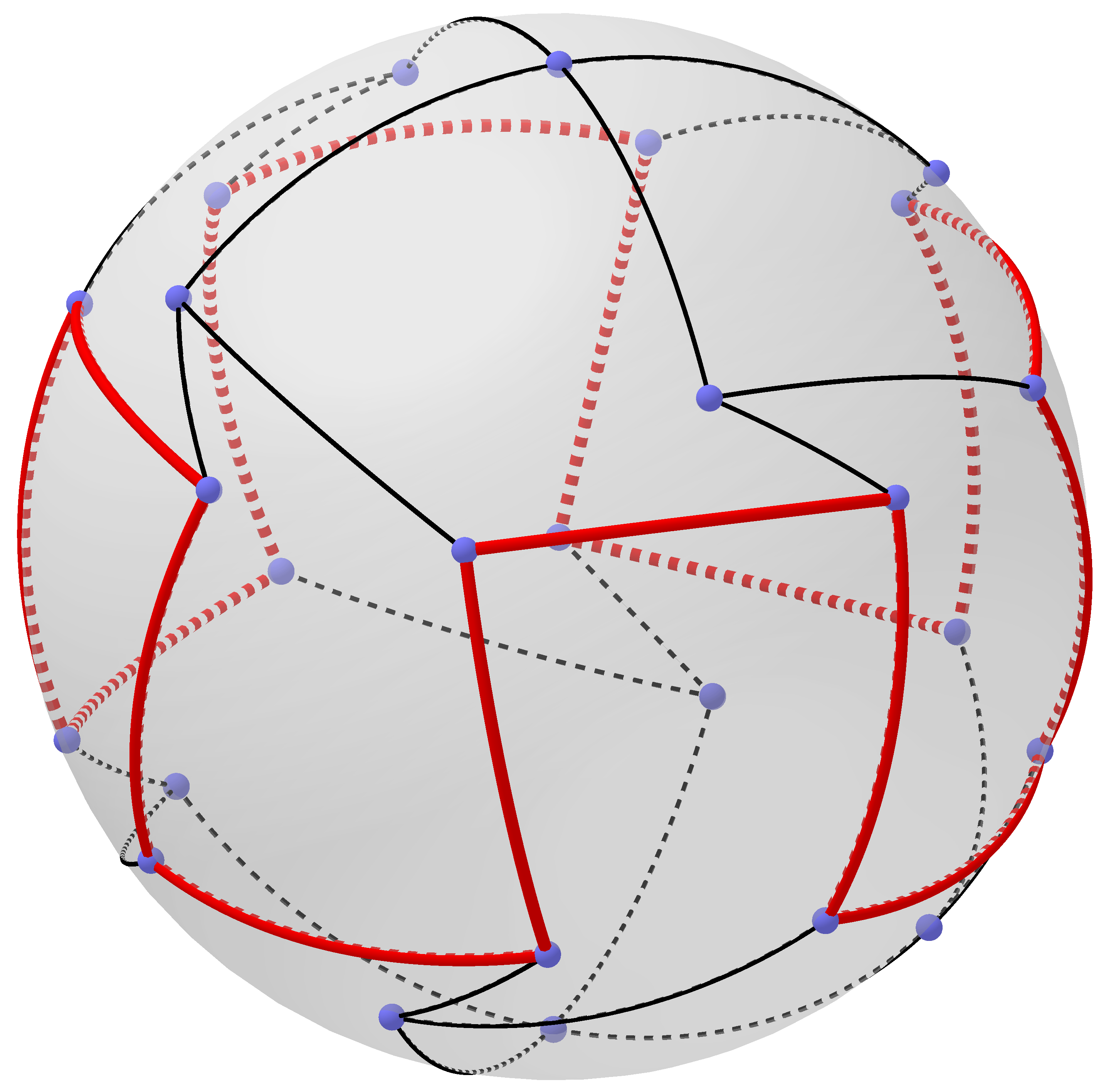} \hspace{60pt}
		\includegraphics[scale=0.14]{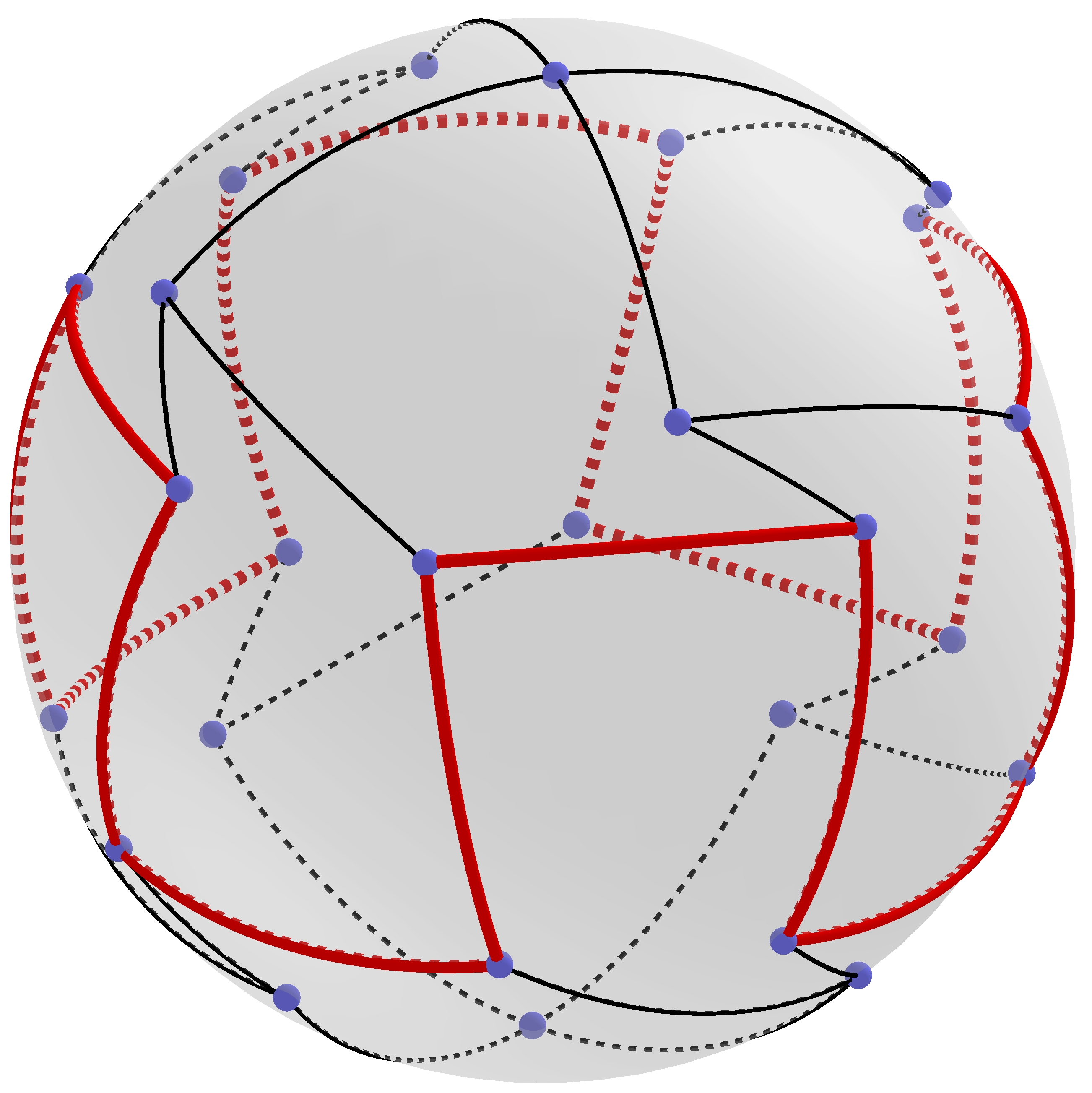}
		\caption{A special flip  modification when $f=16$ and  $\aaa=\frac{\pi}{2}=\ccc$.}
		\label{s16}
	\end{figure}

\begin{figure}[htp]
		\centering
		\includegraphics[angle=90,scale=0.118]{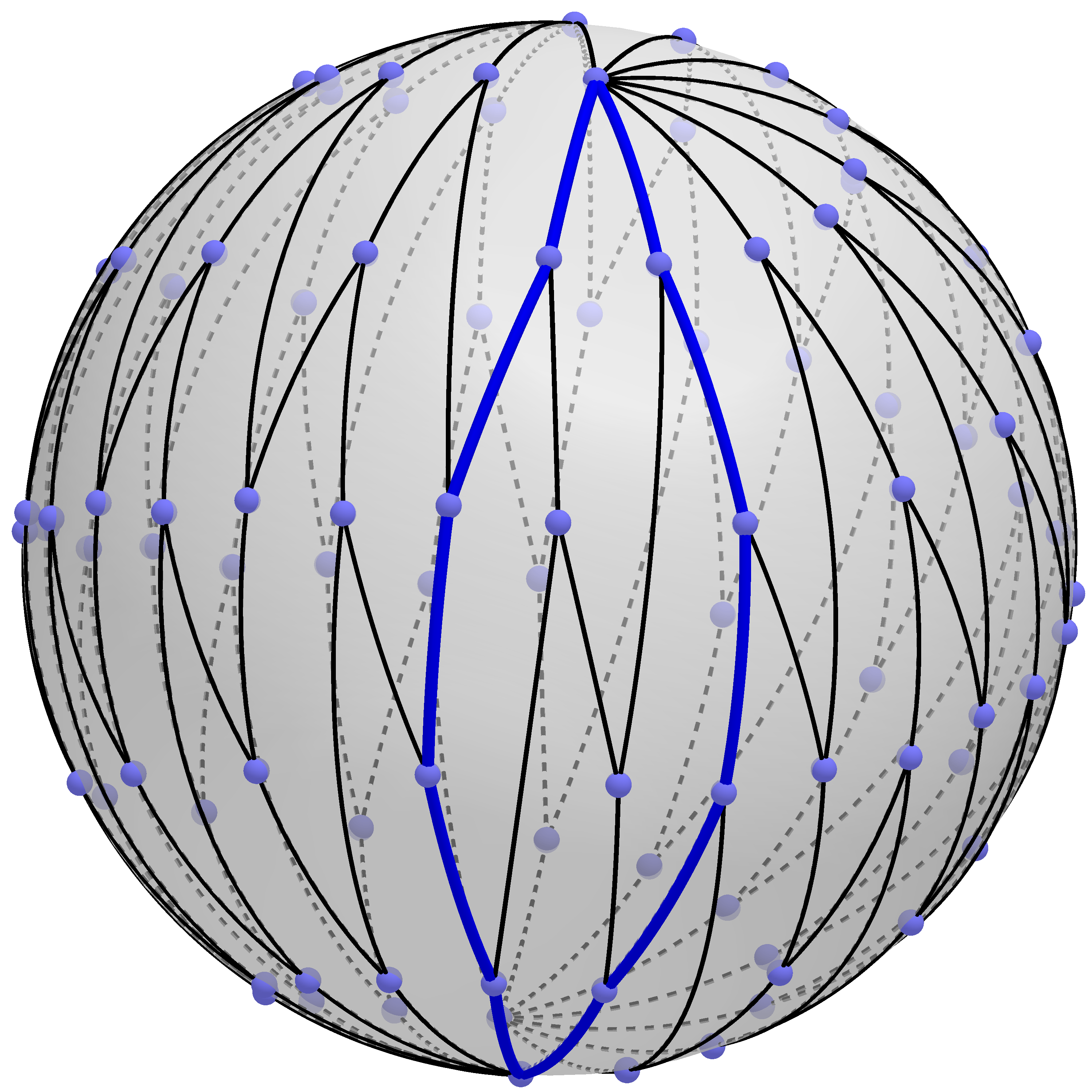}

		\caption{The only UFO for a tiling with AVC $1$ and $60$ tiles.}
		\label{f=60}
	\end{figure}

\begin{figure}[htp]
		\centering
	    \begin{overpic}[scale=0.257]{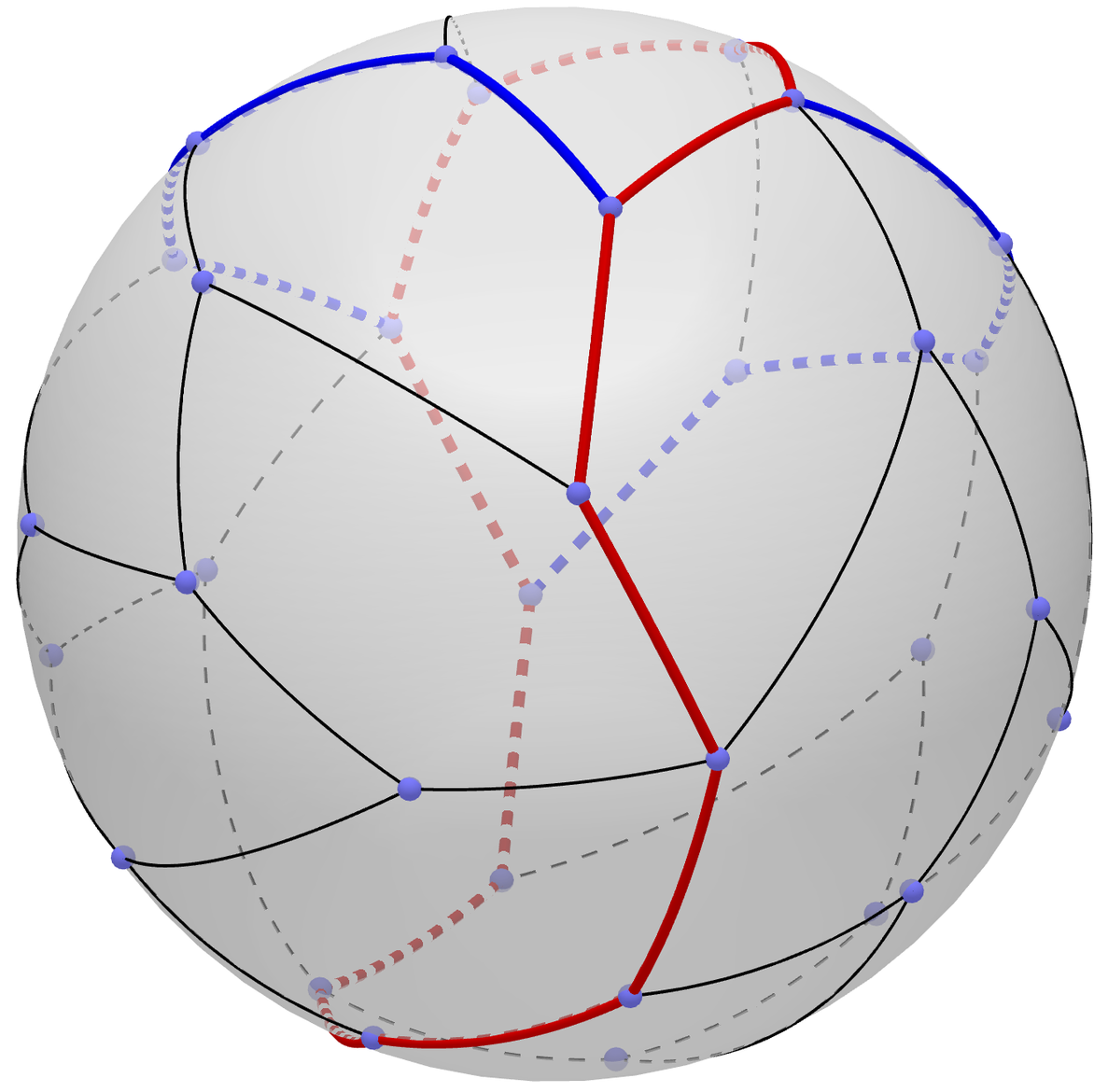}
	    	\put(10,-15){\small AVC $1$, $1$ UFO}
	    \end{overpic}\hspace{10pt}
	    \begin{overpic}[scale=0.142]{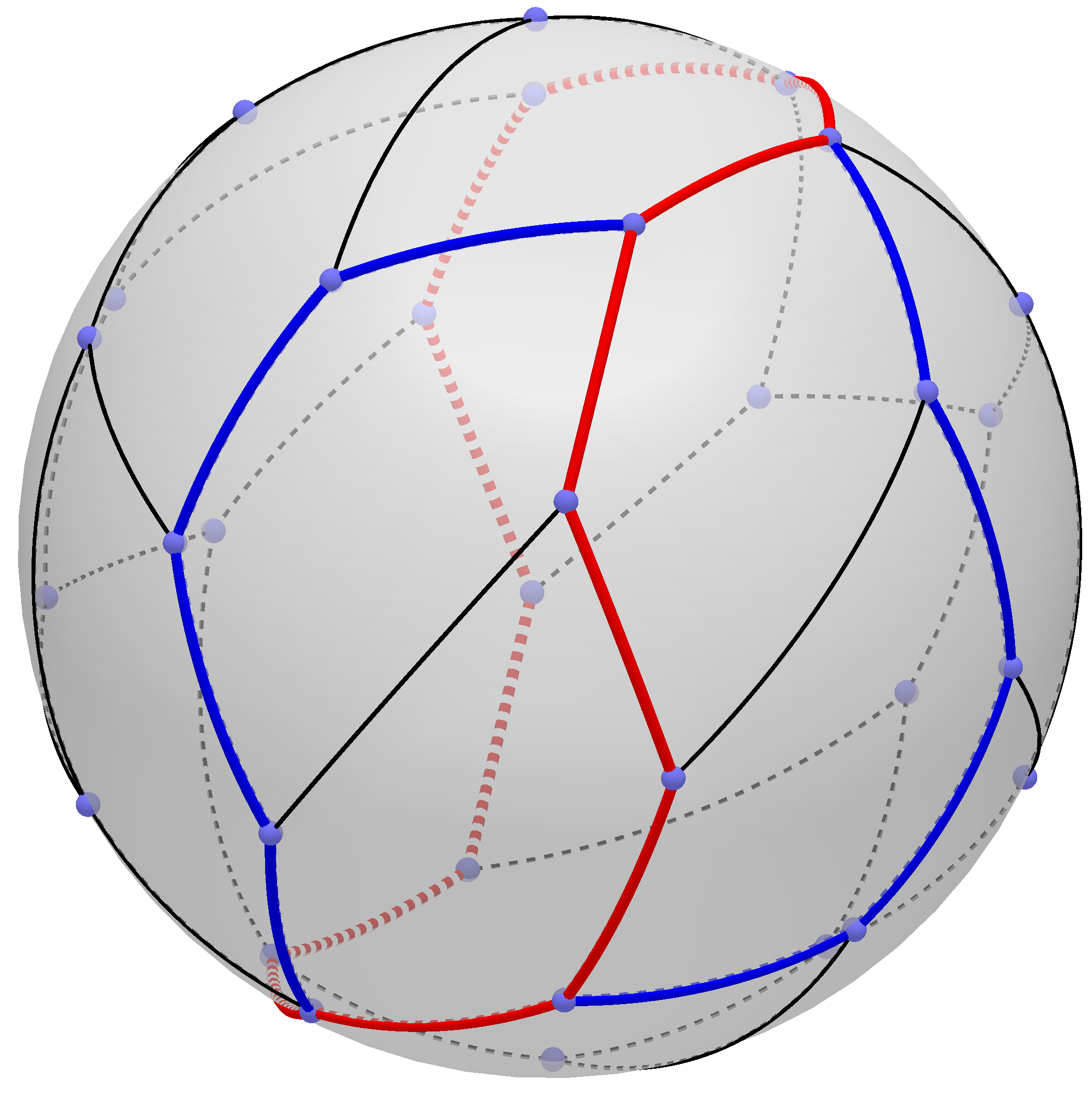}
	    	\put(10,-15){\small AVC $1$, $1$ UFO}
	    \end{overpic}\hspace{11pt}
	    \begin{overpic}[scale=0.1638]{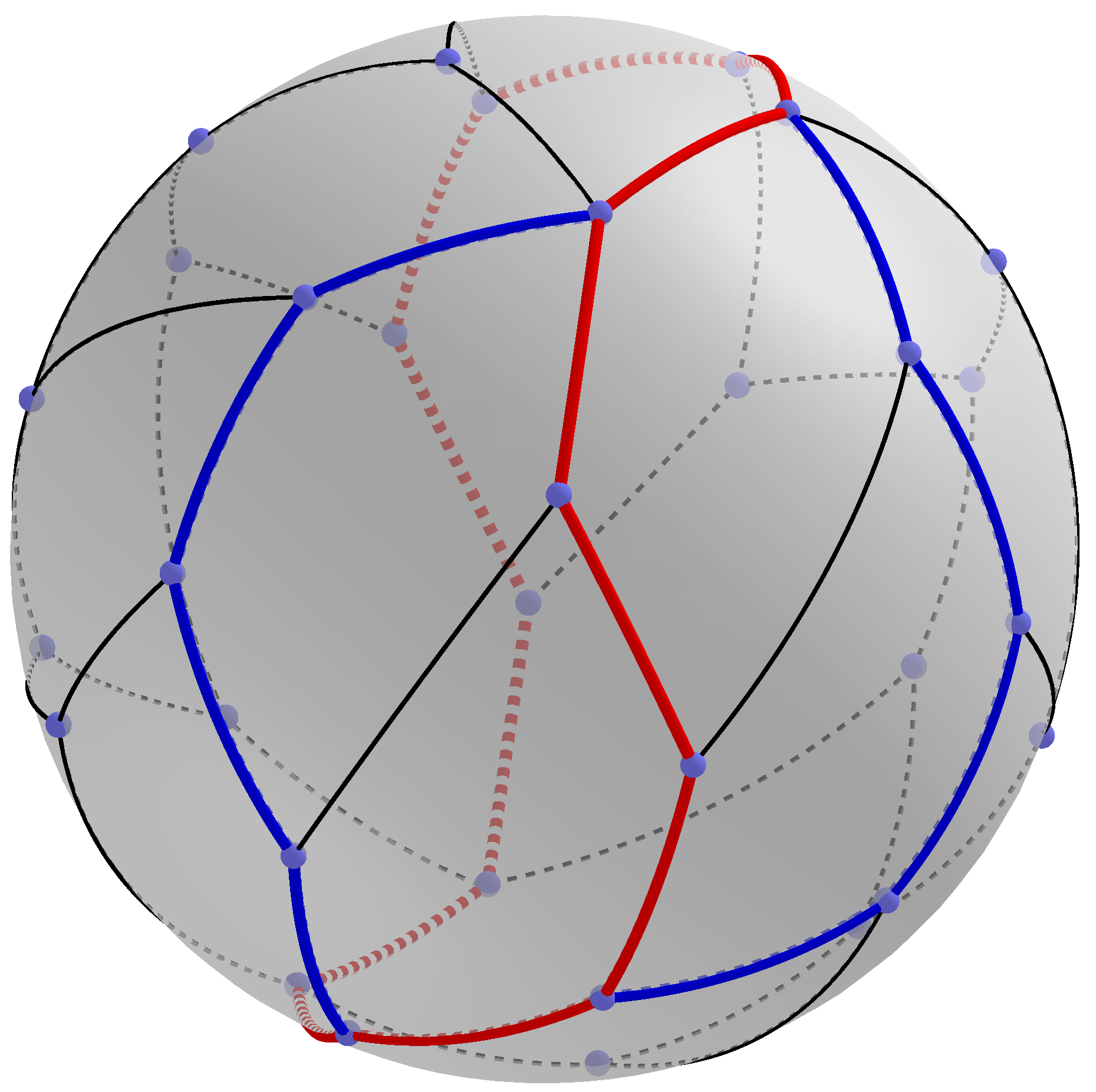}
	    	\put(10,-15){\small AVC $1$, $2$ UFO}
	    \end{overpic}\hspace{11pt}
        \begin{overpic}[scale=0.0546]{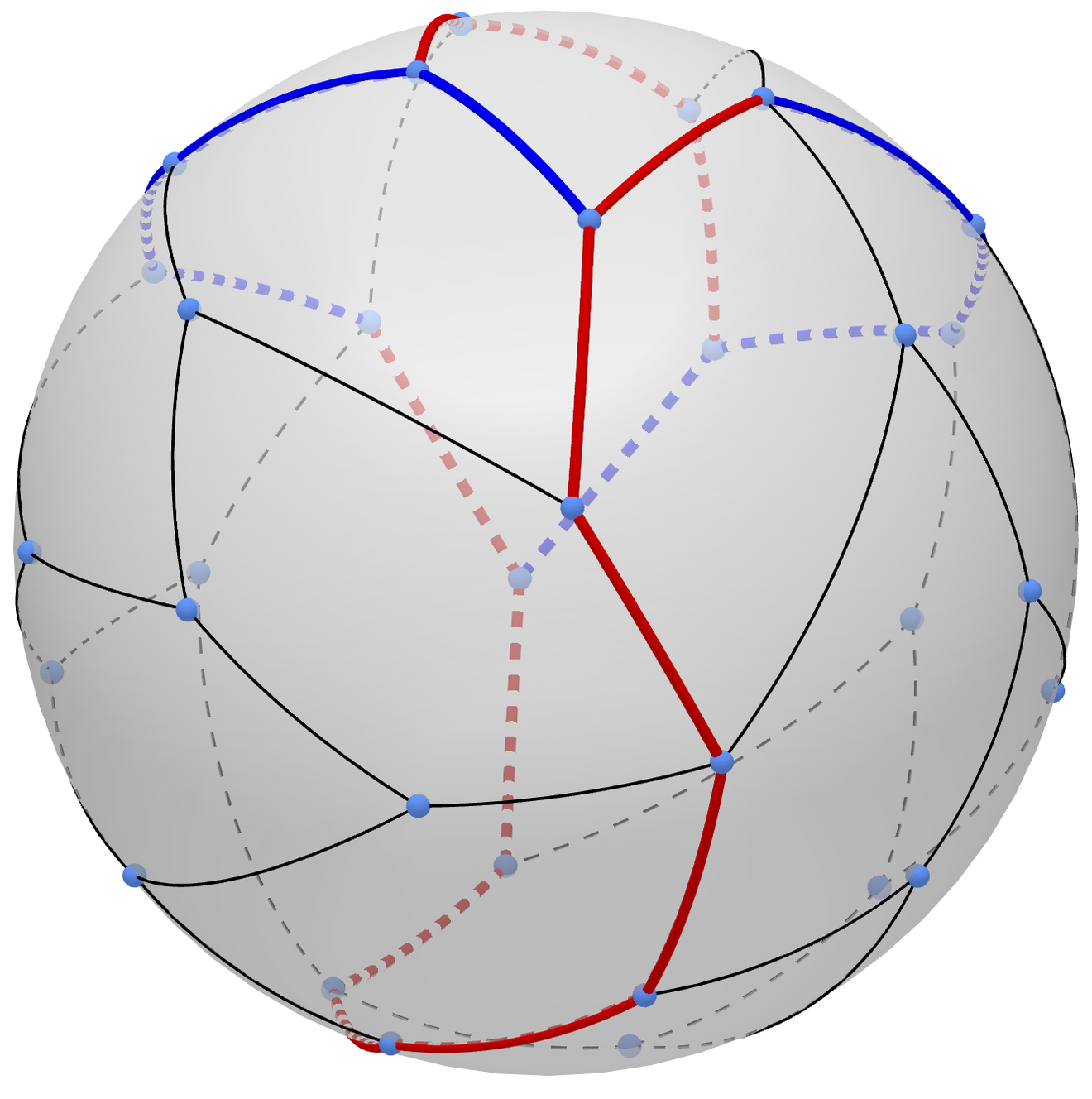}
        	\put(10,-15){\small AVC $2$, $1$ UFO}
        \end{overpic}
		\\\hspace{40pt}

		\begin{overpic}[scale=0.1458]{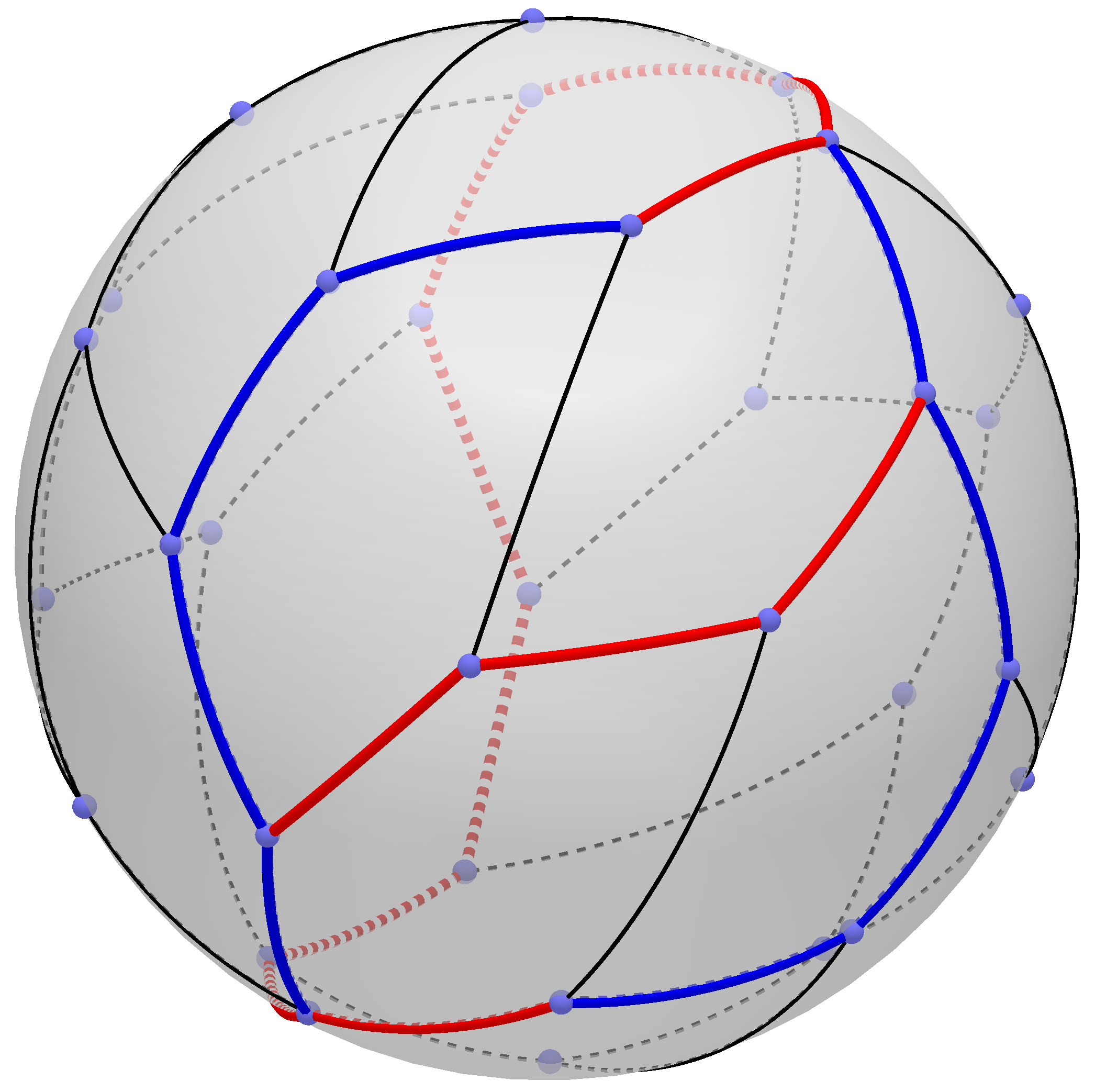}
			\put(10,-15){\small AVC $3$, $3$ UFO}
		\end{overpic}\hspace{10pt}
	    \begin{overpic}[scale=0.1575]{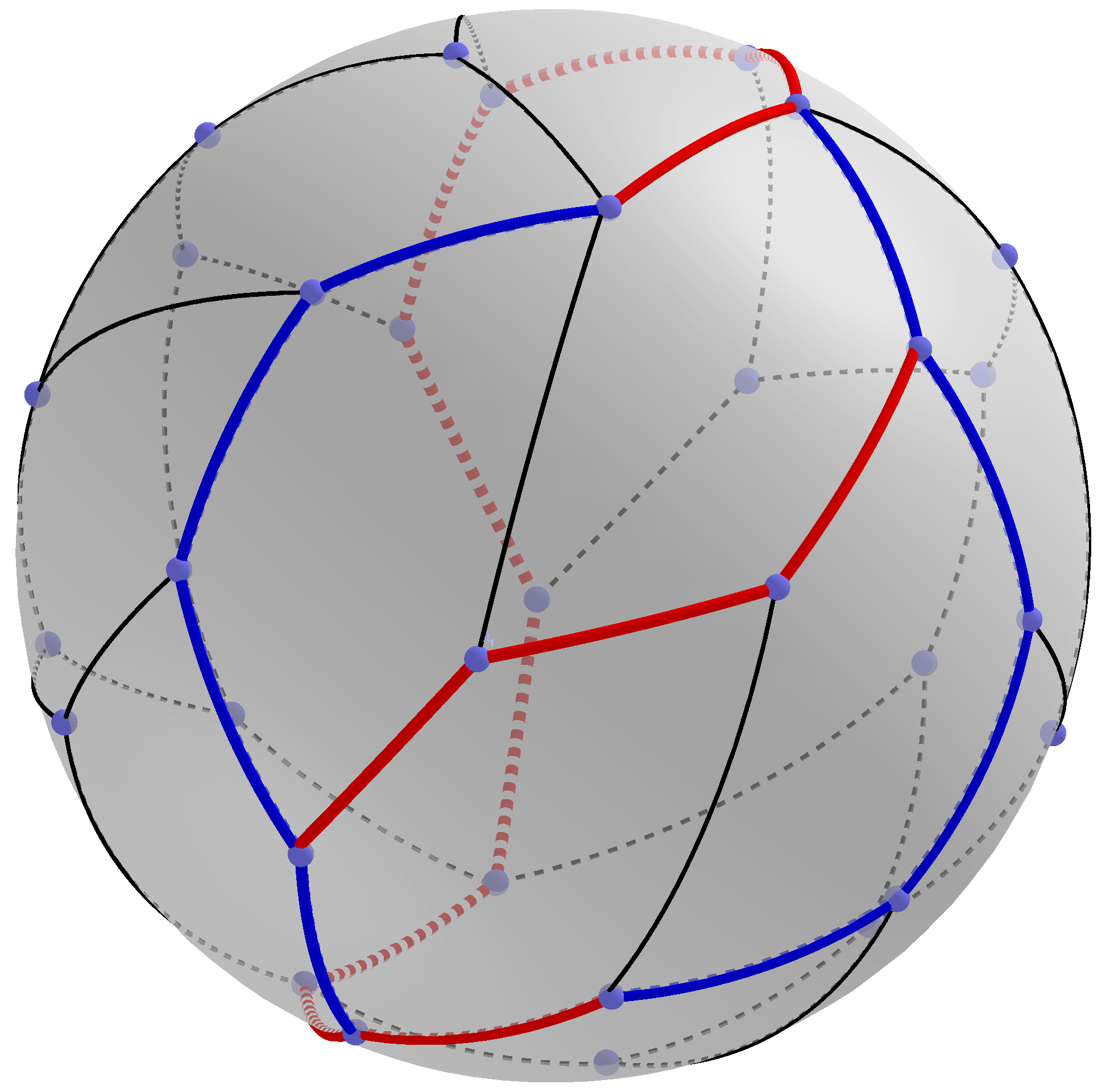}
	    	\put(10,-15){\small AVC $4$, $3$ UFO}
	    \end{overpic}\hspace{10pt}
		\begin{overpic}[scale=0.1395]{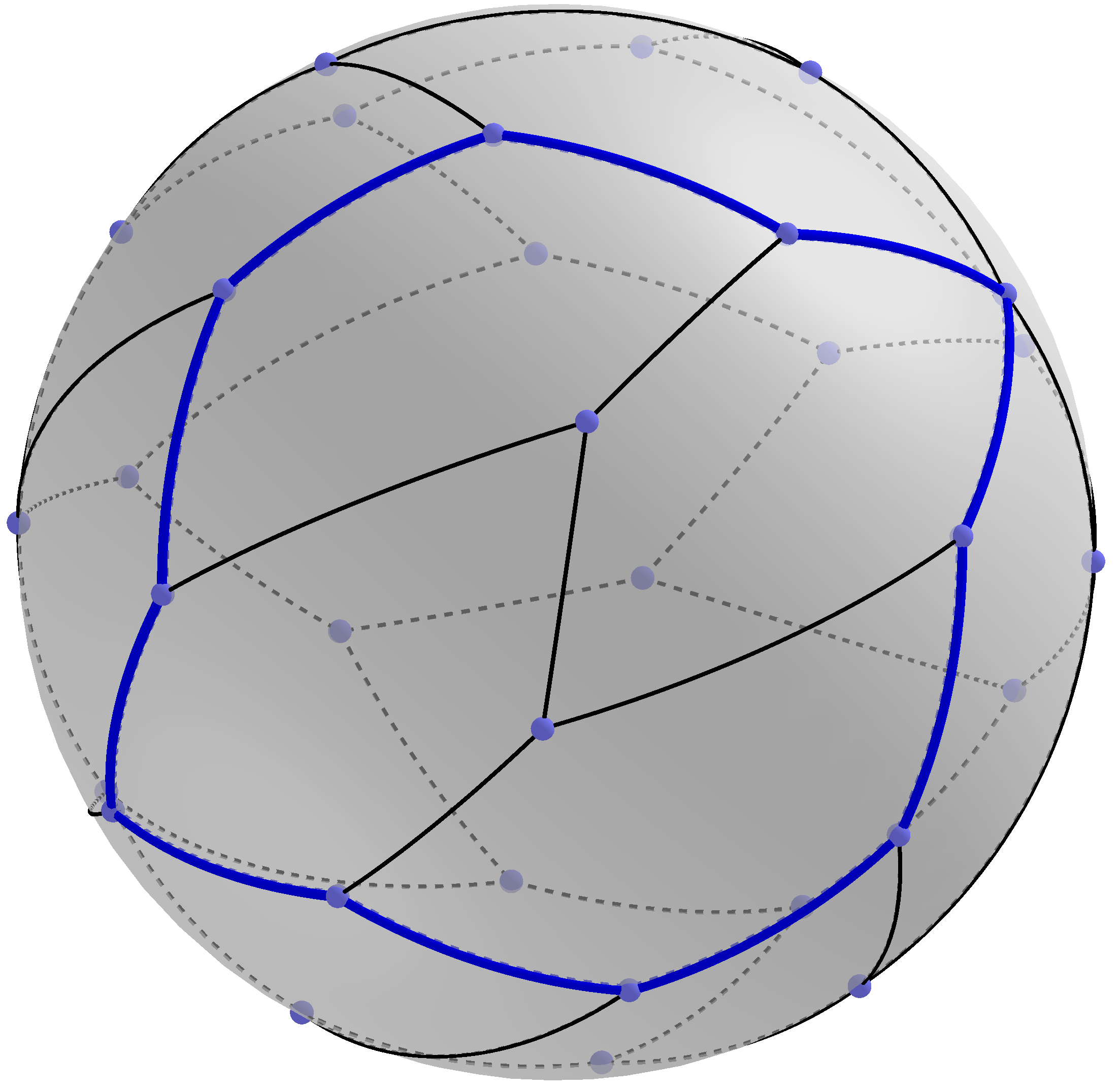}
			\put(10,-15){\small AVC $5$, $3$ UFO}
		\end{overpic}\hspace{10pt}
	    \begin{overpic}[scale=0.24]{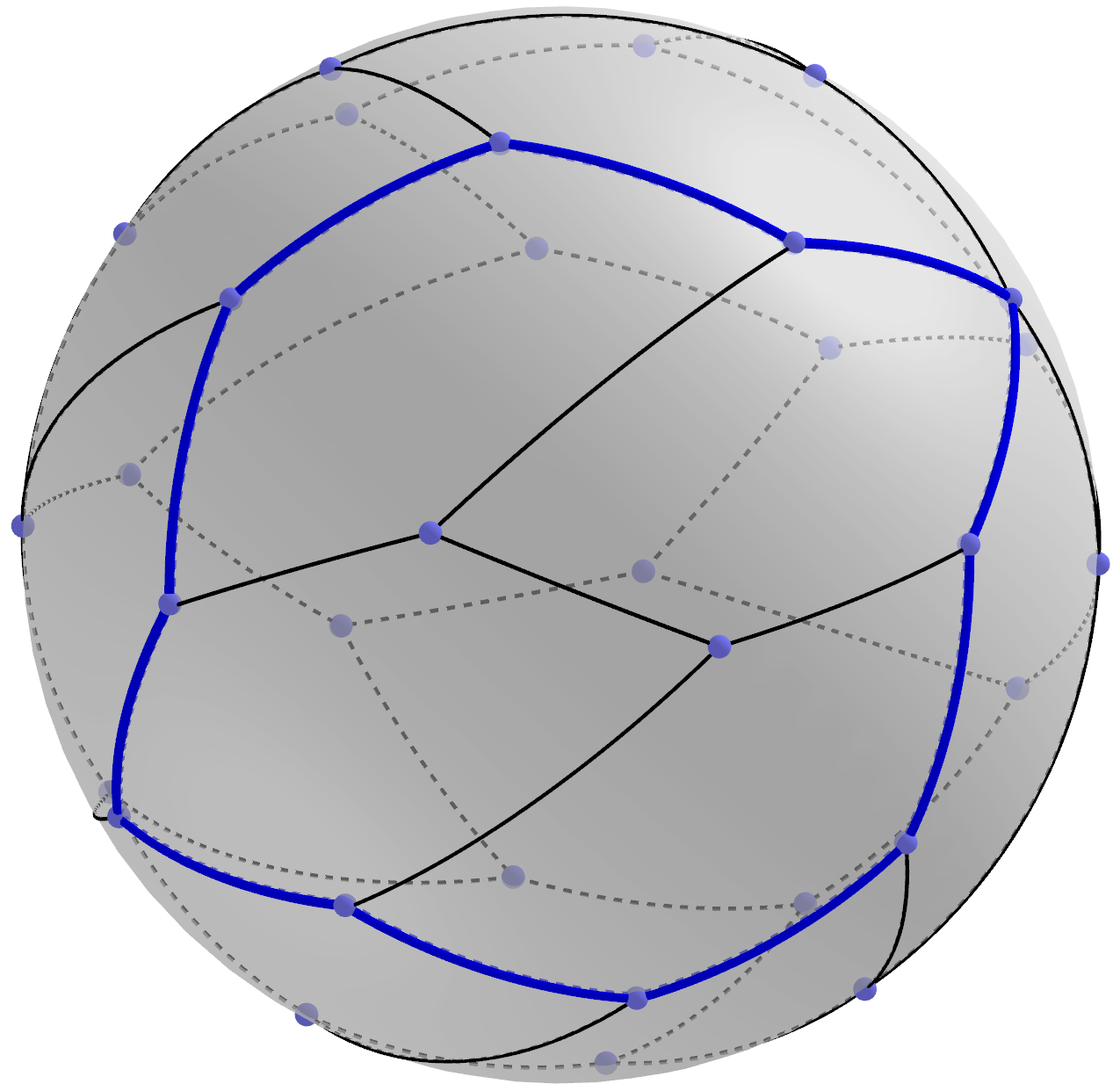}
	    	\put(10,-15){\small AVC $5$, $3$ UFO}
	    \end{overpic}
		\\\hspace{40pt}

		\begin{overpic}[scale=0.18]{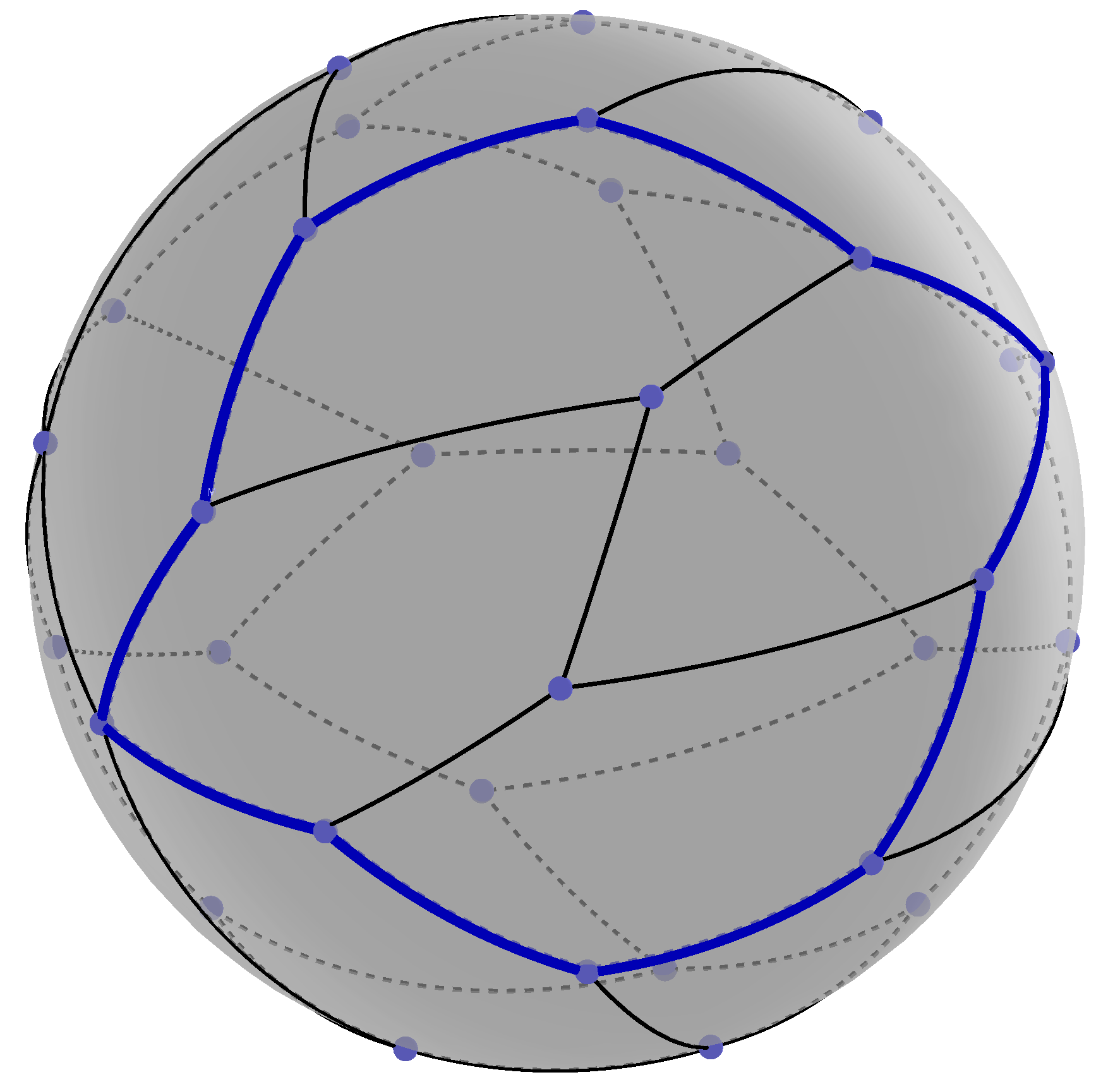}
			\put(10,-15){\small AVC $6$, $4$ UFO}
		\end{overpic}\hspace{10pt}
	    \begin{overpic}[scale=0.179]{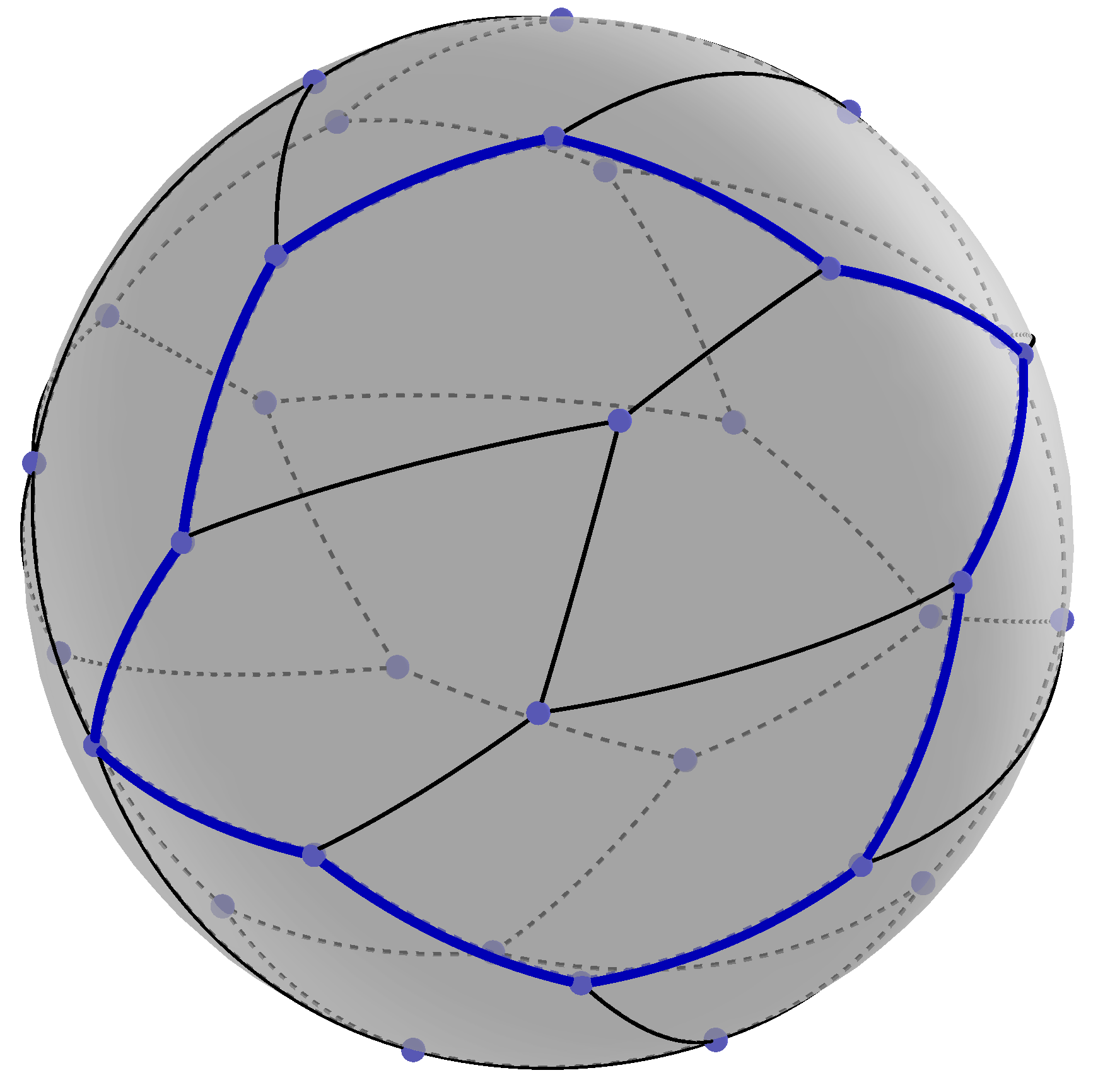}
	    	\put(10,-15){\small AVC $6$, $4$ UFO}
	    \end{overpic}\hspace{10pt}
		\begin{overpic}[scale=0.1503]{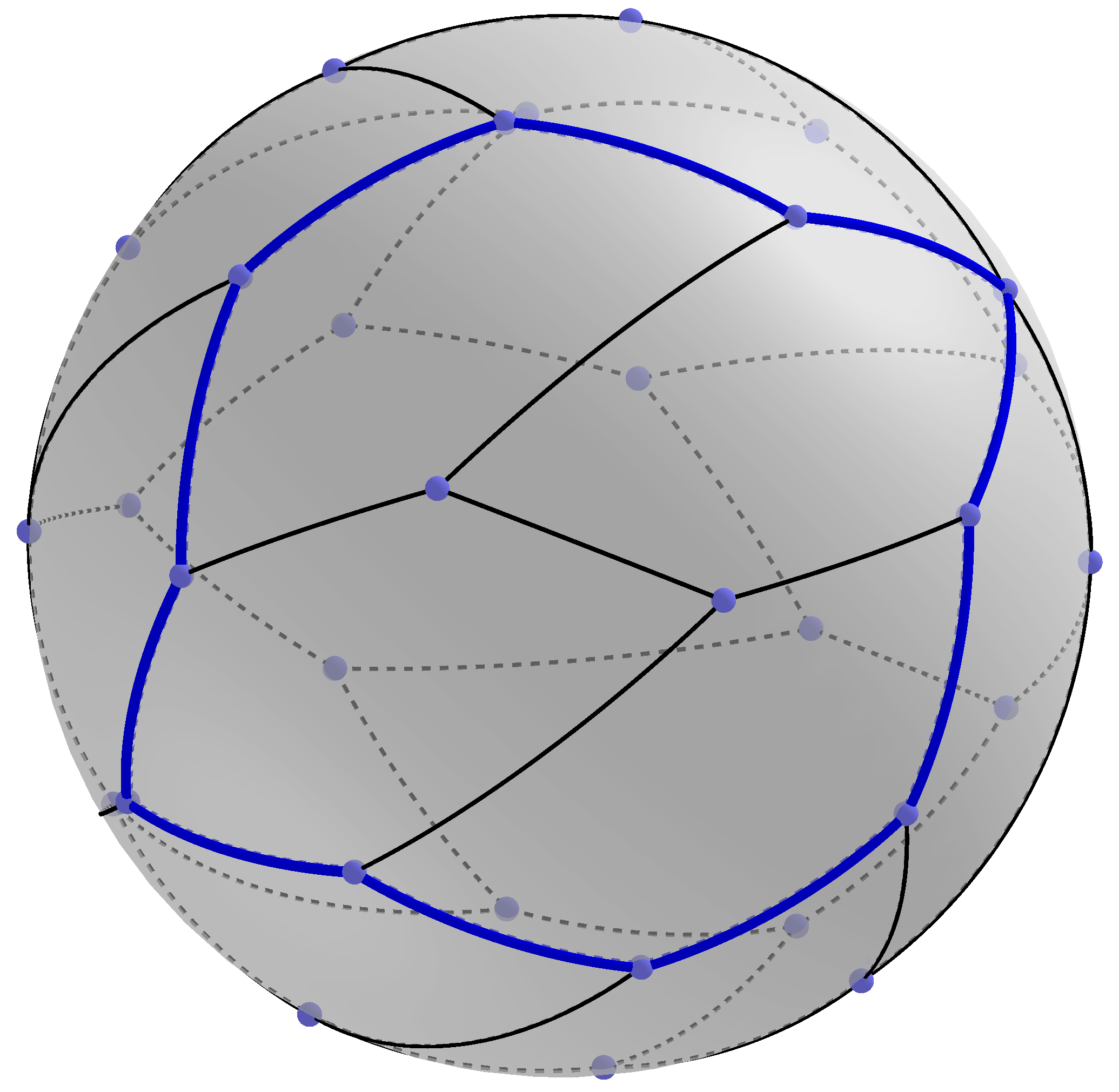}
			\put(10,-15){\small AVC $7$, $4$ UFO}
		\end{overpic}\hspace{10pt}
		\begin{overpic}[scale=0.1458]{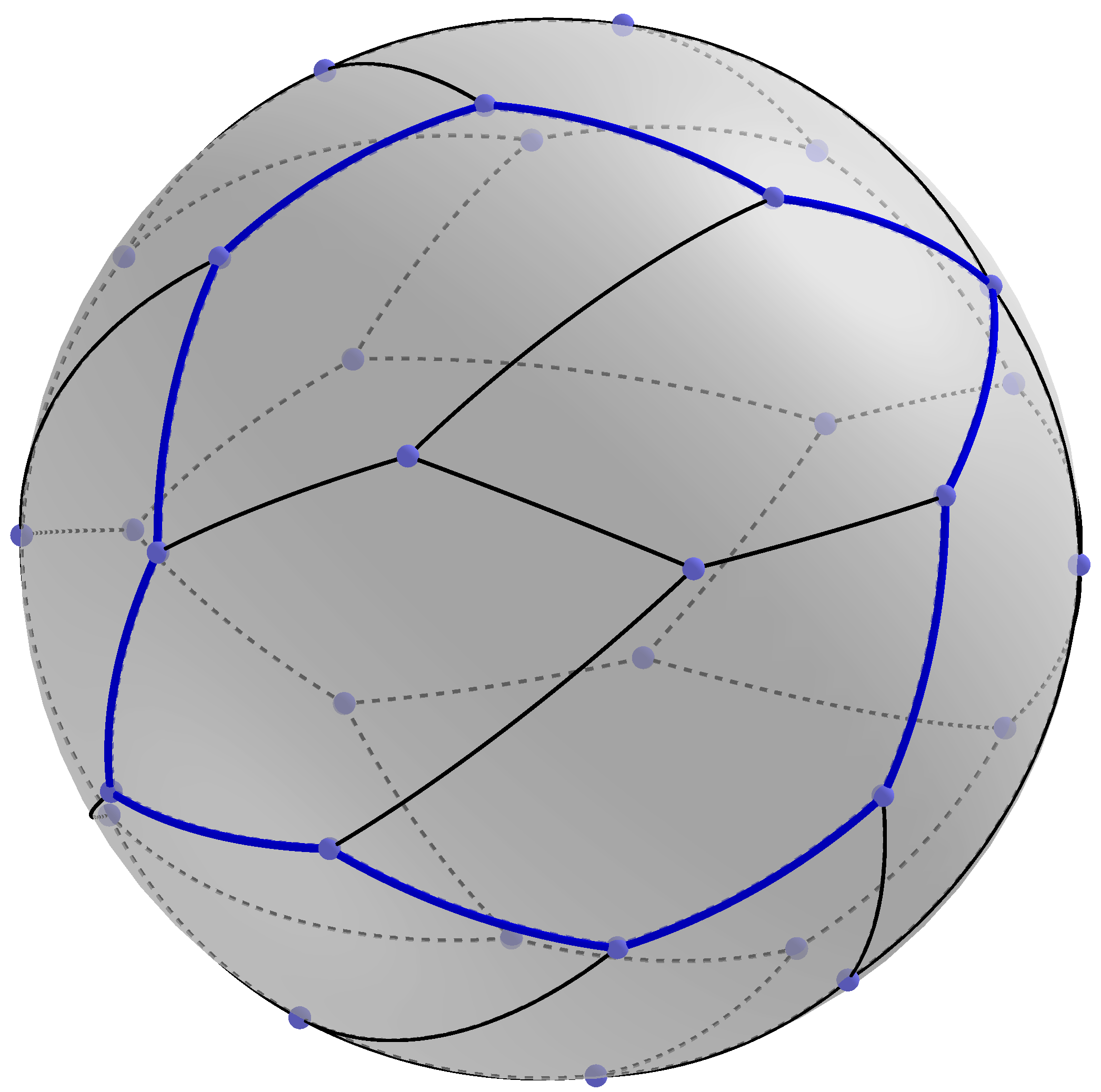}
			\put(10,-15){\small AVC $8$, $4$ UFO}
		\end{overpic}
		\\\hspace{20pt}

		\caption{All modifications for the non-symmetric $3$-layer earth map tiling with $20$ tiles and $\aaa=\bbb=\frac{4\pi}{5}$.}
		\label{various flip}
	\end{figure}

	Henceforth, we only focus on non-symmetric $a^4b$-tilings. We will prove: 

	\begin{theorem} \label{theorem2}
		All non-symmetric $a^4b$-tilings with general angles are:
		\begin{enumerate}
			\item Two $1$-parameter families of pentagonal subdivisions of the octahedron and icosahedron with $24$ tiles $T(24\aaa\ddd\eee,8\ccc^3,6\bbb^4)$ and $60$ tiles $T(60\aaa\ddd\eee,20\ccc^3,12\bbb^5)$ (see the second and third pictures of Figure \ref{poufen} and Table \ref{tab-1}).

			\item A sequence of 1-parameter families of pentagons (whose geometric data are given in Table \ref{tab-2} and shown in  Figure \ref{f=16 3}) each admitting a non-symmetric $3$-layer earth map tiling $T(4m\,\aaa\ddd\eee,2m\,\bbb^2\ccc,2\ccc^m)$ with $4m$ tiles for any $m\ge4$, among which each odd $m=2k+1$ case always admits a standard flip modification: $T((8k+4)\aaa\ddd\eee,4k\,\bbb^2\ccc,4\bbb\ccc^{k+1})$ (see Figure \ref{The earth map}). Moreover,
			\begin{itemize}
				\item For $f=16$, the unique pentagon with angles $(\frac\pi2,\frac{3\pi}{4},\frac\pi2,\ddd,\frac{3\pi}{2}-\ddd)$, $\cos\ddd= \frac{1}{\sqrt[4]{2}}$, $\cos a= \sqrt{\sqrt{2}-1}$, $\cos b= \sqrt{2}\sqrt{\sqrt{2}-1}$, admits a unique flip modification: $T(8\aaa\ddd\eee,8\aaa\bbb^2,8\ccc\ddd\eee,2\ccc^4)$ (see Figure \ref{s16}).
				\item For $f=8k+4$ (i.e. $m=2k+1$) with any $k\ge2$, the unique pentagon determined by $\aaa=\bbb=\pi-\frac{4\pi}{f}$ (in Table \ref{tab-2}), admits a lot more modifications (shown in Table \ref{table-3} and Figure \ref{f=60}, \ref{various flip}).
			\end{itemize}
		\end{enumerate}
	\end{theorem}

	\begin{table}[htp]
		\begin{center}
			\caption{A sequence of $1$-parameter families of pentagons admitting non-symmetric $3$-layer earth map tilings.}\label{tab-2}
			\bgroup
			\resizebox{\linewidth}{!}{
				\begin{tabular}{l}
					\hline
					\multicolumn{1}{c}{$\bbb=\pi-\frac{4\pi}{f}$, $\ccc=\frac{8\pi}{f}$, $f=4m\,(m\ge4)$} \\
					\hline
					\\
					$p:=\cos\aaa,q:=\cos\frac{4\pi}{f},$\\
					$\cos \ddd=\sqrt{\frac{-4 \sqrt{p^{2} (p +1) (p -1)^{2} (2 q -1)^{2} (q -1)^{3} (p q +p -q )}+(8 p^{3}-8 p^{2}-4 p +4) q^{3}-(12 p^{3}-16 p^{2}-4 p +4) q^{2}-(14 p^{2}-p -1) q +4 p^{3}+6 p^{2}-p}{-8 p^{2} q +8 p^{2}+4 q^{2}-4 q +1}}\,$ $(\aaa\in[\frac{\pi}{2},\pi])$, \\
					\\
					$\cos \ddd=\sqrt{\frac{4 \sqrt{p^{2} (p +1) (p -1)^{2} (2 q -1)^{2} (q -1)^{3} (p q +p -q )}+(8 p^{3}-8 p^{2}-4 p +4) q^{3}-(12 p^{3}-16 p^{2}-4 p +4) q^{2}-(14 p^{2}-p -1) q +4 p^{3}+6 p^{2}-p}{-8 p^{2} q +8 p^{2}+4 q^{2}-4 q +1}}\,$   ($\aaa\in(\pi,\frac{3\pi}{2})$; or\\
					$\aaa\in(2 \arctan(\sqrt{-5+4\sqrt{2}}),\frac{\pi}{2})$ for $f=16$, $\aaa\in(\arccos(\frac{\sqrt{2(1-2q)(2q^3-q^2-2q+1+\sqrt{(3q+1)(2q-1)^2(1-q)^3})}}{4q-2}),\frac{\pi}{2})$ for $f\ge20$). \\
					$\eee=2\pi-\aaa-\ddd$, $\cos a= \frac{\sin\frac{\aaa}{2} \sin\frac{2 \pi}{f} \cos\frac{6 \pi}{f}  \cos\ddd-\cos\aaa\sin\frac{4 \pi}{f} \sin(\frac{\aaa}{2}+\ddd) }{\sin\frac{\aaa}{2} \cos\frac{2\pi}{f} \cos\frac{6\pi}{f}  \sin\ddd}$, $\cos b$ is given in Lemma \ref{cosB}.\\
					When $\aaa\in(\arccos(q-2q^2),\pi)$, the pentagon is convex.\\
					$\lim\limits_{f\to\infty}\bbb=\pi,\lim\limits_{f\to\infty}\ccc=0,\lim\limits_{f\to\infty}\ddd=0,\lim\limits_{f\to\infty}\eee=(2\pi-\aaa)\in(\frac{\pi}{2},\frac{3\pi}{2}],\lim\limits_{f\to\infty}b=2\lim\limits_{f\to\infty}a\in[\frac\pi3,\frac\pi2]$.\\
					When $f=16$ and $\aaa\approx0.663\pi$, $\cos a=\cos b\approx0.779$;\\
					When $f=20$ and $\aaa\approx0.628\pi$, $\cos a=\cos b\approx0.776$;\\
					When $f=24$ and $\aaa=0.5\pi$, $\cos a=\cos b\approx0.681$;\\
					When $f=24$ and $\aaa\approx0.522\pi$, $\cos a=\cos b\approx0.706$. (These $4$ equilateral pentagons appeared in \cite{wy3}.)\\
					\hline
			\end{tabular}}
			\egroup
		\end{center}
	\end{table}

	\begin{table}[htp]
		\caption{Various modifications of the non-symmetric $3$-layer earth map tiling with  $\aaa=\bbb=(1-\frac{4}{f})\pi$, $f=8k+4$, $k\ge2$ if both $\aaa\ddd\eee,\bbb\ddd\eee$ appear.}\label{table-3}
		\centering
		\resizebox{\linewidth}{!}{
		\begin{tabular}{c|c}
			\hline
			\multicolumn{2}{c}{Vertices and tilings}  \\
			\hline
			AVC $1$&$\{4\aaa\bbb\ccc,(8k-2)\aaa\ddd\eee,(4k-2)\bbb^2\ccc$, $4\bbb\ddd\eee,2\aaa\ccc^{k+1},2\ccc^{k}\ddd\eee:3\}$ \\[2pt]
			AVC $2$&$T(2\aaa^2\ccc, 6\aaa\bbb\ccc,(8k-6)\aaa\ddd\eee,(4k-4)\bbb^2\ccc,6\bbb\ddd\eee,4\ccc^{k}\ddd\eee)$ \\[2pt]
			AVC $3$&$T(2\aaa^2\ccc,4\aaa\bbb\ccc,(8k-6)\aaa\ddd\eee,(4k-4)\bbb^2\ccc,8\bbb\ddd\eee,2\aaa\ccc^{k+1},2\ccc^{k}\ddd\eee)$ \\[2pt]
			AVC $4$&$T(2\aaa^2\ccc,5\aaa\bbb\ccc,(8k-6)\aaa\ddd\eee,(4k-4)\bbb^2\ccc,7\bbb\ddd\eee,\aaa\ccc^{k+1},3\ccc^{k}\ddd\eee)$ \\[2pt]
			AVC $5$&$\{2\aaa^2\ccc,6\aaa\bbb\ccc,(8k-7)\aaa\ddd\eee,(4k-5)\bbb^2\ccc,8\bbb\ddd\eee,\aaa\ccc^{k+1},3\ccc^{k}\ddd\eee:2k-2\}$ \\[2pt]
			AVC $6$&$\{4\aaa^2\ccc,6\aaa\bbb\ccc,(8k-10)\aaa\ddd\eee,(4k-6)\bbb^2\ccc,10\bbb\ddd\eee,4\ccc^{k}\ddd\eee:2\}$ \\[2pt]
			AVC $7$&$\{4\aaa^2\ccc,7\aaa\bbb\ccc,(8k-11)\aaa\ddd\eee,(4k-7)\bbb^2\ccc,11\bbb\ddd\eee,4\ccc^{k}\ddd\eee:3k-5\}$ \\[2pt]
			AVC $8$&$\{4\aaa^2\ccc,8\aaa\bbb\ccc,(8k-12)\,\aaa\ddd\eee,(4k-8)\bbb^2\ccc,12\bbb\ddd\eee$, $4\ccc^{k}\ddd\eee:(2(k-1)(k-2)+\lfloor\frac{k}{2}\rfloor(k-\lfloor\frac{k}{2}\rfloor))\}$ \\
			\hline
			\end{tabular}}
	\end{table}


    The geometric data in Table \ref{tab-2} depends on the angle parameter $\alpha$,  and Figure \ref{f=16 3} shows some examples in the family with $f=16$: $\aaa$ is close to the lower bound $2 \arctan(\sqrt{-5+4\sqrt{2}})\approx0.43\pi$ in the first picture, is $\pi$ in the second, and is close to the upper bound $\frac{3\pi}{2}$ in the third. Note that the pentagon in Figure \ref{s16} corresponds to the parameter $\alpha=\frac{\pi}{2}$ in this family.

   All tilings in Table \ref{table-3} contain at least one special block called a UFO\label{UFO}. Here we use UFO (Unidentified Flying Object) to denote the symmetric $4$-tile patch bounded by the blue lines in Figure \ref{f=60}, since it looks like a flying saucer or spaceship. Any UFO contains two $\bbb\ddd\eee$ in the interior and is centrally symmetric together with two mirror symmetries. Figure \ref{various flip} shows all tilings with $f=20$ and $8$ types of AVC in Table \ref{table-3}, where the number of UFOs is also counted for each tiling by playing with the 3D pictures in GeoGebra$^\circledR$. Note that only AVC $5,7,8$ in Table \ref{table-3} will admit more different tilings as $f=8k+4$ becomes bigger. However there are still only $4$ UFOs for all these tilings with AVC $7,8$, and only $3$ UFOs for all these tilings with AVC $5$. For all other AVC in Table \ref{table-3}, the number of different tilings and the number of different UFOs are the same as Case $f=20$ in Figure \ref{various flip}.

	A sequel \cite{hllw} of this paper will handle $a^4b$-pentagons with all angles being rational in degrees by solving some trigonometric Diophantine equations, to complete our full classification of edge-to-edge tilings of the sphere by congruent pentagons. It turns out that there are no more new prototiles, nor new tilings.
	Therefore, this paper, its sequel \cite{hllw}, and the earlier papers \cite{wy1,wy2,wy3} together give all edge-to-edge spherical tilings by congruent pentagons. This is in sharp contrast to the quadrilateral case \cite{lw2}, where there is an abundance of tilings of quadrilaterals with all angles being rational.

	\subsubsection*{Outline of the paper}

	The basic facts of tilings and basic techniques for complete classification are presented in Section \ref{basic_facts}, including some technical results specific to $a^4b$-tilings. We call a vertex with three $a$-edges an \textbf{$a^3$-vertex}, and a vertex with two $a$-edges and one $b$-edge an \textbf{$a^2b$-vertex}. Tilings with two different types of $a^2b$-vertices are studied in Section \ref{two $a^2b$-vertices}. Tilings with only one type of $a^2b$-vertex are studied in Section \ref{ad2,bd2 or be2 as the unique $a^2b$-vertex}, \ref{bde as the unique $a^2b$-vertex} and \ref{ade as the unique $a^2b$-vertex}. The non-edge-to-edge quadrilateral tilings induced from degenerate pentagonal tilings in this paper are summarized in the last section.

    \bigskip
    \subsection*{Acknowledgments}
    The results of this paper were reported at the 11th Shanghai Conference on Combinatorics (11SHCC) on May 24 - 28, 2024. The authors would like to thank Shanghai Jiao Tong University and the organizers for great hospitality.

	\section{Basic Facts}\label{basic_facts}
	For the basics of tilings, we adopt the definitions and conventions in the books by Gr\"{u}nbaum and Shephard \cite{gs}, and by Adams \cite{a1}. Angles are in units of $\pi$ radians, except in Lemmas \ref{formula} and \ref{cosB}, where they are in radians. All side lengths are in radians. So the sum of all angles at a vertex is $2$. We present some basic facts and techniques in this section.

	Let $v,e,f$ be the numbers of vertices, edges, and tiles in a pentagonal tiling. Let $v_k$ be the number of vertices of degree $k$. We have the following basic formulas in \cite{wy1}.
	\begin{align}
		f
		&=12+2\sum_{k\ge 4}(k-3)v_k=12+2v_4+4v_5+6v_6+\cdots, \label{vcountf} \\
		v_3
		&=20+\sum_{k\ge 4}(3k-10)v_k=20+2v_4+5v_5+8v_6+\cdots. \label{vcountv}
	\end{align}
	So $f\ge 12$ and $v_3\ge f+8 \ge 20$. The equality \eqref{vcountv} shows that most vertices have degree $3$.

	\begin{lemma}[{\cite[Lemma 4]{wy1}}] \label{anglesum}
		If all tiles in a tiling of the sphere by $f$ pentagons have the same five angles $\aaa,\bbb,\ccc,\ddd,\eee$, then
		\[
		\aaa+\bbb+\ccc+\ddd+\eee = 3+\tfrac{4}{f} ,
		\]
		ranging in $(3,\frac{10}3]$. In particular no vertex contains all five angles.
	\end{lemma}

	\begin{lemma}[{\cite[Lemma 1]{wy1}}] \label{base_tile}
		Any pentagonal tiling of the sphere has a ``special" tile, such that four vertices have degree $3$ and the fifth vertex has degree $3$, $4$ or $5$.
	\end{lemma}

	\begin{lemma}[Parity Lemma, {\cite[Lemma 10]{wy2}}]\label{parity}
		In an $a^4b$-tiling, the total number of $ab$-angles (i.e. $\delta$ and $\epsilon$) at any vertex is even.
	\end{lemma}

	In a tiling of the sphere by $f$ congruent tiles, each angle of the tile appears $f$ times in total. If one vertex has more $\aaa$ than $\bbb$, there must exist another vertex with more $\bbb$ than $\aaa$.
	Such global counting induces many interesting and useful results.

	\begin{lemma}[Balance Lemma, {\cite[Lemma 6]{lw3}}]\label{balance}
		In an $a^4b$-tiling, if either $\ddd^2\cdots$ or $\eee^2\cdots$ is not a vertex, then any vertex either has no $\ddd,\eee$, or is of the form $\ddd\eee\cdots$ with no more $\ddd,\eee$ in the remainder.
	\end{lemma}

	\begin{lemma}\label{l2}
		In an $a^4b$-tiling, an $a^3$-vertex and an $a^2b$-vertex must both appear.
	\end{lemma}

	\begin{proof}
		By Lemma \ref{base_tile}, we know that either $\ddd\cdots$ or $\eee\cdots$ is an $a^2b$-vertex. If there is no $a^3$-vertex, then every degree 3 vertex is an $a^2b$-vertex. By the equations (\ref{vcountf}) and (\ref{vcountv}), we have $\#\ddd+\#\eee\ge2v_3>2f$ in the tiling, contradicting $\#\ddd+\#\eee=2f$.
	\end{proof}

	The very useful tool \textit{adjacent angle deduction} (abbreviated as \textbf{AAD}) has been introduced in \cite[Section 2.5]{wy1}. We give a quick review here using  Figure \ref{anglecombo}. Let ``$\thin$'' denote an a-edge and ``$\thick$'' denote a b-edge. Then we indicate the arrangements of angles and edges by denoting the vertices as $\thin\bbb\thin\bbb\thin\ddd\thick\ddd\thin$. The notation can be reversed, such as $\thin\bbb\thin\bbb\thin\ddd\thick\ddd\thin=\thin\ddd\thick\ddd\thin\bbb\thin\bbb\thin$; and it can be rotated, such as $\thin\bbb\thin\bbb\thin\ddd\thick\ddd\thin=\thin\bbb\thin\ddd\thick\ddd\thin\bbb\thin=\thick\ddd\thin\bbb\thin\bbb\thin\ddd\thick$. We also denote the first  vertex in Figure \ref{anglecombo} as $\ddd\thick\ddd\cdots$, $\thin\ddd\thick\ddd\thin\cdots,\bbb\thin\ddd\cdots$, $\thin\ddd\thick\ddd\thin\bbb\thin\cdots$, and denote the consecutive angle segments as $\ddd\thick\ddd$, $\thin\ddd\thick\ddd\thin,\bbb\thin\ddd$, $\thin\ddd\thick\ddd\thin\bbb\thin$.

	\begin{figure}[htp]
		\centering
		\begin{tikzpicture}[>=latex,scale=1]

			\begin{scope}[xshift=-4cm]

				\draw
				(0,0) -- (0,0.8) -- (-0.6,1.2) -- (-1.2,0.6) -- (-0.8,0) -- (0,0);
				\draw
				(0,0) -- (0,0.8) -- (0.6,1.2) -- (1.2,0.6) -- (0.8,0) -- (0,0);
				\draw
				(0,0) -- (0,-0.8) -- (0.6,-1.2) -- (1.2,-0.6) -- (0.8,0) -- (0,0);
				\draw
				(0,0) -- (0,-0.8) -- (-0.6,-1.2) -- (-1.2,-0.6) -- (-0.8,0) -- (0,0);

				\draw[line width=1.5]
				(-1.2,0.6) -- (-0.8,0);
				\draw[line width=1.5]
				(0,0) -- (0.8,0);
				\draw[line width=1.5]
				(-1.2,-0.6) -- (-0.8,0);

				\node at (-0.7,0.2) {\small $\delta$};
				\node at (-1,0.6) {\small $\epsilon$};
				\node at (-0.6,0.95) {\small $\gamma$};
				\node at (-0.2,0.7) {\small $\alpha$};
				\node at (-0.2,0.2) {\small $\beta$};

				\node at (0.7,0.2) {\small $\eee$};
				\node at (1,0.6) {\small $\ccc$};
				\node at (0.6,0.95) {\small $\aaa$};
				\node at (0.2,0.7) {\small $\bbb$};
				\node at (0.2,0.2) {\small $\ddd$};

				\node at (0.7,-0.2) {\small $\eee$};
				\node at (1,-0.6) {\small $\ccc$};
				\node at (0.6,-0.95) {\small $\aaa$};
				\node at (0.2,-0.7) {\small $\bbb$};
				\node at (0.2,-0.2) {\small $\ddd$};

				\node at (-0.7,-0.2) {\small $\delta$};
				\node at (-1,-0.6) {\small $\epsilon$};
				\node at (-0.6,-0.95) {\small $\gamma$};
				\node at (-0.2,-0.7) {\small $\alpha$};
				\node at (-0.2,-0.2) {\small $\beta$};

			\end{scope}

			\begin{scope}[xshift=0cm]

				\draw
				(0,0) -- (0,0.8) -- (-0.6,1.2) -- (-1.2,0.6) -- (-0.8,0) -- (0,0);
				\draw
				(0,0) -- (0,0.8) -- (0.6,1.2) -- (1.2,0.6) -- (0.8,0) -- (0,0);
				\draw
				(0,0) -- (0,-0.8) -- (0.6,-1.2) -- (1.2,-0.6) -- (0.8,0) -- (0,0);
				\draw
				(0,0) -- (0,-0.8) -- (-0.6,-1.2) -- (-1.2,-0.6) -- (-0.8,0) -- (0,0);

				\draw[line width=1.5]
				(-1.2,0.6) -- (-0.8,0);
				\draw[line width=1.5]
				(0,0) -- (0.8,0);
				\draw[line width=1.5]
				(0,-0.8) -- (-0.6,-1.2);

				\node at (-0.7,0.2) {\small $\ddd$};
				\node at (-1,0.6) {\small $\eee$};
				\node at (-0.6,0.95) {\small $\ccc$};
				\node at (-0.2,0.7) {\small $\aaa$};
				\node at (-0.2,0.2) {\small $\beta$};

				\node at (0.7,0.2) {\small $\eee$};
				\node at (1,0.6) {\small $\ccc$};
				\node at (0.6,0.95) {\small $\aaa$};
				\node at (0.2,0.7) {\small $\bbb$};
				\node at (0.2,0.2) {\small $\ddd$};

				\node at (0.7,-0.2) {\small $\eee$};
				\node at (1,-0.6) {\small $\ccc$};
				\node at (0.6,-0.95) {\small $\aaa$};
				\node at (0.2,-0.7) {\small $\bbb$};
				\node at (0.2,-0.2) {\small $\ddd$};

				\node at (-0.7,-0.2) {\small $\aaa$};
				\node at (-1,-0.6) {\small $\ccc$};
				\node at (-0.6,-0.95) {\small $\eee$};
				\node at (-0.2,-0.7) {\small $\ddd$};
				\node at (-0.2,-0.2) {\small $\beta$};

			\end{scope}

			\begin{scope}[xshift=4cm]

				\draw
				(0,0) -- (0,0.8) -- (-0.6,1.2) -- (-1.2,0.6) -- (-0.8,0) -- (0,0);
				\draw
				(0,0) -- (0,0.8) -- (0.6,1.2) -- (1.2,0.6) -- (0.8,0) -- (0,0);
				\draw
				(0,0) -- (0,-0.8) -- (0.6,-1.2) -- (1.2,-0.6) -- (0.8,0) -- (0,0);
				\draw
				(0,0) -- (0,-0.8) -- (-0.6,-1.2) -- (-1.2,-0.6) -- (-0.8,0) -- (0,0);

				\draw[line width=1.5]
				(0,0.8) -- (-0.6,1.2);
				\draw[line width=1.5]
				(0,0) -- (0.8,0);
				\draw[line width=1.5]
				(0,-0.8) -- (-0.6,-1.2);

				\node at (-0.7,0.2) {\small $\aaa$};
				\node at (-1,0.6) {\small $\ccc$};
				\node at (-0.6,0.95) {\small $\eee$};
				\node at (-0.2,0.7) {\small $\ddd$};
				\node at (-0.2,0.2) {\small $\beta$};

				\node at (0.7,0.2) {\small $\eee$};
				\node at (1,0.6) {\small $\ccc$};
				\node at (0.6,0.95) {\small $\aaa$};
				\node at (0.2,0.7) {\small $\bbb$};
				\node at (0.2,0.2) {\small $\ddd$};

				\node at (0.7,-0.2) {\small $\eee$};
				\node at (1,-0.6) {\small $\ccc$};
				\node at (0.6,-0.95) {\small $\aaa$};
				\node at (0.2,-0.7) {\small $\bbb$};
				\node at (0.2,-0.2) {\small $\ddd$};

				\node at (-0.7,-0.2) {\small $\aaa$};
				\node at (-1,-0.6) {\small $\ccc$};
				\node at (-0.6,-0.95) {\small $\eee$};
				\node at (-0.2,-0.7) {\small $\ddd$};
				\node at (-0.2,-0.2) {\small $\beta$};

			\end{scope}

		\end{tikzpicture}
		\caption{Different Adjacent Angle Deductions of $\bbb^2\ddd^2$.} \label{anglecombo}
	\end{figure}

	The pictures of Figure \ref{anglecombo} have the same vertex $\thin\bbb\thin\bbb\thin\ddd\thick\ddd\thin$, but different arrangements of the four tiles. To indicate the difference, we write ${}^{\lambda}\theta^{\mu}$ to mean $\lambda,\mu$ are the two angles adjacent to $\theta$ in the quadrilateral. The first picture has the AAD $\thin^{\aaa}\bbb^{\ddd}\thin^{\ddd}\bbb^{\aaa}\thin^{\bbb}\ddd^{\eee}\thick^{\eee}\ddd^{\bbb}\thin$. The second and third have the AAD $\thin^{\aaa}\bbb^{\ddd}\thin^{\aaa}\bbb^{\ddd}\thin^{\bbb}\ddd^{\eee}\thick^{\eee}\ddd^{\bbb}\thin$ and $\thin^{\ddd}\bbb^{\aaa}\thin^{\aaa}\bbb^{\ddd}\thin^{\bbb}\ddd^{\eee}\thick^{\eee}\ddd^{\bbb}\thin$ respectively. The following useful lemma is from \cite[Lemma 10]{wy1}.

	\begin{lemma}[{\cite[Lemma 8]{wy1}}]\label{hdeg}
		Suppose an angle $\theta$ does not appear at degree $3$ vertices in a tiling of the sphere by pentagons with the same angle combination. Then one of $\{\theta^{4}, \theta^{3}\rho, \theta^{5}\}$ is a vertex, where $\rho\ne\theta$.
	\end{lemma}

	\begin{lemma}[{\cite[Lemma 2.4]{lw1}}] \label{geometry2}
		In an $a^4b$-tiling, the following inequalities hold.
		\begin{enumerate}
			\item If $\delta,\epsilon<1$, then $\alpha+2\beta>1$ and $\alpha+2\gamma>1$.
			\item If $\beta,\delta<1$, then $2\alpha+\gamma>1$ and $\gamma+2\epsilon>1$.
			\item If $\gamma,\epsilon<1$, then $2\alpha+\beta>1$ and $\beta+2\delta>1$.
		\end{enumerate}
	\end{lemma}

	\begin{lemma}[{\cite[Lemma 2]{wy2}}] \label{geometry1}
		In an $a^4b$-tiling, $\beta>\gamma$ is equivalent to $\delta<\epsilon$, and $\beta<\gamma$ is equivalent to $\delta>\epsilon$.
	\end{lemma}

	\begin{lemma}[{\cite[Lemma 18]{clm2}}]\label{formula}
		If the angles $\alpha,\beta,\gamma,\delta,\epsilon$ and the edge length $a$ of an almost equilateral pentagon satisfy $a\not\in {\bb Z}\pi$, then
		\begin{align}
			&
			[
			((1-\cos\beta)\sin(\delta-\tfrac{1}{2}\alpha)
			-(1-\cos\gamma)\sin(\epsilon-\tfrac{1}{2}\alpha))\sin\tfrac{1}{2}(\delta-\epsilon) \nonumber \\
			&-(1-\cos(\beta-\gamma))\sin\tfrac{1}{2}\alpha\sin\tfrac{1}{2}(\delta+\epsilon)]\cos\tfrac{1}{2}(\delta+\epsilon-\alpha)=0, \label{coolsaet_eq1} \\
			& \sin\tfrac{1}{2}\alpha\sin\tfrac{1}{2}(\beta-\gamma)
			(\sin\tfrac{1}{2}\beta\sin\delta\cos a
			-\cos\tfrac{1}{2}\beta\cos\delta) \nonumber \\
			& +\sin\tfrac{1}{2}\gamma
			\sin\tfrac{1}{2}(\delta-\epsilon)
			\cos\tfrac{1}{2}(\delta+\epsilon-\alpha)=0, \label{coolsaet_eq2} \\
			& \sin\tfrac{1}{2}\alpha\sin\tfrac{1}{2}(\beta-\gamma)
			(\sin\tfrac{1}{2}\gamma\sin\epsilon\cos a
			-\cos\tfrac{1}{2}\gamma\cos\epsilon) \nonumber \\
			& +\sin\tfrac{1}{2}\beta
			\sin\tfrac{1}{2}(\delta-\epsilon)
			\cos\tfrac{1}{2}(\delta+\epsilon-\alpha)=0. \label{coolsaet_eq3}
		\end{align}
		Conversely, if $\alpha,\beta,\gamma,\delta,\epsilon,a$ satisfy the three equalities, and $\alpha,\beta-\gamma\not\in 2{\bb Z}\pi$, and one of $\delta,\epsilon\not\in {\bb Z}\pi$, then there is an almost equilateral pentagon with the given $\alpha,\beta,\gamma,\delta,\epsilon,a$.
	\end{lemma}

	\begin{lemma}\label{cosB}
	Given $\aaa$, $\bbb$, $\ccc$ and $a$, the following formula solves $b$:
	\begin{eqnarray}
		\cos b &=& B_4 \cos^4 a +B_3 \cos^3 a+B_2 \cos^2 a+B_1 \cos a+B_0 \label{cosb}
	\end{eqnarray}
	\begin{eqnarray*}
		B_4 &=& (1-\cos\alpha)(1-\cos\beta)(1-\cos\gamma) \\
		B_3 &=& -B_1= \sin\alpha \left( \sin(\beta+\gamma) -\sin\beta -\sin\gamma \right) \\
		B_2 &=& \cos \alpha \left( \cos (\beta -\gamma) -\cos\beta  -\cos \gamma+1 \right) -2 \cos \beta \cos \gamma +\cos \beta +\cos\gamma \\
		B_0 &=& \cos\beta \cos \gamma - \cos\alpha \sin\beta \sin\gamma .
	\end{eqnarray*}
	\end{lemma}
    The proof is similar to {\cite[Lemma 11]{wy1}}.

    \vspace{6pt}

	A vertex $\aaa^{n_1}\bbb^{n_2}\ccc^{n_3}\ddd^{n_4}\eee^{n_5}$ can be effectively represented by its vector type $\boldsymbol{n}=(n_1, n_2, n_3, n_4, n_5)\in \mathbb{N}^5$. For convenience, both representations will be alternately used later.

	\begin{lemma}[Irrational Angle Lemma]\label{irrational}
		Given three different vertices $\boldsymbol{k},\boldsymbol{l},\boldsymbol{m}\in \mathbb{N}^5$ in an $a^4b$-tiling with any irrational angle and $f$ tiles, if $\boldsymbol{k},\boldsymbol{l},\boldsymbol{m}$ and $\boldsymbol{u}:=(1\, 1 \, 1 \, 1 \, 1)$ are linearly independent, then any other vertex $\boldsymbol{n}$ belongs to
		\[
		\mathbb{N}^5\cap\mathrm{3\,\,dimensional\,\,affine\,\,subspace\,\,containing}\left\{\boldsymbol{k},\boldsymbol{l},\boldsymbol{m	},\frac{2f}{3f+4}\boldsymbol{u}\right\}.
		\]
		In particular, the square matrix of $\boldsymbol{k},\boldsymbol{l},\boldsymbol{m}, \boldsymbol{u},\bn$ has zero determinant.
	\end{lemma}
	\begin{proof}
		By Lemma \ref{anglesum} and the vertices $\boldsymbol{k},\boldsymbol{l},\boldsymbol{m}, \boldsymbol{n}$, the angles satisfy a linear system of equations with the augmented matrix
		\[C=\left(\begin{matrix}
			1 & 1 & 1 & 1 & 1 & 3+\tfrac4f \\
			k_1 & k_2 & k_3 & k_4 & k_5 & 2\\
			l_1 & l_2 & l_3 & l_4 & l_5 & 2\\
			m_1 & m_2 & m_3 & m_4 & m_5 & 2\\
			n_1 & n_2 & n_3 & n_4 & n_5 & 2
		\end{matrix}\right).
		\]
		If $\bn$ is not a linear combination of $\boldsymbol{k},\boldsymbol{l},\boldsymbol{m}, \boldsymbol{u}$, then the system has a unique rational solution, a contradiction. Therefore, the last row of $C$ must be some linear combination of the other $4$ rows, which implies that $\boldsymbol{n}$ lies in a three-dimensional affine subspace containing $4$ points $\boldsymbol{k},\boldsymbol{l},\boldsymbol{m},\frac{2f}{3f+4}\boldsymbol{u}$.
	\end{proof}

	 Since there is no $b^2$-angle, any degree $3$ vertex is either an $a^3$-vertex: \[\{\aaa^{3}, \aaa^{2}\bbb, \aaa^{2}\ccc, \aaa\bbb\ccc, \aaa\bbb^{2}, \aaa\ccc^{2},\bbb^{3}, \bbb^{2}\ccc, \bbb\ccc^{2}, \ccc^{3}\}\] or an $a^2b$-vertex: \[\{\aaa\ddd^2, \aaa\ddd\eee, \aaa\eee^2, \bbb\ddd^2, \bbb\ddd\eee,\bbb\eee^2, \ccc\ddd^2, \ccc\ddd\eee, \ccc\eee^2\}.\]

	\begin{lemma}\label{three degree $3$ vertices with one $b$-edge}
		A non-symmetric $a^4b$-tiling with general angles can have at most two different types of $a^2b$-vertex.
	\end{lemma}

	\begin{proof}
		For $a^2b$-vertices with no $\aaa$, Lemma 5 in \cite{wy2} implies that there are at most $2$ different types of such vertices, and that the only possible combinations of two such vertices are $\{\bbb\ddd^2,\ccc\eee^2\}$, $\{\bbb\ddd\eee,\ccc\eee^2\}$, or $\{\bbb\ddd^2,\ccc\ddd\eee\}$. Up to the symmetry of exchanging $\bbb\leftrightarrow\ccc$, $\ddd\leftrightarrow\eee$, it is enough to consider $\{\bbb\ddd^2,\ccc\eee^2\}$ or $\{\bbb\ddd\eee,\ccc\eee^2\}$.

		If an $a^2b$-vertex has $\aaa$, then it is $\aaa\ddd^2$, $\aaa\ddd\eee$ or $\aaa\eee^2$. These three vertices cannot appear simultaneously (otherwise $\ddd=\eee$ and the pentagon is symmetric). We only need to show that the following six combinations of three $a^2b$-vertices are impossible:
		\[
		1.\, \aaa\ddd\eee,\bbb\ddd^2,\ccc\eee^2\quad\quad\quad\quad 2.\,\aaa\ddd\eee,\bbb\ddd\eee,\ccc\eee^2 \quad\quad\quad\quad 3.\,\aaa\ddd^2,\bbb\ddd^2,\ccc\eee^2 \quad\quad\quad\quad
		\]
		\[
		4.\, \aaa\ddd^2,\bbb\ddd\eee,\ccc\eee^2\quad\quad\quad\quad 5.\,\aaa\eee^2,\bbb\ddd^2,\ccc\eee^2 \quad\quad\quad\quad 6.\,\aaa\eee^2,\bbb\ddd\eee,\ccc\eee^2 \quad\quad\quad\quad
		\]

		\begin{table*}[htp]
			\centering
			\caption{All possible cases with three $a^2b$-vertex types.}\label{1}
			~\\
			\resizebox{\textwidth}{12mm}{\begin{tabular}{|c|c|c|c|c|}

					\hline
					Vertex& Angles & Irrational Angle Lemma & AAD & Contradiction\\
					\hline\hline
					$\aaa\ddd\eee,\bbb\ddd^2,\ccc\eee^2$&$(\frac12+\frac2f,2-2\ddd,-1+2\ddd+\frac4f,\ddd,\frac32-\ddd-\frac2f)$
					&$2n_2+n_5=2n_3+n_4$& $\thin^{\ddd}\bbb^{\aaa}\thin^{\bbb}\ddd^{\eee}\thick^{\eee}\ddd^{\bbb}\thin$ &$\aaa\bbb\cdots=\aaa\bbb\ccc\Rightarrow \ddd=\eee$\\
					\hline
					$\aaa\ddd\eee,\bbb\ddd\eee,\ccc\eee^2$&$(1-\frac23\ddd+\frac4{3f},1-\frac23\ddd+\frac4{3f},\frac23\ddd+\frac8{3f},\ddd,1-\frac13\ddd-\frac4{3f})$
					&$2n_1+2n_2+n_5=2n_3+3n_4$& \multirow{2}{*}{$\thin^{\eee}\ccc^{\aaa}\thin^{\ccc}\eee^{\ddd}\thick^{\ddd}\eee^{\ccc}\thin$} &No $\aaa\ccc\cdots$ by the\\
					\cline{1-3}
					$\aaa\ddd^2,\bbb\ddd^2,\ccc\eee^2$&$(\frac23\eee+\frac8{3f},\frac23\eee+\frac8{3f},2-2\eee,1-\frac13\eee-\frac4{3f},\eee)$
					&$2n_1+2n_2+3n_5=6n_3+n_4$& & irrational angle lemma\\
					\hline
					$\aaa\ddd^2,\bbb\ddd\eee,\ccc\eee^2$&$(2-2\ddd,\frac12+\frac2f,-1+2\ddd+\frac4f,\ddd,\frac32-\ddd-\frac2f)$
					&$2n_1+n_5=2n_3+n_4$& $\thin^{\ccc}\aaa^{\bbb}\thin^{\bbb}\ddd^{\eee}\thick^{\eee}\ddd^{\bbb}\thin$ &No $\bbb^2\cdots$\\
					\hline
					$\aaa\eee^2,\bbb\ddd^2,\ccc\eee^2$&$(\frac23\ddd+\frac8{3f},2-2\ddd,\frac23\ddd+\frac8{3f},\ddd,1-\frac13\ddd-\frac4{3f})$
					&$2n_1+2n_3+3n_4=6n_2+n_5$& $\thin^{\bbb}\aaa^{\ccc}\thin^{\ccc}\eee^{\ddd}\thick^{\ddd}\eee^{\ccc}\thin$ & No $\ccc^2\cdots$\\
					\hline
					$\aaa\eee^2,\bbb\ddd\eee,\ccc\eee^2$&$(\frac12+\frac2f,\frac54-\ddd+\frac1f,\frac12+\frac2f,\ddd,\frac34-\frac1f)$
					&$n_2=n_4$& $\thin^{\eee}\ccc^{\aaa}\thin^{\ccc}\eee^{\ddd}\thick^{\ddd}\eee^{\ccc}\thin$ & No $\ddd^2\cdots$\\
					\hline
			\end{tabular}}
		\end{table*}

		\subsubsection*{Case 1. $\{\aaa\ddd\eee,\bbb\ddd^2,\ccc\eee^2\}$}
		We will explain the first row of Table \ref{1}. The AAD of $\thin^{\ddd}\bbb^{\aaa}\thin^{\bbb}\ddd^{\eee}\thick^{\eee}\ddd^{\bbb}\thin$ gives $\aaa\bbb\cdots$, denoted by $\bn$ with $n_1\ge1,n_2\ge1$. By the irrational angle lemma, we have $2n_2+n_5=2n_3+n_4 \ge 2$.  Let $R(\aaa\bbb\cdots)$ denote the remainder or $``\cdots"$ part of the vertex $\aaa\bbb\cdots$. By the angle expressions in Table \ref{1}, we get $R(\aaa\bbb\cdots)<2\ddd,\ddd+\eee,2\eee,2\ccc$. By the parity lemma, we must have $n_4=n_5=0$. By the parity lemma, we have $n_3=1$, i.e. $\aaa\bbb\cdots=\aaa\bbb\ccc$. This new vertex implies $f=12$. By the angle equation \eqref{coolsaet_eq1} in Lemma \ref{formula}, we get $\ddd=\eee=\frac23$, a contradiction.

		The other five cases also lead to contradictions similarly following Table \ref{1}.
	\end{proof}

	In fact, the above proof also induces the following two lemmas.

	\begin{lemma}\label{two degree $3$ vertices with one $b$-edge}
		In a non-symmetric $a^4b$-tiling, if there are only two $a^2b$-vertex types, then it must be one of the following (up to the exchange $\bbb\leftrightarrow\ccc$, $\ddd\leftrightarrow\eee$):
		\[
		1.\,\bbb\ddd^2,\ccc\eee^2 \quad\,\,\, 2.\,\bbb\ddd\eee,\ccc\eee^2  \quad\,\,\, 3.\, \aaa\ddd^2,\bbb\ddd^2 \quad\,\,\, 4.\,\aaa\ddd^2,\bbb\ddd\eee \quad\,\,\, 5.\,\aaa\ddd^2,\bbb\eee^2 \quad\,\,\, 6.\, \aaa\ddd^2,\ccc\ddd^2
		\]
		\[
		7.\, \aaa\ddd^2,\ccc\ddd\eee \quad\,\,\,  8.\, \aaa\ddd^2,\ccc\eee^2  \quad\,\,\, 9.\,\aaa\ddd\eee,\bbb\ddd^2 \quad\,\, 10.\,\aaa\ddd\eee,\bbb\eee^2 \quad 11.\,\aaa\ddd\eee,\bbb\ddd\eee \quad\quad\quad\quad\quad\,\,\,\,
		\]
	\end{lemma}

	\begin{lemma}\label{one degree $3$ vertices with one $b$-edge}
		In a non-symmetric $a^4b$-tiling, if there is only one $a^2b$-vertex type, then it must be one of the following (up to the exchange $\bbb\leftrightarrow\ccc$, $\ddd\leftrightarrow\eee$):
		\[
		1.\,\aaa\ddd^2 \quad\quad 2.\, \bbb\ddd^2\quad\quad 3.\, \ccc\ddd^2\quad\quad 4.\,\aaa\ddd\eee \quad\quad 5.\,\bbb\ddd\eee
		\]
	\end{lemma}

	\section{Tilings with two $a^2b$-vertex types}
	\label{two $a^2b$-vertices}

	To facilitate discussion, we denote by $T_i$ the tile labeled $i$, by $\theta_i$ the angle $\theta$ in $T_i$. When we say a tile is {\em determined}, we mean that we know all the edges and angles of the tile.

	For the eleven pairs of $a^2b$-vertices in Lemma \ref{two degree $3$ vertices with one $b$-edge}, this section will show that only $\{\aaa\ddd\eee,\bbb\ddd\eee\}$ admits tilings.

	\begin{proposition}\label{ad2,ce2}
		There is no non-symmetric $a^4b$-tiling with general angles, such that both members of any one of the first ten pairs in Lemma \ref{two degree $3$ vertices with one $b$-edge} are vertices.
	\end{proposition}

	\begin{proof}
		For each pair of $a^2b$-vertices, we will first apply Lemma \ref{hdeg} to get a new vertex, then apply the irrational angle lemma to control the AVC, and finally try to get a contradiction by the AAD. We summarize such computations and deductions in Table \ref{1.1} and \ref{1.2}.

		For example, when $\{\aaa\ddd^{2},\ccc\eee^2\}$ are vertices, if $\bbb$ does not appear in degree 3 vertices, by Lemma \ref{hdeg}, one of $\{\aaa\bbb^{3}, \bbb^{3}\ccc, \bbb^{4}, \bbb^{5}\}$ must appear. If $\bbb$ appears in a degree 3 vertex, by Lemma \ref{three degree $3$ vertices with one $b$-edge}, this vertex must be an $a^3$-vertex i.e. one of $\{\aaa^{2}\bbb, \aaa\bbb^{2}, \aaa\bbb\ccc, \bbb^{3}, \bbb^{2}\ccc, \bbb\ccc^{2}\}$. Applying the irrational angle lemma to each case, we get the first ten rows of Table \ref{1.2} for the pair $\{\aaa\ddd^{2},\ccc\eee^2\}$. We show the case \{$\aaa\ddd^{2},\ccc\eee^{2},\bbb\ccc^{2}$\} in detail.

		By the irrational angle lemma, we have $6n_1+2n_3=4n_2+3n_4+n_5$. The angles are
		\[\aaa=2-2\ddd,\,\bbb=-\frac23+\frac{4\ddd}{3}+\frac{16}{3f},\,\ccc=\frac43-\frac{2\ddd}{3}-\frac{8}{3f},\,\eee=\frac13+\frac{\ddd}{3}+\frac{4}{3f}\]

		The AAD of $\thin^{\ccc}\aaa^{\bbb}\thin^{\bbb}\ddd^{\eee}\thick^{\eee}\ddd^{\bbb}\thin$ gives $\bbb^2\cdots$, denoted by $\bn$ with $n_2\ge2$. By $R(\bbb^2\cdots)<\aaa+\ccc,4\ccc$ and  $6n_1+2n_3=4n_2+3n_4+n_5\ge8$, we get $n_1\ge2, n_3=0$. Then $4n_2+3n_4+n_5=6n_1\ge12$. So we get $(n_2-2)+n_4+n_5\ge1$, contradicting
		\[
		2\aaa+3\bbb>2,2\aaa+2\bbb+\ddd>2 \,\,\,\text{and}\,\,\,  2\aaa+2\bbb+\eee>2.
		\]
		Therefore, Case \{$\aaa\ddd^{2},\ccc\eee^{2},\bbb\ccc^{2}$\} admits no tiling.

		\vspace{6pt}

		When the new vertex is $\aaa\bbb\ccc$, $\bbb^{3}$, $\bbb^{4}$ or $\bbb^{5}$, we actually deduce a combinatorial tiling of the sphere. However, there is no solution to \eqref{coolsaet_eq1}, \eqref{coolsaet_eq2}, \eqref{coolsaet_eq3} in Lemma \ref{formula} (see also \cite[Proposition 20]{clm2}), i.e. there are no such pentagons geometrically.

		All other pairs are discussed similarly by Table \ref{1.1} and \ref{1.2}.
	\end{proof}

		\begin{table*}[htp]
			\centering
			\caption{All possible cases with two $a^2b$-vertex types, part one.}\label{1.1}
			~\\
			\resizebox{\linewidth}{!}{\begin{tabular}{|c|c|c|c|c|c|}
					\hline
					\multicolumn{2}{|c|}{Vertex} & Angles & Irrational Angle Lemma & AAD & Contradiction \\
					\hline\hline
					\multirow{10}{*}{$\bbb\ddd^2,\ccc\eee^{2}$}&$\aaa^{2}\bbb$&$(\ddd,2-2\ddd,\frac8f,\ddd,1-\frac4f)$
					&$n_1+n_4=2n_2$&\multirow{3}{*}{ $\thin^{\eee}\ccc^{\aaa}\thin^{\ccc}\eee^{\ddd}\thick^{\ddd}\eee^{\ccc}\thin$}&No $\aaa\ccc\cdots$ by \\
					\cline{2-4}
					&$\aaa\bbb^{2}$&$(-2+4\ddd,2-2\ddd,4-6\ddd+\frac8f,\ddd,-1+3\ddd-\frac4f)$
					&$4n_1+n_4+3n_5=2n_2+6n_3$&&the irrational \\
					\cline{2-4}
					&$\aaa^{3}\bbb$&$(\frac{2}3\ddd,2-2\ddd,\frac{2}3\ddd+\frac8f,\ddd,1-\frac{1}3\ddd-\frac4f)$
					&$2n_1+2n_3+3n_4=6n_2+n_5$& &angle lemma\\
					\cline{2-6}
					&$\aaa^{2}\ccc$&$(\eee,\frac8f,2-2\eee,1-\frac4f,\eee)$
					&$n_1+n_5=2n_3$& \multirow{3}{*}{ $\thin^{\ddd}\bbb^{\aaa}\thin^{\bbb}\ddd^{\eee}\thick^{\eee}\ddd^{\bbb}\thin$} &\multirow{3}{*}{No $\aaa\bbb\cdots$}\\
					\cline{2-4}
					&$\aaa\ccc^{2}$&$(-\frac23+\frac{4}3\ddd+\frac{16}{3f},2-2\ddd,\frac43-\frac{2}3\ddd-\frac8{3f},\ddd,\frac13+\frac{1}3\ddd+\frac4{3f})$
					&$4n_1+3n_4+n_5=6n_2+2n_3$& &\\
					\cline{2-4}
					&$\aaa^{3}\ccc$&$(2-2\ddd-\frac8f,2-2\ddd,-4+6\ddd+\frac{24}f,\ddd,3-3\ddd-\frac{12}f)$
					&$2n_1+2n_2+3n_5=6n_3+n_4$& &\\
					\cline{2-6}
					&$\aaa\bbb\ccc$&$(\frac8f,2-2\ddd,2\ddd-\frac8f,\ddd,1-\ddd+\frac4f)$
					&\multirow{4}{*}{$2n_2+n_5=2n_3+n_4$}& &No solution to\\
					\cline{2-3}
					&$\aaa^{3}$&$(\frac23,2-2\ddd,-\frac43+2\ddd+\frac8f,\ddd,\frac53-\ddd-\frac4f)$
					& && \eqref{coolsaet_eq1}, \eqref{coolsaet_eq2}, \eqref{coolsaet_eq3}\\
					\cline{2-3}
					&$\aaa^{4}$&$(\frac12,2-2\ddd,-1+2\ddd+\frac8f,\ddd,\frac32-\ddd-\frac4f)$
					&& & in Lemma \ref{formula}\\
					\cline{2-3}
					&$\aaa^{5}$&$(\frac25,2-2\ddd,-\frac45+2\ddd+\frac8f,\ddd,\frac75-\ddd-\frac4f)$
					&& &\\
					\hline
					\multirow{10}{*}{$\bbb\ddd\eee,\ccc\eee^{2}$}&$\aaa\bbb^{2}$&$(-1+2\eee+\frac4f,\frac32-\eee-\frac2f,2-2\eee,\frac12+\frac2f,\eee)$
					&$2n_1+n_5=n_2+2n_3$&\multirow{10}{*}{ $\thin^{\eee}\ccc^{\aaa}\thin^{\ccc}\eee^{\ddd}\thick^{\ddd}\eee^{\ccc}\thin$}&No $\aaa\ccc\cdots$\\
					\cline{2-4}
					\cline{6-6}
					&$\aaa^{2}\bbb$& $(\frac13+\frac{2\ddd}{3}-\frac4{3f},\frac43-\frac{4}3\ddd+\frac8{3f},\frac23-\frac{2}3\ddd+\frac{16}{3f},\ddd,\frac23+\frac{\ddd}3-\frac8{3f})$
					&$2n_1+3n_4+n_5=4n_2+2n_3$& &\\
					\cline{2-4}
					&$\aaa\bbb\ccc$&$(1-2\ddd+\frac{12}f,1-\frac4f,2\ddd-\frac8f,\ddd,1-\ddd+\frac4f)$
					& $2n_1+n_5=2n_3+n_4$&& \\
					\cline{2-4}
					&$\aaa^{3}\bbb$&$(\frac15+\frac{2}5\ddd-\frac4{5f},\frac75-\frac{6}5\ddd+\frac{12}{5f},\frac45-\frac{2}5\ddd+\frac{24}{5f},\ddd,\frac35+\frac{1}5\ddd-\frac{12}{5f})$
					&$2n_1+5n_4+n_5=6n_2+2n_3$& &\\
					\cline{2-4}
					&$\aaa^{2}\ccc$&$(1-\frac4f,1-\ddd+\frac4f,\frac8f,\ddd,1-\frac4f)$
					&\multirow{6}{*}{$n_2=n_4$}& &\multirow{4}{*}{No $\ddd^{2}\cdots$}\\
					\cline{2-3}
					&$\aaa\ccc^{2}$& $(\frac8f,\frac32-\ddd-\frac2f,1-\frac4f,\ddd,\frac12+\frac2f)$
					&& &\\
					\cline{2-3}
					&$\aaa^{3}\ccc$&$(\frac12-\frac2f,\frac54-\ddd+\frac3f,\frac12+\frac6f,\ddd,\frac34-\frac3f)$
					&& &\\
					\cline{2-3}
					&$\aaa^{3}$&$(\frac23,\frac76-\ddd+\frac2f,\frac13+\frac4f,\ddd,\frac56-\frac2f)$
					&& &\\
					\cline{2-3}
					&$\aaa^{4}$&$(\frac12,\frac54-\ddd+\frac2f,\frac12+\frac4f,\ddd,\frac34-\frac2f)$
					&& &\\
					\cline{2-3}
					&$\aaa^{5}$&$(\frac25,\frac{13}{10}-\ddd+\frac2f,\frac35+\frac4f,\ddd,\frac{7}{10}-\frac2f)$
					&& &\\
					\hline
					\multirow{10}{*}{$\aaa\ddd^2,\bbb\ddd^{2}$}&$\aaa^2\ccc$&\multirow{3}{*}{$(2-2\ddd,2-2\ddd,-2+4\ddd,\ddd,1-\ddd+\frac4f)$}
					&\multirow{3}{*}{$2n_1+2n_2+n_5=4n_3+n_4$}&\multirow{20}{*}{ $\thin^{\ccc}\aaa^{\bbb}\thin^{\bbb}\ddd^{\eee}\thick^{\eee}\ddd^{\bbb}\thin$} &\multirow{20}{*}{No $\eee^2\cdots$}\\
					\cline{2-2}
					&$\aaa\bbb\ccc$&&& &\\
					\cline{2-2}
					&$\bbb^{2}\ccc$&&&&\\
					\cline{2-4}
					&$\aaa\ccc^{2}$&\multirow{2}{*}{$(2-2\ddd,2-2\ddd,\ddd,\ddd,-1+2\ddd+\frac4f)$}
					&\multirow{2}{*}{$2n_1+2n_2=n_3+n_4+2n_5$}& &\\
					\cline{2-2}
					&$\bbb\ccc^{2}$&&& &\\
					\cline{2-4}
					&$\aaa\ccc^{3}$&\multirow{2}{*}{$(2-2\ddd,2-2\ddd,\frac{2\ddd}3,\ddd,-1+\frac{7\ddd}3+\frac4f)$}
					&\multirow{2}{*}{$6n_1+6n_2=2n_3+3n_4+7n_5$}& & \\
					\cline{2-2}
					&$\bbb\ccc^{3}$&
					&& & \\
					\cline{2-4}
					&$\ccc^{3}$&$(2-2\ddd,2-2\ddd,\frac23,\ddd,-\frac53+3\ddd+\frac4f)$
					&\multirow{3}{*}{$2n_1+2n_2=n_4+3n_5$}& & \\
					\cline{2-3}
					&$\ccc^{4}$&$(2-2\ddd,2-2\ddd,\frac12,\ddd,-\frac32+3\ddd+\frac4f)$
					&& & \\
					\cline{2-3}
					&$\ccc^{5}$&$(2-2\ddd,2-2\ddd,\frac25,\ddd,-\frac75+3\ddd+\frac4f)$
					&& & \\
					\cline{1-4}
					\multirow{10}{*}{$\aaa\ddd^2,\bbb\ddd\eee$}&$\aaa\bbb\ccc$&$(2-2\ddd,1-\frac4f,-1+2\ddd+\frac4f,\ddd,1-\ddd+\frac4f)$
					&$2n_1+n_5=2n_3+n_4$&&\\
					\cline{2-4}
					&$\bbb^{2}\ccc$&$(2-2\ddd,\frac32-\ddd-\frac2f,-1+2\ddd+\frac4f,\ddd,\frac12+\frac2f)$
					&$2n_1+n_2=2n_3+n_4$& &\\
					\cline{2-4}
					&$\bbb\ccc^{2}$&$(2-2\ddd,4-4\ddd-\frac8f,-1+2\ddd+\frac4f,\ddd,-2+3\ddd+\frac8f)$
					&$2n_1+4n_2=2n_3+n_4+3n_5$& &\\
					\cline{2-4}
					&$\bbb\ccc^{3}$&$(2-2\ddd,5-6\ddd-\frac{12}f,-1+2\ddd+\frac4f,\ddd,-3+5\ddd+\frac{12}f)$
					&$2n_1+6n_2=2n_3+n_4+5n_5$& &\\
					\cline{2-4}
					&$\aaa^2\ccc$&$(1-\frac4f,\frac32-\eee-\frac2f,\frac8f,\frac12+\frac2f,\eee)$&\multirow{7}{*}{$n_2=n_5$}& &\\
					\cline{2-3}
					&$\aaa\ccc^{2}$&$(\frac8f,1-\eee+\frac4f,1-\frac4f,1-\frac4f,\eee)$&& &\\
					\cline{2-3}
					&$\aaa\ccc^{3}$&$(\frac12+\frac6f,\frac54-\eee+\frac3f,\frac12-\frac2f,\frac34-\frac3f,\eee)$&& &\\
					\cline{2-3}
					&$\ccc^{3}$&$(\frac13+\frac4f,\frac76-\eee+\frac2f,\frac23,\frac56-\frac2f,\eee)$&& &\\
					\cline{2-3}
					&$\ccc^{4}$&$(\frac12+\frac4f,\frac54-\eee+\frac2f,\frac12,\frac34-\frac2f,\eee)$&& &\\
					\cline{2-3}
					&$\ccc^{5}$&$(\frac35+\frac4f,\frac{13}{10}-\eee+\frac2f,\frac25,\frac{7}{10}-\frac2f,\eee)$&&&\\
					\hline
					\multirow{10}{*}{$\aaa\ddd^2,\bbb\eee^2$}&$\aaa\ccc^{2}$&$(2-2\ddd,\frac8f,\ddd,\ddd,1-\frac4f)$
					&$2n_1=n_3+n_4$& \multirow{2}{*}{$\thin^{\aaa}\bbb^{\ddd}\thin^{\ccc}\eee^{\ddd}\thick^{\ddd}\eee^{\ccc}\thin$} & \multirow{2}{*}{No $\ccc\ddd\cdots$}\\
					\cline{2-4}
					&$\bbb^{2}\ccc$&$(2-2\ddd,\frac43-\frac{2}3\ddd-\frac8{3f},-\frac23+\frac{4}3\ddd+\frac{16}{3f},\ddd,\frac13+\frac{1}3\ddd+\frac4{3f})$
					&$6n_1+2n_2=4n_3+3n_4+n_5$& &\\
					\cline{2-6}
					&$\ccc^{3}$&$(2-2\ddd,-\frac43+2\ddd+\frac8f,\frac23,\ddd,\frac53-\ddd-\frac4f)$
					&\multirow{2}{*}{$2n_1+n_5=2n_2+n_4$}& $\thin^{\eee}\ccc^{\aaa}\thin^{\eee}\ccc^{\aaa}\thin\cdots$& No $\aaa\eee\cdots$\\
					\cline{2-3}
					\cline{5-6}
					&$\aaa\bbb\ccc$&$(2-2\ddd,2\ddd-\frac8f,\frac8f,\ddd,1-\ddd+\frac4f)$
					&&\multirow{27}{*}{ $\thin^{\ccc}\aaa^{\bbb}\thin^{\bbb}\ddd^{\eee}\thick^{\eee}\ddd^{\bbb}\thin$}&\multirow{8}{*}{No $\bbb^2\cdots$}\\
					\cline{2-4}
					&$\aaa^2\ccc$&$(2-2\ddd,4-6\ddd+\frac8f,-2+4\ddd,\ddd,-1+3\ddd-\frac4f)$
					&$2n_1+6n_2=4n_3+n_4+3n_5$& &\\
					\cline{2-4}
					&$\bbb\ccc^{2}$&$(\frac8f,2-2\eee,\eee,1-\frac4f,\eee)$
					&$2n_2=n_3+n_5$& &\\
					\cline{2-4}
					&$\aaa\ccc^{3}$&$(2-2\ddd,\frac{2}3\ddd+\frac8f,\frac{2}3\ddd,\ddd,1-\frac{1}3\ddd-\frac4f)$
					&$6n_1+n_5=2n_2+2n_3+3n_4$& & \\
					\cline{2-4}
					&$\bbb\ccc^{3}$&$(2-2\ddd,-4+6\ddd+\frac{24}f,2-2\ddd-\frac8f,\ddd,3-3\ddd-\frac{12}f)$
					&$2n_1+2n_3+3n_5=6n_2+n_4$& &\\
					\cline{2-4}
					&$\ccc^{4}$&$(2-2\ddd,-1+2\ddd+\frac8f,\frac12,\ddd,\frac32-\ddd-\frac4f)$&\multirow{2}{*}{$2n_1+n_5=2n_2+n_4$}& &\\
					\cline{2-3}
					&$\ccc^{5}$&$(2-2\ddd,-\frac45+2\ddd+\frac8f,\frac25,\ddd,\frac75-\ddd-\frac4f)$
					&&  & \\
					\cline{1-4}
					\cline{6-6}
					\multirow{10}{*}{$\aaa\ddd^2,\ccc\ddd^2$}&$\aaa^{2}\bbb$&\multirow{3}{*}{$(2-2\ddd,-2+4\ddd,2-2\ddd,\ddd,1-\ddd+\frac4f)$}
					&\multirow{3}{*}{$2n_1+2n_3+n_5=4n_2+n_4$}&&\multirow{10}{*}{No $\eee^2\cdots$}\\
					\cline{2-2}
					&$\aaa\bbb\ccc$&&& & \\
					\cline{2-2}
					&$\bbb\ccc^{2}$&&& & \\
					\cline{2-4}
					&$\aaa\bbb^{2}$&\multirow{2}{*}{$(2-2\ddd,\ddd,2-2\ddd,\ddd,-1+2\ddd+\frac4f)$}
					&\multirow{2}{*}{$2n_1+2n_3=n_2+n_4+2n_5$}& &\\
					\cline{2-2}
					&$\bbb^{2}\ccc$&&& & \\
					\cline{2-4}
					&$\aaa\bbb^{3}$&\multirow{2}{*}{$(2-2\ddd,\frac{2}3\ddd,2-2\ddd,\ddd,-1+\frac{7}3\ddd+\frac4f)$}
					&\multirow{2}{*}{$6n_1+6n_3=2n_2+3n_4+7n_5$}& & \\
					\cline{2-2}
					&$\bbb^{3}\ccc$&&& & \\
					\cline{2-4}
					&$\bbb^{3}$&$(2-2\ddd,\frac23,2-2\ddd,\ddd,-\frac53+3\ddd+\frac4f)$
					&\multirow{3}{*}{$2n_1+2n_3=n_4+3n_5$}& & \\
					\cline{2-3}
					&$\bbb^{4}$&$(2-2\ddd,\frac12,2-2\ddd,\ddd,-\frac32+3\ddd+\frac4f)$
					&& & \\
					\cline{2-3}
					&$\bbb^{5}$&$(2-2\ddd,\frac25,2-2\ddd,\ddd,-\frac75+3\ddd+\frac4f)$
					&& & \\
					\cline{1-4}
					\cline{6-6}
					\multirow{10}{*}{$\aaa\ddd^2,\ccc\ddd\eee$}&$\aaa\bbb\ccc$&$(2-2\ddd,-1+2\ddd+\frac{4}f,1-\frac{4}f,\ddd,1-\ddd+\frac{4}f)$
					&$2n_1+n_5=2n_2+n_4$&&\multirow{2}{*}{No $\bbb^2\cdots$}\\
					\cline{2-4}
					&$\bbb\ccc^{2}$&$(2-2\ddd,-1+2\ddd+\frac{4}f,\frac{3}2-\ddd-\frac{2}f,\ddd,\frac{1}2+\frac{2}f)$
					&$2n_1+n_3=2n_2+n_4$& & \\
					\cline{2-4}
					\cline{6-6}
					&$\bbb^{2}\ccc$&$(2-2\ddd,-1+2\ddd+\frac4f,4-4\ddd-\frac{8}f,\ddd,-2+3\ddd+\frac{8}f)$
					&$2n_1+4n_3=2n_2+n_4+3n_5$& &\multirow{8}{*}{No $\eee^2\cdots$}\\
					\cline{2-4}
					&$\bbb^{3}\ccc$&$(2-2\ddd,-1+2\ddd+\frac4f,5-6\ddd-\frac{12}f,\ddd,-3+5\ddd+\frac{12}f)$
					&$2n_1+6n_3=2n_2+n_4+5n_5$& &\\
					\cline{2-4}
					&$\aaa^{2}\bbb$&$(1-\frac{4}f,\frac{8}f,\frac{3}2-\eee-\frac{2}f,\frac{1}2+\frac{2}f,\eee)$&\multirow{6}{*}{$n_3=n_5$}& &\\
					\cline{2-3}
					&$\aaa\bbb^{2}$&$(\frac{8}f,1-\frac{4}f,1-\eee+\frac{4}f,1-\frac{4}f,\eee)$&& &\\
					\cline{2-3}
					&$\aaa\bbb^{3}$&$(\frac12+\frac6f,\frac12-\frac2f,\frac54-\eee+\frac3f,\frac34-\frac3f,\eee)$&& &\\
					\cline{2-3}
					&$\bbb^{3}$&$(\frac{1}3+\frac{4}f,\frac{2}3,\frac{7}6-\eee+\frac{2}f,\frac{5}6-\frac{2}f,\eee)$&& &\\
					\cline{2-3}
					&$\bbb^{4}$&$(\frac12+\frac4f,\frac12,\frac54-\eee+\frac2f,\frac34-\frac2f,\eee)$&& &\\
					\cline{2-3}
					&$\bbb^{5}$&$(\frac35+\frac4f,\frac25,\frac{13}{10}-\eee+\frac2f,\frac7{10}-\frac2f,\eee)$
					&&  & \\
					\hline
			\end{tabular}}
		\end{table*}

	\begin{table*}[htp]
		\centering
		\caption{All possible cases with two $a^2b$-vertex types, part two.}\label{1.2}
		~\\
		\resizebox{\linewidth}{!}{\begin{tabular}{|c|c|c|c|c|c|}
				\hline
				\multicolumn{2}{|c|}{Vertex} & Angles & Irrational Angle Lemma & AAD & Contradiction \\
				\hline\hline
				\multirow{10}{*}{$\aaa\ddd^2,\ccc\eee^{2}$}&$\aaa^{2}\bbb$&$(2-2\ddd,-2+4\ddd,4-6\ddd+\frac{8}f,\ddd,-1+3\ddd-\frac{4}f)$
				&$2n_1+6n_3=4n_2+n_4+3n_5$& \multirow{2}{*}{$\thin^{\ccc}\aaa^{\bbb}\thin^{\bbb}\ddd^{\eee}\thick^{\eee}\ddd^{\bbb}\thin$} &\multirow{2}{*}{No $\bbb^2\cdots$}\\
				\cline{2-4}
				&$\bbb\ccc^{2}$&$(2-2\ddd,-\frac{2}3+\frac{4\ddd}3+\frac{16}{3f},\frac{4}3-\frac{2\ddd}3-\frac{8}{3f},\ddd,\frac{1}3+\frac{\ddd}3+\frac{4}{3f})$
				&$6n_1+2n_3=4n_2+3n_4+n_5$&  & \\
				\cline{2-6}
				&$\aaa\bbb^{2}$&$(2-2\ddd,\ddd,\frac{8}f,\ddd,1-\frac{4}f)$
				&$2n_1=n_2+n_4$&\multirow{4}{*}{$\thin^{\eee}\ccc^{\aaa}\thin^{\ccc}\eee^{\ddd}\thick^{\ddd}\eee^{\ccc}\thin$}  &\multirow{4}{*}{No $\aaa\ccc\cdots$} \\
				\cline{2-4}
				&$\bbb^{2}\ccc$&$(\frac{8}f,\eee,2-2\eee,1-\frac{4}f,\eee)$
				&$n_2+n_5=2n_3$& & \\
				\cline{2-4}
				&$\aaa\bbb^{3}$&$(2-2\ddd,\frac{2\ddd}3,\frac{2\ddd}3+\frac{8}f,\ddd,1-\frac{\ddd}3-\frac{4}f)$
				&$6n_1+n_5=2n_2+2n_3+3n_4$& & \\
				\cline{2-4}
				&$\bbb^{3}\ccc$&$(2-2\ddd,2-2\ddd-\frac{8}f,-4+6\ddd+\frac{24}f,\ddd,3-3\ddd-\frac{12}f)$
				&$2n_1+2n_2+3n_5=6n_3+n_4$& &\\
				\cline{2-6}
				& $\aaa\bbb\ccc$&$(2-2\ddd,\frac{8}f,2\ddd-\frac{8}f,\ddd,1-\ddd+\frac{4}f)$
				&\multirow{4}{*}{$2n_1+n_5=2n_3+n_4$}&&No solution to  \\
				\cline{2-3}
				&$\bbb^{3}$&$(2-2\ddd,\frac{2}3,-\frac{4}3+2\ddd+\frac{8}f,\ddd,\frac{5}3-\ddd-\frac{4}f)$
				&&&\eqref{coolsaet_eq1}, \eqref{coolsaet_eq2}, \eqref{coolsaet_eq3}\\
				\cline{2-3}
				&$\bbb^{4}$&$(2-2\ddd,\frac{1}2,-1+2\ddd+\frac{8}f,\ddd,\frac{3}2-\ddd-\frac{4}f)$&&&in Lemma \ref{formula}\\
				\cline{2-3}
				& $\bbb^{5}$&$(2-2\ddd,\frac{2}5,-\frac45+2\ddd+\frac8f,\ddd,\frac{7}5-\ddd-\frac{4}f)$ &&&\\
				\hline
				\multirow{10}{*}{$\aaa\ddd\eee,\bbb\ddd^2$}&$\aaa^2\ccc$&$(\frac{3}2-\ddd-\frac{2}f,2-2\ddd,-1+2\ddd+\frac{4}f,\ddd,\frac{1}2+\frac{2}f)$
				&$n_1+2n_2=2n_3+n_4$& \multirow{10}{*}{ $\thin^{\ddd}\bbb^{\aaa}\thin^{\bbb}\ddd^{\eee}\thick^{\eee}\ddd^{\bbb}\thin$}& No $\aaa\bbb\cdots$\\
				\cline{2-4}
				\cline{6-6}
				&$\aaa\bbb\ccc$&$(1-\frac{4}f,2-2\ddd,-1+2\ddd+\frac{4}f,\ddd,1-\ddd+\frac{4}f)$
				&$2n_2+n_5=2n_3+n_4$& &\multirow{9}{*}{No $\eee^2\cdots$}\\
				\cline{2-4}
				&$\aaa\ccc^{2}$&$(4-4\ddd-\frac{8}f,2-2\ddd,-1+2\ddd+\frac{4}f,\ddd,-2+3\ddd+\frac{8}f)$
				&$4n_1+2n_2=2n_3+n_4+3n_5$& & \\
				\cline{2-4}
				&$\aaa\ccc^{3}$&$(5-6\ddd-\frac{12}f,2-2\ddd,-1+2\ddd+\frac{4}f,\ddd,-3+5\ddd+\frac{12}f)$
				&$6n_1+2n_2=2n_3+n_4+5n_5$& & \\
				\cline{2-4}	&$\bbb^{2}\ccc$&$(\frac{3}2-\eee-\frac{2}f,1-\frac{4}f,\frac{8}f,\frac{1}2+\frac{2}f,\eee)$&\multirow{6}{*}{$n_1=n_5$}& &\\
				\cline{2-3}
				&$\bbb\ccc^{2}$&$(1-\eee+\frac{4}f,\frac{8}f,1-\frac{4}f,1-\frac{4}f,\eee)$
				&& &\\
				\cline{2-3}
				&$\bbb\ccc^{3}$&$(\frac{5}4-\eee+\frac{3}f,\frac{1}2+\frac{6}f,\frac{1}2-\frac{2}f,\frac{3}4-\frac{3}f,\eee)$
				&& &\\
				\cline{2-3}
				&$\ccc^{3}$&$(\frac{7}6-\eee+\frac{2}f,\frac{1}3+\frac{4}f,\frac{2}3,\frac{5}6-\frac{2}f,\eee)$
				&& &\\
				\cline{2-3}
				&$\ccc^{4}$&$(\frac{5}4-\eee+\frac{2}f,\frac{1}2+\frac{4}f,\frac{1}2,\frac{3}4-\frac{2}f,\eee)$
				&& &\\
				\cline{2-3}
				&$\ccc^{5}$&$(\frac{13}{10}-\eee+\frac{2}f,\frac{3}5+\frac{4}f,\frac{2}5,\frac{7}{10}-\frac{2}f,\eee)$
				&&  & \\
				\hline
				\multirow{10}{*}{$\aaa\ddd\eee,\bbb\eee^2$}&$\aaa^2\ccc$&$(\frac{3}2-\eee-\frac{2}f,2-2\eee,-1+2\eee+\frac{4}f,\frac{1}2+\frac{2}f,\eee)$
				&$n_1+2n_2=2n_3+n_5$&\multirow{10}{*}{ $\thin^{\aaa}\bbb^{\ddd}\thin^{\ccc}\eee^{\ddd}\thick^{\ddd}\eee^{\ccc}\thin$}& No $\ccc\ddd\cdots$\\
				\cline{2-4}
				\cline{6-6}
				&$\aaa\bbb\ccc$&$(1-\frac{4}f,2\ddd-\frac{8}f,1-2\ddd+\frac{12}f,\ddd,1-\ddd+\frac{4}f)$
				&$2n_2+n_4=2n_3+n_5$& & \multirow{9}{*}{No $\ddd^2\cdots$}\\
				\cline{2-4}
				&$\aaa\ccc^{2}$&$(\frac{4}3-\frac{4}3\ddd+\frac{8}{3f},\frac{2}3-\frac{2}3\ddd+\frac{16}{3f},\frac{1}3+\frac{2}3\ddd-\frac{4}{3f},\ddd,\frac{2}3+\frac{1}3\ddd-\frac{8}{3f})$
				&$4n_1+2n_2=2n_3+3n_4+n_5$& & \\
				\cline{2-4}
				&$\aaa\ccc^{3}$&$(\frac{7}5-\frac{6}5\ddd+\frac{12}{5f},\frac{4}{5}-\frac{2}5\ddd+\frac{24}{5f},\frac{1}5+\frac{2}5\ddd-\frac{4}{5f},\ddd,\frac{3}5+\frac{1}5\ddd-\frac{12}{5f})$
				&$6n_1+2n_2=2n_3+5n_4+n_5$& & \\
				\cline{2-4}
				&$\bbb^{2}\ccc$&$(\frac{3}2-\ddd-\frac{2}f,1-\frac{4}f,\frac{8}f,\ddd,\frac{1}2+\frac{2}f)$
				&\multirow{6}{*}{$n_1=n_4$}& &\\
				\cline{2-3}
				&$\bbb\ccc^{2}$&$(1-\ddd+\frac{4}f,\frac{8}f,1-\frac{4}f,\ddd,1-\frac{4}f)$
				&& &\\
				\cline{2-3}
				&$\bbb\ccc^{3}$&$(\frac{5}4-\ddd+\frac{3}f,\frac{1}2+\frac{6}f,\frac{1}2-\frac{2}f,\ddd,\frac{3}4-\frac{3}f)$
				&& &\\
				\cline{2-3}
				&$\ccc^{3}$&$(\frac{7}6-\ddd+\frac{2}f,\frac{1}3+\frac{4}f,\frac{2}3,\ddd,\frac{5}6-\frac{2}f)$
				&& &\\
				\cline{2-3}
				&$\ccc^{4}$&$(\frac{5}4-\ddd+\frac{2}f,\frac{1}2+\frac{4}f,\frac{1}2,\ddd,\frac{3}4-\frac{2}f)$
				&& &\\
				\cline{2-3}
				&$\ccc^{5}$&$(\frac{13}{10}-\ddd+\frac{2}f,\frac{3}5+\frac{4}f,\frac{2}5,\ddd,\frac{7}{10}-\frac{2}f)$
				&&  & \\
				\hline
		\end{tabular}}
	\end{table*}

	\begin{proposition}\label{ade,bde}
		All non-symmetric $a^4b$-tilings with general angles, such that $\aaa\ddd\eee$, $\bbb\ddd\eee$ appear as vertices, are the following:
		\begin{itemize}
			\item A unique tiling $T(8\aaa\ccc^2,8\aaa\ddd\eee,8\bbb\ddd\eee,2\bbb^4)$ (see the second picture of Figure \ref{s16}) with $f=16$, by the unique pentagon with angles $(\frac12,\frac12,\frac34,\frac32-\eee,\eee)$, $\cos(\eee\pi)=1/\sqrt[4]{2}$, and with sides $\cos a= \sqrt{\sqrt{2}-1}$, $\cos b= \sqrt{2}\sqrt{\sqrt{2}-1}$.
			\item Many different tilings (see Table \ref{table-3} and Figure \ref{various flip}) with $f=8k+4$ for each $k\ge2$, by the unique pentagon with $\aaa=\bbb=1-\frac4f$ in Table \ref{tab-2}.
		\end{itemize}
	\end{proposition}

	\begin{proof}
		If $\ccc$ does not appear at degree $3$ vertices, by Lemma \ref{hdeg}, one of $\{\aaa\ccc^{3}, \bbb\ccc^{3},  \ccc^{4},\ccc^{5}\}$ must appear. If $\ccc$ appears in some $a^3$-vertex, one of $\{\aaa^2\ccc, \aaa\bbb\ccc, \aaa\ccc^{2}, \bbb^{2}\ccc, \bbb\ccc^{2},\ccc^{3}\}$ must appear. In all cases, the irrational angle lemma implies $n_4=n_5$ and we get Table \ref{ade bde-}. We discuss each case as follows.

		\begin{table*}[htp]
			\centering
			\caption{All cases when $\aaa\ddd\eee$ and $\bbb\ddd\eee$ appear.}\label{ade+bde}
			~\\
			\resizebox{\linewidth}{!}{\begin{tabular}{|c|c|c|c|c|}
					\hline
					\multicolumn{2}{|c|}{Vertex} & Angles & Irrational Angle Lemma & Conclusion \\
					\hline\hline
					\multirow{10}{*}{$\aaa\ddd\eee,\bbb\ddd\eee$}&$\aaa^2\ccc$&\multirow{3}{*}{$(1-\frac4f,1-\frac4f,\frac8f,\ddd,1-\ddd+\frac4f)$}&\multirow{10}{*}{$n_4=n_5$}&\multirow{3}{*}{See page \pageref{ade bde 2ac}}\\
					\cline{2-2}
					&$\aaa\bbb\ccc$&&&\\
					\cline{2-2}
					&$\bbb^2\ccc$&&&\\
					\cline{2-3}
					\cline{5-5}
					&$\aaa\ccc^2$&\multirow{2}{*}{$(\frac8f,\frac8f,1-\frac4f,\ddd,2-\ddd-\frac8f)$}&&\multirow{2}{*}{See page \pageref{ade bde a2c}}\\
					\cline{2-2}
					&$\bbb\ccc^2$&&&\\
					\cline{2-3}
					\cline{5-5}
					&$\ccc^3$&$(\tfrac13+\tfrac4f,\tfrac13+\tfrac4f,\tfrac23,\ddd,\tfrac53-\ddd-\tfrac4f)$ &  & \multirow{4}{*}{See Table \ref{Tab-3.1}}\\
					\cline{2-3}
					&$\aaa\ccc^3$&\multirow{2}{*}{$(\tfrac12+\tfrac6f,\tfrac12+\tfrac6f,\tfrac12-\tfrac2f,\ddd,\tfrac32-\ddd-\tfrac6f)$}&&\\
					\cline{2-2}
					&$\bbb\ccc^3$&&&  \\
					\cline{2-3}
					&$\ccc^4$&$(\tfrac12+\tfrac4f,\tfrac12+\tfrac4f,\tfrac12,\ddd,\tfrac32-\ddd-\tfrac4f)$&& \\
					\cline{2-3}
					\cline{5-5}
					&$\ccc^5$&$(\frac35+\frac4f,\frac35+\frac4f,\frac25,\ddd,\frac75-\ddd-\frac4f)$&&See page \pageref{ade bde 5c}\\
					\hline
			\end{tabular}}\label{ade bde-}
		\end{table*}

		\subsection*{Case. $\{\aaa\ddd\eee, \bbb\ddd\eee, \ccc^5\}$}\label{ade bde 5c}
		If $\aaa^{n_1}\bbb^{n_2}\ccc^{n_3}\ddd^{n_4}\eee^{n_5}$ is a vertex, then we have
		\[
		(\tfrac{3}{5}+\tfrac{4}{f})n_1+(\tfrac{3}{5}+\tfrac{4}{f})n_2+\tfrac{2}{5}n_3+(\tfrac{7}{5}-\tfrac{4}{f})n_4=2.
		\]
		 Since $\ccc^5$ is a vertex, we get $f>12$. We also have $\aaa=\bbb > \frac{3}{5},\ccc=\frac{2}{5},\ddd+\eee>\frac{16}{15}$.  This implies $n_1 \le 3,n_2\le 3,n_3 \le 5,n_4=n_5\le 1$. We substitute the finitely many combinations of exponents satisfying the bounds into the equation above and solve for $f$. By the angle values and the parity lemma, we get all possible AVC in Table \ref{Tab-2.1}. Its first row  ``$f=\text{all}$'' means that the angle combinations can be vertices for any $f$; all other rows are mutually exclusive.  All possible tilings based on the AVC of Table \ref{Tab-2.1} are deduced as follows. Note that the $\text{AVC}\sub\{\aaa\ddd\eee,\bbb\ddd\eee,\ccc^5\}$ in the first row admits no solution satisfying the balance lemma.

		\begin{table}[htp]
			\centering
			\caption{AVC for Case $\{\aaa\ddd\eee,\bbb\ddd\eee,\ccc^5\}$.}\label{Tab-2.1}
			\begin{minipage}{0.4\textwidth}
				\begin{tabular}{c|c}
					\hline
					$f$&vertex\\
					\hline
					\hline
					all&$\aaa\ddd\eee,\bbb\ddd\eee,\ccc^5$\\
					$20$&$\aaa^2\ccc$, $\aaa\bbb\ccc$, $\bbb^2\ccc$, $\aaa\ccc^3$, $\bbb\ccc^3$, $\ccc^2\ddd\eee$ \\
					$60$&$\aaa^3$, $\aaa^2\bbb$, $\aaa\bbb^2$, $\bbb^3$ \\
				\end{tabular}
			\end{minipage}
		\end{table}

		\subsubsection*{Table \ref{Tab-2.1}, $f=20$.}
		We have the unique AAD $\ccc^5=\thin^{\aaa}\ccc^{\eee}\thin^{\aaa}\ccc^{\eee}\thin^{\aaa}\ccc^{\eee}\thin\cdots$. By $\aaa\eee\cdots=\aaa\ddd\eee$ and $\bbb^2\cdots=\bbb^2\ccc$, we can determine $16$ tiles in Figure \ref{ade}. This tiling is the non-symmetric $3$-layer earth map tiling $T(20\aaa\ddd\eee,10\bbb^2\ccc,2\ccc^5)$. But there is no $\bbb\ddd\eee$ in this tiling, so it actually belongs to Section \ref{ade as the unique $a^2b$-vertex}.

		\subsubsection*{Table \ref{Tab-2.1}, $f=60$.}
		We have the unique AAD for $\bbb\ddd\eee=\thin\bbb\thin^{\ccc}\eee^{\ddd}\thick^{\eee}\ddd^{\bbb}\thin$. This gives a vertex $\aaa\ccc\cdots$ or $\ccc\ddd\cdots$, contradicting the AVC.

		\vspace{6pt}

		In summary, Case $\{\aaa\ddd\eee, \bbb\ddd\eee, \ccc^5\}$ admits no tiling.

		Cases $\{\aaa\ddd\eee, \bbb\ddd\eee, \ccc^3\}$, $\{\aaa\ddd\eee, \bbb\ddd\eee, \ccc^4\}$ and $\{\aaa\ddd\eee, \bbb\ddd\eee, \aaa\ccc^3(\bbb\ccc^3) \}$ are similar to Case \{$\aaa\ddd\eee, \bbb\ddd\eee, \ccc^5$\} and admit no tiling. We only give the AVC in Table \ref{Tab-3.1}, and skip the details.

		\begin{table}[htp]
			\centering
			\caption{AVC for Cases $\{\aaa\ddd\eee, \bbb\ddd\eee, \ccc^3\}$, $\{\aaa\ddd\eee, \bbb\ddd\eee, \ccc^4\}$ and $\{\aaa\ddd\eee, \bbb\ddd\eee, \aaa\ccc^3(\bbb\ccc^3) \}$.}\label{Tab-3.1}
			\begin{minipage}{1\textwidth}
				\resizebox{\linewidth}{!}{\begin{tabular}{c|c||c|c||c|c}
					\hline
					\multicolumn{2}{c||}{AVC for $\{\aaa\ddd\eee, \bbb\ddd\eee, \ccc^3\}$}&\multicolumn{2}{c||}{AVC for $\{\aaa\ddd\eee, \bbb\ddd\eee, \ccc^4\}$}&\multicolumn{2}{c}{AVC for $\{\aaa\ddd\eee, \bbb\ddd\eee, \aaa\ccc^3(\bbb\ccc^3) \}$}\\
					\hline
					$f$&vertex&$f$&vertex&$f$&vertex\\
					\hline
					\hline
					all&$\aaa\ddd\eee,\bbb\ddd\eee,\ccc^3$&all&$\aaa\ddd\eee,\bbb\ddd\eee,\ccc^4$&all&$\aaa\ddd\eee,\bbb\ddd\eee,\aaa\ccc^3,\bbb\ccc^3$\\
					$24$&$\aaa^4,\aaa^3\bbb,\aaa^2\bbb^2,\aaa\bbb^3,\bbb^4$&$16$&$\aaa^2\ccc,\aaa\bbb\ccc,\bbb^2\ccc$&$20$&$\aaa^2\ccc,\aaa\bbb\ccc,\aaa\ccc^3,\bbb^2\ccc,\bbb\ccc^3,\ccc^2\ddd\eee,\ccc^5$\\
					$36$&$\aaa^3\ccc,\aaa^2\bbb\ccc,\aaa\bbb^2\ccc,\bbb^3\ccc$ &$24$&$\aaa^3,\aaa^2\bbb,\aaa\bbb^2,\bbb^3$&$36$&$\aaa^3,\aaa^2\bbb,\aaa\bbb^2,\aaa\ccc^3,\bbb^3,\bbb\ccc^3$\\
					$60$&$\aaa^2\bbb^3,\aaa^3\bbb^2$ &&&&\\
				\end{tabular}}
			\end{minipage}
		\end{table}

		\subsection*{Case. $\{\aaa\ddd\eee, \bbb\ddd\eee, \aaa\ccc^2(\bbb\ccc^2)\}$}\label{ade bde a2c}
		The angle sums of $\aaa\ddd\eee, \bbb\ddd\eee, \aaa\ccc^2(\bbb\ccc^2)$ imply $\aaa=\bbb=\tfrac8f,\ccc=1-\frac4f,\ddd+\eee=2-\frac8f$. By $\bbb\neq\ccc$, we get $f>12$ and $\bbb<\ccc$.  By Lemma \ref{geometry1}, we have $\ddd>\eee$. Then $\eee<1$. By Lemma \ref{geometry2}, we know $\tfrac{24}f=2\aaa+\bbb>1$. This implies $f<24$. We give the AVC in Table \ref{Tab-2.2}. All possible tilings based on the AVC of Table \ref{Tab-2.2} are deduced as follows. Note that the $\text{AVC}\sub\{\aaa\ddd\eee,\bbb\ddd\eee,\aaa\ccc^2,\bbb\ccc^2\}$ in the first row admits no solution satisfying the balance lemma.

		\begin{table}[htp]
			\centering
			\caption{AVC for  Case $\{\aaa\ddd\eee,\bbb\ddd\eee,\aaa\ccc^2(\bbb\ccc^2)\}$.}\label{Tab-2.2}
			\begin{minipage}{0.7\textwidth}
				\begin{tabular}{c|c}
					\hline
					$f$&vertex\\
					\hline
					\hline
					all&$\aaa\ddd\eee,\bbb\ddd\eee,\aaa\ccc^2,\bbb\ccc^2$\\
					$16$&$\aaa^4,\aaa^3\bbb,\aaa^2\bbb^2,\aaa\bbb^3,\bbb^4$ \\
					$20$&$\aaa^3\ccc,\aaa^2\bbb\ccc,\aaa\bbb^2\ccc,\bbb^3\ccc,\aaa^5,\aaa^4\bbb,\aaa^3\bbb^2,\aaa^2\bbb^3,\aaa\bbb^4,\bbb^5$ \\
				\end{tabular}
			\end{minipage}
		\end{table}

		\subsubsection*{Table \ref{Tab-2.2}, $f=16$.}
		By the AVC, we know $\bbb\ccc\cdots$, $\ccc\ddd\cdots$, $\ccc\eee\cdots$ and $\eee^2\cdots$ are not vertices. If $\aaa^2\bbb^2$ is a vertex, we have the AAD for $\aaa^2\bbb^2=\thin^{\bbb}\aaa^{\ccc}\thin^{\ccc}\aaa^{\bbb}\thin^{\aaa}\bbb^{\ddd}\thin^{\aaa}\bbb^{\ddd}\thin$ or $\thin^{\ccc}\aaa^{\bbb}\thin^{\ddd}\bbb^{\aaa}\thin^{\ccc}\aaa^{\bbb}\thin^{\ddd}\bbb^{\aaa}\thin$. The first case gives a vertex $\ccc^2\cdots$. By $\aaa\ccc^2$, we get $\ccc^2\cdots=\thin\aaa\thin^{\eee}\ccc^{\aaa}\thin^{\aaa}\ccc^{\eee}\thin$. This implies $\ccc\eee\cdots$ is a vertex, contradicting the AVC. The second case gives a vertex $\aaa\ccc\cdots$. By $\aaa\ccc^2$, we get $\aaa\ccc\cdots=\thin^{\ccc}\aaa^{\bbb}\thin^{\aaa}\ccc^{\eee}\thin\ccc\thin$. This implies $\ccc\eee\cdots$ is a vertex, contradicting the AVC. Similarly, we know $\aaa\bbb^3$, $\aaa^3\bbb$, and $\aaa^4$ are not vertices.

		Therefore $\text{AVC}\sub\{\aaa\ddd\eee,\bbb\ddd\eee,\aaa\ccc^2,\bbb\ccc^2,\bbb^4\}$. We have
		\[(\aaa,\bbb,\ccc,\ddd,\eee)=(\tfrac12,\tfrac12,\tfrac34,\tfrac32-\eee,\eee).\]

		By the equations (\ref{coolsaet_eq1}), (\ref{coolsaet_eq2}) and (\ref{cosb}), we get:
	\begin{equation}
		\cos(\eee\pi) = \frac{1}{\sqrt[4]{2}},\quad \cos a = {\textstyle \sqrt{\sqrt{2}-1} }, \quad \cos b = {\textstyle \sqrt{2}\sqrt{\sqrt{2}-1} } \label{coseee}
	\end{equation}

		In Figure \ref{figure 5}, the unique AAD
		\begin{equation*}
			\bbb^4=\thin^{\aaa}\bbb^{\ddd}\thin^{\aaa}\bbb^{\ddd}\thin^{\aaa}\bbb^{\ddd}\thin^{\aaa}\bbb^{\ddd}\thin
		\end{equation*}
		determines $T_1,T_2,T_3,T_4$. Then $\aaa_2\ddd_1\cdots=\aaa_2\ddd_1\eee_5$ determines $T_5$. Similarly, we can determine $T_6,T_7,T_8$. So $\ddd_6\eee_2\cdots=\aaa_{10}\ddd_6\eee_2$ or $\bbb_{10}\ddd_6\eee_2$, which induces two tilings in Figure \ref{figure 5} respectively. It is clear that flipping the lower half of one tiling gives the other tiling.  Their 3D pictures are given in Figure \ref{s16}. The first tiling of Figure \ref{figure 5} is the non-symmetric $3$-layer earth map tiling $T\{16\aaa\ddd\eee,8\bbb\ccc^2,2\bbb^4\}$. But there is no $\bbb\ddd\eee$ in this tiling, so it actually belongs to Section \ref{ade as the unique $a^2b$-vertex}.
		\begin{figure}[htp]
			\centering
			\begin{tikzpicture}[>=latex,scale=0.33]

				\foreach \b in {0,1,2,3}
				{
					\begin{scope}[xshift=4*\b cm]
						\fill[gray!50]  (0-5,0+3)--(0-5,-2+3)--(1-5,-3+3)--(3-5,-3+3)--(4-5,-2+3)--(4-5,0+3)--(0-5,0+3);
						\fill[gray!50]  (4-5,-2+3)--(5-5,-3+3)--(5-5,-5+3)--(3-5,-5+3)--(3-5,-3+3)--(4-5,-2+3);

						\draw (-1-4,3)--(-1-4,1)--(0-4,0)--(2-4,0)--(3-4,1)--(3-4,3)
						(-3,-5)--(-3,-3)--(-2,-2)--(0,-2)--(1,-3)--(1,-5)
						(-2,-2)--(-2,0)--(-1,1)--(0,0)--(0,-2)
						(0-4,0)--(0-4,-2)--(1-4,-3)--(2-4,-2)--(2-4,0);

						\draw[line width=1.5] (-1,1)--(-2,0);
						\draw[line width=1.5] (0-4,-2)--(1-4,-3);

						\node at (1-4,3){\small $\bbb$};
						\node at (-0.4-4,1.3){\small $\aaa$};
						\node at (0.3-4,0.5){\small $\ccc$};
						\node at (1.7-4,0.5){\small $\eee$};
						\node at (2.5-4,1.3){\small $\ddd$};

						\node at (-1,-4.8){\small $\bbb$};
						\node at (-2.4,-3.3){\small $\aaa$};
						\node at (-2,-2.5){\small $\ccc$};
						\node at (0,-2.5){\small $\eee$};
						\node at (0.5,-3.3){\small $\ddd$};

						\node at (-1.7,-1.6){\small $\bbb$};
						\node at (-1.7,-0.2){\small $\ddd$};
						\node at (-1,0.4){\small $\eee$};
						\node at (-0.3,-0.2){\small $\ccc$};
						\node at (-0.3,-1.7){\small $\aaa$};

						\node at (0.4-4,-0.4){\small $\bbb$};
						\node at (0.4-4,-1.7){\small $\ddd$};
						\node at (1-4,-2.4){\small $\eee$};
						\node at (1.6-4,-1.8){\small $\ccc$};
						\node at (1.6-4,-0.5){\small $\aaa$};
					\end{scope}
				}

				\draw[line width=1.5] (0-4+16,-2)--(1-4+16,-3);

				\node[draw,shape=circle, inner sep=0.5] at (1-4,1.6) {\small $1$};
				\node[draw,shape=circle, inner sep=0.5] at (1,1.6) {\small $2$};
				\node[draw,shape=circle, inner sep=0.5] at (1+4,1.6) {\small $3$};
				\node[draw,shape=circle, inner sep=0.5] at (1+8,1.6) {\small $4$};

				\node[draw,shape=circle, inner sep=0.5] at (-1,-0.8) {\small $5$};
				\node[draw,shape=circle, inner sep=0.5] at (-1+4,-0.8) {\small $6$};
				\node[draw,shape=circle, inner sep=0.5] at (-1+8,-0.8) {\small $7$};
				\node[draw,shape=circle, inner sep=0.5] at (-1+12,-0.8) {\small $8$};

				\node[draw,shape=circle, inner sep=0.5] at (1-4,-1.1) {\small $9$};
				\node[draw,shape=circle, inner sep=0.5] at (1,-1.1) {\small $10$};
				\node[draw,shape=circle, inner sep=0.5] at (1+4,-1.1) {\small $11$};
				\node[draw,shape=circle, inner sep=0.5] at (1+8,-1.1) {\small $12$};

				\node[draw,shape=circle, inner sep=0.5] at (-1,-3.5) {\small $13$};
				\node[draw,shape=circle, inner sep=0.5] at (-1+4,-3.5) {\small $14$};
				\node[draw,shape=circle, inner sep=0.5] at (-1+8,-3.5) {\small $15$};
				\node[draw,shape=circle, inner sep=0.5] at (-1+12,-3.5) {\small $16$};
			\end{tikzpicture}
			\begin{tikzpicture}[>=latex,scale=0.33]

				\foreach \b in {0,1,2,3}
				{
					\begin{scope}[xshift=4*\b cm]
						\fill[gray!50]  (0-5,0+3)--(0-5,-2+3)--(1-5,-3+3)--(1-5,-5+3)--(2-5,-6+3)--(2-5,-8+3)--(6-5,-8+3)--(6-5,-6+3)--(5-5,-5+3)--(5-5,-3+3)--(4-5,-2+3)--(4-5,0+3)--(0-5,0+3);

						\draw (-1-4,3)--(-1-4,1)--(0-4,0)--(2-4,0)--(3-4,1)--(3-4,3)
						(-3,-5)--(-3,-3)--(-2,-2)--(0,-2)--(1,-3)--(1,-5)
						(-2,-2)--(-2,0)--(-1,1)--(0,0)--(0,-2)
						(0-4,0)--(0-4,-2)--(1-4,-3)--(2-4,-2)--(2-4,0);

						\draw[line width=1.5] (-1,1)--(-2,0);
						\draw[line width=1.5] (2-4,-2)--(1-4,-3);

						\node at (1-4,3){\small $\bbb$};
						\node at (-0.4-4,1.3){\small $\aaa$};
						\node at (0.3-4,0.5){\small $\ccc$};
						\node at (1.7-4,0.5){\small $\eee$};
						\node at (2.5-4,1.3){\small $\ddd$};

						\node at (-1,-4.8){\small $\bbb$};
						\node at (-2.4,-3.3){\small $\ddd$};
						\node at (-2,-2.5){\small $\eee$};
						\node at (0,-2.5){\small $\ccc$};
						\node at (0.5,-3.3){\small $\aaa$};

						\node at (-1.7,-1.6){\small $\bbb$};
						\node at (-1.7,-0.2){\small $\ddd$};
						\node at (-1,0.4){\small $\eee$};
						\node at (-0.3,-0.2){\small $\ccc$};
						\node at (-0.3,-1.7){\small $\aaa$};

						\node at (0.4-4,-0.4){\small $\aaa$};
						\node at (0.4-4,-1.9){\small $\ccc$};
						\node at (1-4,-2.4){\small $\eee$};
						\node at (1.6-4,-1.7){\small $\ddd$};
						\node at (1.6-4,-0.5){\small $\bbb$};

					\end{scope}
				}

				\node[draw,shape=circle, inner sep=0.5] at (1-4,1.6) {\footnotesize $1$};
				\node[draw,shape=circle, inner sep=0.5] at (1,1.6) {\footnotesize $2$};
				\node[draw,shape=circle, inner sep=0.5] at (1+4,1.6) {\footnotesize $3$};
				\node[draw,shape=circle, inner sep=0.5] at (1+8,1.6) {\footnotesize $4$};

				\node[draw,shape=circle, inner sep=0.5] at (-1,-0.8) {\footnotesize $5$};
				\node[draw,shape=circle, inner sep=0.5] at (-1+4,-0.8) {\footnotesize $6$};
				\node[draw,shape=circle, inner sep=0.5] at (-1+8,-0.8) {\footnotesize $7$};
				\node[draw,shape=circle, inner sep=0.5] at (-1+12,-0.8) {\footnotesize $8$};

				\node[draw,shape=circle, inner sep=0.5] at (1-4,-1.1) {\small $9$};
				\node[draw,shape=circle, inner sep=0.5] at (1,-1.1) {\small $10$};
				\node[draw,shape=circle, inner sep=0.5] at (1+4,-1.1) {\small $11$};
				\node[draw,shape=circle, inner sep=0.5] at (1+8,-1.1) {\small $12$};

				\node[draw,shape=circle, inner sep=0.5] at (-1,-3.5) {\small $13$};
				\node[draw,shape=circle, inner sep=0.5] at (-1+4,-3.5) {\small $14$};
				\node[draw,shape=circle, inner sep=0.5] at (-1+8,-3.5) {\small $15$};
				\node[draw,shape=circle, inner sep=0.5] at (-1+12,-3.5) {\small $16$};
			\end{tikzpicture}
			\caption{For $f=16$, $T\{16\aaa\ddd\eee,8\bbb\ccc^2,2\bbb^4\}$ and $T\{8\aaa\ddd\eee,8\bbb\ddd\eee,8\aaa\ccc^2,2\bbb^4\}$.} \label{figure 5}
		\end{figure}

		\subsubsection*{Table \ref{Tab-2.2}, $f=20$.}
		By the AVC, we know $\ccc^2\cdots$, $\ccc\ddd\cdots$, $\ccc\eee\cdots$ and $\ddd^2\cdots$ are not vertices. If $\bbb^3\ccc$ is a vertex, we have the unique AAD for $\bbb^3\ccc=\thin^{\aaa}\bbb^{\ddd}\thin^{\aaa}\bbb^{\ddd}\thin^{\aaa}\bbb^{\ddd}\thin^{\aaa}\ccc^{\eee}\thin$. This gives a vertex $\aaa\ddd\cdots$. By $\aaa\ddd\eee$, we get $\aaa\ddd\cdots=\thin^{\ccc}\aaa^{\bbb}\thin^{\bbb}\ddd^{\eee}\thick^{\ddd}\eee^{\ccc}\thin$. This implies $\ccc^2\cdots$ is a vertex, a contradiction.

		If $\aaa^2\bbb\ccc$ is a vertex, we have the unique AAD for $\aaa^2\bbb\ccc=\thin^{\ccc}\aaa^{\bbb}\thin^{\ddd}\bbb^{\aaa}\thin^{\ddd}\bbb^{\aaa}\thin^{\eee}\ccc^{\aaa}\thin$, $\thin^{\bbb}\aaa^{\ccc}\thin^{\aaa}\bbb^{\ddd}\thin^{\aaa}\bbb^{\ddd}\thin^{\aaa}\ccc^{\eee}\thin$, $\thin^{\ccc}\aaa^{\bbb}\thin^{\ddd}\bbb^{\aaa}\thin^{\eee}\ccc^{\aaa}\thin^{\ddd}\bbb^{\aaa}\thin$ or $\thin^{\bbb}\aaa^{\ccc}\thin^{\aaa}\bbb^{\ddd}\thin^{\aaa}\ccc^{\eee}\thin^{\aaa}\bbb^{\ddd}\thin$. The first and second cases give a vertex $\aaa\ddd\cdots$. By $\aaa\ddd\eee$, we get $\aaa\ddd\cdots=\thin^{\ccc}\aaa^{\bbb}\thin^{\bbb}\ddd^{\eee}\thick^{\ddd}\eee^{\ccc}\thin$. This implies $\ccc^2\cdots$ is a vertex, a contradiction. The third and fourth cases give a vertex $\bbb\ddd\cdots$. By $\bbb\ddd\eee$, we get $\bbb\ddd\cdots=\thin^{\ddd}\bbb^{\aaa}\thin^{\bbb}\ddd^{\eee}\thick^{\ddd}\eee^{\ccc}\thin$. This implies $\ccc\ddd\cdots$ is a vertex, a contradiction.

		If $\aaa\bbb^2\ccc$ is a vertex, we have the unique AAD for $\aaa\bbb^2\ccc=\thin^{\bbb}\aaa^{\ccc}\thin^{\bbb}\aaa^{\ccc}\thin^{\aaa}\bbb^{\ddd}\thin^{\aaa}\ccc^{\eee}\thin$, $\thin^{\ccc}\aaa^{\bbb}\thin^{\ccc}\aaa^{\bbb}\thin^{\ddd}\bbb^{\aaa}\thin^{\eee}\ccc^{\aaa}\thin$, $\thin^{\ccc}\aaa^{\bbb}\thin^{\ddd}\bbb^{\aaa}\thin^{\ccc}\aaa^{\bbb}\thin^{\eee}\ccc^{\aaa}\thin$ or $\thin^{\bbb}\aaa^{\ccc}\thin^{\aaa}\bbb^{\ddd}\thin^{\bbb}\aaa^{\ccc}\thin^{\aaa}\ccc^{\eee}\thin$. The first case gives a vertex $\bbb\ccc\cdots$. By $\aaa\bbb^2\ccc$, we get $\bbb\ccc\cdots=\thin^{\bbb}\aaa^{\ccc}\thin^{\ccc}\aaa^{\bbb}\thin^{\ddd}\bbb^{\aaa}\thin^{\aaa}\ccc^{\eee}\thin$. This implies $\ccc^2\cdots$ is a vertex, a contradiction. The second, third and fourth cases give a vertex $\bbb\ddd\cdots$. By $\bbb\ddd\eee$, we get $\bbb\ddd\cdots=\thin^{\ddd}\bbb^{\aaa}\thin^{\bbb}\ddd^{\eee}\thick^{\ddd}\eee^{\ccc}\thin$. This implies $\ccc\ddd\cdots$ is a vertex, a contradiction.

		By the previous discussion, we know $\bbb\ccc\cdots$ is not a vertex.

		If $\aaa^3\ccc$ is a vertex, we have the AAD for $\aaa^3\ccc=\thin^{\aaa}\ccc^{\eee}\thin^{\bbb}\aaa^{\ccc}\thin\aaa\thin$. This implies $\bbb\ccc\cdots$ or $\ccc^2\cdots$ is a vertex. This is a contradiction. Similarly, we know $\aaa^5,\aaa^4\bbb,\aaa^3\bbb^2,\aaa^2\bbb^3,\aaa\bbb^4,\bbb^5$ are not vertices.

	\subsection*{Case. $\{\aaa\ddd\eee, \bbb\ddd\eee, \aaa^2\ccc(\aaa\bbb\ccc,\bbb^2\ccc)\}$}\label{ade bde 2ac}
		The angle sums of $\alpha\delta\epsilon,\bbb\ddd\eee,\aaa^2\ccc(\aaa\bbb\ccc,\bbb^2\ccc)$ and the angle sum for the pentagon imply
		\[
		\alpha=\beta=1-\tfrac{4}{f},\;
		\gamma=\tfrac{8}{f},\;
		\delta+\epsilon=1+\tfrac{4}{f}.
		\]

		We have $1>\beta>\gamma$. By Lemma \ref{geometry1}, this implies $\delta<\epsilon$. By the angle values and $n_4=n_5$, we have
		\begin{align}
			\text{AVC}
			\sub\{\aaa\ddd\eee,\bbb\ddd\eee,\aaa^2\ccc,\aaa\bbb\ccc,\bbb^2\ccc,\aaa\ccc^{\frac{f+4}8},\bbb\ccc^{\frac{f+4}8},\ccc^{\frac{f-4}8}\ddd\eee,\ccc^{\frac{f}4}\}.
		\end{align}\label{ade bde}

	    By the balance lemma,  $\aaa\ccc^{\frac{f+4}8},\bbb\ccc^{\frac{f+4}8},\ccc^{\frac{f-4}8}\ddd\eee$ or $\ccc^{\frac{f}4}$ must appear.

	    \subsubsection*{Subcase. $\ccc^{\frac{f}4}$ is a vertex.}\label{c 1}
	    In Figure \ref{ade}, the unique AAD of  $\ccc^{\frac{f}{4}}=\thin^{\aaa}\ccc_{1}^{\eee}\thin^{\aaa}\ccc_{1'}^{\eee}\thin\cdots$ determines $\frac{f}{4}$ tiles $T_1,T_{1'},\cdots$. Then $\aaa_1\eee\cdots=\aaa_1\ddd_2\eee$ determines $T_2$; $\aaa_{1'}\eee\cdots=\aaa_{1'}\ddd_{2'}\eee$ determines $T_{2'}$. Then $\bbb_1\bbb_2\cdots=\bbb_1\bbb_2\ccc_3$ and $\ddd_1\eee_{2'}\cdots=\aaa_3\ddd_1\eee_{2'}$ determine $T_3$; $\aaa_2\eee_3\cdots=\aaa_2\ddd_4\eee_3$ determines $T_4$. The tiles $T_1,T_2,T_3,T_4$ together form a \textbf{time zone}. Similarly, we can determine $T_{1'},T_{2'},T_{3'},T_{4'}$. By repeating the process, we get the non-symmetric $3$-layer  earth map tiling $T(f\aaa\ddd\eee,\frac{f}2\bbb^2\ccc,2\ccc^{\frac{f}4})$. But there is no $\bbb\ddd\eee$ in this tiling, so it actually belongs to Section \ref{ade as the unique $a^2b$-vertex}.
	    \begin{figure}[htp]
	    	\centering
	    	\begin{tikzpicture}[>=latex,scale=0.37]
	    		\foreach \a in {0,1,2}
	    		{
	    			\begin{scope}[xshift=4*\a cm]

	    				\fill[gray!50]  (-3,-5)--(-3,-3)--(-2,-2)--(0,-2)--(1,-3)--(1,-5);
	    				\fill[gray!50]  (0,0)--(0,-2)--(1,-3)--(2,-2)--(2,0);

	    				\draw (-1,3)--(-1,1)--(0,0)--(2,0)--(3,1)--(3,3)
	    				(-3,-5)--(-3,-3)--(-2,-2)--(0,-2)--(1,-3)--(1,-5)
	    				(-2,-2)--(-2,0)--(-1,1)--(0,0)--(0,-2)
	    				(0,0)--(0,-2)--(1,-3)--(2,-2)--(2,0);

	    				\draw[line width=1.5] (-1,1)--(-2,0);
	    				\draw[line width=1.5] (0,-2)--(1,-3);

	    				\node at (1,3){\small $\ccc$};
	    				\node at (-0.4,1.3){\small $\aaa$};
	    				\node at (-1.5,1.3){\small $\eee$};
	    				\node at (0.3,0.5){\small $\bbb$};
	    				\node at (1.7,0.5){\small $\ddd$};
	    				\node at (2.5,1.3){\small $\eee$};

	    				\node at (-1,-4.8){\small $\ccc$};
	    				\node at (-2.4,-3.3){\small $\aaa$};
	    				\node at (-2,-2.5){\small $\bbb$};
	    				\node at (0,-2.5){\small $\ddd$};
	    				\node at (0.5,-3.3){\small $\eee$};

	    				\node at (-1.7,-1.6){\small $\ccc$};
	    				\node at (-1.7,-0.2){\small $\eee$};
	    				\node at (-1,0.4){\small $\ddd$};
	    				\node at (-0.3,-0.2){\small $\bbb$};
	    				\node at (-0.3,-1.7){\small $\aaa$};

	    				\node at (0.4,-0.4){\small $\ccc$};
	    				\node at (0.4,-1.7){\small $\eee$};
	    				\node at (1,-2.4){\small $\ddd$};
	    				\node at (1.6,-1.7){\small $\bbb$};
	    				\node at (1.6,-0.3){\small $\aaa$};
	    			\end{scope}
	    		}

	    		\fill (12,-1) circle (0.1);
	    		\fill (12.6,-1) circle (0.1);
	    		\fill (13.2,-1) circle (0.1);

	    		\node[draw,shape=circle, inner sep=0.5] at (-3+4,3-1) {\small $1$};
	    		\node[draw,shape=circle, inner sep=0.5] at (-3+4+4,3-1) {\small $1'$};
	    		\node[draw,shape=circle, inner sep=0.5] at (-4+1+4,-1) {\small $3$};
	    		\node[draw,shape=circle, inner sep=0.5] at (-4+1+4+4,-1) {\small $3'$};
	    		\node[draw,shape=circle, inner sep=0.5] at (-4+1+6,-1) {\small $2'$};
	    		\node[draw,shape=circle, inner sep=0.5] at (-4+1+2,-1) {\small $2$};
	    		\node[draw,shape=circle, inner sep=0.5] at (1-2,-3) {\small $4$};
	    		\node[draw,shape=circle, inner sep=0.5] at (1+2,-3) {\small $4'$};

	    		\draw[line width=1.5] (11,1)--(10,0);
	    	\end{tikzpicture}
	    	\caption{The non-symmetric $3$-layer earth map tiling $T(f \, \aaa\ddd\eee, \frac{f}2 \, \bbb^2\ccc, 2\ccc^{\frac{f}4})$.} \label{ade}
	    \end{figure}

	    \subsubsection*{Subcase. $\bbb\ccc^{\frac{f+4}{8}}$ is a vertex.}\label{c 2}
	    In the first row of Figure \ref{9.1}, the unique AAD of  $\bbb\ccc^{\frac{f+4}{8}}=\thin^{\aaa}\ccc_{1}^{\eee}\thin^{\aaa}\ccc_{2}^{\eee}\thin\cdots\thin^{\aaa}\ccc_{3}^{\eee}\thin^{\aaa}\bbb_4^{\ddd}\thin$ determines $\frac{f+12}{8}$ tiles $T_1,T_{2},\cdots,T_3,T_{4}$. Then $R(\bbb_4\ccc_1\cdots)=\ccc^{\frac{f-4}8}$ and such $\ccc^{\frac{f-4}8}$ determines $\frac{f-4}8$ time zones ($\frac{f}2-2$ tiles). Then $\aaa_4\eee_3\cdots=\aaa_4\ddd_5\eee_3$ determines $T_{5}$. These timezones together with $T_1$ and $T_5$ have $\frac{f}{2}$ tiles, and form half of the sphere with boundary vertices $A_1,\cdots,A_{12}$.  We have $\bbb_5\ccc_4\cdots=\bbb^2\ccc$ or $\bbb\ccc^{\frac{f+4}{8}}$.

	    In the left of the second row of Figure \ref{9.1}, $\bbb_5\ccc_4\cdots=\bbb_5\bbb_6\ccc_4$. Then $\aaa_5\cdots=\aaa_5\ddd_6\eee_7$ determines $T_6,T_7$. Then $R(\bbb\ccc_5\cdots)=\ccc^{\frac{f-4}8}$ and such $\ccc^{\frac{f-4}8}$ determines $\frac{f-4}8$ time zones ($\frac{f}2-2$ tiles). Then $\aaa\ddd\cdots=\aaa\ddd\eee_8$ determines $T_{8}$. These timezones together with $T_4$ and $T_8$ also have $\frac{f}{2}$ tiles, and form exactly the other half of the sphere. Its boundary vertices $A_1$, $\cdots$, $A_{12}$ can be glued with the $A_1$, $\cdots$, $A_{12}$ in the first row and form a tiling.

	    In the right of the second row of Figure \ref{9.1}, $\bbb_5\ccc_4\cdots=\bbb\ccc^{\frac{f+4}{8}}$. Then $R(\bbb_5\ccc_4\cdots)=\ccc^{\frac{f-4}8}$ and such $\ccc^{\frac{f-4}8}$ determines $\frac{f-4}8$ time zones ($\frac{f}2-2$ tiles). Then $\aaa_5\eee\cdots=\aaa_5\ddd_6\eee$ determines $T_{6}$. These timezones together with $T_4$ and $T_6$ also have $\frac{f}{2}$ tiles, and form exactly the other half of the sphere. Its boundary vertices $A_1$, $\cdots$, $A_{12}$ can be glued with the $A_1$, $\cdots$, $A_{12}$ in the first row and form another tiling.

	    These two tilings turn out to be equivalent after a rotation, and give the standard flip modification $T(f\,\aaa\ddd\eee,\frac{f-4}{2}\bbb^2\ccc,4\bbb\ccc^{\frac{f+4}{8}})$ (see Figure \ref{ade 2bc filp} and \ref{The earth map}) of the non-symmetric $3$-layer earth map tiling. But there is no $\bbb\ddd\eee$ in this tiling, so it actually belongs to Section \ref{ade as the unique $a^2b$-vertex}.

	    \begin{figure}[htp]
	    	\centering
	    	\begin{tikzpicture}[>=latex,scale=0.26]
	    		\begin{scope}[xshift=0cm,yshift=0cm]
	    			\fill[gray!50]  (-1,-5) -- (-1,-3) -- (0,-2) -- (2,-2) -- (3,-3) -- (3,-5);
	    			\fill[gray!50]  (-1+8,-5) -- (-1+8,-3) -- (0+8,-2) -- (2+8,-2) -- (3+8,-3) -- (3+8,-5);
	    			\fill[gray!50]  (-2+4,0) -- (0+4,0) -- (-0+4,-2) -- (-1+4,-3) -- (-2+4,-2) -- (-2+4,0);
	    			\fill[gray!50]  (-2+12,0) -- (0+12,0) -- (-0+12,-2) -- (-1+12,-3) -- (-2+12,-2) -- (-2+12,0);
	    			\fill[gray!50]  (-3+16,3) -- (-3+16,1) -- (-2+16,0) -- (0+16,0) -- (1+16,1) -- (1+16,3);

	    			\draw
	    			(-3,3) -- (-3,1) -- (-2,0) -- (0,0) -- (1,1) -- (1,3);
	    			\draw
	    			(-2+4,0) -- (0+4,0) -- (-0+4,-2) -- (-1+4,-3) -- (-2+4,-2) -- (-2+4,0);
	    			\draw
	    			(-3+4,1) -- (-2+4,0) -- (-2+4,-2) -- (-4+4,-2) -- (-4+4,-0) -- (-3+4,1);
	    			\draw
	    			(-1,-5) -- (-1,-3) -- (0,-2) -- (2,-2) -- (3,-3) -- (3,-5);
	    			\draw
	    			(-3+4,3) -- (-3+4,1) -- (-2+4,0) -- (0+4,0) -- (1+4,1) -- (1+4,3);
	    			\draw
	    			(-2+12,0) -- (0+12,0) -- (-0+12,-2) -- (-1+12,-3) -- (-2+12,-2) -- (-2+12,0);
	    			\draw
	    			(-3+12,1) -- (-2+12,0) -- (-2+12,-2) -- (-4+12,-2) -- (-4+12,-0) -- (-3+12,1);
	    			\draw
	    			(-1+8,-5) -- (-1+8,-3) -- (0+8,-2) -- (2+8,-2) -- (3+8,-3) -- (3+8,-5);
	    			\draw
	    			(-3+12,3) -- (-3+12,1) -- (-2+12,0) -- (0+12,0) -- (1+12,1) -- (1+12,3);
	    			\draw
	    			(-3+16,1) -- (-2+16,0) -- (-2+16,-2) -- (-4+16,-2) -- (-4+16,-0) -- (-3+16,1);
	    			\draw
	    			(-3+16,3) -- (-3+16,1) -- (-2+16,0) -- (0+16,0) -- (1+16,1) -- (1+16,3);

	    			\draw[line width=1.5]
	    			(0,0) -- (1,1);
	    			\draw[line width=1.5]
	    			(0+4,0) -- (1+4,1);
	    			\draw[line width=1.5]
	    			(0+8,0) -- (1+8,1);
	    			\draw[line width=1.5]
	    			(0+12,0) -- (1+12,1);
	    			\draw[line width=1.5]
	    			(0+16,0) -- (1+16,1);

	    			\draw[line width=1.5]
	    			(2,-2) -- (3,-3);
	    			\draw[line width=1.5]
	    			(2+8,-2) -- (3+8,-3);

	    			\fill (5.5,-1) circle (0.1); \fill (6,-1) circle (0.1); \fill (6.5,-1) circle (0.1);

	    			\node at (-3-1,1) {\scriptsize $A_2$};
	    			\node at (-3-1,3) {\scriptsize $A_1$};
	    			\node at (13+1,1+0.5) {\scriptsize $A_{12}$};
	    			\node at (14+1,0+0.5) {\scriptsize $A_{11}$};
	    			\node at (14+1,-4+2) {\scriptsize $A_{10}$};
	    			\node at (12+0.7,-4+1) {\scriptsize $A_{9}$};
	    			\node at (14+1+0.7-4,-4+0.5) {\scriptsize $A_{8}$};
	    			\node at (14+1+0.7-4,-5) {\scriptsize $A_{7}$};
	    			\node at (-1-1,-3) {\scriptsize $A_6$};
	    			\node at (-3+2.1,-2+0.3) {\scriptsize $A_{5}$};
	    			\node at (-3+2.1,0-0.5) {\scriptsize $A_{4}$};
	    			\node at (-3+0.3,-0.5) {\scriptsize $A_{3}$};

	    			\fill (-3,1) circle (0.15);
	    			\fill (-3,3) circle (0.15);
	    			\fill (1+12,1) circle (0.15);
	    			\fill (-2+16,0) circle (0.15);
	    			\fill (0+14,-2) circle (0.15);
	    			\fill (0+12,-2) circle (0.15);
	    			\fill (3+8,-3) circle (0.15);
	    			\fill (3+8,-5) circle (0.15);
	    			\fill (-1,-3) circle (0.15);
	    			\fill (0,-2) circle (0.15);
	    			\fill (0,0) circle (0.15);
	    			\fill (-2,0) circle (0.15);

	    			\node[draw,shape=circle, inner sep=0.5] at (-3+2,3-1) {\small $1$};
	    			\node[draw,shape=circle, inner sep=0.5] at (-3+6,3-1) {\small $2$};
	    			\node[draw,shape=circle, inner sep=0.5] at (-3+14,3-1) {\small $3$};
	    			\node[draw,shape=circle, inner sep=0.5] at (-3+18.4,3-1) {\small $4$};
	    			\node[draw,shape=circle, inner sep=0.5] at (-4+1+16,-1) {\small $5$};

	    			\draw[line width=3pt, ->]
	    			(-3,-5) -- (-5,-7);
	    			\draw[line width=3pt, ->]
	    			(13,-5) -- (15,-7);
	    		\end{scope}

	    		\begin{scope}[xshift=12cm,yshift=-10cm]
	    			\fill[gray!50]  (-1,-5) -- (-1,-3) -- (0,-2) -- (2,-2) -- (3,-3) -- (3,-5);
	    			\fill[gray!50]  (-1+8,-5) -- (-1+8,-3) -- (0+8,-2) -- (2+8,-2) -- (3+8,-3) -- (3+8,-5);
	    			\fill[gray!50]  (-2+4,0) -- (0+4,0) -- (-0+4,-2) -- (-1+4,-3) -- (-2+4,-2) -- (-2+4,0);
	    			\fill[gray!50]  (-2+12,0) -- (0+12,0) -- (-0+12,-2) -- (-1+12,-3) -- (-2+12,-2) -- (-2+12,0);

	    			\draw
	    			(-3,3) -- (-3,1) -- (-2,0) -- (0,0) -- (1,1) -- (1,3);
	    			\draw
	    			(-2+4,0) -- (0+4,0) -- (-0+4,-2) -- (-1+4,-3) -- (-2+4,-2) -- (-2+4,0);
	    			\draw
	    			(-3+4,1) -- (-2+4,0) -- (-2+4,-2) -- (-4+4,-2) -- (-4+4,-0) -- (-3+4,1);
	    			\draw
	    			(-1,-5) -- (-1,-3) -- (0,-2) -- (2,-2) -- (3,-3) -- (3,-5);
	    			\draw
	    			(-3+4,3) -- (-3+4,1) -- (-2+4,0) -- (0+4,0) -- (1+4,1) -- (1+4,3);
	    			\draw
	    			(-2+12,0) -- (0+12,0) -- (-0+12,-2) -- (-1+12,-3) -- (-2+12,-2) -- (-2+12,0);
	    			\draw
	    			(-3+12,1) -- (-2+12,0) -- (-2+12,-2) -- (-4+12,-2) -- (-4+12,-0) -- (-3+12,1);
	    			\draw
	    			(-1+8,-5) -- (-1+8,-3) -- (0+8,-2) -- (2+8,-2) -- (3+8,-3) -- (3+8,-5);
	    			\draw
	    			(-3+12,3) -- (-3+12,1) -- (-2+12,0) -- (0+12,0) -- (1+12,1) -- (1+12,3);
	    			\draw
	    			(-3+16,1) -- (-2+16,0) -- (-2+16,-2) -- (-4+16,-2) -- (-4+16,-0) -- (-3+16,1);

	    			\draw[line width=1.5]
	    			(0,0) -- (1,1);
	    			\draw[line width=1.5]
	    			(0+4,0) -- (1+4,1);
	    			\draw[line width=1.5]
	    			(0+8,0) -- (1+8,1);
	    			\draw[line width=1.5]
	    			(0+12,0) -- (1+12,1);

	    			\draw[line width=1.5]
	    			(2,-2) -- (3,-3);
	    			\draw[line width=1.5]
	    			(2+8,-2) -- (3+8,-3);

	    			\fill (5.5,-1) circle (0.1); \fill (6,-1) circle (0.1); \fill (6.5,-1) circle (0.1);

	    			\node at (-3-1,1) {\scriptsize $A_{12}$};
	    			\node at (-3-1,3) {\scriptsize $A_{11}$};
	    			\node at (13+1,1+0.5) {\scriptsize $A_{10}$};
	    			\node at (14+1,0+0.5) {\scriptsize $A_{9}$};
	    			\node at (14+1,-4+2) {\scriptsize $A_{8}$};
	    			\node at (12+0.7,-4+1) {\scriptsize $A_{7}$};
	    			\node at (14+1+0.7-4,-4+0.5) {\scriptsize $A_{6}$};
	    			\node at (14+1+0.7-4,-5) {\scriptsize $A_{5}$};
	    			\node at (-1-1,-3) {\scriptsize $A_4$};
	    			\node at (-3+2.1,-2+0.3) {\scriptsize $A_{3}$};
	    			\node at (-3+2.1,0-0.5) {\scriptsize $A_{2}$};
	    			\node at (-3+0.3,-0.5) {\scriptsize $A_{1}$};

	    			\fill (-3,1) circle (0.15);
	    			\fill (-3,3) circle (0.15);
	    			\fill (1+12,1) circle (0.15);
	    			\fill (-2+16,0) circle (0.15);
	    			\fill (0+14,-2) circle (0.15);
	    			\fill (0+12,-2) circle (0.15);
	    			\fill (3+8,-3) circle (0.15);
	    			\fill (3+8,-5) circle (0.15);
	    			\fill (-1,-3) circle (0.15);
	    			\fill (0,-2) circle (0.15);
	    			\fill (0,0) circle (0.15);
	    			\fill (-2,0) circle (0.15);

	    			\node[draw,shape=circle, inner sep=0.5] at (-3+2,3-1) {\small $4$};
	    			\node[draw,shape=circle, inner sep=0.5] at (-4+1+16,-1) {\small $6$};
	    		\end{scope}

	    		\begin{scope}[xshift=-12cm,yshift=-10cm]
	    			\fill[gray!50]  (-1,-5) -- (-1,-3) -- (0,-2) -- (2,-2) -- (3,-3) -- (3,-5);
	    			\fill[gray!50]  (-1+8,-5) -- (-1+8,-3) -- (0+8,-2) -- (2+8,-2) -- (3+8,-3) -- (3+8,-5);
	    			\fill[gray!50]  (-2+4,0) -- (0+4,0) -- (-0+4,-2) -- (-1+4,-3) -- (-2+4,-2) -- (-2+4,0);
	    			\fill[gray!50]  (-2+12,0) -- (0+12,0) -- (-0+12,-2) -- (-1+12,-3) -- (-2+12,-2) -- (-2+12,0);

	    			\draw
	    			(-3,3) -- (-3,1) -- (-2,0) -- (0,0) -- (1,1) -- (1,3);
	    			\draw
	    			(-2+4,0) -- (0+4,0) -- (-0+4,-2) -- (-1+4,-3) -- (-2+4,-2) -- (-2+4,0);
	    			\draw
	    			(-3+4,1) -- (-2+4,0) -- (-2+4,-2) -- (-4+4,-2) -- (-4+4,-0) -- (-3+4,1);
	    			\draw
	    			(-1,-5) -- (-1,-3) -- (0,-2) -- (2,-2) -- (3,-3) -- (3,-5);
	    			\draw
	    			(-3+4,3) -- (-3+4,1) -- (-2+4,0) -- (0+4,0) -- (1+4,1) -- (1+4,3);
	    			\draw
	    			(-2+12,0) -- (0+12,0) -- (-0+12,-2) -- (-1+12,-3) -- (-2+12,-2) -- (-2+12,0);
	    			\draw
	    			(-3+12,1) -- (-2+12,0) -- (-2+12,-2) -- (-4+12,-2) -- (-4+12,-0) -- (-3+12,1);
	    			\draw
	    			(-1+8,-5) -- (-1+8,-3) -- (0+8,-2) -- (2+8,-2) -- (3+8,-3) -- (3+8,-5);
	    			\draw
	    			(-3+12,3) -- (-3+12,1) -- (-2+12,0) -- (0+12,0) -- (1+12,1) -- (1+12,3);
	    			\draw
	    			(-3+16,1) -- (-2+16,0) -- (-2+16,-2) -- (-4+16,-2) -- (-4+16,-0) -- (-3+16,1);

	    			\draw[line width=1.5]
	    			(0,0) -- (1,1);
	    			\draw[line width=1.5]
	    			(0+4,0) -- (1+4,1);
	    			\draw[line width=1.5]
	    			(0+8,0) -- (1+8,1);
	    			\draw[line width=1.5]
	    			(0+12,0) -- (1+12,1);

	    			\draw[line width=1.5]
	    			(2,-2) -- (3,-3);
	    			\draw[line width=1.5]
	    			(2+8,-2) -- (3+8,-3);

	    			\fill (5.5,-1) circle (0.1); \fill (6,-1) circle (0.1); \fill (6.5,-1) circle (0.1);

	    			\node at (-3-1,1) {\scriptsize $A_{4}$};
	    			\node at (-3-1+16+2,3) {\scriptsize $A_{3}$};
	    			\node at (13+1,1+0.5) {\scriptsize $A_{2}$};
	    			\node at (14+1,0+0.5) {\scriptsize $A_{1}$};
	    			\node at (14+1,-4+2) {\scriptsize $A_{12}$};
	    			\node at (12+0.7,-4+1+0.2) {\scriptsize $A_{11}$};
	    			\node at (14+1+1-4,-4+0.5) {\scriptsize $A_{10}$};
	    			\node at (14+1+0.7-4-12-2,-5) {\scriptsize $A_{9}$};
	    			\node at (-1-1,-3) {\scriptsize $A_8$};
	    			\node at (-3+2.1,-2+0.3) {\scriptsize $A_{7}$};
	    			\node at (-3+2.1,0-0.5) {\scriptsize $A_{6}$};
	    			\node at (-3+0.3,-0.5) {\scriptsize $A_{5}$};

	    			\fill (-3,1) circle (0.15);
	    			\fill (-3+16,3) circle (0.15);
	    			\fill (1+12,1) circle (0.15);
	    			\fill (-2+16,0) circle (0.15);
	    			\fill (0+14,-2) circle (0.15);
	    			\fill (0+12,-2) circle (0.15);
	    			\fill (3+8,-3) circle (0.15);
	    			\fill (3+8-12,-5) circle (0.15);
	    			\fill (-1,-3) circle (0.15);
	    			\fill (0,-2) circle (0.15);
	    			\fill (0,0) circle (0.15);
	    			\fill (-2,0) circle (0.15);

	    			\node[draw,shape=circle, inner sep=0.5] at (-3+2,3-1) {\small $8$};
	    			\node[draw,shape=circle, inner sep=0.5] at (1+8,-4) {\small $7$};
	    			\node[draw,shape=circle, inner sep=0.5] at (-4+1+14,-1) {\small $6$};
	    			\node[draw,shape=circle, inner sep=0.5] at (-4+1+16,-1) {\small $4$};
	    		\end{scope}
	    	\end{tikzpicture}
	    	\caption{The standard flip modification of the $3$-layer earth map tiling.}
	    	\label{9.1}
	    \end{figure}
        \begin{figure}[htp]
        	\centering
        	\begin{tikzpicture}[>=latex,scale=0.27]
        		\begin{scope}
        			\draw (0,5.34)--(3.2,4.31)--(5.17,1.59)--(5.17,-1.77)--(3.2,-4.49)--(0,-5.53)--(-3.2,-4.49)--(-5.17,-1.77)--(-5.17,1.59)--(-3.2,4.31)--(0,5.34)
        			(3.2,4.31)--(2.22,2.96)
        			(5.17,-1.77)--(3.43,-1.21)
        			(-5.17,-1.77)--(-3.52,-1.24)
        			(-3.2,4.31)--(-2.23,2.97)
        			(-5.17,1.59)--(-6.51,2.02)
        			(-3.2,-4.49)--(-4.25,-5.94)
        			(3.2,-4.49)--(4.34,-6.06)
        			(5.17,1.59)--(6.59,2.05)
        			(5.03,6.83)--(-5.53,-7.67);

        			\draw[line width=1.5] (-3.2,-4.49)--(-4.25,-5.94);
        			\draw[line width=1.5] (3.2,4.31)--(2.22,2.96);
        			\draw[line width=1.5] (-5.17,1.59)--(-3.2,4.31);
        			\draw[line width=1.5] (5.17,-1.77)--(3.2,-4.49);

        			\node at (0,4.5){\small $\ccc^q$};
        			\node at (2.2,3.96){\small $\eee$};
        			\node at (3.3,3.16){\small $\ddd$};
        			\node at (4.24,1.3){\small $\bbb$};
        			\node at (4.67,-1.04){\small $\aaa$};
        			\node at (4.2,-2.17){\small $\ddd$};
        			\node at (3.0,-3.9){\small $\eee$};
        			\node at (0.5,-4.5){\small $\ccc^{q+1}$};
        			\node at (-2.8,-4){\small $\aaa$};
        			\node at (-4.1,-2.07){\small $\bbb$};
        			\node at (-4.3,-0.97){\small $\ccc$};
        			\node at (-4.57,1.5){\small $\eee$};
        			\node at (-3.2,3.41){\small $\ddd$};
        			\node at (-2.3,4.01){\small $\aaa$};

        			\node at (0.5,6.04){\small $\ccc^{q+1}$};
        			\node at (3.5,4.81){\small $\aaa$};
        			\node at (5.47,2.39){\small $\bbb$};
        			\node at (5.67,1.09){\small $\ccc$};
        			\node at (5.77,-1.87){\small $\eee$};
        			\node at (4.0,-4.29){\small $\ddd$};
        			\node at (3.1,-5.19){\small $\aaa$};
        			\node at (0,-6.13){\small $\ccc^q$};
        			\node at (-2.9,-5.05){\small $\eee$};
        			\node at (-3.95,-4.49){\small $\ddd$};
        			\node at (-5.67,-1.87){\small $\bbb$};
        			\node at (-5.67,1.1){\small $\aaa$};
        			\node at (-5.47,2.29){\small $\ddd$};
        			\node at (-3.5,4.71){\small $\eee$};
        			\draw[line width=3pt, ->]
        			(8.5,0) -- (12,0);
        		\end{scope}
        		\begin{scope}[xshift=20cm]
        			\draw (0,5.34)--(3.2,4.31)--(5.17,1.59)--(5.17,-1.77)--(3.2,-4.49)--(0,-5.53)--(-3.2,-4.49)--(-5.17,-1.77)--(-5.17,1.59)--(-3.2,4.31)--(0,5.34)
        			(3.2,4.31)--(2.22,2.96)
        			(5.17,-1.77)--(3.43,-1.21)
        			(-3.2,4.31)--(-2.23,2.97)
        			(-5.17,1.59)--(-6.51,2.02)
        			(-3.2,-4.49)--(-4.25,-5.94)
        			(3.2,-4.49)--(4.34,-6.06)
        			(5.17,1.59)--(6.59,2.05)
        			(0,-5.5)--(0,-3.5);

        			\draw[line width=1.5] (-3.2,-4.49)--(-4.25,-5.94);
        			\draw[line width=1.5] (3.2,4.31)--(2.22,2.96);
        			\draw[line width=1.5] (-5.17,1.59)--(-3.2,4.31);
        			\draw[line width=1.5] (5.17,-1.77)--(3.2,-4.49);

        			\node at (0,4.5){\small $\bbb$};
        			\node at (2.2,3.96){\small $\ddd$};
        			\node at (3.3,3.16){\small $\eee$};
        			\node at (4.24,1.3){\small $\ccc^q$};
        			\node at (4.67,-1.04){\small $\aaa$};
        			\node at (4.2,-2.17){\small $\ddd$};
        			\node at (3.0,-3.9){\small $\eee$};
        			\node at (-0.6,-4.5){\small $\bbb$};
        			\node at (0.6,-4.5){\small $\ccc$};
        			\node at (-2.8,-4){\small $\aaa$};
        			\node at (-3.8,-1.4){\small $\ccc^{q+1}$};
        			\node at (-4.57,1.5){\small $\eee$};
        			\node at (-3.2,3.41){\small $\ddd$};
        			\node at (-2.3,4.01){\small $\aaa$};

        			\node at (0.5,6.04){\small $\ccc^{q+1}$};
        			\node at (3.5,4.81){\small $\aaa$};
        			\node at (5.47,2.39){\small $\bbb$};
        			\node at (5.67,1.09){\small $\ccc$};
        			\node at (5.77,-1.87){\small $\eee$};
        			\node at (4.0,-4.29){\small $\ddd$};
        			\node at (3.1,-5.19){\small $\aaa$};
        			\node at (0,-6.13){\small $\ccc^q$};
        			\node at (-2.9,-5.05){\small $\eee$};
        			\node at (-3.95,-4.49){\small $\ddd$};
        			\node at (-5.67,-1.87){\small $\bbb$};
        			\node at (-5.67,1.1){\small $\aaa$};
        			\node at (-5.47,2.29){\small $\ddd$};
        			\node at (-3.5,4.71){\small $\eee$};

        			\node at (0,0){\small flipped};
        		\end{scope}
        	\end{tikzpicture}
        	\caption{The loop of edges dividing the sphere into two identical halves.} \label{ade 2bc filp}
        \end{figure}

        \vspace{3pt}

	    Henceforth, we can assume
		\[
		\text{AVC}\subset
		\{\aaa\ddd\eee,\bbb\ddd\eee,\aaa^2\ccc,\aaa\bbb\ccc,\bbb^2\ccc,\aaa\ccc^{\frac{f+4}8},\ccc^{\frac{f-4}8}\ddd\eee\}.
		\]

		\subsubsection*{Subcase. $\aaa\ccc^{\frac{f+4}8}$ is a vertex.}
		By the AVC, we get that $\eee^2\cdots$ is not a vertex. Then we have
	\begin{eqnarray*}
		\aaa\ccc^{\frac{f+4}8} &=& \thin^{\bbb}\aaa^{\ccc}\thin^{\eee}\ccc^{\aaa}\thin\cdots\thin^{\eee}\ccc^{\aaa}\thin, \\
		 &\text{or}& \thin^{\bbb}\aaa^{\ccc}\thin^{\eee}\ccc^{\aaa}\thin\cdots\thin^{\eee}\ccc^{\aaa}\thin^{\aaa}\ccc^{\eee}\thin, \\
		 &\text{or}& \thin^{\bbb}\aaa^{\ccc}\thin^{\eee}\ccc^{\aaa}\thin\cdots\thin^{\eee}\ccc^{\aaa}\thin^{\aaa}\ccc^{\eee}\thin\cdots\thin^{\aaa}\ccc^{\eee}\thin, \\
		 &\text{or}& \thin^{\bbb}\aaa^{\ccc}\thin^{\aaa}\ccc^{\eee}\thin\cdots\thin^{\aaa}\ccc^{\eee}\thin.
	\end{eqnarray*}

	\vspace{6pt}

    $\bullet$	In Figure \ref{ac1}, we have $\aaa\ccc^{\frac{f+4}8}=\thin^{\bbb}\aaa_1^{\ccc}\thin^{\eee}\ccc_2^{\aaa}\thin\cdots\thin^{\eee}\ccc_4^{\aaa}\thin$, which determines $T_1$, $T_2$, $T_3$, $\cdots$, $T_4$ ($\frac{f+4}8$ tiles). Then $\aaa_2\eee_3\cdots=\aaa_2\ddd_5\eee_3$ determines $T_5$; $\bbb_2\bbb_5\cdots=\bbb_2\bbb_5\ccc_6$ determines $T_6$; $\aaa_5\eee_6\cdots=\aaa_5\ddd_7\eee_6$ determines $T_7$; Similarly, $R(\aaa_1\ccc_2\ccc_4\cdots)=\ccc^{\frac{f-12}8}$ and such $\ccc^{\frac{f-12}8}$ determines $\frac{f-12}8$ time zones ($\frac{3f-36}8$ tiles). Then $\aaa_6\ddd_2\cdots=\aaa_6\ddd_2\eee_8$ determines $T_8$.

	By $\ddd_8\eee_2\cdots=\thin^{\aaa}\ccc^{\eee}\thin\cdots\thin^{\aaa}\ccc^{\eee}\thin^{\bbb}\ddd^{\eee}\thick^{\ddd}\eee^{\ccc}\thin$, we can determine $T_{9},\cdots,T_{10}$, altogether determining $\frac{f-4}8$ tiles. Then $\aaa_2\eee_9\cdots=\aaa_2\ddd_{11}\eee_9$ determines $T_{11}$; $\bbb_{9}\bbb_{11}\cdots=\bbb_{9}\bbb_{11}\ccc_{12}$ determines $T_{12}$; $\aaa_{11}\eee_{12}\cdots=\aaa_{11}\ddd_{13}\eee_{12}$ determines $T_{13}$. Then $R(\thin^{\bbb}\ddd_8^{\eee}\thick^\ddd\eee_2^\ccc$ $\thin^\aaa\ccc_{1}^{\eee}\thin\cdots)=\ccc^{\frac{f-4}8}$ and such $\ccc^{\frac{f-4}8}$ determines $\frac{f-4}8$ time zones ($(\frac{f}2-2)$ tiles). Then $\aaa_4\ccc_{17}\cdots=\aaa_{4}\bbb_{1}\cdots=\aaa_{4}\bbb_{1}\ccc_{17}$ determines $T_{17}$; $\bbb_4\eee_{17}\cdots=\bbb_4\ddd_{18}\eee_{17}$ determines $T_{18}$.

	More repetitions give the unique tiling with the AVC
	\{$4\aaa\bbb\ccc,(f-6)\aaa\ddd\eee,(\frac{f}2-4)\bbb^2\ccc,4\bbb\ddd\eee,2\aaa\ccc^{\frac{f+4}8},2\ccc^{\frac{f-4}8}\ddd\eee$\}. In Figure \ref{ac1}, $T_{4}$, $T_{17}$, $T_{18}$, $T_{19}$ form a UFO (see Page \pageref{UFO} for its definition and Figure \ref{f=60} and \ref{various flip}).

		\begin{figure}[htp]
			\centering
			\begin{tikzpicture}[>=latex,scale=0.28]


				\begin{scope}[xshift=-9.5cm]
					\fill[gray!50]  (-3-1,3) -- (-3-1,1) -- (-2-1,0) -- (0-1,0) -- (1-1,1) -- (1-1,3);
					\fill[gray!50]  (-3+4-1,3) -- (-3+4-1,1) -- (-2+4-1,0) -- (0+4-1,0) -- (1+4-1,1) -- (1+4-1,3);
					\fill[gray!50]  (-3+12-1,3) -- (-3+12-1,1) -- (-2+12-1,0) -- (0+12-1,0) -- (1+12-1,1) -- (1+12-1,3);

					\fill[gray!50]  (-3-1,1) -- (-2-1,0) -- (-2-1,-2) -- (-4-1,-2) -- (-4-1,-0) -- (-3-1,1);
					\fill[gray!50]  (-3+4-1,1) -- (-2+4-1,0) -- (-2+4-1,-2) -- (-4+4-1,-2) -- (-4+4-1,-0) -- (-3+4-1,1);
					\fill[gray!50]  (-3+11,1) -- (-2+11,0) -- (-2+11,-2) -- (-4+11,-2) -- (-4+11,-0) -- (-3+11,1);
					\fill[gray!50]  (-3+16-1,1) -- (-2+16-1,0) -- (-2+16-1,-2) -- (-4+16-1,-2) -- (-4+16-1,-0) -- (-3+16-1,1);
					\fill[gray!50]  (9,0) -- (11,0) -- (11,-2) -- (10,-3) -- (9,-2) -- (9,0);

					\draw
					(-3-4-1,3) -- (-3-4-1,1) -- (-2-4-1,0) -- (-0.5-4-1,0) -- (1-4-1,1) -- (1-4-1,3);
					\draw
					(-3-1,3) -- (-3-1,1) -- (-2-1,0) -- (0-1,0) -- (1-1,1) -- (1-1,3);
					\draw
					(-2-1,0) -- (0-1,0) -- (-0-1,-2) -- (-1-1,-3) -- (-2-1,-2) -- (-2-1,0);
					\draw
					(-3-1,1) -- (-2-1,0) -- (-2-1,-2) -- (-4-1,-2) -- (-4-1,-0) -- (-3-1,1);
					\draw
					(-1-1,-5) -- (-1-1,-3) -- (0-1,-2) -- (2-1,-2) -- (3-1,-3) -- (3-1,-5);
					\draw
					(-3+4-1,3) -- (-3+4-1,1) -- (-2+4-1,0) -- (0+4-1,0) -- (1+4-1,1) -- (1+4-1,3);
					\draw
					(-2+4-1,0) -- (0+4-1,0) -- (-0+4-1,-2) -- (-1+4-1,-3) -- (-2+4-1,-2) -- (-2+4-1,0);
					\draw
					(-3+4-1,1) -- (-2+4-1,0) -- (-2+4-1,-2) -- (-4+4-1,-2) -- (-4+4-1,-0) -- (-3+4-1,1);
					\draw
					(-1+4-1,-5) -- (-1+4-1,-3) -- (0+4-1,-2) -- (2+4-1,-2) -- (3+4-1,-3) -- (3+4-1,-5);
					\draw
					(-3+11,1) -- (-2+11,0) -- (-2+11,-2) -- (-4+11,-2) -- (-4+11,-0) -- (-3+11,1);
					\draw
					(-3+12-1,3) -- (-3+12-1,1) -- (-2+12-1,0) -- (0+12-1,0) -- (1+12-1,1) -- (1+12-1,3);
					\draw
					(-2+12-1,0) -- (0+12-1,0) -- (-0+12-1,-2) -- (-1+12-1,-3) -- (-2+12-1,-2) -- (-2+12-1,0);
					\draw
					(-3+16-1,1) -- (-2+16-1,0) -- (-2+16-1,-2) -- (-4+16-1,-2) -- (-4+16-1,-0) -- (-3+16-1,1);
					\draw
					(-1+12-1,-5) -- (-1+12-1,-3) -- (0+12-1,-2) -- (2+12-1,-2);

					\draw[line width=1.5]
					(-2-4-1,0) -- (-0.5-4-1,0);
					\draw[line width=1.5]
					(-3-1,1) -- (-2-1,0);
					\draw[line width=1.5]
					(-3+4-1,1) -- (-2+4-1,0);
					\draw[line width=1.5]
					(-3+12-1,1) -- (-2+12-1,0);

					\draw[line width=1.5]
					(-1-1,-3) -- (0-1,-2);
					\draw[line width=1.5]
					(-1+4-1,-3) -- (0+4-1,-2);

					\draw[line width=1.5]
					(-4+16-1,-2) -- (-4+16-1,-0);

					\fill (4.5,-1) circle (0.1); \fill (5,-1) circle (0.1); \fill (5.5,-1) circle (0.1);

					\node at (-3-1-1,3) {\scriptsize $A_1$};
					\node at (-3-1-1,1+0.5) {\scriptsize $A_2$};
					\node at (-4-0.7-1,0-0.5) {\scriptsize $A_3$};
					\node at (-4-0.7-1,-2-0.5) {\scriptsize $A_4$};
					\node at (-2-0.7-1,-2-0.5) {\scriptsize $A_5$};
					\node at (-1-1-1,-3-0.5) {\scriptsize $A_6$};
					\node at (-1-1-1,-5) {\scriptsize $A_7$};
					\node at (11+0.7-1,-3-0.7) {\scriptsize $A_8$};
					\node at (12+0.7-1,-2-0.5) {\scriptsize $A_9$};
					\node at (14+0.7-1,-2-0.2) {\scriptsize $A_{10}$};
					\node at (14+1-1,0) {\scriptsize $A_{11}$};
					\node at (13+1-1,1+0.5) {\scriptsize $A_{12}$};

					\fill (-3-1,3) circle (0.15);
					\fill (-3-1,1) circle (0.15);
					\fill (-4-1,-0) circle (0.15);
					\fill (-4-1,-2) circle (0.15);
					\fill (-2-1,-2) circle (0.15);
					\fill (-1-1,-3) circle (0.15);
					\fill (-1-1,-5) circle (0.15);
					\fill (-1+12-1,-3) circle (0.15);
					\fill (-0+12-1,-2) circle (0.15);
					\fill (2+12-1,-2) circle (0.15);
					\fill (2+12-1,0) circle (0.15);
					\fill (-3+16-1,1) circle (0.15);

					\node[draw,shape=circle, inner sep=0.5] at (-3+2-1,3-1) {\small $2$};
					\node[draw,shape=circle, inner sep=0.5] at (-4+1-1,-1) {\small $8$};
					\node[draw,shape=circle, inner sep=0.5] at (-4+3-1,-1) {\small $6$};
					\node[draw,shape=circle, inner sep=0.5] at (-4+5-1,-1) {\small $5$};
					\node[draw,shape=circle, inner sep=0.5] at (1-1,-4) {\small $7$};
					\node[draw,shape=circle, inner sep=0.5] at (-3+6-1,3-1) {\small $3$};
					\node[draw,shape=circle, inner sep=0.5] at (-3+6+8-1,3-1) {\small $4$};
					\node[draw,shape=circle, inner sep=0.5] at (-4+3+12-1,-1) {\small $18$};
					\node[draw,shape=circle, inner sep=0.5] at (8,-1) {\small $19$};
					\node[draw,shape=circle, inner sep=0.5] at (-4+5+12-1,-1) {\small $17$};
				\end{scope}

				\begin{scope}[xshift=9.5cm]
					\fill[gray!50]  (-1,-5) -- (-1,-3) -- (0,-2) -- (2,-2) -- (3,-3) -- (3,-5);
					\fill[gray!50]  (-1+4,-5) -- (-1+4,-3) -- (0+4,-2) -- (2+4,-2) -- (3+4,-3) -- (3+4,-5);
					\fill[gray!50]  (-1+12,-5) -- (-1+12,-3) -- (0+12,-2) -- (2+12,-2) -- (3+12,-3) -- (3+12,-5);

					\fill[gray!50]  (-2+4,0) -- (0+4,0) -- (-0+4,-2) -- (-1+4,-3) -- (-2+4,-2) -- (-2+4,0);
					\fill[gray!50]  (-2+8,0) -- (0+8,0) -- (-0+8,-2) -- (-1+8,-3) -- (-2+8,-2) -- (-2+8,0);
					\fill[gray!50]  (-2+16,0) -- (0+16,0) -- (-0+16,-2) -- (-1+16,-3) -- (-2+16,-2) -- (-2+16,0);

					\draw
					(-3,3) -- (-3,1) -- (-2,0) -- (0,0) -- (1,1) -- (1,3);
					\draw
					(-2+4,0) -- (0+4,0) -- (-0+4,-2) -- (-1+4,-3) -- (-2+4,-2) -- (-2+4,0);
					\draw
					(-3+4,1) -- (-2+4,0) -- (-2+4,-2) -- (-4+4,-2) -- (-4+4,-0) -- (-3+4,1);
					\draw
					(-1,-5) -- (-1,-3) -- (0,-2) -- (2,-2) -- (3,-3) -- (3,-5);
					\draw
					(-3+4,3) -- (-3+4,1) -- (-2+4,0) -- (0+4,0) -- (1+4,1) -- (1+4,3);
					\draw
					(-2+8,0) -- (0+8,0) -- (-0+8,-2) -- (-1+8,-3) -- (-2+8,-2) -- (-2+8,0);
					\draw
					(-3+8,1) -- (-2+8,0) -- (-2+8,-2) -- (-4+8,-2) -- (-4+8,-0) -- (-3+8,1);
					\draw
					(-1+4,-5) -- (-1+4,-3) -- (0+4,-2) -- (2+4,-2) -- (3+4,-3) -- (3+4,-5);
					\draw
					(-3+12,3) -- (-3+12,1) -- (-2+12,0) -- (0+12,0) -- (1+12,1) -- (1+12,3);
					\draw
					(-2+16,0) -- (0+16,0) -- (-0+16,-2) -- (-1+16,-3) -- (-2+16,-2) -- (-2+16,0);
					\draw
					(-3+16,1) -- (-2+16,0) -- (-2+16,-2) -- (-4+16,-2) -- (-4+16,-0) -- (-3+16,1);
					\draw
					(-1+12,-5) -- (-1+12,-3) -- (0+12,-2) -- (2+12,-2) -- (3+12,-3) -- (3+12,-5);

					\draw[line width=1.5]
					(0,0) -- (1,1);
					\draw[line width=1.5]
					(0+4,0) -- (1+4,1);
					\draw[line width=1.5]
					(0+12,0) -- (1+12,1);

					\draw[line width=1.5]
					(2,-2) -- (3,-3);
					\draw[line width=1.5]
					(2+4,-2) -- (3+4,-3);
					\draw[line width=1.5]
					(2+12,-2) -- (3+12,-3);

					\fill (9.5,-1) circle (0.1); \fill (10,-1) circle (0.1); \fill (10.5,-1) circle (0.1);

					\node at (-3-1,1) {\scriptsize $A_1$};
					\node at (-3-1,3) {\scriptsize $A_2$};
					\node at (13+0.7,1+0.5) {\scriptsize $A_{3}$};
					\node at (14+0.7,0+0.5) {\scriptsize $A_{4}$};
					\node at (14+2+0.7,0+0.4) {\scriptsize $A_{5}$};
					\node at (14+2+0.7,-4+2) {\scriptsize $A_{6}$};
					\node at (14+1+0.7,-4+0.5) {\scriptsize $A_{7}$};
					\node at (14+1+0.7,-5) {\scriptsize $A_{8}$};
					\node at (-1-1,-3) {\scriptsize $A_9$};
					\node at (-3+2.1,-2+0.3) {\scriptsize $A_{10}$};
					\node at (-3+2.1,0-0.5) {\scriptsize $A_{11}$};
					\node at (-3+0.5,-0.5) {\scriptsize $A_{12}$};

					\fill (-3,1) circle (0.15);
					\fill (-3,3) circle (0.15);
					\fill (1+12,1) circle (0.15);
					\fill (-2+16,0) circle (0.15);
					\fill (0+16,0) circle (0.15);
					\fill (0+16,-2) circle (0.15);
					\fill (3+12,-3) circle (0.15);
					\fill (3+12,-5) circle (0.15);
					\fill (-1,-3) circle (0.15);
					\fill (0,-2) circle (0.15);
					\fill (0,0) circle (0.15);
					\fill (-2,0) circle (0.15);

					\node[draw,shape=circle, inner sep=0.5] at (-3+2,3-1) {\small $1$};
					\node[draw,shape=circle, inner sep=0.5] at (-4+3+4,-1) {\small $12$};
					\node[draw,shape=circle, inner sep=0.5] at (-4+5,-1) {\small $11$};
					\node[draw,shape=circle, inner sep=0.5] at (1,-4) {\small $13$};
					\node[draw,shape=circle, inner sep=0.5] at (-3+6,3-1) {\small $9$};
					\node[draw,shape=circle, inner sep=0.5] at (-3+6+8,3-1) {\small $10$};
					\node[draw,shape=circle, inner sep=0.5] at (-4+3+12+4,-1) {\small $15$};
					\node[draw,shape=circle, inner sep=0.5] at (-4+5+12,-1) {\small $14$};
					\node[draw,shape=circle, inner sep=0.5] at (1+12,-4) {\small $16$};
				\end{scope}
			\end{tikzpicture}
			\caption{$\aaa\ccc^{\frac{f+4}8}=\thin^{\bbb}\aaa^{\ccc}\thin^{\eee}\ccc^{\aaa}\thin\cdots\thin^{\eee}\ccc^{\aaa}\thin$.}
			\label{ac1}
		\end{figure}

		When $f=20$, Figure \ref{ac1} gives the second tiling with AVC $1$ in Figure \ref{various flip}.

		\vspace{6pt}

		 $\bullet$ In Figure \ref{ac2}, we have $\aaa\ccc^{\frac{f+4}8}=\thin\aaa_1\thin^{\eee}\ccc_2^{\aaa}\thin\cdots\thin^{\eee}\ccc_3^{\aaa}\thin^{\aaa}\ccc_4^{\eee}\thin$, which determines $T_1,\cdots,T_4$ ($\frac{f+12}8$ tiles). Just like in Figure \ref{ac1}, we can also determine all tiles except $T_6$, $T_7$, $T_8$, $T_9$, which form a UFO and can be flipped.

		 One tiling is $T(2\aaa^2\ccc,4\aaa\bbb\ccc
		 ,(f-10)\aaa\ddd\eee,(\frac{f}2-6)\bbb^2\ccc,8\bbb\ddd\eee,2\aaa\ccc^{\frac{f+4}8},2\ccc^{\frac{f-4}8}\ddd\eee)$.
		 Figure \ref{ac2} shows a tiling with three UFOs: \{$T_{4}$, $T_8$, $T_9$, $T_{10}$\}, \{$T_6$, $T_7$, $T_8$, $T_9$\} and \{$T_{1}$, $T_{4}$, $T_{10}$, $T_{11}$\}.

		 The other tiling has the AVC $\{2\aaa^2\ccc,6\aaa\bbb\ccc
		 ,(f-11)\aaa\ddd\eee,(\frac{f}2-7)\bbb^2\ccc,8\bbb\ddd\eee,\aaa\ccc^{\frac{f+4}8},3\ccc^{\frac{f-4}8}\ddd\eee\}$.
		 Flipping \{$T_6$, $T_7$, $T_8$, $T_9$\} in Figure \ref{ac2} gives three UFOs: \{$T_{3}$, $T_{5}$, $T_{8}$, $T_{9}$\}, \{$T_{6}$, $T_{7}$, $T_{8}$, $T_{9}$\} and \{$T_{1}$, $T_{4}$, $T_{10}$, $T_{11}$\}.

		\begin{figure}[htp]
			\centering
			\begin{tikzpicture}[>=latex,scale=0.28]


				\begin{scope}[xshift=-9.5cm]
					\fill[gray!50]  (-3-1,3) -- (-3-1,1) -- (-2-1,0) -- (0-1,0) -- (1-1,1) -- (1-1,3);
					\fill[gray!50]  (-3+8-1,3) -- (-3+8-1,1) -- (-2+8-1,0) -- (0+8-1,0) -- (1+8-1,1) -- (1+8-1,3);

					\fill[gray!50]  (-3-1,1) -- (-2-1,0) -- (-2-1,-2) -- (-4-1,-2) -- (-4-1,-0) -- (-3-1,1);
					\fill[gray!50]  (4,1) -- (5,0) -- (5,-2) -- (3,-2) -- (3,-0) -- (4,1);

					\draw
					(-3-1,3) -- (-3-1,1) -- (-2-1,0) -- (0-1,0) -- (1-1,1) -- (1-1,3);
					\draw
					(-2-1,0) -- (0-1,0) -- (-0-1,-2) -- (-1-1,-3) -- (-2-1,-2) -- (-2-1,0);
					\draw
					(-3-1,1) -- (-2-1,0) -- (-2-1,-2) -- (-4-1,-2) -- (-4-1,-0) -- (-3-1,1);
					\draw
					(-1-1,-5) -- (-1-1,-3) -- (0-1,-2) -- (2-1,-2) -- (3-1,-3) -- (3-1,-5);
					\draw
					(-2+4-1,-2) -- (-4+4-1,-2) -- (-4+4-1,-0) -- (-3+4-1,1);
					\draw
					(-3+8-1,3) -- (-3+8-1,1) -- (-2+8-1,0) -- (0+8-1,0) -- (1+8-1,1) -- (1+8-1,3);
					\draw
					(-2+8-1,0) -- (0+8-1,0) -- (-0+8-1,-2) -- (-1+8-1,-3) -- (-2+8-1,-2) -- (-2+8-1,0);
					\draw
					(-3+12-1,1) -- (-2+12-1,0) -- (-2+12-1,-2) -- (-4+12-1,-2) -- (-4+12-1,-0) -- (-3+12-1,1);
					\draw
					(-1+8-1,-5) -- (-1+8-1,-3) -- (0+8-1,-2) -- (2+8-1,-2);
					\draw
					(-3+12-1,3) -- (-3+12-1,1) -- (-2+12-1,0) -- (0+12-1,0) -- (1+12-1,1) -- (1+12-1,3);
					\draw
					(-2+12-1,0) -- (0+12-1,0) -- (-0+12-1,-2) -- (-1+12-1,-3) -- (-2+12-1,-2) -- (-2+12-1,0);
					\draw
					(-3+16-1,1) -- (-2+16-1,0) -- (-2+16-1,-2) -- (-4+16-1,-2) -- (-4+16-1,-0) -- (-3+16-1,1);
					\draw
					(-1+12-1,-5) -- (-1+12-1,-3) -- (0+12-1,-2) -- (2+12-1,-2);
					\draw
					(4,1) -- (3,0) -- (3,-2) -- (5,-2);

					\draw[line width=1.5]
					(-3-1,1) -- (-2-1,0);
					\draw[line width=1.5]
					(-3+8-1,1) -- (-2+8-1,0);

					\draw[line width=1.5]
					(-3+16-1,1) -- (-4+16-1,-0);

					\draw[line width=1.5]
					(-1-1,-3) -- (0-1,-2);
					\draw[line width=1.5]
					(-1+8-1,-3) -- (0+8-1,-2);

					\draw[line width=1.5]
					(-4+14-1,-2) -- (-4+14-1,-0);

					\fill (0.5,-1) circle (0.1); \fill (1,-1) circle (0.1); \fill (1.5,-1) circle (0.1);

					\node at (-3-1-1,3) {\scriptsize $A_1$};
					\node at (-3-1-1,1+0.5) {\scriptsize $A_2$};
					\node at (-4-0.7-1,0-0.5) {\scriptsize $A_3$};
 					\node at (-4-0.7-1,-2-0.5) {\scriptsize $A_4$};
					\node at (-2-0.7-1,-2-0.5) {\scriptsize $A_5$};
					\node at (-1-1-1,-3-0.5) {\scriptsize $A_6$};
					\node at (-1-1-1,-5) {\scriptsize $A_7$};
					\node at (11+0.7-1,-3-0.7) {\scriptsize $A_8$};
					\node at (12+0.7-1,-2-0.5) {\scriptsize $A_9$};
					\node at (14+1-1,-2-0.2) {\scriptsize $A_{10}$};
					\node at (14+1-1,0) {\scriptsize $A_{11}$};
					\node at (13+1-1,1+0.5) {\scriptsize $A_{12}$};

					\fill (-3-1,3) circle (0.15);
					\fill (-3-1,1) circle (0.15);
					\fill (-4-1,-0) circle (0.15);
					\fill (-4-1,-2) circle (0.15);
					\fill (-2-1,-2) circle (0.15);
					\fill (-1-1,-3) circle (0.15);
					\fill (-1-1,-5) circle (0.15);
					\fill (-1+12-1,-3) circle (0.15);
					\fill (-0+12-1,-2) circle (0.15);
					\fill (2+12-1,-2) circle (0.15);
					\fill (2+12-1,0) circle (0.15);
					\fill (-3+16-1,1) circle (0.15);

					\node[draw,shape=circle, inner sep=0.5] at (-3+2-1,3-1) {\small $2$};
					\node[draw,shape=circle, inner sep=0.5] at (-4+8,-1) {\small $5$};
					\node[draw,shape=circle, inner sep=0.5] at (-3+6+4-1,3-1) {\small $3$};
					\node[draw,shape=circle, inner sep=0.5] at (-3+6+8-1,3-1) {\small $4$};
					\node[draw,shape=circle, inner sep=0.5] at (-4+3+12-1,-1) {\small $9$};
					\node[draw,shape=circle, inner sep=0.5] at (-4+5+12-1,-1) {\small $10$};
					\node[draw,shape=circle, inner sep=0.1] at (-4+3+12-4-1,-1) {\small $6$};
					\node[draw,shape=circle, inner sep=0.5] at (-4+5+12-4-1,-1) {\small $8$};
					\node[draw,shape=circle, inner sep=0.5] at (1+8-1,-4) {\small $7$};
				\end{scope}

				\begin{scope}[xshift=9.5cm]
					\fill[gray!50]  (-1,-5) -- (-1,-3) -- (0,-2) -- (2,-2) -- (3,-3) -- (3,-5);
					\fill[gray!50]  (-1+4,-5) -- (-1+4,-3) -- (0+4,-2) -- (2+4,-2) -- (3+4,-3) -- (3+4,-5);
					\fill[gray!50]  (-1+12,-5) -- (-1+12,-3) -- (0+12,-2) -- (2+12,-2) -- (3+12,-3) -- (3+12,-5);

					\fill[gray!50]  (-2+4,0) -- (0+4,0) -- (-0+4,-2) -- (-1+4,-3) -- (-2+4,-2) -- (-2+4,0);
					\fill[gray!50]  (-2+8,0) -- (0+8,0) -- (-0+8,-2) -- (-1+8,-3) -- (-2+8,-2) -- (-2+8,0);
					\fill[gray!50]  (-2+16,0) -- (0+16,0) -- (-0+16,-2) -- (-1+16,-3) -- (-2+16,-2) -- (-2+16,0);

					\draw
					(-3,3) -- (-3,1) -- (-2,0) -- (0,0) -- (1,1) -- (1,3);
					\draw
					(-2+4,0) -- (0+4,0) -- (-0+4,-2) -- (-1+4,-3) -- (-2+4,-2) -- (-2+4,0);
					\draw
					(-3+4,1) -- (-2+4,0) -- (-2+4,-2) -- (-4+4,-2) -- (-4+4,-0) -- (-3+4,1);
					\draw
					(-1,-5) -- (-1,-3) -- (0,-2) -- (2,-2) -- (3,-3) -- (3,-5);
					\draw
					(-3+4,3) -- (-3+4,1) -- (-2+4,0) -- (0+4,0) -- (1+4,1) -- (1+4,3);
					\draw
					(-2+8,0) -- (0+8,0) -- (-0+8,-2) -- (-1+8,-3) -- (-2+8,-2) -- (-2+8,0);
					\draw
					(-3+8,1) -- (-2+8,0) -- (-2+8,-2) -- (-4+8,-2) -- (-4+8,-0) -- (-3+8,1);
					\draw
					(-1+4,-5) -- (-1+4,-3) -- (0+4,-2) -- (2+4,-2) -- (3+4,-3) -- (3+4,-5);
					\draw
					(-3+12,3) -- (-3+12,1) -- (-2+12,0) -- (0+12,0) -- (1+12,1) -- (1+12,3);
					\draw
					(-2+16,0) -- (0+16,0) -- (-0+16,-2) -- (-1+16,-3) -- (-2+16,-2) -- (-2+16,0);
					\draw
					(-3+16,1) -- (-2+16,0) -- (-2+16,-2) -- (-4+16,-2) -- (-4+16,-0) -- (-3+16,1);
					\draw
					(-1+12,-5) -- (-1+12,-3) -- (0+12,-2) -- (2+12,-2) -- (3+12,-3) -- (3+12,-5);

					\draw[line width=1.5]
					(0,0) -- (1,1);
					\draw[line width=1.5]
					(0+4,0) -- (1+4,1);
					\draw[line width=1.5]
					(0+12,0) -- (1+12,1);

					\draw[line width=1.5]
					(2,-2) -- (3,-3);
					\draw[line width=1.5]
					(2+4,-2) -- (3+4,-3);
					\draw[line width=1.5]
					(2+12,-2) -- (3+12,-3);

					\fill (9.5,-1) circle (0.1); \fill (10,-1) circle (0.1); \fill (10.5,-1) circle (0.1);

					\node at (-3-1,1) {\scriptsize $A_1$};
					\node at (-3-1,3) {\scriptsize $A_2$};
					\node at (13+0.7,1+0.5) {\scriptsize $A_{3}$};
					\node at (14+0.7,0+0.5) {\scriptsize $A_{4}$};
					\node at (14+2+0.7,0+0.4) {\scriptsize $A_{5}$};
					\node at (14+2+0.7,-4+2) {\scriptsize $A_{6}$};
					\node at (14+1+0.7,-4+0.5) {\scriptsize $A_{7}$};
					\node at (14+1+0.7,-5) {\scriptsize $A_{8}$};
					\node at (-1-1,-3) {\scriptsize $A_9$};
					\node at (-3+2.1,-2+0.3) {\scriptsize $A_{10}$};
					\node at (-3+2.1,0-0.5) {\scriptsize $A_{11}$};
					\node at (-3+0.3,-0.5) {\scriptsize $A_{12}$};

					\fill (-3,1) circle (0.15);
					\fill (-3,3) circle (0.15);
					\fill (1+12,1) circle (0.15);
					\fill (-2+16,0) circle (0.15);
					\fill (0+16,0) circle (0.15);
					\fill (0+16,-2) circle (0.15);
					\fill (3+12,-3) circle (0.15);
					\fill (3+12,-5) circle (0.15);
					\fill (-1,-3) circle (0.15);
					\fill (0,-2) circle (0.15);
					\fill (0,0) circle (0.15);
					\fill (-2,0) circle (0.15);

					\node[draw,shape=circle, inner sep=0.5] at (-3+2,3-1) {\small $1$};
					\node[draw,shape=circle, inner sep=0.5] at (-4+5,-1) {\small $11$};
				\end{scope}
			\end{tikzpicture}
			\caption{$\aaa\ccc^{\frac{f+4}8}=\thin^{\bbb}\aaa^{\ccc}\thin^{\eee}\ccc^{\aaa}\thin\cdots\thin^{\eee}\ccc^{\aaa}\thin^{\aaa}\ccc^{\eee}\thin$.}
			\label{ac2}
		\end{figure}

	    When $f=20$, Figure \ref{ac2} gives the tiling with AVC $3$ and the first tiling with AVC $5$ in Figure \ref{various flip}.

		\vspace{6pt}

		$\bullet$ In Figure \ref{ac3}, we have $\aaa\ccc^{\frac{f+4}8}=\thin\aaa_1\thin^{\eee}\ccc_2^{\aaa}\thin\cdots\thin^{\eee}\ccc_3^{\aaa}\thin^{\aaa}\ccc_4^{\eee}\thin\cdots\thin^{\aaa}\ccc_5^{\eee}\thin$, which determines $T_{1},T_2,\cdots,T_3,T_4,\cdots,T_5$ ($\frac{f+12}8$ tiles). Note that the position of $\ccc_4$ can be any place not next to $\aaa_1$. Then just like in Figure \ref{ac1}, we can also determine all tiles except $T_{7}$, $T_{8}$, $T_{9}$, $T_{10}$, which form a UFO and can be flipped. So we get many tilings with the same AVC $\{2\aaa^2\ccc,6\aaa\bbb\ccc
		,(f-11)\aaa\ddd\eee,(\frac{f}2-7)\bbb^2\ccc,8\bbb\ddd\eee,\aaa\ccc^{\frac{f+4}8},3\ccc^{\frac{f-4}8}\ddd\eee\}$ and $T(2\aaa^2\ccc,5\aaa\bbb\ccc
		,(f-10)\aaa\ddd\eee,(\frac{f}2-6)\bbb^2\ccc,7\bbb\ddd\eee,\aaa\ccc^{\frac{f+4}8},3\ccc^{\frac{f-4}8}\ddd\eee)$. Figure \ref{ac3} gives a tiling with three UFO: \{$T_{7}$, $T_{8}$, $T_{9}$, $T_{10}$\}, \{$T_{4}$, $T_{9}$, $T_{10}$, $T_{11}$\} and \{$T_{1}$, $T_{5}$, $T_{14}$, $T_{15}$\}. After flipping the UFO \{$T_{7}$, $T_{8}$, $T_{9}$, $T_{10}$\}, the new tiling still has three UFO and the other two UFO are: \{$T_{3}$, $T_{6}$, $T_{9}$, $T_{10}$\} and \{$T_{1}$, $T_{5}$, $T_{14}$, $T_{15}$\}.

		\begin{figure}[htp]
			\centering
			\begin{tikzpicture}[>=latex,scale=0.28]


				\begin{scope}[xshift=0cm,yshift=10cm]
					\fill[gray!50]  (-3,3) -- (-3,1) -- (-2,0) -- (0,0) -- (1,1) -- (1,3);
					\fill[gray!50]  (-3+8,3) -- (-3+8,1) -- (-2+8,0) -- (0+8,0) -- (1+8,1) -- (1+8,3);

					\fill[gray!50]  (-1+20,-5) -- (-1+20,-3) -- (0+20,-2) -- (2+20,-2) -- (3+20,-3) -- (3+20,-5);

					\fill[gray!50]  (-3,1) -- (-2,0) -- (-2,-2) -- (-4,-2) -- (-4,-0) -- (-3,1);
					\fill[gray!50]  (-3+8,1) -- (-2+8,0) -- (-2+8,-2) -- (-4+8,-2) -- (-4+8,-0) -- (-3+8,1);
					\fill[gray!50]  (-2+16,0) -- (0+16,0) -- (-0+16,-2) -- (-1+16,-3) -- (-2+16,-2) -- (-2+16,0);
					\fill[gray!50]  (-2+24,0) -- (0+24,0) -- (-0+24,-2) -- (-1+24,-3) -- (-2+24,-2) -- (-2+24,0);

					\fill[gray!50]  (-1+12,-5) -- (-1+12,-3) -- (0+12,-2) -- (2+12,-2) -- (3+12,-3) -- (3+12,-5);

					\draw
					(-3,3) -- (-3,1) -- (-2,0) -- (0,0) -- (1,1) -- (1,3);
					\draw
					(-2,0) -- (0,0) -- (-0,-2) -- (-1,-3) -- (-2,-2) -- (-2,0);
					\draw
					(-3,1) -- (-2,0) -- (-2,-2) -- (-4,-2) -- (-4,-0) -- (-3,1);
					\draw
					(-1,-5) -- (-1,-3) -- (0,-2) -- (2,-2) -- (3,-3) -- (3,-5);
					\draw
					(5,1) -- (4,0) -- (4,-2) -- (6,-2);
					\draw
					(-3+4+4,3) -- (-3+4+4,1) -- (-2+4+4,0) -- (0+4+4,0) -- (1+4+4,1) -- (1+4+4,3);
					\draw
					(-2+4+4,0) -- (0+4+4,0) -- (-0+4+4,-2) -- (-1+4+4,-3) -- (-2+4+4,-2) -- (-2+4+4,0);
					\draw
					(-1+4+4,-5) -- (-1+4+4,-3) -- (0+4+4,-2) -- (2+4+4,-2) -- (3+4+4,-3) -- (3+4+4,-5);
					\draw
					(-3+8+4,1) -- (-2+8+4,0) -- (-2+8+4,-2) -- (-4+8+4,-2) -- (-4+8+4,-0) -- (-3+8+4,1);
					\draw
					(-3+12,3) -- (-3+12,1) -- (-2+12,0) -- (0+12,0) -- (1+12,1) -- (1+12,3);
					\draw
					(-2+12,0) -- (0+12,0) -- (-0+12,-2) -- (-1+12,-3) -- (-2+12,-2) -- (-2+12,0);
					\draw
					(-3+16,1) -- (-2+16,0) -- (-2+16,-2) -- (-4+16,-2) -- (-4+16,-0) -- (-3+16,1);
					\draw
					(-1+12,-5) -- (-1+12,-3) -- (0+12,-2) -- (2+12,-2);
					\draw
					(-3+16,3) -- (-3+16,1) -- (-2+16,0) -- (0+16,0) -- (1+16,1) -- (1+16,3);
					\draw
					(-2+16,0) -- (0+16,0) -- (-0+16,-2) -- (-1+16,-3) -- (-2+16,-2) -- (-2+16,0);
					\draw
					(-3+20,1) -- (-2+20,0) -- (-2+20,-2) -- (-4+20,-2) -- (-4+20,-0) -- (-3+20,1);
					\draw
					(-1+16,-5) -- (-1+16,-3) -- (0+16,-2) -- (2+16,-2);
					\draw
					(-3+24,3) -- (-3+24,1) -- (-2+24,0) -- (0+24,0) -- (1+24,1) -- (1+24,3);
					\draw
					(-2+24,0) -- (0+24,0) -- (-0+24,-2) -- (-1+24,-3) -- (-2+24,-2) -- (-2+24,0);
					\draw
					(-3+28,1) -- (-2+28,0) -- (-2+28,-2) -- (-4+28,-2) -- (-4+28,-0) -- (-3+28,1);
					\draw
					(-1+24,-5) -- (-1+24,-3) -- (0+24,-2) -- (2+24,-2);

					\draw
					(-1+20,-5) -- (-1+20,-3) -- (0+20,-2) -- (2+20,-2);

					\draw[line width=1.5]
					(-3,1) -- (-2,0);
   				    \draw[line width=1.5]
					(-3+8,1) -- (-2+8,0);
					\draw[line width=1.5]
					(-3+16,1) -- (-4+16,-0);
					\draw[line width=1.5]
					(-3+16+4,1) -- (-4+20,-0);
					\draw[line width=1.5]
					(-3+16+4+8,1) -- (-4+20+8,-0);

					\draw[line width=1.5]
					(-1+24-8,-3) -- (-2+24-8,-2);
					\draw[line width=1.5]
					(-1+24,-3) -- (-2+24,-2);

					\draw[line width=1.5]
					(-1,-3) -- (0,-2);
					\draw[line width=1.5]
					(-1+8,-3) -- (0+8,-2);

					\draw[line width=1.5]
					(-4+18-4,-2) -- (-4+18-4,-0);

					\fill (1.5,-1) circle (0.1); \fill (2,-1) circle (0.1); \fill (2.5,-1) circle (0.1); \fill (19.5,-1) circle (0.1); \fill (20,-1) circle (0.1); \fill (20.5,-1) circle (0.1);

					\node at (-3-1,3) {\scriptsize $A_1$};
					\node at (-3-1,1+0.5) {\scriptsize $A_2$};
					\node at (-4-0.7,0-0.5) {\scriptsize $A_3$};
					\node at (-4-0.7,-2-0.5) {\scriptsize $A_4$};
					\node at (-2-0.7,-2-0.5) {\scriptsize $A_5$};
					\node at (-1-1,-3-0.5) {\scriptsize $A_6$};
					\node at (-1-1,-5) {\scriptsize $A_7$};
					\node at (11+4+0.7+8,-3-0.7) {\scriptsize $A_8$};
					\node at (12+4+0.7+8,-2-0.5) {\scriptsize $A_9$};
					\node at (14+4+1+8,-2-0.2) {\scriptsize $A_{10}$};
					\node at (14+4+1+8,0) {\scriptsize $A_{11}$};
					\node at (13+4+1+8,1+0.5) {\scriptsize $A_{12}$};

					\fill (-3,3) circle (0.15);
					\fill (-3,1) circle (0.15);
					\fill (-4,-0) circle (0.15);
					\fill (-4,-2) circle (0.15);
					\fill (-2,-2) circle (0.15);
					\fill (-1,-3) circle (0.15);
					\fill (-1,-5) circle (0.15);
					\fill (-1+12+4+8,-3) circle (0.15);
					\fill (-0+12+4+8,-2) circle (0.15);
					\fill (2+12+4+8,-2) circle (0.15);
					\fill (2+12+4+8,0) circle (0.15);
					\fill (-3+16+4+8,1) circle (0.15);

					\node[draw,shape=circle, inner sep=0.5] at (-3+2,3-1) {\small $2$};
					\node[draw,shape=circle, inner sep=0.5] at (-4+9,-1) {\small $6$};
					\node[draw,shape=circle, inner sep=0.5] at (-3+6+4,3-1) {\small $3$};
					\node[draw,shape=circle, inner sep=0.5] at (-3+6+8,3-1) {\small $4$};
					\node[draw,shape=circle, inner sep=0.5] at (-3+6+16+4,3-1) {\small $5$};
					\node[draw,shape=circle, inner sep=0.5] at (-4+3+12,-1) {\small $10$};
					\node[draw,shape=circle, inner sep=0.5] at (-4+5+12,-1) {\small $11$};
					\node[draw,shape=circle, inner sep=0.1] at (-4+3+12-4,-1) {\small $7$};
					\node[draw,shape=circle, inner sep=0.5] at (-4+5+12-4,-1) {\small $9$};
					\node[draw,shape=circle, inner sep=0.5] at (1+8,-4) {\small $8$};
					\node[draw,shape=circle, inner sep=0.5] at (17+8,-1) {\small $14$};
					\node[draw,shape=circle, inner sep=0.5] at (1+8+4,-4) {\small $13$};
					\node[draw,shape=circle, inner sep=0.1] at (-4+3+12+4,-1) {\small $12$};
				\end{scope}

				\begin{scope}[xshift=5cm,yshift=0cm]
					\fill[gray!50]  (-1,-5) -- (-1,-3) -- (0,-2) -- (2,-2) -- (3,-3) -- (3,-5);
					\fill[gray!50]  (-1+4,-5) -- (-1+4,-3) -- (0+4,-2) -- (2+4,-2) -- (3+4,-3) -- (3+4,-5);
					\fill[gray!50]  (-1+12,-5) -- (-1+12,-3) -- (0+12,-2) -- (2+12,-2) -- (3+12,-3) -- (3+12,-5);

					\fill[gray!50]  (-2+4,0) -- (0+4,0) -- (-0+4,-2) -- (-1+4,-3) -- (-2+4,-2) -- (-2+4,0);
					\fill[gray!50]  (-2+8,0) -- (0+8,0) -- (-0+8,-2) -- (-1+8,-3) -- (-2+8,-2) -- (-2+8,0);
					\fill[gray!50]  (-2+16,0) -- (0+16,0) -- (-0+16,-2) -- (-1+16,-3) -- (-2+16,-2) -- (-2+16,0);

					\draw
					(-3,3) -- (-3,1) -- (-2,0) -- (0,0) -- (1,1) -- (1,3);
					\draw
					(-2+4,0) -- (0+4,0) -- (-0+4,-2) -- (-1+4,-3) -- (-2+4,-2) -- (-2+4,0);
					\draw
					(-3+4,1) -- (-2+4,0) -- (-2+4,-2) -- (-4+4,-2) -- (-4+4,-0) -- (-3+4,1);
					\draw
					(-1,-5) -- (-1,-3) -- (0,-2) -- (2,-2) -- (3,-3) -- (3,-5);
					\draw
					(-3+4,3) -- (-3+4,1) -- (-2+4,0) -- (0+4,0) -- (1+4,1) -- (1+4,3);
					\draw
					(-2+8,0) -- (0+8,0) -- (-0+8,-2) -- (-1+8,-3) -- (-2+8,-2) -- (-2+8,0);
					\draw
					(-3+8,1) -- (-2+8,0) -- (-2+8,-2) -- (-4+8,-2) -- (-4+8,-0) -- (-3+8,1);
					\draw
					(-1+4,-5) -- (-1+4,-3) -- (0+4,-2) -- (2+4,-2) -- (3+4,-3) -- (3+4,-5);
					\draw
					(-3+12,3) -- (-3+12,1) -- (-2+12,0) -- (0+12,0) -- (1+12,1) -- (1+12,3);
					\draw
					(-2+16,0) -- (0+16,0) -- (-0+16,-2) -- (-1+16,-3) -- (-2+16,-2) -- (-2+16,0);
					\draw
					(-3+16,1) -- (-2+16,0) -- (-2+16,-2) -- (-4+16,-2) -- (-4+16,-0) -- (-3+16,1);
					\draw
					(-1+12,-5) -- (-1+12,-3) -- (0+12,-2) -- (2+12,-2) -- (3+12,-3) -- (3+12,-5);

					\draw[line width=1.5]
					(0,0) -- (1,1);
					\draw[line width=1.5]
					(0+4,0) -- (1+4,1);
					\draw[line width=1.5]
					(0+12,0) -- (1+12,1);

					\draw[line width=1.5]
					(2,-2) -- (3,-3);
					\draw[line width=1.5]
					(2+4,-2) -- (3+4,-3);
					\draw[line width=1.5]
					(2+12,-2) -- (3+12,-3);

					\fill (9.5,-1) circle (0.1); \fill (10,-1) circle (0.1); \fill (10.5,-1) circle (0.1);

					\node at (-3-1,1) {\scriptsize $A_1$};
					\node at (-3-1,3) {\scriptsize $A_2$};
					\node at (13+0.7,1+0.5) {\scriptsize $A_{3}$};
					\node at (14+0.7,0+0.5) {\scriptsize $A_{4}$};
					\node at (14+2+0.7,0+0.4) {\scriptsize $A_{5}$};
					\node at (14+2+0.7,-4+2) {\scriptsize $A_{6}$};
					\node at (14+1+0.7,-4+0.5) {\scriptsize $A_{7}$};
					\node at (14+1+0.7,-5) {\scriptsize $A_{8}$};
					\node at (-1-1,-3) {\scriptsize $A_9$};
					\node at (-3+2.1,-2+0.3) {\scriptsize $A_{10}$};
					\node at (-3+2.1,0-0.5) {\scriptsize $A_{11}$};
					\node at (-3+0.3,-0.5) {\scriptsize $A_{12}$};

					\fill (-3,1) circle (0.15);
					\fill (-3,3) circle (0.15);
					\fill (1+12,1) circle (0.15);
					\fill (-2+16,0) circle (0.15);
					\fill (0+16,0) circle (0.15);
					\fill (0+16,-2) circle (0.15);
					\fill (3+12,-3) circle (0.15);
					\fill (3+12,-5) circle (0.15);
					\fill (-1,-3) circle (0.15);
					\fill (0,-2) circle (0.15);
					\fill (0,0) circle (0.15);
					\fill (-2,0) circle (0.15);

					\node[draw,shape=circle, inner sep=0.5] at (-3+2,3-1) {\small $1$};
					\node[draw,shape=circle, inner sep=0.5] at (1,-1
					) {\small $15$};
				\end{scope}
			\end{tikzpicture}
			\caption{$\aaa\ccc^{\frac{f+4}8}=\thin^{\bbb}\aaa^{\ccc}\thin^{\eee}\ccc^{\aaa}\thin\cdots\thin^{\eee}\ccc^{\aaa}\thin^{\aaa}\ccc^{\eee}\thin\cdots\thin^{\aaa}\ccc^{\eee}\thin$.}
			\label{ac3}
		\end{figure}

	When $f=20$, Figure \ref{ac3} gives the tiling with AVC $4$ and the second tiling with AVC $5$ in Figure \ref{various flip}.

	\vspace{6pt}

	 $\bullet$ In Figure  \ref{ac4}, we have $\aaa\ccc^{\frac{f+4}8}=\thin\aaa_1\thin^{\aaa}\ccc_2^{\eee}\thin^{\aaa}\ccc_3^{\eee}\thin\cdots\thin^{\aaa}\ccc_4^{\eee}\thin$, which determines $T_2,T_3,\cdots,T_4$ ($\frac{f+4}8$ tiles). Just like in Figure \ref{ac1}, we can also determine $T_5, \cdots, T_{9}$. We have $\aaa_2\ccc_1\cdots=\aaa^2\ccc,\aaa\bbb\ccc$ or $\aaa\ccc^{\frac{f+4}{8}}$.

	 When $\aaa_2\ccc_1\cdots=\aaa^2\ccc$, if $\aaa\aaa_2\ccc_1=\thin^{\bbb}\aaa^{\ccc}\thin^{\eee}\ccc_1^{\aaa}\thin^{\ccc}\aaa_2^{\bbb_2}\thin$, this gives a vertex $\bbb_2^{\aaa}\thin^{\aaa}\bbb=\thin^{\ddd}\bbb^{\aaa}\thin^{\aaa}\bbb^{\ddd}\thin^{\aaa}\ccc^{\eee}\thin$. This gives a vertex $\ddd\thin\eee\cdots$, contradicting the AVC. Therefore $\aaa\aaa_2\ccc_1=\thin^{\ccc}\aaa^{\bbb}\thin^{\eee}\ccc_1^{\aaa}\thin^{\ccc}\aaa_2^{\bbb_2}\thin$. This gives a vertex $\bbb_2^{\aaa}\thin^{\aaa}\ccc\cdots=\thin\aaa\thin^{\ddd}\bbb^{\aaa}\thin^{\aaa}\ccc^{\eee}\thin$. Then we have the AAD $\ddd_2\eee_5\cdots=\thin^{\ddd}\bbb^{\aaa}\thin^{\bbb}\ddd^{\eee}\thick^{\ddd}\eee^{\ccc}\thin$ or $\ccc^{\frac{f-4}8}\ddd\eee=\thin^{\eee}\ccc^{\aaa}\thin\cdots\thin^{\eee}\ccc^{\aaa}\thin^{\bbb}\ddd^{\eee}\thick^{\ddd}\eee^{\ccc}\thin$. Hence we get $\bbb_7\ccc_5\cdots=\bbb\ccc\ddd\cdots$ or $\bbb\ccc\eee\cdots$, contradicting the AVC.

		\begin{figure}[htp]
			\centering
			\begin{tikzpicture}[>=latex,scale=0.28]


				\begin{scope}[xshift=0cm]
					\fill[gray!50]  (-1+4,-5) -- (-1+4,-3) -- (0+4,-2) -- (2+4,-2) -- (3+4,-3) -- (3+4,-5);
					\fill[gray!50]  (-1,-5) -- (-1,-3) -- (0,-2) -- (2,-2) -- (3,-3) -- (3,-5);
					\fill[gray!50]  (-2+4,0) -- (0+4,0) -- (-0+4,-2) -- (-1+4,-3) -- (-2+4,-2) -- (-2+4,0);
					\fill[gray!50]  (-2+8,0) -- (0+8,0) -- (-0+8,-2) -- (-1+8,-3) -- (-2+8,-2) -- (-2+8,0);

					\draw
					(-3-4,3) -- (-3-4,1) -- (-2-4,0) -- (-0.5-4,0) -- (1-4,1) -- (1-4,3);
					\draw
					(-3,3) -- (-3,1) -- (-2,0) -- (0,0) -- (1,1) -- (1,3);
					\draw
					(-1,-5) -- (-1,-3) -- (0,-2) -- (2,-2) -- (3,-3) -- (3,-5);
					\draw
					(-3+4,3) -- (-3+4,1) -- (-2+4,0) -- (0+4,0) -- (1+4,1) -- (1+4,3);
					\draw
					(-2+4,0) -- (0+4,0) -- (-0+4,-2) -- (-1+4,-3) -- (-2+4,-2) -- (-2+4,0);
					\draw
					(-3+4,1) -- (-2+4,0) -- (-2+4,-2) -- (-4+4,-2) -- (-4+4,-0) -- (-3+4,1);
					\draw
					(-1+4,-5) -- (-1+4,-3) -- (0+4,-2) -- (2+4,-2) -- (3+4,-3) -- (3+4,-5);
					\draw
					(-3+8,1) -- (-2+8,0) -- (-2+8,-2) -- (-4+8,-2) -- (-4+8,-0) -- (-3+8,1);
					\draw
					(-2+8,0) -- (0+8,0) -- (-0+8,-2) -- (-1+8,-3) -- (-2+8,-2) -- (-2+8,0);
					\draw
					(-3+12,3) -- (-3+12,1) -- (-2+12,0) -- (0+12,0) -- (1+12,1) -- (1+12,3);
					\draw
					(-3+16,1) -- (-2+16,0) -- (-2+16,-2) -- (-4+16,-2) -- (-4+16,-0) -- (-3+16,1);
					\draw
					(-1+12,-5) -- (-1+12,-3) -- (0+12,-2) -- (2+12,-2);
					\draw
					(-2+12+4,0) -- (0+12+4,0) -- (-0+12+4,-2) -- (-1+12+4,-3) -- (-2+12+4,-2) -- (-2+12+4,0);

					\draw[line width=1.5]
					(-3+4,1) -- (-2+2,0);
					\draw[line width=1.5]
					(-3+4+4,1) -- (-2+4+2,0);
					\draw[line width=1.5]
					(-3+12+4,1) -- (-2+12+2,0);
					\draw[line width=1.5]
					(-2-4,0) -- (-0.5-4,0);

					\draw[line width=1.5]
					(-1+4,-3) -- (0+2,-2);
					\draw[line width=1.5]
					(-1+4+4,-3) -- (0+4+2,-2);
					\draw[line width=1.5]
					(-2+12+4,0) -- (0+12+4,0);

					\fill (9.5,-1) circle (0.1); \fill (10,-1) circle (0.1); \fill (10.5,-1) circle (0.1);

					\fill (-3,1) circle (0.15);
					\fill (-2,0) circle (0.15);
					\fill (0,0) circle (0.15);
					\fill (0,-2) circle (0.15);

					\node[draw,shape=circle, inner sep=0.5] at (-3+2-4,3-1) {\small $1$};
					\node[draw,shape=circle, inner sep=0.5] at (-3+2,3-1) {\small $2$};
					\node[draw,shape=circle, inner sep=0.5] at (-4+3+4,-1) {\small $6$};
					\node[draw,shape=circle, inner sep=0.5] at (-4+5,-1) {\small $5$};
					\node[draw,shape=circle, inner sep=0.5] at (1,-4) {\small $7$};
					\node[draw,shape=circle, inner sep=0.5] at (-3+6,3-1) {\small $3$};
					\node[draw,shape=circle, inner sep=0.5] at (-3+6+8,3-1) {\small $4$};
					\node[draw,shape=circle, inner sep=0.5] at (-4+3+12+4,-1) {\scriptsize $9$};
					\node[draw,shape=circle, inner sep=0.5] at (-4+5+12,-1) {\small $8$};
				\end{scope}
			\end{tikzpicture}
			\caption{$\aaa\ccc^{\frac{f+4}8}=\thin^{\bbb}\aaa^{\ccc}\thin^{\aaa}\ccc^{\eee}\thin\cdots\thin^{\aaa}\ccc^{\eee}\thin$ and $\aaa_2\ccc_1\cdots=\aaa^2\ccc$.}
			\label{ac4}
		\end{figure}

		When $\aaa_2\ccc_1\cdots=\aaa\bbb\ccc$, just like in Figure \ref{ac1}, we can also determine all tiles in Figure \ref{ac5} which gives a tiling with the AVC \{$4\aaa\bbb\ccc,(f-6)\aaa\ddd\eee,(\frac{f}2-4)\bbb^2\ccc,4\bbb\ddd\eee,2\aaa\ccc^{\frac{f+4}8},2\ccc^{\frac{f-4}8}\ddd\eee$\}.
		It has one UFO: \{$T_{1}$, $T_{4}$, $T_{8}$, $T_{9}$\}.

		\begin{figure}[htp]
			\centering
			\begin{tikzpicture}[>=latex,scale=0.28]


				\begin{scope}[xshift=-12cm]
					\fill[gray!50]  (-1+4,-5) -- (-1+4,-3) -- (0+4,-2) -- (2+4,-2) -- (3+4,-3) -- (3+4,-5);
					\fill[gray!50]  (-1,-5) -- (-1,-3) -- (0,-2) -- (2,-2) -- (3,-3) -- (3,-5);

					\fill[gray!50]  (-2+4,0) -- (0+4,0) -- (-0+4,-2) -- (-1+4,-3) -- (-2+4,-2) -- (-2+4,0);
					\fill[gray!50]  (-2+8,0) -- (0+8,0) -- (-0+8,-2) -- (-1+8,-3) -- (-2+8,-2) -- (-2+8,0);
					\draw
					(-3,3) -- (-3,1) -- (-2,0) -- (0,0) -- (1,1) -- (1,3);
					\draw
					(-1,-5) -- (-1,-3) -- (0,-2) -- (2,-2) -- (3,-3) -- (3,-5);
					\draw
					(-3+4,3) -- (-3+4,1) -- (-2+4,0) -- (0+4,0) -- (1+4,1) -- (1+4,3);
					\draw
					(-2+4,0) -- (0+4,0) -- (-0+4,-2) -- (-1+4,-3) -- (-2+4,-2) -- (-2+4,0);
					\draw
					(-3+4,1) -- (-2+4,0) -- (-2+4,-2) -- (-4+4,-2) -- (-4+4,-0) -- (-3+4,1);
					\draw
					(-1+4,-5) -- (-1+4,-3) -- (0+4,-2) -- (2+4,-2) -- (3+4,-3) -- (3+4,-5);
					\draw
					(-3+8,1) -- (-2+8,0) -- (-2+8,-2) -- (-4+8,-2) -- (-4+8,-0) -- (-3+8,1);
					\draw
					(-2+8,0) -- (0+8,0) -- (-0+8,-2) -- (-1+8,-3) -- (-2+8,-2) -- (-2+8,0);
					\draw
					(-3+12,3) -- (-3+12,1) -- (-2+12,0) -- (0+12,0) -- (1+12,1) -- (1+12,3);
					\draw
					(-3+16,1) -- (-2+16,0) -- (-2+16,-2) -- (-4+16,-2) -- (-4+16,-0) -- (-3+16,1);
					\draw
					(-1+12,-5) -- (-1+12,-3) -- (0+12,-2) -- (2+12,-2);

					\draw[line width=1.5]
					(-3+4,1) -- (-2+2,0);
					\draw[line width=1.5]
					(-3+4+4,1) -- (-2+4+2,0);
					\draw[line width=1.5]
					(-3+12+4,1) -- (-2+12+2,0);

					\draw[line width=1.5]
					(-1+4,-3) -- (0+2,-2);
					\draw[line width=1.5]
					(-1+4+4,-3) -- (0+4+2,-2);

					\fill (9.5,-1) circle (0.1); \fill (10,-1) circle (0.1); \fill (10.5,-1) circle (0.1);

					\node at (-3-1,3) {\scriptsize $A_1$};
					\node at (-3-1,1+0.5) {\scriptsize $A_2$};
					\node at (-4-0.7+2,0-0.5) {\scriptsize $A_3$};
					\node at (-4-0.7+4,0-0.5) {\scriptsize $A_4$};
					\node at (-2-0.7+2,-1-0.5) {\scriptsize $A_5$};
					\node at (-1-1,-3-0.5) {\scriptsize $A_6$};
					\node at (-1-1,-5) {\scriptsize $A_7$};
					\node at (11+0.7,-3-0.7) {\scriptsize $A_8$};
					\node at (12+0.7,-2-0.5) {\scriptsize $A_9$};
					\node at (14+1,-2-0.2) {\scriptsize $A_{10}$};
					\node at (14+1,0) {\scriptsize $A_{11}$};
					\node at (13+1,1+0.5) {\scriptsize $A_{12}$};

					\fill (-3,3) circle (0.15);
					\fill (-3,1) circle (0.15);
					\fill (-4+2,-0) circle (0.15);
					\fill (-4+4,-0) circle (0.15);
					\fill (-2+2,-2) circle (0.15);
					\fill (-1,-3) circle (0.15);
					\fill (-1,-5) circle (0.15);
					\fill (-1+12,-3) circle (0.15);
					\fill (-0+12,-2) circle (0.15);
					\fill (2+12,-2) circle (0.15);
					\fill (2+12,0) circle (0.15);
					\fill (-3+16,1) circle (0.15);

					\node[draw,shape=circle, inner sep=0.5] at (-3+2,3-1) {\small $2$};
					\node[draw,shape=circle, inner sep=0.5] at (-4+3+4,-1) {\small $6$};
					\node[draw,shape=circle, inner sep=0.5] at (-4+5,-1) {\small $5$};
					\node[draw,shape=circle, inner sep=0.5] at (1,-4) {\small $7$};
					\node[draw,shape=circle, inner sep=0.5] at (-3+6,3-1) {\small $3$};
					\node[draw,shape=circle, inner sep=0.5] at (-3+6+8,3-1) {\small $4$};
					\node[draw,shape=circle, inner sep=0.5] at (-4+5+12,-1) {\small $8$};
				\end{scope}

				\begin{scope}[xshift=12cm]
					\fill[gray!50]  (-3,3) -- (-3,1) -- (-2,0) -- (0,0) -- (1,1) -- (1,3);
					\fill[gray!50]  (-3+4,3) -- (-3+4,1) -- (-2+4,0) -- (0+4,0) -- (1+4,1) -- (1+4,3);
					\fill[gray!50]  (-3+12-1,3) -- (-3+12-1,1) -- (-2+12-1,0) -- (0+12-1,0) -- (1+12-1,1) -- (1+12-1,3);

					\fill[gray!50]  (-3,1) -- (-2,0) -- (-2,-2) -- (-4,-2) -- (-4,-0) -- (-3,1);
					\fill[gray!50]  (-3+4,1) -- (-2+4,0) -- (-2+4,-2) -- (-4+4,-2) -- (-4+4,-0) -- (-3+4,1);
					\fill[gray!50]  (-3+8,1) -- (-2+8,0) -- (-2+8,-2) -- (-4+8,-2) -- (-4+8,-0) -- (-3+8,1);

					\draw
					(-2-4,0) -- (0-4,0) -- (-0-4,-2) -- (-1-4,-3) -- (-2-4,-2) -- (-2-4,0);
					\draw
					(-3,1) -- (-2,0) -- (-2,-2) -- (-4,-2) -- (-4,-0) -- (-3,1);
					\draw
					(-1-4,-5) -- (-1-4,-3) -- (0-4,-2) -- (2-4,-2) -- (3-4,-3) -- (3-4,-5);
					\draw
					(-3,3) -- (-3,1) -- (-2,0) -- (0,0) -- (1,1) -- (1,3);
					\draw
					(-2+4,0) -- (0+4,0) -- (-0+4,-2) -- (-1+4,-3) -- (-2+4,-2) -- (-2+4,0);
					\draw
					(-3+4,1) -- (-2+4,0) -- (-2+4,-2) -- (-4+4,-2) -- (-4+4,-0) -- (-3+4,1);
					\draw
					(-1,-5) -- (-1,-3) -- (0,-2) -- (2,-2) -- (3,-3) -- (3,-5);
					\draw
					(-3+4,3) -- (-3+4,1) -- (-2+4,0) -- (0+4,0) -- (1+4,1) -- (1+4,3);
					\draw
					(-3+8,1) -- (-2+8,0) -- (-2+8,-2) -- (-4+8,-2) -- (-4+8,-0) -- (-3+8,1);
					\draw
					(-1+4,-5) -- (-1+4,-3) -- (0+4,-2) -- (2+4,-2) -- (3+4,-3) -- (3+4,-5);
					\draw
					(-3+12-1,3) -- (-3+12-1,1) -- (-2+12-1,0) -- (0+12-1,0) -- (1+12-1,1) -- (1+12-1,3);
					\draw
					(-2+16-4-1,0) -- (0+16-4-1,0) -- (-0+16-4-1,-2) -- (-1+16-4-1,-3) -- (-2+16-4-1,-2) -- (-2+16-4-1,0);
					\draw
					(-1+12-1,-5) -- (-1+12-1,-3) -- (0+12-1,-2) -- (2+12-1,-2) -- (3+12-1,-3) -- (3+12-1,-5);

					\draw[line width=1.5]
					(-3,1) -- (-2,0);
					\draw[line width=1.5]
					(-3+4,1) -- (-2+4,0);
					\draw[line width=1.5]
					(-3+8,1) -- (-2+8,0);
					\draw[line width=1.5]
					(-3+12-1,1) -- (-2+12-1,0);

					\draw[line width=1.5]
					(-1,-3) -- (0,-2);
					\draw[line width=1.5]
					(-1+4,-3) -- (0+4,-2);
					\draw[line width=1.5]
					(-1-4,-3) -- (0-4,-2);
					\draw[line width=1.5]
					(-1+4+8-1,-3) -- (0+4+8-1,-2);

					\fill (7.2,-1) circle (0.1); \fill (7.7,-1) circle (0.1); \fill (8.2,-1) circle (0.1);

					\node at (-3-1,1+0.3) {\scriptsize $A_3$};
					\node at (-3-1,3) {\scriptsize $A_4$};
					\node at (13+0.7-1,1+0.5) {\scriptsize $A_{5}$};
					\node at (14+0.7-2-1,0-0.4) {\scriptsize $A_{6}$};
					\node at (14+0.7-2-1,0+0.5-2) {\scriptsize $A_{7}$};
					\node at (14+2+0.7-2-1,-4+2+0.3) {\scriptsize $A_{8}$};
					\node at (14+1+1-1,-4+0.5+0.3) {\scriptsize $A_{9}$};
					\node at (14+1+1-1,-5) {\scriptsize $A_{10}$};
					\node at (-1-1-4,-3-0.2) {\scriptsize $A_{11}$};
					\node at (-3-4,-0.5-2+0.5) {\scriptsize $A_{12}$};
					\node at (-3-4,0.5) {\scriptsize $A_{1}$};
					\node at (-3+0.3-2,0.5) {\scriptsize $A_{2}$};

					\fill (-3,3-2) circle (0.15);
					\fill (-3,3) circle (0.15);
					\fill (12,1) circle (0.15);
					\fill (-2+16-3,0) circle (0.15);
					\fill (-2+16-2-1,-2) circle (0.15);
					\fill (0+16-2-1,-2) circle (0.15);
					\fill (3+12-1,-3) circle (0.15);
					\fill (3+12-1,-5) circle (0.15);
					\fill (-1-4,-3) circle (0.15);
					\fill (0-6,-2) circle (0.15);
					\fill (0-4,0) circle (0.15);
					\fill (-2-4,0) circle (0.15);

					\node[draw,shape=circle, inner sep=0.5] at (-2-3,-1) {\small $1$};
					\node[draw,shape=circle, inner sep=0.5] at (-3,-4) {\small $9$};
				\end{scope}
			\end{tikzpicture}
			\caption{$\aaa\ccc^{\frac{f+4}8}=\thin^{\bbb}\aaa^{\ccc}\thin^{\aaa}\ccc^{\eee}\thin\cdots\thin^{\aaa}\ccc^{\eee}\thin$ and  $\aaa_2\ccc_1\cdots=\aaa\bbb\ccc$.}
			\label{ac5}
		\end{figure}

		When $f=20$, Figure \ref{ac5} gives the first tiling with AVC $1$ in Figure \ref{various flip}.

		When $\aaa_2\ccc_1\cdots=\aaa\ccc^{\frac{f+4}8}$, just like in Figure \ref{ac1}, we can also determine all tiles in Figure \ref{cde5} which gives a tiling with the AVC \{$4\aaa\bbb\ccc,(f-6)\aaa\ddd\eee,(\frac{f}2-4)\bbb^2\ccc,4\bbb\ddd\eee,2\aaa\ccc^{\frac{f+4}8},2\ccc^{\frac{f-4}8}\ddd\eee$\}.
		It has two UFOs: \{$T_{1}$, $T_{4}$, $T_{8}$, $T_{9}$\} and \{$T_{2}$, $T_{5}$, $T_{10}$, $T_{11}$\}.

		\begin{figure}[htp]
			\centering
			\begin{tikzpicture}[>=latex,scale=0.28]


				\begin{scope}[xshift=-10.5cm]
					\fill[gray!50]  (-1+4,-5) -- (-1+4,-3) -- (0+4,-2) -- (2+4,-2) -- (3+4,-3) -- (3+4,-5);
					\fill[gray!50]  (-1,-5) -- (-1,-3) -- (0,-2) -- (2,-2) -- (3,-3) -- (3,-5);

					\fill[gray!50]  (-2+4,0) -- (0+4,0) -- (-0+4,-2) -- (-1+4,-3) -- (-2+4,-2) -- (-2+4,0);
					\fill[gray!50]  (-2+8,0) -- (0+8,0) -- (-0+8,-2) -- (-1+8,-3) -- (-2+8,-2) -- (-2+8,0);

					\draw
					(-3,3) -- (-3,1) -- (-2,0) -- (0,0) -- (1,1) -- (1,3);
					\draw
					(-1,-5) -- (-1,-3) -- (0,-2) -- (2,-2) -- (3,-3) -- (3,-5);
					\draw
					(-3+4,3) -- (-3+4,1) -- (-2+4,0) -- (0+4,0) -- (1+4,1) -- (1+4,3);
					\draw
					(-2+4,0) -- (0+4,0) -- (-0+4,-2) -- (-1+4,-3) -- (-2+4,-2) -- (-2+4,0);
					\draw
					(-3+4,1) -- (-2+4,0) -- (-2+4,-2) -- (-4+4,-2) -- (-4+4,-0) -- (-3+4,1);
					\draw
					(-1+4,-5) -- (-1+4,-3) -- (0+4,-2) -- (2+4,-2) -- (3+4,-3) -- (3+4,-5);
					\draw
					(-3+8,1) -- (-2+8,0) -- (-2+8,-2) -- (-4+8,-2) -- (-4+8,-0) -- (-3+8,1);
					\draw
					(-2+8,0) -- (0+8,0) -- (-0+8,-2) -- (-1+8,-3) -- (-2+8,-2) -- (-2+8,0);
					\draw
					(-3+12,3) -- (-3+12,1) -- (-2+12,0) -- (0+12,0) -- (1+12,1) -- (1+12,3);
					\draw
					(-3+16,1) -- (-2+16,0) -- (-2+16,-2) -- (-4+16,-2) -- (-4+16,-0) -- (-3+16,1);
					\draw
					(-1+12,-5) -- (-1+12,-3) -- (0+12,-2) -- (2+12,-2);

					\draw[line width=1.5]
					(-3+4,1) -- (-2+2,0);
					\draw[line width=1.5]
					(-3+4+4,1) -- (-2+4+2,0);
					\draw[line width=1.5]
					(-3+12+4,1) -- (-2+12+2,0);

					\draw[line width=1.5]
					(-1+4,-3) -- (0+2,-2);
					\draw[line width=1.5]
					(-1+4+4,-3) -- (0+4+2,-2);

					\fill (9.5,-1) circle (0.1); \fill (10,-1) circle (0.1); \fill (10.5,-1) circle (0.1);

					\node at (-3-1,3) {\scriptsize $A_1$};
					\node at (-3-1,1+0.5) {\scriptsize $A_2$};
					\node at (-4-0.7+2,0-0.5) {\scriptsize $A_3$};
					\node at (-4-0.7+4,0-0.5) {\scriptsize $A_4$};
					\node at (-2-0.7+2,-1-0.5) {\scriptsize $A_5$};
					\node at (-1-1,-3-0.5) {\scriptsize $A_6$};
					\node at (-1-1,-5) {\scriptsize $A_7$};
					\node at (11+0.7,-3-0.7) {\scriptsize $A_8$};
					\node at (12+0.7,-2-0.5) {\scriptsize $A_9$};
					\node at (14+1,-2-0.2) {\scriptsize $A_{10}$};
					\node at (14+1,0) {\scriptsize $A_{11}$};
					\node at (13+1,1+0.5) {\scriptsize $A_{12}$};

					\fill (-3,3) circle (0.15);
					\fill (-3,1) circle (0.15);
					\fill (-4+2,-0) circle (0.15);
					\fill (-4+4,-0) circle (0.15);
					\fill (-2+2,-2) circle (0.15);
					\fill (-1,-3) circle (0.15);
					\fill (-1,-5) circle (0.15);
					\fill (-1+12,-3) circle (0.15);
					\fill (-0+12,-2) circle (0.15);
					\fill (2+12,-2) circle (0.15);
					\fill (2+12,0) circle (0.15);
					\fill (-3+16,1) circle (0.15);

					\node[draw,shape=circle, inner sep=0.5] at (-3+2,3-1) {\small $2$};
					\node[draw,shape=circle, inner sep=0.5] at (-4+3+4,-1) {\small $6$};
					\node[draw,shape=circle, inner sep=0.5] at (-4+5,-1) {\small $5$};
					\node[draw,shape=circle, inner sep=0.5] at (1,-4) {\small $7$};
					\node[draw,shape=circle, inner sep=0.5] at (-3+6,3-1) {\small $3$};
					\node[draw,shape=circle, inner sep=0.5] at (-3+6+8,3-1) {\small $4$};
					\node[draw,shape=circle, inner sep=0.5] at (-4+5+12,-1) {\small $8$};
				\end{scope}

				\begin{scope}[xshift=10.5cm]
					\fill[gray!50]  (-1,-5) -- (-1,-3) -- (0,-2) -- (2,-2) -- (3,-3) -- (3,-5);
					\fill[gray!50]  (-1+4,-5) -- (-1+4,-3) -- (0+4,-2) -- (2+4,-2) -- (3+4,-3) -- (3+4,-5);

					\fill[gray!50]  (-2+4,0) -- (0+4,0) -- (-0+4,-2) -- (-1+4,-3) -- (-2+4,-2) -- (-2+4,0);
					\fill[gray!50]  (-2+8,0) -- (0+8,0) -- (-0+8,-2) -- (-1+8,-3) -- (-2+8,-2) -- (-2+8,0);

					\draw
					(-3,3) -- (-3,1) -- (-2,0) -- (0,0) -- (1,1) -- (1,3);
					\draw
					(-2+4,0) -- (0+4,0) -- (-0+4,-2) -- (-1+4,-3) -- (-2+4,-2) -- (-2+4,0);
					\draw
					(-3+4,1) -- (-2+4,0) -- (-2+4,-2) -- (-4+4,-2) -- (-4+4,-0) -- (-3+4,1);
					\draw
					(-1,-5) -- (-1,-3) -- (0,-2) -- (2,-2) -- (3,-3) -- (3,-5);
					\draw
					(-3+4,3) -- (-3+4,1) -- (-2+4,0) -- (0+4,0) -- (1+4,1) -- (1+4,3);
					\draw
					(-2+8,0) -- (0+8,0) -- (-0+8,-2) -- (-1+8,-3) -- (-2+8,-2) -- (-2+8,0);
					\draw
					(-3+8,1) -- (-2+8,0) -- (-2+8,-2) -- (-4+8,-2) -- (-4+8,-0) -- (-3+8,1);
					\draw
					(-1+4,-5) -- (-1+4,-3) -- (0+4,-2) -- (2+4,-2) -- (3+4,-3) -- (3+4,-5);
					\draw
					(-3+12,3) -- (-3+12,1) -- (-2+12,0) -- (0+12,0) -- (1+12,1) -- (1+12,3);
					\draw
					(-3+16,1) -- (-2+16,0) -- (-2+16,-2) -- (-4+16,-2) -- (-4+16,-0);
					\draw
					(-1+12,-5) -- (-1+12,-3) -- (0+12,-2) -- (2+12,-2) ;

					\draw[line width=1.5]
					(0,0) -- (1,1);
					\draw[line width=1.5]
					(0+4,0) -- (1+4,1);
					\draw[line width=1.5]
					(0+12,0) -- (1+12,1);

					\draw[line width=1.5]
					(2,-2) -- (3,-3);
					\draw[line width=1.5]
					(2+4,-2) -- (3+4,-3);

					\fill (9.5,-1) circle (0.1); \fill (10,-1) circle (0.1); \fill (10.5,-1) circle (0.1);

					\node at (-3-1,1) {\scriptsize $A_1$};
					\node at (-3-1,3) {\scriptsize $A_2$};
					\node at (13+0.7,1+0.5) {\scriptsize $A_{3}$};
					\node at (14+0.7,0+0.5) {\scriptsize $A_{4}$};
					\node at (14+0.7,-2) {\scriptsize $A_{5}$};
					\node at (12.5,-2.7) {\scriptsize $A_{6}$};
					\node at (11+0.7,-4+0.5) {\scriptsize $A_{7}$};
					\node at (11+0.7,-5) {\scriptsize $A_{8}$};
					\node at (-1-1,-3) {\scriptsize $A_9$};
					\node at (-3+2.1,-2+0.3) {\scriptsize $A_{10}$};
					\node at (-3+2.1,0-0.5) {\scriptsize $A_{11}$};
					\node at (-3+0.5,-0.5) {\scriptsize $A_{12}$};

					\fill (-3,1) circle (0.15);
					\fill (-3,3) circle (0.15);
					\fill (1+12,1) circle (0.15);
					\fill (14,0) circle (0.15);
					\fill (0+14,-2) circle (0.15);
					\fill (0+12,-2) circle (0.15);
					\fill (11,-3) circle (0.15);
					\fill (11,-5) circle (0.15);
					\fill (-1,-3) circle (0.15);
					\fill (0,-2) circle (0.15);
					\fill (0,0) circle (0.15);
					\fill (-2,0) circle (0.15);

					\node[draw,shape=circle, inner sep=0.5] at (-3+2,3-1) {\small $1$};
					\node[draw,shape=circle, inner sep=0.5] at (-3+2+2,3-1-3) {\small $9$};
					\node[draw,shape=circle, inner sep=0.5] at (-3+2+12,2) {\small $10$};
					\node[draw,shape=circle, inner sep=0.5] at (-3+2+14,-1) {\small $11$};
				\end{scope}
			\end{tikzpicture}
			\caption{$\aaa\ccc^{\frac{f+4}8}=\thin^{\bbb}\aaa^{\ccc}\thin^{\aaa}\ccc^{\eee}\thin\cdots\thin^{\aaa}\ccc^{\eee}\thin$ and $\aaa_2\ccc_1\cdots=\aaa\ccc^{\frac{f+4}8}$.}
			\label{cde5}
		\end{figure}

	    When $f=20$, Figure \ref{cde5} gives the third tiling with AVC $1$ in Figure \ref{various flip}.

		Henceforth, we can assume
		\[
		\text{AVC}\subset
		\{\aaa\ddd\eee,\bbb\ddd\eee,\aaa^2\ccc,\aaa\bbb\ccc,\bbb^2\ccc,\ccc^{\frac{f-4}8}\ddd\eee\}.
		\]

		\subsubsection*{Subcase. $\ccc^{\frac{f-4}8}\ddd\eee$ is a vertex.}
		By the AVC, we get that $\eee^2\cdots$ is not a vertex. Then we have four possible AAD for $\ccc^{\frac{f-4}8}\ddd\eee$:
	\begin{eqnarray*}
		\ccc^{\frac{f-4}8}\ddd\eee &=& \thin^{\aaa}\ccc^{\eee}\thin\cdots\thin^{\aaa}\ccc^{\eee}\thin^{\bbb}\ddd^{\eee}\thick^{\ddd}\eee^{\ccc}\thin, \\
		 &\text{or}& \thin^{\eee}\ccc^{\aaa}\thin\cdots\thin^{\eee}\ccc^{\aaa}\thin^{\bbb}\ddd^{\eee}\thick^{\ddd}\eee^{\ccc}\thin,\\
		 &\text{or}& \thin^{\eee}\ccc^{\aaa}\thin\cdots\thin^{\eee}\ccc^{\aaa}\thin^{\aaa}\ccc^{\eee}\thin^{\bbb}\ddd^{\eee}\thick^{\ddd}\eee^{\ccc}\thin, \\
		 &\text{or}& \thin^{\eee}\ccc^{\aaa}\thin\cdots\thin^{\eee}\ccc^{\aaa}\thin^{\aaa}\ccc^{\eee}\thin\cdots\thin^{\aaa}\ccc^{\eee}\thin^{\bbb}\ddd^{\eee}\thick^{\ddd}\eee^{\ccc}\thin.
	\end{eqnarray*}

	$\bullet$ In Figure \ref{cde4}, we have $\ccc^{\frac{f-4}8}\ddd\eee=\thin^{\aaa}\ccc_3^{\eee}\thin\cdots\thin^{\aaa}\ccc_4^{\eee}\thin^{\bbb}\ddd_1^{\eee}\thick^{\ddd}\eee_2^{\ccc}\thin$, which determines $T_1,T_2,T_3,\cdots,T_4$ ($\frac{f+12}8$ tiles). Just like in Figure \ref{ac1}, we can also determine all tiles in Figure \ref{cde4} which gives the tiling  $T(2\aaa^2\ccc,6\aaa\bbb\ccc,(f-10)\aaa\ddd\eee,(\frac{f}2-6)\bbb^2\ccc,6\bbb\ddd\eee,4\ccc^{\frac{f-4}8}\ddd\eee)$.
	It has one UFO: \{$T_{5}$, $T_{6}$, $T_{7}$, $T_{8}$\}.

		\begin{figure}[htp]
			\centering
			\begin{tikzpicture}[>=latex,scale=0.28]


				\begin{scope}[xshift=0cm]
					\fill[gray!50]  (-3-4-4,3) -- (-3-4-4,1) -- (-2-4-4,0) -- (-0-4-4,0) -- (1-4-4,1) -- (1-4-4,3);
					\fill[gray!50]  (-3-4,3) -- (-3-4,1) -- (-2-4,0) -- (-0.5-4,0) -- (1-4,1) -- (1-4,3);
					\fill[gray!50]  (-2+12,0) -- (0+12,0) -- (-0+12,-2) -- (-1+12,-3) -- (-2+12,-2) -- (-2+12,0);
					\fill[gray!50]  (-2+4,0) -- (0+4,0) -- (-0+4,-2) -- (-1+4,-3) -- (-2+4,-2) -- (-2+4,0);
					\fill[gray!50]  (-1+8,-5) -- (-1+8,-3) -- (0+8,-2) -- (2+8,-2) -- (3+8,-3) -- (3+8,-5);
					\fill[gray!50]  (-1,-5) -- (-1,-3) -- (0,-2) -- (2,-2) -- (3,-3) -- (3,-5);
					\fill[gray!50]  (-1-4,-5) -- (-1-4,-3) -- (0-4,-2) -- (2-4,-2) -- (3-4,-3) -- (3-4,-5);
					\fill[gray!50]  (-2,0) -- (0,0) -- (-0,-2) -- (-1,-3) -- (-2,-2) -- (-2,0);

					\draw
					(-3-4-4,3) -- (-3-4-4,1) -- (-2-4-4,0) -- (-0-4-4,0) -- (1-4-4,1) -- (1-4-4,3);
					\draw
					(-3-4,3) -- (-3-4,1) -- (-2-4,0) -- (-0.5-4,0) -- (1-4,1) -- (1-4,3);
					\draw
					(-3,3) -- (-3,1) -- (-2,0) -- (0,0) -- (1,1) -- (1,3);
					\draw
					(-1,-5) -- (-1,-3) -- (0,-2) -- (2,-2) -- (3,-3) -- (3,-5);
					\draw
					(-3+8,3) -- (-3+8,1) -- (-2+8,0) -- (0+8,0) -- (1+8,1) -- (1+8,3);
					\draw
					(-2+12,0) -- (0+12,0) -- (-0+12,-2) -- (-1+12,-3) -- (-2+12,-2) -- (-2+12,0);
					\draw
					(-2+4,0) -- (0+4,0) -- (-0+4,-2) -- (-1+4,-3) -- (-2+4,-2) -- (-2+4,0);
					\draw
					(-3+4,1) -- (-2+4,0) -- (-2+4,-2) -- (-4+4,-2) -- (-4+4,-0) -- (-3+4,1);
					\draw
					(-1+8,-5) -- (-1+8,-3) -- (0+8,-2) -- (2+8,-2) -- (3+8,-3) -- (3+8,-5);
					\draw
					(-3+12,1) -- (-2+12,0) -- (-2+12,-2) -- (-4+12,-2) -- (-4+12,-0) -- (-3+12,1);
					\draw
					(-3+12,3) -- (-3+12,1) -- (-2+12,0) -- (0+12,0) -- (1+12,1) -- (1+12,3);
					\draw
					(-3+16,1) -- (-2+16,0) -- (-2+16,-2) -- (-4+16,-2) -- (-4+16,-0) -- (-3+16,1);
					\draw
					(-1+12,-5) -- (-1+12,-3) -- (0+12,-2) -- (2+12,-2);
					\draw
					(-2,0) -- (-2,-2)--(-4,-2)--(-5,-3)--(-5,-5);
					\draw
					(14,0) -- (16,0);
					\draw
					(14,-2) -- (15,-3)--(16,-2);

					\draw[line width=1.5]
					(-2,-2) -- (-1,-3);
					\draw[line width=1.5]
					(16,0) -- (16,-2);
					\draw[line width=1.5]
					(15,-3) -- (15,-5);
					\draw[line width=1.5]
					(-3+4,1) -- (-2+2,0);
					\draw[line width=1.5]
					(-3+4+8,1) -- (-2+2+8,0);
					\draw[line width=1.5]
					(-3+4+12,1) -- (-2+2+12,0);

					\draw[line width=1.5]
					(-3-4,3) -- (-3-4,1);

					\draw[line width=1.5]
					(-1+4+8,-3) -- (-2+4+8,-2);
					\draw[line width=1.5]
					(-1+4,-3) -- (-2+4,-2);

					\fill (5,-1) circle (0.1); \fill (5.5,-1) circle (0.1); \fill (6,-1) circle (0.1);

					\node at (-3-1,3) {\scriptsize $A_1$};
					\node at (-3-1,1+0.5) {\scriptsize $A_2$};
					\node at (-2-0.7,-0.1) {\scriptsize $A_3$};
					\node at (-2.7,-1.4) {\scriptsize $A_4$};
					\node at (-4.2,-1.4) {\scriptsize $A_5$};
					\node at (-6,-3-0.5) {\scriptsize $A_6$};
					\node at (-6,-5) {\scriptsize $A_7$};
					\node at (11+4.7,-3-0.7) {\scriptsize $A_8$};
					\node at (16+0.7,-2-0.5) {\scriptsize $A_9$};
					\node at (16.7,0+0.5) {\scriptsize $A_{10}$};
					\node at (14+1,0+0.5) {\scriptsize $A_{11}$};
					\node at (13+1,1+0.5) {\scriptsize $A_{12}$};

					\fill (-3,3) circle (0.15);
					\fill (-3,1) circle (0.15);
					\fill (-4+2,-0) circle (0.15);
					\fill (-2,-2) circle (0.15);
					\fill (-2-2,-2) circle (0.15);
					\fill (-5,-3) circle (0.15);
					\fill (-5,-5) circle (0.15);
					\fill (-1+16,-3) circle (0.15);
					\fill (16,-2) circle (0.15);
					\fill (2+14,0) circle (0.15);
					\fill (2+12,0) circle (0.15);
					\fill (-3+16,1) circle (0.15);

					\node[draw,shape=circle, inner sep=0.5] at (-3+2-8,3-1) {\small $1$};
					\node[draw,shape=circle, inner sep=0.5] at (-3+2-4,3-1) {\small $2$};
					\node[draw,shape=circle, inner sep=0.5] at (-3+2,3-1) {\small $3$};
					\node[draw,shape=circle, inner sep=0.5] at (11,3-1) {\small $4$};
					\node[draw,shape=circle, inner sep=0.5] at (13,-4) {\small $5$};
					\node[draw,shape=circle, inner sep=0.5] at (15,-1) {\small $6$};
				\end{scope}
			    \begin{scope}[yshift=-10cm]
			    	\fill[gray!50]  (-2+12,0) -- (0+12,0) -- (-0+12,-2) -- (-1+12,-3) -- (-2+12,-2) -- (-2+12,0);
			    	\fill[gray!50]  (-2+4,0) -- (0+4,0) -- (-0+4,-2) -- (-1+4,-3) -- (-2+4,-2) -- (-2+4,0);
			    	\fill[gray!50]  (-1+8,-5) -- (-1+8,-3) -- (0+8,-2) -- (2+8,-2) -- (3+8,-3) -- (3+8,-5);
			    	\fill[gray!50]  (-1,-5) -- (-1,-3) -- (0,-2) -- (2,-2) -- (3,-3) -- (3,-5);
			    	\fill[gray!50]  (-1-4,-5) -- (-1-4,-3) -- (0-4,-2) -- (2-4,-2) -- (3-4,-3) -- (3-4,-5);
			    	\fill[gray!50]  (-2,0) -- (0,0) -- (-0,-2) -- (-1,-3) -- (-2,-2) -- (-2,0);

			    	\draw
			    	(-3,3) -- (-3,1) -- (-2,0) -- (0,0) -- (1,1) -- (1,3);
			    	\draw
			    	(-1,-5) -- (-1,-3) -- (0,-2) -- (2,-2) -- (3,-3) -- (3,-5);
			    	\draw
			    	(-3+8,3) -- (-3+8,1) -- (-2+8,0) -- (0+8,0) -- (1+8,1) -- (1+8,3);
			    	\draw
			    	(-2+12,0) -- (0+12,0) -- (-0+12,-2) -- (-1+12,-3) -- (-2+12,-2) -- (-2+12,0);
			    	\draw
			    	(-2+4,0) -- (0+4,0) -- (-0+4,-2) -- (-1+4,-3) -- (-2+4,-2) -- (-2+4,0);
			    	\draw
			    	(-3+4,1) -- (-2+4,0) -- (-2+4,-2) -- (-4+4,-2) -- (-4+4,-0) -- (-3+4,1);
			    	\draw
			    	(-1+8,-5) -- (-1+8,-3) -- (0+8,-2) -- (2+8,-2) -- (3+8,-3) -- (3+8,-5);
			    	\draw
			    	(-3+12,1) -- (-2+12,0) -- (-2+12,-2) -- (-4+12,-2) -- (-4+12,-0) -- (-3+12,1);
			    	\draw
			    	(-3+12,3) -- (-3+12,1) -- (-2+12,0) -- (0+12,0) -- (1+12,1) -- (1+12,3);
			    	\draw
			    	(-3+16,1) -- (-2+16,0) -- (-2+16,-2) -- (-4+16,-2) -- (-4+16,-0) -- (-3+16,1);
			    	\draw
			    	(-1+12,-5) -- (-1+12,-3) -- (0+12,-2) -- (2+12,-2);
			    	\draw
			    	(-2,0) -- (-2,-2)--(-4,-2)--(-5,-3)--(-5,-5);
			    	\draw
			    	(14,0) -- (16,0);
			    	\draw
			    	(14,-2) -- (15,-3)--(16,-2);

			    	\draw[line width=1.5]
			    	(-2,-2) -- (-1,-3);
			    	\draw[line width=1.5]
			    	(16,0) -- (16,-2);
			    	\draw[line width=1.5]
			    	(15,-3) -- (15,-5);
			    	\draw[line width=1.5]
			    	(-3+4,1) -- (-2+2,0);
			    	\draw[line width=1.5]
			    	(-3+4+8,1) -- (-2+2+8,0);
			    	\draw[line width=1.5]
			    	(-3+4+12,1) -- (-2+2+12,0);

			    	\draw[line width=1.5]
			    	(-1+4+8,-3) -- (-2+4+8,-2);
			    	\draw[line width=1.5]
			    	(-1+4,-3) -- (-2+4,-2);

			    	\fill (5,-1) circle (0.1); \fill (5.5,-1) circle (0.1); \fill (6,-1) circle (0.1);

			    	\node at (-3-1,3) {\scriptsize $A_4$};
			    	\node at (-3-1,1+0.5) {\scriptsize $A_3$};
			    	\node at (-2-0.7,-0.1) {\scriptsize $A_2$};
			    	\node at (-2.7,-1.4) {\scriptsize $A_1$};
			    	\node at (-4.2,-1.4) {\scriptsize $A_{12}$};
			    	\node at (-6,-3-0.5) {\scriptsize $A_{11}$};
			    	\node at (-6,-5) {\scriptsize $A_{10}$};
			    	\node at (11+4.7,-3-0.7) {\scriptsize $A_9$};
			    	\node at (16+0.7,-2-0.5) {\scriptsize $A_8$};
			    	\node at (16.7,0+0.5) {\scriptsize $A_{7}$};
			    	\node at (14+1,0+0.5) {\scriptsize $A_{6}$};
			    	\node at (13+1,1+0.5) {\scriptsize $A_{5}$};

			    	\fill (-3,3) circle (0.15);
			    	\fill (-3,1) circle (0.15);
			    	\fill (-4+2,-0) circle (0.15);
			    	\fill (-2,-2) circle (0.15);
			    	\fill (-2-2,-2) circle (0.15);
			    	\fill (-5,-3) circle (0.15);
			    	\fill (-5,-5) circle (0.15);
			    	\fill (-1+16,-3) circle (0.15);
			    	\fill (16,-2) circle (0.15);
			    	\fill (2+14,0) circle (0.15);
			    	\fill (2+12,0) circle (0.15);
			    	\fill (-3+16,1) circle (0.15);

			    	\node[draw,shape=circle, inner sep=0.5] at (-1,-1) {\small $2$};
			    	\node[draw,shape=circle, inner sep=0.5] at (-3,-4) {\small $1$};
			    	\node[draw,shape=circle, inner sep=0.5] at (13,-4) {\small $7$};
			    	\node[draw,shape=circle, inner sep=0.5] at (15,-1) {\small $8$};
			    \end{scope}
			\end{tikzpicture}
			\caption{$\ccc^{\frac{f-4}8}\ddd\eee=\thin^{\aaa}\ccc^{\eee}\thin\cdots\thin^{\aaa}\ccc^{\eee}\thin^{\bbb}\ddd^{\eee}\thick^{\ddd}\eee^{\ccc}\thin$.}
			\label{cde4}
		\end{figure}

	    When $f=20$, Figure \ref{cde4} gives the tiling with AVC $2$ in Figure \ref{various flip}.

		$\bullet$ For the other AAD of $\ccc^{\frac{f-4}{8}}\ddd\eee$, half of the sphere can be determined through AAD. As shown in Figures \ref{cde1}, \ref{cde3}, \ref{cde2}, we get three different half spheres HS1, HS2, HS3 from $\thin\ddd_1\thick\eee_2\thin^{\eee}\ccc_3^{\aaa}\thin\cdots=\thin\ddd_1\thick\eee_2\thin^{\eee}\ccc_3^{\aaa}\thin\cdots\thin^{\eee}\ccc_4^{\aaa}\thin$, or
		$\thin\ddd_1\thick\eee_2\thin^{\eee}\ccc_3^{\aaa}\thin\cdots\thin^{\eee}\ccc_4^{\aaa}\thin^{\aaa}\ccc_5^{\eee}\thin$, or $\thin\ddd_1\thick\eee_2\thin^{\eee}\ccc_3^{\aaa}\thin\cdots\thin^{\eee}\ccc_4^{\aaa}\thin^{\aaa}\ccc_5^{\eee}\thin\cdots\thin^{\aaa}\ccc_6^{\eee}\thin$ respectively, but with the same boundary and the same labeling of $A_1$, $\cdots$, $A_{12}$.  Similarly, the remaining half of the sphere can also be determined by the AAD.

		\begin{figure}[htp]
			\centering
			\begin{tikzpicture}[>=latex,scale=0.28]


				\begin{scope}[xshift=0cm]
					\fill[gray!50]  (-3-4-4,3) -- (-3-4-4,1) -- (-2-4-4,0) -- (-0-4-4,0) -- (1-4-4,1) -- (1-4-4,3);
					\fill[gray!50]  (-3-4,3) -- (-3-4,1) -- (-2-4,0) -- (-0.5-4,0) -- (1-4,1) -- (1-4,3);
					\fill[gray!50]  (-3,3) -- (-3,1) -- (-2,0) -- (0,0) -- (1,1) -- (1,3);
					\fill[gray!50]  (-3+4+4,3) -- (-3+4+4,1) -- (-2+4+4,0) -- (0+4+4,0) -- (1+4+4,1) -- (1+4+4,3);
					\fill[gray!50]  (-3+12,3) -- (-3+12,1) -- (-2+12,0) -- (0+12,0) -- (1+12,1) -- (1+12,3);

					\fill[gray!50]  (-3,1) -- (-2,0) -- (-2,-2) -- (-4,-2) -- (-4,-0) -- (-3,1);
					\fill[gray!50]  (-3+8,1) -- (-2+8,0) -- (-2+8,-2) -- (-4+8,-2) -- (-4+8,-0) -- (-3+8,1);
					\fill[gray!50]  (-3+12,1) -- (-2+12,0) -- (-2+12,-2) -- (-4+12,-2) -- (-4+12,-0) -- (-3+12,1);
					\fill[gray!50]  (-3+16,1) -- (-2+16,0) -- (-2+16,-2) -- (-4+16,-2) -- (-4+16,-0) -- (-3+16,1);

					\fill[gray!50]  (-2+12,0) -- (0+12,0) -- (-0+12,-2) -- (-1+12,-3) -- (-2+12,-2) -- (-2+12,0);
					\fill[gray!50]  (-2+16,0) -- (0+16,0) -- (-0+16,-2) -- (-1+16,-3) -- (-2+16,-2) -- (-2+16,0);

					\fill[gray!50]  (-1+12,-5) -- (-1+12,-3) -- (0+12,-2) -- (2+12,-2) -- (3+12,-3) -- (3+12,-5);

					\draw
					(-3-4-4,3) -- (-3-4-4,1) -- (-2-4-4,0) -- (-0-4-4,0) -- (1-4-4,1) -- (1-4-4,3);
					\draw
					(-3-4,3) -- (-3-4,1) -- (-2-4,0) -- (-0.5-4,0) -- (1-4,1) -- (1-4,3);
					\draw
					(-3,3) -- (-3,1) -- (-2,0) -- (0,0) -- (1,1) -- (1,3);
					\draw
					(-2,0) -- (0,0) -- (-0,-2) -- (-1,-3) -- (-2,-2) -- (-2,0);
					\draw
					(-3,1) -- (-2,0) -- (-2,-2) -- (-4,-2) -- (-4,-0) -- (-3,1);
					\draw
					(-1,-5) -- (-1,-3) -- (0,-2) -- (2,-2) -- (3,-3) -- (3,-5);
					\draw
					(-3+8,3) -- (-3+8,1) -- (-2+8,0) -- (0+8,0) -- (1+8,1) -- (1+8,3);
					\draw
					(-2+8,0) -- (0+8,0) -- (-0+8,-2) -- (-1+8,-3) -- (-2+8,-2) -- (-2+8,0);
					\draw
					(-3+8,1) -- (-2+8,0) -- (-2+8,-2) -- (-4+8,-2) -- (-4+8,-0) -- (-3+8,1);
					\draw
					(-1+8,-5) -- (-1+8,-3) -- (0+8,-2) -- (2+8,-2) -- (3+8,-3) -- (3+8,-5);
					\draw
					(-3+12,3) -- (-3+12,1) -- (-2+12,0) -- (0+12,0) -- (1+12,1) -- (1+12,3);
					\draw
					(-2+12,0) -- (0+12,0) -- (-0+12,-2) -- (-1+12,-3) -- (-2+12,-2) -- (-2+12,0);
					\draw
					(-3+16,1) -- (-2+16,0) -- (-2+16,-2) -- (-4+16,-2) -- (-4+16,-0) -- (-3+16,1);
					\draw
					(-1+12,-5) -- (-1+12,-3) -- (0+12,-2) -- (2+12,-2);
					\draw
					(-2+12+4,0) -- (0+12+4,0) -- (-0+12+4,-2) -- (-1+12+4,-3) -- (-2+12+4,-2) -- (-2+12+4,0);
					\draw
					(-1+12+4,-5) -- (-1+12+4,-3);
					\draw
					(-3+12,1) -- (-2+12,0) -- (-2+12,-2) -- (-4+12,-2) -- (-4+12,-0) -- (-3+12,1);

					\draw[line width=1.5]
					(-3,1) -- (-2,0);
					\draw[line width=1.5]
					(-3+8,1) -- (-2+8,0);
					\draw[line width=1.5]
					(-3+12,1) -- (-2+12,0);

					\draw[line width=1.5]
					(-3-4,3) -- (-3-4,1);

					\draw[line width=1.5]
					(-1,-3) -- (0,-2);
					\draw[line width=1.5]
					(-1+8,-3) -- (0+8,-2);

					\draw[line width=1.5]
					(-4+16,-2) -- (-4+16,-0);

					\draw[line width=1.5]
					(-1+12+4,-3) -- (-2+12+4,-2);

					\fill (5.5-4,-1) circle (0.1); \fill (6-4,-1) circle (0.1); \fill (6.5-4,-1) circle (0.1);

					\node at (-3-1-0.8,3+0.8) {\scriptsize $A_1(A_2)$};
					\node at (-3-1-0.8,1+0.2) {\scriptsize $A_2(A_1)$};
					\node at (-4-1-1,0-0.5) {\scriptsize $A_3(A_{12})$};
					\node at (-4-1-1,-2) {\scriptsize $A_4(A_{11})$};
					\node at (-2-0.7-0.8,-2-0.5-0.3) {\scriptsize $A_5(A_{10})$};
					\node at (-1-1-1,-3-0.5-0.5) {\scriptsize $A_6(A_9)$};
					\node at (-1-1-1,-5-0.5) {\scriptsize $A_7(A_8)$};
					\node at (11+0.7+4+1,-3-0.7) {\scriptsize $A_8(A_7)$};
					\node at (12+0.7+4+1,-2-0.5) {\scriptsize $A_9(A_6)$};
					\node at (12+0.7+4+1,0-0.5) {\scriptsize $A_{10}(A_5)$};
					\node at (14+1+1,0+0.5) {\scriptsize $A_{11}(A_4)$};
					\node at (13+1+1,1+0.5+0.3) {\scriptsize $A_{12}(A_3)$};

					\fill (-3,3) circle (0.15);
					\fill (-3,1) circle (0.15);
					\fill (-4,-0) circle (0.15);
					\fill (-4,-2) circle (0.15);
					\fill (-2,-2) circle (0.15);
					\fill (-1,-3) circle (0.15);
					\fill (-1,-5) circle (0.15);
					\fill (-1+12+4,-3) circle (0.15);
					\fill (-0+12+4,-2) circle (0.15);
					\fill (2+12+2,0) circle (0.15);
					\fill (2+12,0) circle (0.15);
					\fill (-3+16,1) circle (0.15);

					\node[draw,shape=circle, inner sep=0.5] at (-3+2-8,3-1) {\small $1$};
					\node[draw,shape=circle, inner sep=0.5] at (-3+2-4,3-1) {\small $2$};
					\node[draw,shape=circle, inner sep=0.5] at (-3+2,3-1) {\small $3$};
					\node[draw,shape=circle, inner sep=0.5] at (-3+2+4+8,3-1) {\small $4$};
					\node[draw,shape=circle, inner sep=0.5] at (-4+3+12,-1) {\small $6$};
					\node[draw,shape=circle, inner sep=0.5] at (-4+3+12+4,-1) {\small $8$};
					\node[draw,shape=circle, inner sep=0.5] at (-4+5+12,-1) {\small $7$};
					\node[draw,shape=circle, inner sep=0.5] at (1+8+4,-4) {\small $9$};
					\node[draw,shape=circle, inner sep=0.5] at (-4+5+8,-1) {\small $5$};
				\end{scope}
			\end{tikzpicture}
			\caption{The half sphere HS1 from  $\ccc^{\frac{f-4}8}\ddd\eee=\thin^{\eee}\ccc^{\aaa}\thin\cdots\thin^{\eee}\ccc^{\aaa}\thin^{\bbb}\ddd^{\eee}\thick^{\ddd}\eee^{\ccc}\thin$.}
			\label{cde1}
		\end{figure}

		\begin{figure}[htp]
			\centering
			\begin{tikzpicture}[>=latex,scale=0.28]


				\begin{scope}[xshift=0cm]
					\fill[gray!50]  (-3-4-4,3) -- (-3-4-4,1) -- (-2-4-4,0) -- (-0-4-4,0) -- (1-4-4,1) -- (1-4-4,3);
					\fill[gray!50]  (-3-4,3) -- (-3-4,1) -- (-2-4,0) -- (-0.5-4,0) -- (1-4,1) -- (1-4,3);
					\fill[gray!50]  (-3,3) -- (-3,1) -- (-2,0) -- (0,0) -- (1,1) -- (1,3);
					\fill[gray!50]  (-3+8-1+1,3) -- (-3+8-1+1,1) -- (-2+8-1+1,0) -- (0+8-1+1,0) -- (1+8-1+1,1) -- (1+8-1+1,3);
					\fill[gray!50]  (-1+12-1+1,-5) -- (-1+12-1+1,-3) -- (0+12-1+1,-2) -- (2+12-1+1,-2) -- (3+12-1+1,-3) -- (3+12-1+1,-5);

					\fill[gray!50]  (-3,1) -- (-2,0) -- (-2,-2) -- (-4,-2) -- (-4,-0) -- (-3,1);
					\fill[gray!50]  (-3+4+3+1,1) -- (-2+4+3+1,0) -- (-2+4+3+1,-2) -- (-4+4+3+1,-2) -- (-4+4+3+1,-0) -- (-3+4+3+1,1);
					\fill[gray!50]  (-2+16-1+1,0) -- (0+16-1+1,0) -- (-0+16-1+1,-2) -- (-1+16-1+1,-3) -- (-2+16-1+1,-2) -- (-2+16-1+1,0);

					\draw
					(-3-4-4,3) -- (-3-4-4,1) -- (-2-4-4,0) -- (-0-4-4,0) -- (1-4-4,1) -- (1-4-4,3);
					\draw
					(-3-4,3) -- (-3-4,1) -- (-2-4,0) -- (-0.5-4,0) -- (1-4,1) -- (1-4,3);
					\draw
					(-3,3) -- (-3,1) -- (-2,0) -- (0,0) -- (1,1) -- (1,3);
					\draw
					(-2,0) -- (0,0) -- (-0,-2) -- (-1,-3) -- (-2,-2) -- (-2,0);
					\draw
					(-3,1) -- (-2,0) -- (-2,-2) -- (-4,-2) -- (-4,-0) -- (-3,1);
					\draw
					(-1,-5) -- (-1,-3) -- (0,-2) -- (2,-2) -- (3,-3) -- (3,-5);
					\draw
					(4+1,1) -- (3+1,0) -- (3+1,-2)--(5+1,-2);
					\draw
					(-3+4+4-1+1,3) -- (-3+4+4-1+1,1) -- (-2+4+4-1+1,0) -- (0+4+4-1+1,0) -- (1+4+4-1+1,1) -- (1+4+4-1+1,3);
					\draw
					(-2+4+4-1+1,0) -- (0+4+4-1+1,0) -- (-0+4+4-1+1,-2) -- (-1+4+4-1+1,-3) -- (-2+4+4-1+1,-2) -- (-2+4+4-1+1,0);

					\draw
					(-1+4+4-1+1,-5) -- (-1+4+4-1+1,-3) -- (0+4+4-1+1,-2) -- (2+4+4-1+1,-2) -- (3+4+4-1+1,-3) -- (3+4+4-1+1,-5);
					\draw
					(-3+8+4-1+1,1) -- (-2+8+4-1+1,0) -- (-2+8+4-1+1,-2) -- (-4+8+4-1+1,-2) -- (-4+8+4-1+1,-0) -- (-3+8+4-1+1,1);
					\draw
					(-3+12-1+1,3) -- (-3+12-1+1,1) -- (-2+12-1+1,0) -- (0+12-1+1,0) -- (1+12-1+1,1) -- (1+12-1+1,3);
					\draw
					(-2+12-1+1,0) -- (0+12-1+1,0) -- (-0+12-1+1,-2) -- (-1+12-1+1,-3) -- (-2+12-1+1,-2) -- (-2+12-1+1,0);
					\draw
					(-3+16-1+1,1) -- (-2+16-1+1,0) -- (-2+16-1+1,-2) -- (-4+16-1+1,-2) -- (-4+16-1+1,-0) -- (-3+16-1+1,1);
					\draw
					(-1+12-1+1,-5) -- (-1+12-1+1,-3) -- (0+12-1+1,-2) -- (2+12-1+1,-2);
					\draw
					(-2+12+4-1+1,0) -- (0+12+4-1+1,0) -- (-0+12+4-1+1,-2) -- (-1+12+4-1+1,-3) -- (-2+12+4-1+1,-2) -- (-2+12+4-1+1,0);
					\draw
					(-1+12+4-1+1,-5) -- (-1+12+4-1+1,-3);

					\draw[line width=1.5]
					(-3,1) -- (-2,0);
					\draw[line width=1.5]
					(-3+8-1+1,1) -- (-2+8-1+1,0);
					\draw[line width=1.5]
					(-3+12+4-1+1,1) -- (-2+12+2-1+1,0);

					\draw[line width=1.5]
					(-3-4,3) -- (-3-4,1);

					\draw[line width=1.5]
					(-1,-3) -- (0,-2);
					\draw[line width=1.5]
					(-1+8-1+1,-3) -- (0+8-1+1,-2);

					\draw[line width=1.5]
					(-4+16-2-1+1,-2) -- (-4+16-2-1+1,-0);

					\draw[line width=1.5]
					(-1+12+4-1+1,-3) -- (-2+12+4-1+1,-2);

					\fill (1.5,-1) circle (0.1); \fill (2,-1) circle (0.1); \fill (2.5,-1) circle (0.1);

					\node at (-3-1-0.8,3+0.8) {\scriptsize $A_1(A_2)$};
					\node at (-3-1-0.8,1+0.2) {\scriptsize $A_2(A_1)$};
					\node at (-4-1-1,0-0.5) {\scriptsize $A_3(A_{12})$};
					\node at (-4-1-1,-2) {\scriptsize $A_4(A_{11})$};
					\node at (-2-0.7-0.8,-2-0.5-0.3) {\scriptsize $A_5(A_{10})$};
					\node at (-1-1-1,-3-0.5-0.5) {\scriptsize $A_6(A_9)$};
					\node at (-1-1-1,-5-0.5) {\scriptsize $A_7(A_8)$};
					\node at (11+0.7+4+1-1+1,-3-0.7) {\scriptsize $A_8(A_7)$};
					\node at (12+0.7+4+1-1+1,-2-0.5) {\scriptsize $A_9(A_6)$};
					\node at (12+0.7+4+1-1+1,0-0.5) {\scriptsize $A_{10}(A_5)$};
					\node at (14+1+1-1+1,0+0.5) {\scriptsize $A_{11}(A_4)$};
					\node at (13+1+1-1+1,1+0.5+0.3) {\scriptsize $A_{12}(A_3)$};

					\fill (-3,3) circle (0.15);
					\fill (-3,1) circle (0.15);
					\fill (-4,-0) circle (0.15);
					\fill (-4,-2) circle (0.15);
					\fill (-2,-2) circle (0.15);
					\fill (-1,-3) circle (0.15);
					\fill (-1,-5) circle (0.15);
					\fill (-1+12+4-1+1,-3) circle (0.15);
					\fill (-0+12+4-1+1,-2) circle (0.15);
					\fill (2+12+2-1+1,0) circle (0.15);
					\fill (2+12-1+1,0) circle (0.15);
					\fill (-3+16-1+1,1) circle (0.15);

					\node[draw,shape=circle, inner sep=0.5] at (-3+2-8,3-1) {\small $1$};
					\node[draw,shape=circle, inner sep=0.5] at (-3+2-4,3-1) {\small $2$};
					\node[draw,shape=circle, inner sep=0.5] at (-3+2,3-1) {\small $3$};
					\node[draw,shape=circle, inner sep=0.5] at (-3+2+8-1+1,3-1) {\small $4$};
					\node[draw,shape=circle, inner sep=0.5] at (-3+2+4+8-1+1,3-1) {\small $5$};
					\node[draw,shape=circle, inner sep=0.5] at (-4+3+12-1+1,-1) {\small $8$};
					\node[draw,shape=circle, inner sep=0.5] at (-4+5+12-1+1,-1) {\small $10$};
					\node[draw,shape=circle, inner sep=0.5] at (-4+3+12-4-1+1,-1) {\small $6$};
					\node[draw,shape=circle, inner sep=0.5] at (-4+5+12-4-1+1,-1) {\small $7$};
					\node[draw,shape=circle, inner sep=0.5] at (1+8+4-4-1+1,-4) {\small $9$};
				\end{scope}
			\end{tikzpicture}
			\caption{The half sphere HS2 from  $\ccc^{\frac{f-4}8}\ddd\eee=\thin^{\eee}\ccc^{\aaa}\thin\cdots\thin^{\eee}\ccc^{\aaa}\thin^{\aaa}\ccc^{\eee}\thin^{\bbb}\ddd^{\eee}\thick^{\ddd}\eee^{\ccc}\thin$.}
			\label{cde3}
		\end{figure}

		\begin{figure}[htp]
			\centering
			\begin{tikzpicture}[>=latex,scale=0.28]


				\begin{scope}[xshift=0cm]
					\fill[gray!50]  (-3-4-4,3) -- (-3-4-4,1) -- (-2-4-4,0) -- (-0-4-4,0) -- (1-4-4,1) -- (1-4-4,3);
					\fill[gray!50]  (-3-4,3) -- (-3-4,1) -- (-2-4,0) -- (-0.5-4,0) -- (1-4,1) -- (1-4,3);
					\fill[gray!50]  (-3,3) -- (-3,1) -- (-2,0) -- (0,0) -- (1,1) -- (1,3);
					\fill[gray!50]  (-3+8,3) -- (-3+8,1) -- (-2+8,0) -- (0+8,0) -- (1+8,1) -- (1+8,3);
					\fill[gray!50]  (-1+12,-5) -- (-1+12,-3) -- (0+12,-2) -- (2+12,-2) -- (3+12,-3) -- (3+12,-5);
					\fill[gray!50]  (-1+24-4,-5) -- (-1+24-4,-3) -- (0+24-4,-2) -- (2+24-4,-2) -- (3+24-4,-3) -- (3+24-4,-5);

					\fill[gray!50]  (-3,1) -- (-2,0) -- (-2,-2) -- (-4,-2) -- (-4,-0) -- (-3,1);
					\fill[gray!50]  (-3+8,1) -- (-2+8,0) -- (-2+8,-2) -- (-4+8,-2) -- (-4+8,-0) -- (-3+8,1);
					\fill[gray!50]  (-2+16,0) -- (0+16,0) -- (-0+16,-2) -- (-1+16,-3) -- (-2+16,-2) -- (-2+16,0);
					\fill[gray!50]  (-2+28-4,0) -- (0+28-4,0) -- (-0+28-4,-2) -- (-1+28-4,-3) -- (-2+28-4,-2) -- (-2+28-4,0);

					\draw
					(-3-4-4,3) -- (-3-4-4,1) -- (-2-4-4,0) -- (-0-4-4,0) -- (1-4-4,1) -- (1-4-4,3);
					\draw
					(-3-4,3) -- (-3-4,1) -- (-2-4,0) -- (-0.5-4,0) -- (1-4,1) -- (1-4,3);
					\draw
					(-3,3) -- (-3,1) -- (-2,0) -- (0,0) -- (1,1) -- (1,3);
					\draw
					(-2,0) -- (0,0) -- (-0,-2) -- (-1,-3) -- (-2,-2) -- (-2,0);
					\draw
					(-3,1) -- (-2,0) -- (-2,-2) -- (-4,-2) -- (-4,-0) -- (-3,1);
					\draw
					(-1,-5) -- (-1,-3) -- (0,-2) -- (2,-2) -- (3,-3) -- (3,-5);
					\draw
					(5,1) -- (4,0)--(4,-2)--(6,-2);
					\draw
					(-3+4+4,3) -- (-3+4+4,1) -- (-2+4+4,0) -- (0+4+4,0) -- (1+4+4,1) -- (1+4+4,3);
					\draw
					(-2+4+4,0) -- (0+4+4,0) -- (-0+4+4,-2) -- (-1+4+4,-3) -- (-2+4+4,-2) -- (-2+4+4,0);
					\draw
					(-1+4+4,-5) -- (-1+4+4,-3) -- (0+4+4,-2) -- (2+4+4,-2) -- (3+4+4,-3) -- (3+4+4,-5);
					\draw
					(-3+8+4,1) -- (-2+8+4,0) -- (-2+8+4,-2) -- (-4+8+4,-2) -- (-4+8+4,-0) -- (-3+8+4,1);
					\draw
					(-3+12,3) -- (-3+12,1) -- (-2+12,0) -- (0+12,0) -- (1+12,1) -- (1+12,3);
					\draw
					(-2+12,0) -- (0+12,0) -- (-0+12,-2) -- (-1+12,-3) -- (-2+12,-2) -- (-2+12,0);
					\draw
					(-1+12,-5) -- (-1+12,-3) -- (0+12,-2) -- (2+12,-2);
					\draw
					(-3+16,3) -- (-3+16,1) -- (-2+16,0) -- (0+16,0);
					\draw
					(-2+16,0) -- (0+16,0) -- (-0+16,-2) -- (-1+16,-3) -- (-2+16,-2) -- (-2+16,0);
					\draw
					(-4+20,-2) -- (-4+20,-0);
					\draw
					(-1+16,-5) -- (-1+16,-3) -- (0+16,-2);
					\draw
					(-3+24-4,3) -- (-3+24-4,1) -- (-2+24-4,0) -- (0+24-4,0) -- (1+24-4,1) -- (1+24-4,3);
					\draw
					(-2+24-4,0) -- (0+24-4,0) -- (-0+24-4,-2) -- (-1+24-4,-3);
					\draw
					(-3+28-4,1) -- (-2+28-4,0) -- (-2+28-4,-2) -- (-4+28-4,-2) -- (-4+28-4,-0) -- (-3+28-4,1);
					\draw
					(-1+24-4,-5) -- (-1+24-4,-3) -- (0+24-4,-2) -- (2+24-4,-2);
					\draw
					(-1+16+12-4,-5) -- (-1+16+12-4,-3);
					\draw
					(-2+24+4-4,0) -- (0+24+4-4,0) -- (-0+24+4-4,-2) -- (-1+24+4-4,-3) -- (-2+24+4-4,-2) -- (-2+24+4-4,0);

					\draw[line width=1.5]
					(-3,1) -- (-2,0);
					\draw[line width=1.5]
					(-3+8,1) -- (-2+8,0);
					\draw[line width=1.5]
					(-3+16,1) -- (-4+16,-0);
					\draw[line width=1.5]
					(-3+16+4+8-4,1) -- (-4+20+8-4,-0);

					\draw[line width=1.5]
					(-3-4,3) -- (-3-4,1);

					\draw[line width=1.5]
					(-1+24-8,-3) -- (-2+24-8,-2);
					\draw[line width=1.5]
					(-1+24+4-4,-3) -- (-2+24+4-4,-2);

					\draw[line width=1.5]
					(-1,-3) -- (0,-2);
					\draw[line width=1.5]
					(-1+8,-3) -- (0+8,-2);

					\draw[line width=1.5]
					(-4+18-4,-2) -- (-4+18-4,-0);

					\fill (1.5,-1) circle (0.1); \fill (2,-1) circle (0.1); \fill (2.5,-1) circle (0.1); \fill (19.5-2,-1) circle (0.1); \fill (20-2,-1) circle (0.1); \fill (20.5-2,-1) circle (0.1);

					\node at (-3-1-0.8,3+0.8) {\scriptsize $A_1(A_2)$};
					\node at (-3-1-0.8,1+0.2) {\scriptsize $A_2(A_1)$};
					\node at (-4-1-1,0-0.5) {\scriptsize $A_3(A_{12})$};
					\node at (-4-1-1,-2) {\scriptsize $A_4(A_{11})$};
					\node at (-2-0.7-0.8,-2-0.5-0.3) {\scriptsize $A_5(A_{10})$};
					\node at (-1-1-1,-3-0.5-0.5) {\scriptsize $A_6(A_9)$};
					\node at (-1-1-1,-5-0.5) {\scriptsize $A_7(A_8)$};
					\node at (11+0.7+4+1+8+4-4,-3-0.7) {\scriptsize $A_8(A_7)$};
					\node at (12+0.7+4+1+8+4-4,-2-0.5) {\scriptsize $A_9(A_6)$};
					\node at (12+0.7+4+1+8+4-4,0-0.5) {\scriptsize $A_{10}(A_5)$};
					\node at (14+1+1+8+4-4,0+0.5) {\scriptsize $A_{11}(A_4)$};
					\node at (13+1+1+8+4-4,1+0.5+0.3) {\scriptsize $A_{12}(A_3)$};

					\fill (-3,3) circle (0.15);
					\fill (-3,1) circle (0.15);
					\fill (-4,-0) circle (0.15);
					\fill (-4,-2) circle (0.15);
					\fill (-2,-2) circle (0.15);
					\fill (-1,-3) circle (0.15);
					\fill (-1,-5) circle (0.15);
					\fill (-1+12+4+8+4-4,-3) circle (0.15);
					\fill (-0+12+4+8+4-4,-2) circle (0.15);
					\fill (2+12+4+8+2-4,0) circle (0.15);
					\fill (2+12+4+8-4,0) circle (0.15);
					\fill (-3+16+4+8-4,1) circle (0.15);

					\node[draw,shape=circle, inner sep=0.5] at (-3+2-8,3-1) {\small $1$};
					\node[draw,shape=circle, inner sep=0.5] at (-3+2-4,3-1) {\small $2$};
					\node[draw,shape=circle, inner sep=0.5] at (-3+2,3-1) {\small $3$};
					\node[draw,shape=circle, inner sep=0.5] at (-3+6+4,3-1) {\small $4$};
					\node[draw,shape=circle, inner sep=0.5] at (-3+6+8,3-1) {\small $5$};
					\node[draw,shape=circle, inner sep=0.5] at (-3+6+16,3-1) {\small $6$};
					\node[draw,shape=circle, inner sep=0.5] at (-4+3+12,-1) {\small $9$};
					\node[draw,shape=circle, inner sep=0.5] at (-4+5+12,-1) {\small $11$};
					\node[draw,shape=circle, inner sep=0.1] at (-4+3+12-4,-1) {\small $7$};
					\node[draw,shape=circle, inner sep=0.5] at (-4+5+12-4,-1) {\small $8$};
					\node[draw,shape=circle, inner sep=0.5] at (1+8,-4) {\small $10$};
				\end{scope}

			\end{tikzpicture}
			\caption{The half sphere HS3 from $\ccc^{\frac{f-4}8}\ddd\eee=\thin^{\eee}\ccc^{\aaa}\thin\cdots\thin^{\aaa}\ccc^{\eee}\thin\cdots\thin^{\aaa}\ccc^{\eee}\thin^{\bbb}\ddd^{\eee}\thick^{\ddd}\eee^{\ccc}\thin$.}
			\label{cde2}
		\end{figure}

    When $f\ge28$, each of HS1, HS2, HS3 with boundary vertices $A_1$, $\cdots$, $A_{12}$ outside the parenthesis can be glued with any of HS1, HS2, HS3, including itself, by matching the $A_1$, $\cdots$, $A_{12}$ inside the parenthesis.

    There are two UFOs: \{$T_{4},T_{5},T_{6},T_{7}$\} and \{$T_{6},T_{7},T_{8},T_{9}$\} in HS1.	There are two UFOs: \{$T_{6},T_{7},T_{8},T_{9}$\} and \{$T_{5},T_{7},T_{8},T_{10}$\} in HS2.	There are two UFOs: \{$T_{7},T_{8},T_{9},T_{10}$\} and \{$T_{5},T_{8},T_{9},T_{11}$\} in HS3. Whatever tiling is formed by gluing these half spheres, there are always four UFOs in the end.

	When $f=20$, there is no HS3. Let HS2$'$ denote the modification of HS2 by flipping its UFO \{$T_{6}$, $T_{7}$, $T_{8}$, $T_{9}$\}.  Then HS1$+$HS1$=$HS2$'$$+$HS2$'$ gives the first tiling with AVC $6$ in Figure \ref{various flip}; and HS1$+$HS2$=$HS2$+$HS2$'$ gives the tiling with AVC $7$ in Figure \ref{various flip}; and HS1$+$HS2$'$ gives the second tiling with AVC $6$ in Figure \ref{various flip}; and HS2$+$HS2 gives the tiling with AVC $8$ in Figure \ref{various flip}.

	\vspace{6pt}

	In conclusion, Case. $\{\aaa\ddd\eee, \bbb\ddd\eee, \aaa^2\ccc(\aaa\bbb\ccc,\bbb^2\ccc)\}$ admits various modifications of the non-symmetric $3$-layer earth map tiling with $f=8k+4$, $k\ge2$, summarized  in Table \ref{table-3}. In particular, we count the number of different tilings associated to each possible AVC explicitly in this table.
	\end{proof}

	\section{Tilings with a unique $a^2b$-vertex type $\aaa\ddd^2$, $\bbb\ddd^2$ or $\ccc\ddd^2$}
	\label{ad2,bd2 or be2 as the unique $a^2b$-vertex}

	Our strategy in the following sections is to first tile the partial neighborhood of a special tile from Lemma \ref{base_tile}.

	\begin{proposition}\label{bd2}
		There is no non-symmetric $a^4b$-tiling with general angles, such that $\aaa\ddd^2$, $\bbb\ddd^2$ or $\ccc\ddd^2$ is the unique $a^2b$-vertex type.
	\end{proposition}

	\begin{proof}
		Since $\eee$ can not appear in an $a^2b$-vertex (thus not in any degree $3$ vertex), the special tile must have a vertex of degree $4$ or $5$ at $\eee$ and must have degree $3$ vertices at the other four corners (see Figure \ref{fbd2}).

		By Lemma \ref{hdeg} and the parity lemma, we know that $\ddd\eee^{3}$ or $\eee^{4}$ is a vertex.

		If $\aaa\ddd^{2}$ is the unique $a^2b$-vertex type, the partial neighborhood of a special tile is shown in the first and second pictures of Figure \ref{fbd2}. In both pictures $\ddd_1\cdots=\aaa_3\ddd_1\ddd_2$ determines $T_2$. These two pictures show two possible arrangements of $T_3$. In the first picture the degree $3$ vertex $\bbb_1\ccc_3\cdots=\aaa\bbb\ccc$, $\bbb^2\ccc$ or $\bbb\ccc^2$. In the second picture the degree $3$ vertex $\bbb_1\bbb_3\cdots=\aaa\bbb^2$, $\bbb^3$ or $\bbb^2\ccc$. So three vertices $\aaa\ddd^2$ + one of \{$\ddd\eee^{3}, \eee^{4}$\} + one of \{$\aaa\bbb^{2}, \aaa\bbb\ccc, \bbb^3, \bbb^{2}\ccc, \bbb\ccc^2$\} must appear.

		If $\bbb\ddd^{2}$ is the unique $a^2b$-vertex type, the partial neighborhood of a special tile is shown in the third and fourth pictures of Figure \ref{fbd2}. In both pictures $\ddd_1\cdots=\bbb_3\ddd_1\ddd_2$ determines $T_2$. These two pictures show two possible arrangements of $T_3$. Then we have the degree $3$ vertex $\bbb_1\ddd_3\cdots=\bbb_1\ddd_3\ddd_4$, which determines $T_4$. In the third and fourth pictures the degree $3$ vertex $\aaa\bbb\cdots=\aaa^2\bbb$, $\aaa\bbb^2$ or $\aaa\bbb\ccc$. So three vertices $\bbb\ddd^2$ + one of \{$\ddd\eee^{3}, \eee^{4}$\} + one of \{$\aaa^{2}\bbb, \aaa\bbb\ccc, \aaa\bbb^{2}$\} must appear.

		If $\ccc\ddd^{2}$ is the unique $a^2b$-vertex type, the partial neighborhood of a special tile is shown in the fifth picture of Figure \ref{fbd2}. In this picture $\ddd_1\cdots=\ccc_3\ddd_1\ddd_2$ determines $T_2$, $T_3$. The degree $3$ vertex $\aaa_3\bbb_1\cdots=\aaa^2\bbb$, $\aaa\bbb^2$ or $\aaa\bbb\ccc$. So three vertices $\ccc\ddd^{2}$ + one of \{$\ddd\eee^{3}, \eee^{4}$\} + one of \{$\aaa^2\bbb, \aaa\bbb^2, \aaa\bbb\ccc$\} must appear.

		\begin{figure}[htp]
			\centering
			\begin{tikzpicture}[>=latex]


				\begin{scope}[xshift=-4cm,yshift=3 cm,scale=0.9]

					\fill (0,0.7) circle (0.1);

					\foreach \x in {1,2,3}
					\draw[rotate=72*\x]
					(18:0.7) -- (90:0.7);

					\foreach \x in {0,1,2}
					\draw[dotted,rotate=-72*\x]
					(18:0.7) -- (18:1.3) -- (-18:1.7) -- (-54:1.3) -- (-54:0.7);

					\draw
					(18:0.7) -- (18:1.3) -- (-18:1.7)
					(90:0.7) -- (18:0.7)
					(-18:1.7) -- (-54:1.3)
					(-54:1.3) -- (-54:0.7)
					(18:0.7) -- (-54:0.7)
					(18:1.3) -- (40:1.7) -- (60:1.5) -- (90:0.7);

					\draw[dotted]
					(90:0.7) -- (120:1.5) -- (140:1.7) -- (162:1.3);

					\draw[line width=1.5]
					(90:0.7) -- (18:0.7)
					(-18:1.7) -- (-54:1.3);

					\node at (90:0.4) {\small $\eee$}; 
					\node at (162:0.45) {\small $\ccc$};
					\node at (18:0.45) {\small $\ddd$};
					\node at (234:0.45) {\small $\aaa$};
					\node at (-54:0.45) {\small $\bbb$};

					\node at (32:0.8) {\small $\ddd$}; 
					\node at (70:0.75) {\small $\eee$};
					\node at (26:1.2) {\small $\bbb$};
					\node at (39:1.5) {\small $\aaa$};
					\node at (55:1.3) {\small $\ccc$};

					\node at (4:0.8) {\small $\aaa$}; 
					\node at (-40:0.8) {\small $\ccc$};
					\node at (-44:1.15) {\small $\eee$};
					\node at (-18:1.45) {\small $\ddd$};
					\node at (10:1.15) {\small $\bbb$};

					\node[draw,shape=circle, inner sep=0.5] at (0,0) {\small $1$};
					\node[draw,shape=circle, inner sep=0.5] at (45:1.15) {\small $2$};\node[draw,shape=circle, inner sep=0.5] at (-18:1) {\small $3$};
				\end{scope}
				\begin{scope}[xshift=0cm,yshift=3 cm,scale=0.9]

					\fill (0,0.7) circle (0.1);

					\foreach \x in {1,2,3}
					\draw[rotate=72*\x]
					(18:0.7) -- (90:0.7);

					\foreach \x in {0,1,2}
					\draw[dotted,rotate=-72*\x]
					(18:0.7) -- (18:1.3) -- (-18:1.7) -- (-54:1.3) -- (-54:0.7);

					\draw
					(18:0.7) -- (18:1.3) -- (-18:1.7)
					(90:0.7) -- (18:0.7)
					(-18:1.7) -- (-54:1.3)
					(-54:1.3) -- (-54:0.7)
					(18:0.7) -- (-54:0.7)
					(18:1.3) -- (40:1.7) -- (60:1.5) -- (90:0.7);

					\draw[dotted]
					(90:0.7) -- (120:1.5) -- (140:1.7) -- (162:1.3);

					\draw[line width=1.5]
					(90:0.7) -- (18:0.7)
					(-18:1.7) -- (-54:1.3);

					\node at (90:0.4) {\small $\eee$}; 
					\node at (162:0.45) {\small $\ccc$};
					\node at (18:0.45) {\small $\ddd$};
					\node at (234:0.45) {\small $\aaa$};
					\node at (-54:0.45) {\small $\bbb$};

					\node at (32:0.8) {\small $\ddd$}; 
					\node at (70:0.75) {\small $\eee$};
					\node at (26:1.2) {\small $\bbb$};
					\node at (39:1.5) {\small $\aaa$};
					\node at (55:1.3) {\small $\ccc$};

					\node at (4:0.8) {\small $\aaa$}; 
					\node at (-40:0.8) {\small $\bbb$};
					\node at (-44:1.15) {\small $\ddd$};
					\node at (-18:1.45) {\small $\eee$};
					\node at (10:1.15) {\small $\ccc$};

					\node[draw,shape=circle, inner sep=0.5] at (0,0) {\small $1$};
					\node[draw,shape=circle, inner sep=0.5] at (45:1.15) {\small $2$};\node[draw,shape=circle, inner sep=0.5] at (-18:1) {\small $3$};
				\end{scope}


				\begin{scope}[xshift=4 cm,yshift=3 cm,scale=0.9]

					\fill (0,0.7) circle (0.1);

					\foreach \x in {1,2,3}
					\draw[rotate=72*\x]
					(18:0.7) -- (90:0.7);

					\foreach \x in {0,1,2}
					\draw[dotted,rotate=-72*\x]
					(18:0.7) -- (18:1.3) -- (-18:1.7) -- (-54:1.3) -- (-54:0.7);

					\draw
					(18:0.7) -- (18:1.3) -- (-18:1.7)
					(90:0.7) -- (18:0.7)
					(-18:1.7) -- (-54:1.3)
					(18:0.7) -- (-54:0.7)
					(18:1.3) -- (40:1.7) -- (60:1.5) -- (90:0.7)
					(-54:0.7) -- (-54:1.3) -- (-90:1.7) -- (234:1.3) -- (234:0.7);

					\draw[dotted]
					(90:0.7) -- (120:1.5) -- (140:1.7) -- (162:1.3);

					\draw[line width=1.5]
					(90:0.7) -- (18:0.7)
					(-54:1.3) -- (-54:0.7);

					\node at (90:0.4) {\small $\eee$}; 
					\node at (162:0.45) {\small $\ccc$};
					\node at (18:0.45) {\small $\ddd$};
					\node at (234:0.45) {\small $\aaa$};
					\node at (-54:0.45) {\small $\bbb$};

					\node at (32:0.8) {\small $\ddd$}; 
					\node at (70:0.75) {\small $\eee$};
					\node at (26:1.2) {\small $\bbb$};
					\node at (39:1.5) {\small $\aaa$};
					\node at (55:1.3) {\small $\ccc$};

					\node at (4:0.8) {\small $\bbb$}; 
					\node at (-40:0.8) {\small $\ddd$};
					\node at (-44:1.15) {\small $\eee$};
					\node at (-18:1.45) {\small $\ccc$};
					\node at (10:1.15) {\small $\aaa$};

					\node at (-68:0.8) {\small $\ddd$}; 
					\node at (-112:0.8) {\small $\bbb$};
					\node at (-116:1.15) {\small $\aaa$};
					\node at (-90:1.45) {\small $\ccc$};
					\node at (-62:1.15) {\small $\eee$};

					\node[draw,shape=circle, inner sep=0.5] at (0,0) {\small $1$};
					\node[draw,shape=circle, inner sep=0.5] at (45:1.15) {\small $2$};\node[draw,shape=circle, inner sep=0.5] at (-18:1) {\small $3$};
					\node[draw,shape=circle, inner sep=0.5] at (-90:1.05) {\small $4$};
				\end{scope}

				\begin{scope}[xshift=-2 cm,yshift=0 cm,scale=0.9]

					\fill (0,0.7) circle (0.1);

					\foreach \x in {1,2,3}
					\draw[rotate=72*\x]
					(18:0.7) -- (90:0.7);

					\foreach \x in {0,1,2}
					\draw[dotted,rotate=-72*\x]
					(18:0.7) -- (18:1.3) -- (-18:1.7) -- (-54:1.3) -- (-54:0.7);

					\draw
					(18:0.7) -- (18:1.3) -- (-18:1.7)
					(90:0.7) -- (18:0.7)
					(-18:1.7) -- (-54:1.3)
					(-54:1.3) -- (-54:0.7)
					(18:0.7) -- (-54:0.7)
					(18:1.3) -- (40:1.7) -- (60:1.5) -- (90:0.7);

					\draw[dotted]
					(90:0.7) -- (120:1.5) -- (140:1.7) -- (162:1.3);

					\draw[line width=1.5]
					(90:0.7) -- (18:0.7)
					(18:1.3) -- (-18:1.7);

					\node at (90:0.4) {\small $\eee$}; 
					\node at (162:0.45) {\small $\ccc$};
					\node at (18:0.45) {\small $\ddd$};
					\node at (234:0.45) {\small $\aaa$};
					\node at (-54:0.45) {\small $\bbb$};

					\node at (32:0.8) {\small $\ddd$}; 
					\node at (70:0.75) {\small $\eee$};
					\node at (26:1.2) {\small $\bbb$};
					\node at (39:1.5) {\small $\aaa$};
					\node at (55:1.3) {\small $\ccc$};

					\node at (4:0.8) {\small $\bbb$}; 
					\node at (-40:0.8) {\small $\aaa$};
					\node at (-44:1.15) {\small $\ccc$};
					\node at (-18:1.45) {\small $\eee$};
					\node at (10:1.15) {\small $\ddd$};

					\node[draw,shape=circle, inner sep=0.5] at (0,0) {\small $1$};
					\node[draw,shape=circle, inner sep=0.5] at (45:1.15) {\small $2$};\node[draw,shape=circle, inner sep=0.5] at (-18:1) {\small $3$};
				\end{scope}
				\begin{scope}[xshift=2cm,yshift=0 cm,scale=0.9]

					\fill (0,0.7) circle (0.1);

					\foreach \x in {1,2,3}
					\draw[rotate=72*\x]
					(18:0.7) -- (90:0.7);

					\foreach \x in {0,1,2}
					\draw[dotted,rotate=-72*\x]
					(18:0.7) -- (18:1.3) -- (-18:1.7) -- (-54:1.3) -- (-54:0.7);

					\draw
					(18:0.7) -- (18:1.3) -- (-18:1.7)
					(90:0.7) -- (18:0.7)
					(-18:1.7) -- (-54:1.3)
					(-54:1.3) -- (-54:0.7)
					(18:0.7) -- (-54:0.7)
					(18:1.3) -- (40:1.7) -- (60:1.5) -- (90:0.7);

					\draw[dotted]
					(90:0.7) -- (120:1.5) -- (140:1.7) -- (162:1.3);

					\draw[line width=1.5]
					(90:0.7) -- (18:0.7)
					(18:1.3) -- (-18:1.7);

					\node at (90:0.4) {\small $\eee$}; 
					\node at (162:0.45) {\small $\ccc$};
					\node at (18:0.45) {\small $\ddd$};
					\node at (234:0.45) {\small $\aaa$};
					\node at (-54:0.45) {\small $\bbb$};

					\node at (32:0.8) {\small $\ddd$}; 
					\node at (70:0.75) {\small $\eee$};
					\node at (26:1.2) {\small $\bbb$};
					\node at (39:1.5) {\small $\aaa$};
					\node at (55:1.3) {\small $\ccc$};

					\node at (4:0.8) {\small $\ccc$}; 
					\node at (-40:0.8) {\small $\aaa$};
					\node at (-44:1.15) {\small $\bbb$};
					\node at (-18:1.45) {\small $\ddd$};
					\node at (10:1.15) {\small $\eee$};

					\node[draw,shape=circle, inner sep=0.5] at (0,0) {\small $1$};
					\node[draw,shape=circle, inner sep=0.5] at (45:1.15) {\small $2$};\node[draw,shape=circle, inner sep=0.5] at (-18:1) {\small $3$};
				\end{scope}
			\end{tikzpicture}
			\caption{$\aaa\ddd^2$, $\bbb\ddd^2$ or $\ccc\ddd^2$ is the unique $a^2b$-vertex type.}
			\label{fbd2}
		\end{figure}

		For each combination of $3$ vertices, we will apply the irrational angle lemma and AAD to get Table \ref{9}. The proof is completed by the contradictions in the last column of Table \ref{9}, similar to Proposition \ref{ad2,ce2}.

\begin{table}[htp]
	\centering
	\caption{All $a^2b$-vertices are $\aaa\ddd^2,\bbb\ddd^2$ or $\ccc\ddd^2$.}\label{9}
	\resizebox{\linewidth}{!}{\begin{tabular}{|c|c|c|c|c|c|c|}
			\hline
			\multicolumn{3}{|c|}{Vertex}& Angles & Irrational Angle Lemma & AAD  & Contradiction\\
			\hline\hline
			\multirow{11}{*}{$\aaa\ddd^2$} & \multirow{6}{*}{$\ddd\eee^3$} & \multirow{2}{*}{$\aaa\bbb^2$}& \multirow{2}{*}{$(2-2\ddd,\ddd,\frac13+\frac{1}3\ddd+\frac4f,\ddd,\frac23-\frac{1}3\ddd)$} & \multirow{2}{*}{$6n_1+n_5=3n_2+n_3+3n_4$}& \multirow{2}{*}{$\thin^{\ccc}\aaa^{\bbb}\thin^{\bbb}\ddd^{\eee}\thick^{\eee}\ddd^{\bbb}\thin$} &No $\bbb\ccc\cdots$ by the\\
			&&&&&& irrational angle lemma\\
			\cline{3-7}
			&& $\aaa\bbb\ccc$ & $(1-\frac{12}f,1-\ccc+\frac{12}f,\ccc,\frac12+\frac6f,\frac12-\frac2f)$&$n_2=n_3$&\multirow{3}{*}{$\thick^{\eee}\ddd^{\bbb}\thin^{\ccc}\eee^{\ddd}\thick^{\ddd}\eee^{\ccc}\thin^{\ccc}\eee^{\ddd}\thick$}  &\multirow{3}{*}{No $\ccc^2\cdots$}\\
			\cline{3-5}
			& & $\bbb^3$ & $(2-2\ddd,\frac23,-\frac13+\frac{4}3\ddd+\frac4f,\ddd,\frac23-\frac{1}3\ddd)$ & $6n_1+n_5=4n_3+3n_4$& & \\
			\cline{3-5}
			& & $\bbb^2\ccc$& $(2-2\ddd,\frac53-\frac{4}3\ddd-\frac4f,-\frac43+\frac{8}3\ddd+\frac8f,\ddd,\frac23-\frac{1}3\ddd)$ & $6n_1+4n_2+n_5=8n_3+3n_4$& & \\
			\cline{3-7}
			& & $\bbb\ccc^2$& $(2-2\ddd,-\frac43+\frac{8}3\ddd+\frac8f,\frac53-\frac{4}3\ddd-\frac4f,\ddd,\frac23-\frac{1}3\ddd)$ & $6n_1+4n_3+n_5=8n_2+3n_4$& \multirow{3}{*}{$\thin^{\ccc}\aaa^{\bbb}\thin^{\bbb}\ddd^{\eee}\thick^{\eee}\ddd^{\bbb}\thin$}&No $\bbb^2\cdots$\\
			\cline{2-5}
			\cline{7-7}
			& \multirow{5}{*}{$\eee^4$} & $\aaa\bbb^2$& $(2-2\ddd,\ddd,\frac12+\frac4f,\ddd,\frac12)$ & $2n_1=n_2+n_4$&&\multirow{2}{*}{No $\bbb\ccc\cdots$}\\
			\cline{3-5}
			& & $\bbb^3$& $(2-2\ddd,\frac23,-\frac16+\ddd+\frac4f,\ddd,\frac12)$ & $2n_1=n_3+n_4$& & \\
			\cline{3-7}
			&& $\aaa\bbb\ccc$& $(1-\frac8f,1-\ccc+\frac8f,\ccc,\frac12+\frac4f,\frac12)$ & $n_2=n_3$& \multirow{2}{*}{$\thick^{\ddd}\eee^{\ccc}\thin^{\ccc}\eee^{\ddd}\thick^{\ddd}\eee^{\ccc}\thin^{\ccc}\eee^{\ddd}\thick$} &\multirow{2}{*}{No $\ccc^2\cdots$}\\
			\cline{3-5}
			& & $\bbb^2\ccc$& $(2-2\ddd,\frac32-\ddd-\frac4f,-1+2\ddd+\frac8f,\ddd,\frac12)$ & $2n_1+n_2=2n_3+n_4$& & \\
			\cline{3-7}
			&& $\bbb\ccc^2$& $(2-2\ddd,-1+2\ddd+\frac8f,\frac32-\ddd-\frac4f,\ddd,\frac12)$ & $2n_1+n_3=2n_2+n_4$& $\thin^{\ccc}\aaa^{\bbb}\thin^{\bbb}\ddd^{\eee}\thick^{\eee}\ddd^{\bbb}\thin$ &No $\bbb^2\cdots$\\
			\cline{3-7}
			\hline
			\multirow{6}{*}{$\bbb\ddd^2$} & \multirow{3}{*}{$\ddd\eee^3$} & $\aaa^2\bbb$& $(\ddd,2-2\ddd,\frac13+\frac{1}3\ddd+\frac4f,\ddd,\frac23-\frac{1}3\ddd)$ &$3n_1+n_3+3n_4=6n_2+n_5$ & \multirow{3}{*}{$\thick^{\eee}\ddd^{\bbb}\thin^{\ccc}\eee^{\ddd}\thick^{\ddd}\eee^{\ccc}\thin^{\ccc}\eee^{\ddd}\thick$}&No $\bbb\ccc\cdots$\\
			\cline{3-5}
			\cline{7-7}
			& & $\aaa\bbb^2$& $(-2+4\ddd,2-2\ddd,\frac73-\frac{8}3\ddd+\frac4f,\ddd,\frac23-\frac{1}3\ddd)$ & $12n_1+3n_4=6n_2+8n_3+n_5$& &\multirow{7}{*}{No $\ccc^2\cdots$}\\
			\cline{3-5}
			& & $\aaa\bbb\ccc$& $(1-\ccc+\frac{12}f,1-\frac{12}f,\ccc,\frac12+\frac6f,\frac12-\frac2f)$ & $n_1=n_3$& &\\
			\cline{3-6}
			\cline{2-2}
			& \multirow{3}{*}{$\eee^4$} & $\aaa^2\bbb$& $(\ddd,2-2\ddd,\frac12+\frac4f,\ddd,\frac12)$ &$n_1+n_4=2n_2$& \multirow{3}{*}{$\thick^{\ddd}\eee^{\ccc}\thin^{\ccc}\eee^{\ddd}\thick^{\ddd}\eee^{\ccc}\thin^{\ccc}\eee^{\ddd}\thick$} &\\
			\cline{3-5}
			& & $\aaa\bbb^2$& $(-2+4\ddd,2-2\ddd,\frac52-3\ddd+\frac4f,\ddd,\frac12)$ & $4n_1+n_4=2n_2+3n_3$&  &\\
			\cline{3-5}
			& & $\aaa\bbb\ccc$& $(1-\ccc+\frac8f,1-\frac8f,\ccc,\frac12+\frac4f,\frac12)$ & $n_1=n_3$& &\\
			\cline{1-6}
			\multirow{6}{*}{$\ccc\ddd^2$} & \multirow{3}{*}{$\ddd\eee^3$} & $\aaa^2\bbb$& $(\frac53-\frac43\ddd-\frac4f,-\frac43+\frac83\ddd+\frac8f,2-2\ddd,\ddd,\frac23-\frac13\ddd)$ & $4n_1+6n_3+n_5=8n_2+3n_4$& \multirow{2}{*}{$\thick^{\eee}\ddd^{\bbb}\thin^{\ccc}\eee^{\ddd}\thick^{\ddd}\eee^{\ccc}\thin^{\ccc}\eee^{\ddd}\thick$} &\\
			\cline{3-5}
			& & $\aaa\bbb^2$ & $(-\frac43+\frac83\ddd+\frac8f,\frac53-\frac43\ddd-\frac4f,2-2\ddd,\ddd,\frac23-\frac13\ddd)$ &$8n_1+3n_4=4n_2+6n_3+n_5$& & \\
			\cline{3-7}
			&& $\aaa\bbb\ccc$& $(1-\bbb+\frac{12}f,\bbb,1-\frac{12}f,\frac12+\frac6f,\frac12-\frac2f)$ &$n_1=n_2$& $\thin^{\aaa}\ccc^{\eee}\thin^{\bbb}\ddd^{\eee}\thick^{\eee}\ddd^{\bbb}\thin$ &No $\bbb\eee\cdots$\\
			\cline{3-7}
			\cline{2-2}
			&\multirow{3}{*}{$\eee^4$}& $\aaa^2\bbb$& $(\frac32-\ddd-\frac4f,-1+2\ddd+\frac8f,2-2\ddd,\ddd,\frac12)$ & $n_1+2n_3=2n_2+n_4$ & \multirow{2}{*}{$\thick^{\ddd}\eee^{\ccc}\thin^{\ccc}\eee^{\ddd}\thick^{\ddd}\eee^{\ccc}\thin^{\ccc}\eee^{\ddd}\thick$}&\multirow{2}{*}{No $\ccc^2\cdots$}\\
			\cline{3-5}
			& & $\aaa\bbb^2$ & $(-1+2\ddd+\frac8f,\frac32-\ddd-\frac4f,2-2\ddd,\ddd,\frac12)$ &$2n_1+n_4=n_2+2n_3$ & & \\
			\cline{3-7}
			&& $\aaa\bbb\ccc$& $(1-\bbb+\frac8f,\bbb,1-\frac8f,\frac12+\frac4f,\frac12)$ &$n_1=n_2$& $\thin^{\aaa}\ccc^{\eee}\thin^{\bbb}\ddd^{\eee}\thick^{\eee}\ddd^{\bbb}\thin$ &No $\bbb\eee\cdots$\\
			\hline
	\end{tabular}}

\end{table}

	\end{proof}

	\section{Tilings with a unique $a^2b$-vertex type $\bbb\ddd\eee$}
	\label{bde as the unique $a^2b$-vertex}
	We divide the discussion into five cases in Figure \ref{a4b_pentagon}, according to the location of the vertex $H$ of a special tile (indicated by $\bullet$). The degree of $H$ is $3,4$ or $5$.

	\begin{figure}[htp]
		\centering
		\begin{tikzpicture}[>=latex,scale=1]


			\foreach \a in {1,2,3,4,5}
			{
				\begin{scope}[xshift=-2cm + 2.5*\a cm]

					\draw
					(-54:0.8) -- (18:0.8) -- (90:0.8) -- (162:0.8) -- (234:0.8);

					\draw[line width=1.5]
					(-54:0.8) -- (-126:0.8);

					\fill (18+72*\a:0.8) circle (0.1);

					\node at (90:0.5) {\small $\alpha$};
					\node at (162:0.5) {\small $\beta$};
					\node at (18:0.5) {\small $\gamma$};
					\node at (234:0.5) {\small $\delta$};
					\node at (-54:0.5) {\small $\epsilon$};

					\node at (0,0) {\a};

				\end{scope}
			}

		\end{tikzpicture}
		\caption{Special tile for the edge combination $a^4b$.}
		\label{a4b_pentagon}
	\end{figure}

	\begin{proposition}\label{bde,a2b}
		There is no non-symmetric $a^4b$-tiling with general angles, such that $\bbb\ddd\eee$ is the unique $a^2b$-vertex type and $\aaa^{2}\bbb$ is a vertex.
	\end{proposition}

	\begin{proof}
		By Lemma \ref{hdeg}, one of $\{\aaa\bbb\ccc,\aaa^2\ccc,\aaa\ccc^{2},  \bbb^{2}\ccc, \bbb\ccc^{2}, \ccc^{3}, \aaa\ccc^{3}, \bbb\ccc^{3},  \ccc^{4}, \ccc^{5}\}$ must appear. In any of these cases, we have $n_4=n_5$ by the irrational angle lemma.

		If there is no degree $\ge4$ vertex $\ddd\eee\cdots$, by $\bbb\ddd\eee$ being the unique $a^2b$-vertex type and $n_4=n_5$, we get $\ddd\cdots=\ddd\eee\cdots=\bbb\ddd\eee$. Then $\#\bbb\ge\#\bbb\ddd\eee+\#\aaa^2\bbb=f+\#\aaa^2\bbb>f$, a contradiction. Therefore, we know degree $\ge4$ vertex $\ddd\eee\cdots$ must appear.

		For each case of $\{\aaa\bbb\ccc,\aaa\ccc^{2},  \bbb^{2}\ccc, \bbb\ccc^{2}, \ccc^{3}, \aaa\ccc^{3}, \bbb\ccc^{3},  \ccc^{4}, \ccc^{5}\}$, we summarize the computations of angle values in Table \ref{table bde+a^2b}, excluding the symmetric ($\bbb=\ccc$) Case $\{\bbb\ddd\eee,\aaa^2\bbb,\aaa^2\ccc\}$.

		\subsection*{Case. $\{\bbb\ddd\eee,\aaa^2\bbb,\ccc^3\}$}\label{bde 2ab 3c}
		The angle sums of $\bbb\ddd\eee,\aaa^2\bbb,\ccc^3$ and the angle sum for the pentagon imply
		\begin{align*}
			\aaa=\tfrac13+\tfrac4f,\;
			\bbb=\tfrac43-\tfrac8f,\;
			\ccc=\tfrac23,\;
			\ddd+\eee=\tfrac23+\tfrac8f.
		\end{align*}
		By the angle values and $n_4=n_5$, we get $\ddd\cdots=\bbb\ddd\eee$, $\aaa^2\ddd\eee$, $\aaa\ccc\ddd\eee$, $\ddd^2\eee^2$ or $\aaa^3\ddd\eee$.

		If $\aaa^2\ddd\eee$ is a vertex, then $\text{AVC}\sub\{\bbb\ddd\eee,\aaa^2\bbb,\ccc^3,\aaa^4,\aaa^2\ddd\eee,\ddd^2\eee^2\}$. We have the AAD $\thick^{\eee}\ddd^{\bbb}\thin\aaa\thin$ at $\aaa^2\ddd\eee$. This gives a vertex $\bbb^2\cdots$ or $\bbb\ccc\cdots$, contradicting the AVC. Similarly, $\ddd^2\eee^2$ and $\aaa^3\ddd\eee$ are not vertices.

		If $\aaa\ccc\ddd\eee$ is a vertex, then $\text{AVC}\sub\{\bbb\ddd\eee,\aaa^2\bbb,\ccc^3, \aaa^3\ccc,\aaa\ccc\ddd\eee\}$. We have the  AAD $\aaa\ccc\ddd\eee=\thin^{\bbb}\ddd^{\eee}\thick^{\ddd}\eee^{\ccc}\thin^{\ccc}\aaa^{\bbb}\thin^{\aaa}\ccc^{\eee}\thin$ or $\thin^{\bbb}\ddd^{\eee}\thick^{\ddd}\eee^{\ccc}\thin^{\ccc}\aaa^{\bbb}\thin^{\eee}\ccc^{\aaa}\thin$. This gives a vertex $\aaa\bbb\cdots=\thin^{\bbb}\aaa^{\ccc}\thin\bbb\thin\aaa\thin$. This implies that $\bbb^2\cdots$ or $\bbb\ccc\cdots$ is a vertex, contradicting the AVC.

		Therefore, we know $\ddd\eee\cdots=\bbb\ddd\eee$, a contradiction. Then $\#\bbb\ge\#\bbb\ddd\eee+\#\aaa^2\bbb>f$, a contradiction.

		\subsection*{Case. $\{\bbb\ddd\eee,\aaa^2\bbb,\aaa\ccc^2\}$}\label{bde 2ab a2c}
		The angle sums of $\bbb\ddd\eee,\aaa^{2}\bbb,\aaa\ccc^2$ and the angle sum for the pentagon imply
		\begin{align*}
			\aaa=\tfrac8f,\;
			\bbb=2-\tfrac{16}f,\;
			\ccc=1-\tfrac4f,\;
			\ddd+\eee=\tfrac{16}f.
		\end{align*}
		By $\bbb\neq\ccc$, we have $f\neq12$. By the angle values and $n_4=n_5$, we get the
		\[
		\text{AVC}\sub\{\bbb\ddd\eee,\aaa^2\bbb,\aaa\ccc^2,\aaa^{\frac{f+4}8}\ccc,\aaa^{\frac{f}4}, \aaa^x\ddd^{\frac{f-4x}8}\eee^{\frac{f-4x}8},\aaa^y\ccc\ddd^{\frac{f-8y+4}{16}}\eee^{\frac{f-8y+4}{16}}\},
		\]
		where $x=0,1,2,\cdots,\lfloor\tfrac{f}4\rfloor$, $y=0,1,2,\cdots,\lfloor\tfrac{f+4}8\rfloor$. Then we know there is no $\bbb^2\cdots$ and $\bbb\ccc\cdots$ vertex, which further implies that (by AAD) the $\text{AVC}\sub\{\bbb\ddd\eee,\aaa^2\bbb,\aaa\ccc^2,\aaa\ccc\ddd\eee,\aaa^{\frac{f+4}8}\ccc,\aaa^{\frac{f}4}\}$.

		If $\aaa\ccc\ddd\eee$ is a vertex, we have the unique AAD $\aaa\ccc\ddd\eee=\thin^{\ccc}\aaa^{\bbb}\thin^{\aaa}\ccc^{\eee}\thin^{\bbb}\ddd^{\eee}\thick^{\ddd}\eee^{\ccc}\thin$. This gives a vertex $\aaa\bbb\cdots=\thin^{\ccc}\aaa^{\bbb}\thin^{\ccc}\aaa^{\bbb}\thin^{\ddd}\bbb^{\aaa}\thin$. This gives a vertex $\bbb\ccc\cdots$, contradicting the AVC. Therefore, we know $\ddd\eee\cdots=\bbb\ddd\eee$, a contradiction.

		\begin{table*}[htp]
			\centering
			\caption{$\bbb\ddd\eee$, $\aaa^2\bbb$ are vertices.}\label{table bde+a^2b}
			~\\
			\resizebox{\linewidth}{!}{\begin{tabular}{|c|c|c|c|c|}
					\hline
					\multicolumn{2}{|c|}{Vertex} & Angles & Irrational Angle Lemma & Contradiction \\
					\hline\hline
					\multirow{9}{*}{$\bbb\ddd\eee,\aaa^{2}\bbb$}&$\aaa\bbb\ccc$&\multirow{2}{*}{$(\frac12+\frac2f,1-\frac4f,\frac12+\frac2f,\ddd,1-\ddd+\frac4f)$}
					&\multirow{9}{*}{$n_4=n_5$} &\multirow{9}{*}{No degree $\ge4$ vertex $\ddd\eee\cdots$}\\
					\cline{2-2}
					&$\bbb\ccc^{2}$&
					&&\\
					\cline{2-3}
					&$\bbb^{2}\ccc$&$(\frac35+\frac4{5f},\frac45-\frac8{5f},\frac25+\frac{16}{5f},\ddd,\frac65-\ddd+\frac8{5f})$
					& &\\
					\cline{2-3}
					&$\aaa\ccc^{3}$&$(\frac12+\frac6f,1-\frac{12}f,\frac12-\frac2f,\ddd,1-\ddd+\frac{12}f)$
					& &\\
					\cline{2-3}
					&$\bbb\ccc^{3}$&$(\frac35+\frac{12}{5f},\frac45-\frac{24}{5f},\frac25+\frac8{5f},\ddd,\frac65-\ddd+\frac{24}{5f})$
					& & \\
					\cline{2-3}
					&$\ccc^{4}$&$(\frac12+\frac4f,1-\frac8f,\frac12,\ddd,1-\ddd+\frac8f)$
					& &\\
					\cline{2-3}
					&$\ccc^{5}$&$(\frac35+\frac4f,\frac45-\frac8f,\frac25,\ddd,\frac65-\ddd+\frac8f)$&&\\
					\cline{2-3}
					&$\ccc^3$&$(\frac13+\frac4f,\frac43-\frac8f,\frac23,\ddd,\frac23-\ddd+\frac8f)$&&\\
					\cline{2-3}
					&$\aaa\ccc^2$&$(\frac8f,2-\frac{16}{f},1-\frac4f,\ddd,-\ddd+\frac{16}{f})$&&\\
					\hline
			\end{tabular}}
		\end{table*}

	\vspace{6pt}

	The other cases lead to contradictions in a similar and much easier way.
	\end{proof}

	\begin{lemma}\label{bde H}
		In a non-symmetric $a^4b$-tiling with general angles, if $\bbb\ddd\eee$ is the unique $a^2b$-vertex type, then $\aaa^3$, $\aaa^2\ccc$, $\aaa\bbb\ccc$ or $\aaa\ccc^2$ must appear. In each case, three linearly independent vertices in Table \ref{10} must appear.
	\end{lemma}

	\begin{proof}
	By Proposition \ref{bde,a2b}, we know that $\aaa^2\bbb$ is not a vertex. We divide the discussion into five cases in Figure \ref{a4b_pentagon}, according to the location of the vertex $H$ of a special tile. Our strategy is to identify three linearly independent vertices.

	\subsection*{Case $H=\aaa_1\cdots$}
	Let the first of Figure \ref{a4b_pentagon} be the center tile $T_1$ in the partial neighborhood in Figure \ref{bde(a)}. Since $\bbb\ddd\eee$ is the unique $a^2b$-vertex type and $H=\aaa_1\cdots$, the degree $3$ vertex $\eee_1\cdots=\bbb\ddd\eee$ determines $T_3$, $T_4$. The two pictures show two possible arrangements of $T_5$. In the right-hand side of Figure \ref{bde(a)}, we have the degree $3$ vertex $\bbb_1\ddd_5\cdots=\bbb_1\ddd_5\eee_6$, which determines $T_6$.
	\begin{figure}[htp]
		\centering
		\begin{tikzpicture}[>=latex,scale=0.9]


			\begin{scope}[xshift=-3cm]

				\fill (0,0.7) circle (0.1);

				\foreach \x in {1,2,3}
				\draw[rotate=72*\x]
				(18:0.7) -- (90:0.7);

				\foreach \x in {0,1,2}
				\draw[dotted,rotate=-72*\x]
				(18:0.7) -- (18:1.3) -- (-18:1.7) -- (-54:1.3) -- (-54:0.7);

				\draw
				(18:0.7) -- (18:1.3) -- (-18:1.7)
				(90:0.7) -- (18:0.7)
				(18:0.7) -- (-54:0.7)
				(-54:0.7) -- (-54:1.3) -- (-90:1.7) -- (234:1.3) -- (234:0.7)
				(-18:1.7) -- (-54:1.3)
				(234:1.3) -- (198:1.7) -- (162:1.3) -- (162:0.7);

				\draw[dotted]
				(90:0.7) -- (120:1.5) -- (140:1.7) -- (162:1.3) 
				(18:1.3) -- (40:1.7) -- (60:1.5) -- (90:0.7); 

				\draw[line width=1.5]
				(234:0.7) -- (-54:0.7)
				(-18:1.7) -- (-54:1.3)
				(234:1.3) -- (198:1.7);

				\node at (90:0.4) {\small $\aaa$}; 
				\node at (162:0.45) {\small $\bbb$};
				\node at (18:0.45) {\small $\ccc$};
				\node at (234:0.45) {\small $\ddd$};
				\node at (-54:0.45) {\small $\eee$};

				\node at (32:0.8) {\small $\bbb$}; 

				\node at (4:0.8) {\small $\aaa$}; 
				\node at (-40:0.8) {\small $\bbb$};
				\node at (-44:1.15) {\small $\ddd$};
				\node at (-18:1.45) {\small $\eee$};
				\node at (10:1.15) {\small $\ccc$};

				\node at (-68:0.8) {\small $\ddd$}; 
				\node at (-112:0.8) {\small $\eee$};
				\node at (-116:1.15) {\small $\ccc$};
				\node at (-90:1.45) {\small $\aaa$};
				\node at (-62:1.15) {\small $\bbb$};

				\node at (176:0.8) {\small $\aaa$}; 
				\node at (220:0.8) {\small $\bbb$};
				\node at (224:1.15) {\small $\ddd$};
				\node at (198:1.45) {\small $\eee$};
				\node at (170:1.2) {\small $\ccc$};

				\node at (146:0.75) {\small $\ccc$}; 

				\node[draw,shape=circle, inner sep=0.5] at (0,0) {\small $1$};
				\node[draw,shape=circle, inner sep=0.5] at (45:1.15) {\small $2$};\node[draw,shape=circle, inner sep=0.5] at (-18:1) {\small $3$};
				\node[draw,shape=circle, inner sep=0.5] at (-90:1.05) {\small $4$};
				\node[draw,shape=circle, inner sep=0.5] at (198:1) {\small $5$};
				\node[draw,shape=circle, inner sep=0.5] at (135:1.07) {\small $6$};
			\end{scope}


			\begin{scope}[xshift=3cm]

				\fill (0,0.7) circle (0.1);

				\foreach \x in {1,2,3}
				\draw[rotate=72*\x]
				(18:0.7) -- (90:0.7);

				\foreach \x in {0,1,2}
				\draw[dotted,rotate=-72*\x]
				(18:0.7) -- (18:1.3) -- (-18:1.7) -- (-54:1.3) -- (-54:0.7);

				\draw
				(18:0.7) -- (18:1.3) -- (-18:1.7)
				(90:0.7) -- (18:0.7)
				(18:0.7) -- (-54:0.7)
				(-54:0.7) -- (-54:1.3) -- (-90:1.7) -- (234:1.3) -- (234:0.7)
				(-18:1.7) -- (-54:1.3)
				(234:1.3) -- (198:1.7) -- (162:1.3) -- (162:0.7)
				(90:0.7) -- (120:1.5) -- (140:1.7) -- (162:1.3);

				\draw[dotted]
				(18:1.3) -- (40:1.7) -- (60:1.5) -- (90:0.7); 

				\draw[line width=1.5]
				(234:0.7) -- (-54:0.7)
				(-18:1.7) -- (-54:1.3)
				(162:1.3) -- (162:0.7);

				\node at (90:0.4) {\small $\aaa$}; 
				\node at (162:0.45) {\small $\bbb$};
				\node at (18:0.45) {\small $\ccc$};
				\node at (234:0.45) {\small $\ddd$};
				\node at (-54:0.45) {\small $\eee$};

				\node at (4:0.8) {\small $\aaa$}; 
				\node at (-40:0.8) {\small $\bbb$};
				\node at (-44:1.15) {\small $\ddd$};
				\node at (-18:1.45) {\small $\eee$};
				\node at (10:1.15) {\small $\ccc$};

				\node at (-68:0.8) {\small $\ddd$}; 
				\node at (-112:0.8) {\small $\eee$};
				\node at (-116:1.15) {\small $\ccc$};
				\node at (-90:1.45) {\small $\aaa$};
				\node at (-62:1.15) {\small $\bbb$};

				\node at (176:0.8) {\small $\ddd$}; 
				\node at (220:0.8) {\small $\bbb$};
				\node at (224:1.15) {\small $\aaa$};
				\node at (198:1.45) {\small $\ccc$};
				\node at (170:1.2) {\small $\eee$};

				\node at (146:0.75) {\small $\eee$}; 
				\node at (112:0.75) {\small $\ccc$};
				\node at (152:1.2) {\small $\ddd$};
				\node at (140:1.5) {\small $\bbb$};
				\node at (125:1.35) {\small $\aaa$};

				\node[draw,shape=circle, inner sep=0.5] at (0,0) {\small $1$};
				\node[draw,shape=circle, inner sep=0.5] at (45:1.15) {\small $2$};\node[draw,shape=circle, inner sep=0.5] at (-18:1) {\small $3$};
				\node[draw,shape=circle, inner sep=0.5] at (-90:1.05) {\small $4$};
				\node[draw,shape=circle, inner sep=0.5] at (198:1) {\small $5$};
				\node[draw,shape=circle, inner sep=0.5] at (135:1.07) {\small $6$};
			\end{scope}

		\end{tikzpicture}
		\caption{$H=\aaa_1\cdots$.}
		\label{bde(a)}
	\end{figure}

	In the left-hand side of Figure \ref{bde(a)}, we have $\aaa_5\bbb_1\cdots=\aaa\bbb^2$ or $\aaa\bbb\ccc$, $\aaa_3\ccc_1\cdots=\aaa^2\ccc$, $\aaa\bbb\ccc$ or $\aaa\ccc^2$. By $\bbb\neq\ccc$, we get $(\aaa_5\bbb_1\cdots,\aaa_3\ccc_1\cdots)=(\aaa\bbb^2,\aaa^2\ccc)$, $(\aaa\bbb\ccc,\aaa^2\ccc)$ or $(\aaa\bbb\ccc,\aaa\bbb\ccc)$. If $\aaa_5\bbb_1\cdots=\aaa_3\ccc_1\cdots=\aaa\bbb\ccc$, then we get
	\begin{align*}
	H
	&=\{\aaa^3, \aaa^4, \aaa^3\bbb, \aaa^3\ccc, \aaa^2\ddd^2, \aaa^2\ddd\eee, \aaa^2\eee^2, \aaa^5, \aaa^4\bbb, \aaa^4\ccc, \aaa^3\bbb^2, \aaa^3\ddd^2, \aaa^3\ddd\eee, \aaa^3\eee^2,\\
	&\quad\,\,\,\aaa^2\bbb\ddd^2,\aaa^2\bbb\eee^2,\aaa^2\ccc\ddd^2,\aaa^2\ccc\eee^2,\aaa\ddd^3\eee,\aaa\ddd\eee^3\}.
    \end{align*}
	In the right-hand side of Figure \ref{bde(a)}, we have $\aaa_3\ccc_1\cdots=\aaa_2\aaa_3\ccc_1$, $\aaa_3\bbb_2\ccc_1$ or $\aaa_3\ccc_1\ccc_2$. Each case induces the following $H$ respectively:
	\begin{align*}
		H
		&=\{\aaa\bbb\ccc, \aaa\ccc^2, \aaa\bbb^2\ccc, \aaa\bbb\ccc^2, \aaa\ccc^3, \aaa\bbb^3\ccc, \aaa\bbb^2\ccc^2, \aaa\bbb\ccc^3, \aaa\bbb\ccc\ddd^2, \aaa\bbb\ccc\eee^2, \aaa\ccc^4, \aaa\ccc^2\ddd^2,  \\
		&\quad\,\,\,\aaa\ccc^2\ddd\eee, \aaa\ccc^2\eee^2\} \\
		&\text{or}\,\,\{\aaa^2\ccc, \aaa^3\ccc, \aaa\ccc\ddd^2, \aaa^4\ccc, \aaa^2\ccc\ddd^2, \aaa^2\ccc\eee^2, \aaa\ccc^2\ddd^2\} \\
		&\text{or}\,\,\{\aaa^2\ccc, \aaa^3\ccc, \aaa^2\bbb\ccc, \aaa\ccc\ddd\eee, \aaa\ccc\eee^2, \aaa^4\ccc, \aaa^3\bbb\ccc, \aaa^2\bbb^2\ccc, \aaa^2\ccc\ddd^2, \aaa^2\ccc\ddd\eee, \aaa^2\ccc\eee^2\}.
	\end{align*}

	\subsection*{Case $H=\bbb_1\cdots$}
	Let the second of Figure \ref{a4b_pentagon} be the center tile $T_1$ in the partial neighborhood in Figure \ref{bde(b)}. Since $\bbb\ddd\eee$ is an $a^2b$-vertex and $H=\bbb_1\cdots$, we have the degree $3$ vertex $\eee_1\cdots=\bbb\ddd\eee$, which determines $T_4$, $T_5$. The four pictures show four possible arrangements of $T_3$.

	\begin{figure}[htp]
		\centering
		\begin{tikzpicture}[>=latex,scale=0.85]


			\begin{scope}[xshift=-5cm]

				\fill (0,0.7) circle (0.1);

				\foreach \x in {1,2,3}
				\draw[rotate=72*\x]
				(18:0.7) -- (90:0.7);

				\foreach \x in {0,1,2}
				\draw[dotted,rotate=-72*\x]
				(18:0.7) -- (18:1.3) -- (-18:1.7) -- (-54:1.3) -- (-54:0.7);

				\draw
				(18:0.7) -- (18:1.3) -- (-18:1.7)
				(90:0.7) -- (18:0.7)
				(18:0.7) -- (-54:0.7)
				(-54:0.7) -- (-54:1.3) -- (-90:1.7) -- (234:1.3) -- (234:0.7)
				(-18:1.7) -- (-54:1.3)
				(234:1.3) -- (198:1.7) -- (162:1.3) -- (162:0.7);

				\draw[dotted]
				(90:0.7) -- (120:1.5) -- (140:1.7) -- (162:1.3) 
				(18:1.3) -- (40:1.7) -- (60:1.5) -- (90:0.7); 

				\draw[line width=1.5]
				(234:0.7) -- (162:0.7)
				(-18:1.7) -- (18:1.3)
				(234:1.3) -- (-90:1.7);

				\node at (90:0.4) {\small $\bbb$}; 
				\node at (162:0.45) {\small $\ddd$};
				\node at (18:0.45) {\small $\aaa$};
				\node at (234:0.45) {\small $\eee$};
				\node at (-54:0.45) {\small $\ccc$};

				\node at (4:0.8) {\small $\bbb$}; 
				\node at (-40:0.8) {\small $\aaa$};
				\node at (-44:1.15) {\small $\ccc$};
				\node at (-18:1.45) {\small $\eee$};
				\node at (10:1.15) {\small $\ddd$};

				\node at (-68:0.8) {\small $\aaa$}; 
				\node at (-112:0.8) {\small $\bbb$};
				\node at (-116:1.15) {\small $\ddd$};
				\node at (-90:1.45) {\small $\eee$};
				\node at (-62:1.15) {\small $\ccc$};

				\node at (176:0.8) {\small $\eee$}; 
				\node at (220:0.8) {\small $\ddd$};
				\node at (224:1.15) {\small $\bbb$};
				\node at (198:1.45) {\small $\aaa$};
				\node at (170:1.2) {\small $\ccc$};

				\node at (146:0.75) {\small $\bbb$}; 

				\node[draw,shape=circle, inner sep=0.5] at (0,0) {\small $1$};
				\node[draw,shape=circle, inner sep=0.5] at (45:1.15) {\small $2$};\node[draw,shape=circle, inner sep=0.5] at (-18:1) {\small $3$};
				\node[draw,shape=circle, inner sep=0.5] at (-90:1.05) {\small $4$};
				\node[draw,shape=circle, inner sep=0.5] at (198:1) {\small $5$};
				\node[draw,shape=circle, inner sep=0.5] at (135:1.07) {\small $6$};
			\end{scope}


			\begin{scope}[xshift=-1.5cm]

				\fill (0,0.7) circle (0.1);

				\foreach \x in {1,2,3}
				\draw[rotate=72*\x]
				(18:0.7) -- (90:0.7);

				\foreach \x in {0,1,2}
				\draw[dotted,rotate=-72*\x]
				(18:0.7) -- (18:1.3) -- (-18:1.7) -- (-54:1.3) -- (-54:0.7);

				\draw
				(18:0.7) -- (18:1.3) -- (-18:1.7)
				(90:0.7) -- (18:0.7)
				(18:0.7) -- (-54:0.7)
				(-54:0.7) -- (-54:1.3) -- (-90:1.7) -- (234:1.3) -- (234:0.7)
				(-18:1.7) -- (-54:1.3)
				(234:1.3) -- (198:1.7) -- (162:1.3) -- (162:0.7)
				(18:1.3) -- (40:1.7) -- (60:1.5) -- (90:0.7);

				\draw[dotted]
				(90:0.7) -- (120:1.5) -- (140:1.7) -- (162:1.3); 

				\draw[line width=1.5]
				(234:0.7) -- (162:0.7)
				(-18:1.7) -- (18:1.3)
				(234:1.3) -- (-90:1.7)
				(40:1.7) -- (60:1.5);

				\node at (90:0.4) {\small $\bbb$}; 
				\node at (162:0.45) {\small $\ddd$};
				\node at (18:0.45) {\small $\aaa$};
				\node at (234:0.45) {\small $\eee$};
				\node at (-54:0.45) {\small $\ccc$};

				\node at (32:0.8) {\small $\aaa$}; 

				\node at (4:0.8) {\small $\ccc$}; 
				\node at (-40:0.8) {\small $\aaa$};
				\node at (-44:1.15) {\small $\bbb$};
				\node at (-18:1.45) {\small $\ddd$};
				\node at (10:1.15) {\small $\eee$};

				\node at (-68:0.8) {\small $\aaa$}; 
				\node at (-112:0.8) {\small $\bbb$};
				\node at (-116:1.15) {\small $\ddd$};
				\node at (-90:1.45) {\small $\eee$};
				\node at (-62:1.15) {\small $\ccc$};

				\node at (176:0.8) {\small $\eee$}; 
				\node at (220:0.8) {\small $\ddd$};
				\node at (224:1.15) {\small $\bbb$};
				\node at (198:1.45) {\small $\aaa$};
				\node at (170:1.2) {\small $\ccc$};

				\node at (146:0.75) {\small $\bbb$}; 

				\node[draw,shape=circle, inner sep=0.5] at (0,0) {\small $1$};
				\node[draw,shape=circle, inner sep=0.5] at (45:1.15) {\small $2$};\node[draw,shape=circle, inner sep=0.5] at (-18:1) {\small $3$};
				\node[draw,shape=circle, inner sep=0.5] at (-90:1.05) {\small $4$};
				\node[draw,shape=circle, inner sep=0.5] at (198:1) {\small $5$};
				\node[draw,shape=circle, inner sep=0.5] at (135:1.07) {\small $6$};
			\end{scope}


			\begin{scope}[xshift=2cm]

				\fill (0,0.7) circle (0.1);

				\foreach \x in {1,2,3}
				\draw[rotate=72*\x]
				(18:0.7) -- (90:0.7);

				\foreach \x in {0,1,2}
				\draw[dotted,rotate=-72*\x]
				(18:0.7) -- (18:1.3) -- (-18:1.7) -- (-54:1.3) -- (-54:0.7);

				\draw
				(18:0.7) -- (18:1.3) -- (-18:1.7)
				(90:0.7) -- (18:0.7)
				(18:0.7) -- (-54:0.7)
				(-54:0.7) -- (-54:1.3) -- (-90:1.7) -- (234:1.3) -- (234:0.7)
				(-18:1.7) -- (-54:1.3)
				(234:1.3) -- (198:1.7) -- (162:1.3) -- (162:0.7);

				\draw[dotted]
				(90:0.7) -- (120:1.5) -- (140:1.7) -- (162:1.3) 
				(18:1.3) -- (40:1.7) -- (60:1.5) -- (90:0.7); 

				\draw[line width=1.5]
				(234:0.7) -- (162:0.7)
				(-18:1.7) -- (-54:1.3)
				(234:1.3) -- (-90:1.7);

				\node at (90:0.4) {\small $\bbb$}; 
				\node at (162:0.45) {\small $\ddd$};
				\node at (18:0.45) {\small $\aaa$};
				\node at (234:0.45) {\small $\eee$};
				\node at (-54:0.45) {\small $\ccc$};

				\node at (4:0.8) {\small $\aaa$}; 
				\node at (-40:0.8) {\small $\bbb$};
				\node at (-44:1.15) {\small $\ddd$};
				\node at (-18:1.45) {\small $\eee$};
				\node at (10:1.15) {\small $\ccc$};

				\node at (-68:0.8) {\small $\aaa$}; 
				\node at (-112:0.8) {\small $\bbb$};
				\node at (-116:1.15) {\small $\ddd$};
				\node at (-90:1.45) {\small $\eee$};
				\node at (-62:1.15) {\small $\ccc$};

				\node at (176:0.8) {\small $\eee$}; 
				\node at (220:0.8) {\small $\ddd$};
				\node at (224:1.15) {\small $\bbb$};
				\node at (198:1.45) {\small $\aaa$};
				\node at (170:1.2) {\small $\ccc$};

				\node at (146:0.75) {\small $\bbb$}; 

				\node[draw,shape=circle, inner sep=0.5] at (0,0) {\small $1$};
				\node[draw,shape=circle, inner sep=0.5] at (45:1.15) {\small $2$};\node[draw,shape=circle, inner sep=0.5] at (-18:1) {\small $3$};
				\node[draw,shape=circle, inner sep=0.5] at (-90:1.05) {\small $4$};
				\node[draw,shape=circle, inner sep=0.5] at (198:1) {\small $5$};
				\node[draw,shape=circle, inner sep=0.5] at (135:1.07) {\small $6$};
			\end{scope}


			\begin{scope}[xshift=5.5cm]

				\fill (0,0.7) circle (0.1);

				\foreach \x in {1,2,3}
				\draw[rotate=72*\x]
				(18:0.7) -- (90:0.7);

				\foreach \x in {0,1,2}
				\draw[dotted,rotate=-72*\x]
				(18:0.7) -- (18:1.3) -- (-18:1.7) -- (-54:1.3) -- (-54:0.7);

				\draw
				(18:0.7) -- (18:1.3) -- (-18:1.7)
				(90:0.7) -- (18:0.7)
				(18:0.7) -- (-54:0.7)
				(-54:0.7) -- (-54:1.3) -- (-90:1.7) -- (234:1.3) -- (234:0.7)
				(-18:1.7) -- (-54:1.3)
				(234:1.3) -- (198:1.7) -- (162:1.3) -- (162:0.7);

				\draw[dotted]
				(90:0.7) -- (120:1.5) -- (140:1.7) -- (162:1.3) 
				(18:1.3) -- (40:1.7) -- (60:1.5) -- (90:0.7); 

				\draw[line width=1.5]
				(234:0.7) -- (162:0.7)
				(-18:1.7) -- (-54:1.3)
				(234:1.3) -- (-90:1.7);

				\node at (90:0.4) {\small $\bbb$}; 
				\node at (162:0.45) {\small $\ddd$};
				\node at (18:0.45) {\small $\aaa$};
				\node at (234:0.45) {\small $\eee$};
				\node at (-54:0.45) {\small $\ccc$};

				\node at (4:0.8) {\small $\aaa$}; 
				\node at (-40:0.8) {\small $\ccc$};
				\node at (-44:1.15) {\small $\eee$};
				\node at (-18:1.45) {\small $\ddd$};
				\node at (10:1.15) {\small $\bbb$};

				\node at (-68:0.8) {\small $\aaa$}; 
				\node at (-112:0.8) {\small $\bbb$};
				\node at (-116:1.15) {\small $\ddd$};
				\node at (-90:1.45) {\small $\eee$};
				\node at (-62:1.15) {\small $\ccc$};

				\node at (176:0.8) {\small $\eee$}; 
				\node at (220:0.8) {\small $\ddd$};
				\node at (224:1.15) {\small $\bbb$};
				\node at (198:1.45) {\small $\aaa$};
				\node at (170:1.2) {\small $\ccc$};

				\node at (146:0.75) {\small $\bbb$}; 

				\node[draw,shape=circle, inner sep=0.5] at (0,0) {\small $1$};
				\node[draw,shape=circle, inner sep=0.5] at (45:1.15) {\small $2$};\node[draw,shape=circle, inner sep=0.5] at (-18:1) {\small $3$};
				\node[draw,shape=circle, inner sep=0.5] at (-90:1.05) {\small $4$};
				\node[draw,shape=circle, inner sep=0.5] at (198:1) {\small $5$};
				\node[draw,shape=circle, inner sep=0.5] at (135:1.07) {\small $6$};
			\end{scope}

		\end{tikzpicture}
		\caption{$H=\bbb_1\cdots$.}
		\label{bde(b)}
	\end{figure}

	In the first picture, we have $\aaa_1\bbb_3\cdots=\aaa\bbb^2$ or $\aaa\bbb\ccc$.
	In the second picture, we have $\aaa_1\ccc_3\cdots=\aaa^2\ccc$, $\aaa\bbb\ccc$ or $\aaa\ccc^2$. If $\aaa_1\ccc_3\cdots=\aaa^2\ccc$, we have
	\begin{align*}
		H=\{
		&\aaa\bbb^2, \aaa\bbb\ccc, \aaa^2\bbb^2, \aaa\bbb^3, \aaa\bbb^2\ccc, \aaa\bbb\ccc^2, \bbb^2\ddd^2, \bbb\ccc\ddd^2, \aaa^2\bbb^3, \aaa\bbb^4, \aaa\bbb^3\ccc, \aaa\bbb^2\ccc^2, \\
		&\aaa\bbb^2\ddd^2, \aaa\bbb^2\eee^2,\aaa\bbb\ccc^3,\aaa\bbb\ccc\ddd^2, \aaa\bbb\ccc\eee^2, \bbb^3\ddd^2, \bbb^2\ccc\ddd^2, \bbb\ccc^2\ddd^2\}.
	\end{align*}
	In the third and fourth pictures, we have $\aaa_1\aaa_3\cdots=\aaa^3$ or $\aaa^2\ccc$.

	\subsection*{Case $H=\ccc_1\cdots$}
	Let the fifth of Figure \ref{a4b_pentagon} be the center tile $T_1$ in the partial neighborhood in Figure \ref{bde(c)}. Since $\bbb\ddd\eee$ is an $a^2b$-vertex and $H=\ccc_1\cdots$, we have the degree $3$ vertex $\eee_1\cdots=\bbb\ddd\eee$, which determines $T_3$. The two pictures show two possible arrangements of $T_4$. In the first picture of Figure \ref{bde(c)}, we have the degree $3$ vertex $\aaa_4\bbb_1\cdots=\aaa_4\bbb_1\bbb_5$ or $\aaa_4\bbb_1\ccc_5$, which determines $\aaa_5$. In the second picture of Figure \ref{bde(c)}, we have the degree $3$ vertex $\bbb_1\ddd_4\cdots=\bbb\ddd\eee$, which determines $T_5$.

	\begin{figure}[htp]
		\centering
		\begin{tikzpicture}[>=latex,scale=0.9]


			\begin{scope}[xshift=-3cm]

				\fill (0,0.7) circle (0.1);

				\foreach \x in {1,2,3}
				\draw[rotate=72*\x]
				(18:0.7) -- (90:0.7);

				\foreach \x in {0,1,2}
				\draw[dotted,rotate=-72*\x]
				(18:0.7) -- (18:1.3) -- (-18:1.7) -- (-54:1.3) -- (-54:0.7);

				\draw
				(18:0.7) -- (18:1.3) -- (-18:1.7)
				(90:0.7) -- (18:0.7)
				(18:0.7) -- (-54:0.7)
				(-54:0.7) -- (-54:1.3) -- (-90:1.7) -- (234:1.3) -- (234:0.7)
				(-18:1.7) -- (-54:1.3)
				(234:1.3) -- (198:1.7) -- (162:1.3) -- (162:0.7);

				\draw[dotted]
				(90:0.7) -- (120:1.5) -- (140:1.7) -- (162:1.3) 
				(18:1.3) -- (40:1.7) -- (60:1.5) -- (90:0.7); 

				\draw[line width=1.5]
				(198:1.7) -- (234:1.3)
				(18:0.7) -- (-54:0.7)
				(-54:1.3) -- (-90:1.7);

				\node at (90:0.4) {\small $\ccc$}; 
				\node at (162:0.45) {\small $\aaa$};
				\node at (18:0.45) {\small $\eee$};
				\node at (234:0.45) {\small $\bbb$};
				\node at (-54:0.45) {\small $\ddd$};

				\node at (32:0.8) {\small $\bbb$}; 

				\node at (4:0.8) {\small $\ddd$}; 
				\node at (-40:0.8) {\small $\eee$};
				\node at (-44:1.15) {\small $\ccc$};
				\node at (-18:1.45) {\small $\aaa$};
				\node at (10:1.15) {\small $\bbb$};

				\node at (-68:0.8) {\small $\bbb$}; 
				\node at (-112:0.8) {\small $\aaa$};
				\node at (-116:1.15) {\small $\ccc$};
				\node at (-90:1.45) {\small $\eee$};
				\node at (-62:1.15) {\small $\ddd$};

				\node at (176:0.8) {\small $\aaa$}; 

				\node[draw,shape=circle, inner sep=0.5] at (0,0) {\small $1$};
				\node[draw,shape=circle, inner sep=0.5] at (45:1.15) {\small $2$};\node[draw,shape=circle, inner sep=0.5] at (-18:1) {\small $3$};
				\node[draw,shape=circle, inner sep=0.5] at (-90:1.05) {\small $4$};
				\node[draw,shape=circle, inner sep=0.5] at (198:1) {\small $5$};
				\node[draw,shape=circle, inner sep=0.5] at (135:1.07) {\small $6$};
			\end{scope}


			\begin{scope}[xshift=3cm]

				\fill (0,0.7) circle (0.1);

				\foreach \x in {1,2,3}
				\draw[rotate=72*\x]
				(18:0.7) -- (90:0.7);

				\foreach \x in {0,1,2}
				\draw[dotted,rotate=-72*\x]
				(18:0.7) -- (18:1.3) -- (-18:1.7) -- (-54:1.3) -- (-54:0.7);

				\draw
				(18:0.7) -- (18:1.3) -- (-18:1.7)
				(90:0.7) -- (18:0.7)
				(18:0.7) -- (-54:0.7)
				(-54:0.7) -- (-54:1.3) -- (-90:1.7) -- (234:1.3) -- (234:0.7)
				(-18:1.7) -- (-54:1.3)
				(234:1.3) -- (198:1.7) -- (162:1.3) -- (162:0.7);

				\draw[dotted]
				(90:0.7) -- (120:1.5) -- (140:1.7) -- (162:1.3) 
				(18:1.3) -- (40:1.7) -- (60:1.5) -- (90:0.7); 

				\draw[line width=1.5]
				(234:1.3) -- (234:0.7)
				(18:0.7) -- (-54:0.7);

				\node at (90:0.4) {\small $\ccc$}; 
				\node at (162:0.45) {\small $\aaa$};
				\node at (18:0.45) {\small $\eee$};
				\node at (234:0.45) {\small $\bbb$};
				\node at (-54:0.45) {\small $\ddd$};

				\node at (32:0.8) {\small $\bbb$}; 

				\node at (4:0.8) {\small $\ddd$}; 
				\node at (-40:0.8) {\small $\eee$};
				\node at (-44:1.15) {\small $\ccc$};
				\node at (-18:1.45) {\small $\aaa$};
				\node at (10:1.15) {\small $\bbb$};

				\node at (-68:0.8) {\small $\bbb$}; 
				\node at (-112:0.8) {\small $\ddd$};
				\node at (-116:1.15) {\small $\eee$};
				\node at (-90:1.45) {\small $\ccc$};
				\node at (-62:1.15) {\small $\aaa$};

				\node at (176:0.8) {\small $\ccc$}; 
				\node at (220:0.8) {\small $\eee$};
				\node at (224:1.15) {\small $\ddd$};
				\node at (198:1.45) {\small $\bbb$};
				\node at (170:1.2) {\small $\aaa$};

				\node[draw,shape=circle, inner sep=0.5] at (0,0) {\small $1$};
				\node[draw,shape=circle, inner sep=0.5] at (45:1.15) {\small $2$};\node[draw,shape=circle, inner sep=0.5] at (-18:1) {\small $3$};
				\node[draw,shape=circle, inner sep=0.5] at (-90:1.05) {\small $4$};
				\node[draw,shape=circle, inner sep=0.5] at (198:1) {\small $5$};
				\node[draw,shape=circle, inner sep=0.5] at (135:1.07) {\small $6$};
			\end{scope}

		\end{tikzpicture}
		\caption{$H=\ccc_1\cdots$.}
		\label{bde(c)}
	\end{figure}

 	In the first picture of Figure \ref{bde(c)}, we have $\aaa_1\aaa_5\cdots=\aaa^3$ or $\aaa^2\ccc$. In the second picture of Figure \ref{bde(c)}, we have $\aaa_1\ccc_5\cdots=\aaa^2\ccc$, $\aaa\bbb\ccc$ or $\aaa\ccc^2$. Each case induces the following $H$ respectively:
	\begin{align*}
		H
		&=\{\aaa\bbb\ccc, \aaa\ccc^2, \aaa\bbb^2\ccc, \aaa\bbb\ccc^2, \aaa\ccc^3, \bbb\ccc\ddd^2, \ccc^2\ddd^2, \ccc^2\ddd\eee, \aaa\bbb^3\ccc, \aaa\bbb^2\ccc^2, \aaa\bbb\ccc^3, \aaa\bbb\ccc\ddd^2, \\
		&\quad\,\,\,\aaa\bbb\ccc\eee^2, \aaa\ccc^4, \aaa\ccc^2\ddd^2,\aaa\ccc^2\ddd\eee,\aaa\ccc^2\eee^2, \bbb^2\ccc\ddd^2,\bbb\ccc^2\ddd^2, \ccc^3\ddd^2, \ccc^3\ddd\eee\} \\
		&\text{or}\,\,\{\aaa^2\ccc, \aaa^3\ccc, \aaa\ccc\ddd^2, \aaa^4\ccc, \aaa^2\ccc\ddd^2,\aaa^2\ccc\eee^2, \aaa\ccc^2\ddd^2, \ccc\ddd^4, \ccc\ddd^3\eee\} \\
		&\text{or}\,\,\{\aaa^2\ccc, \aaa^3\ccc, \aaa^2\bbb\ccc, \aaa\ccc\ddd^2,\aaa\ccc\ddd\eee, \aaa\ccc\eee^2, \aaa^4\ccc, \aaa^3\bbb\ccc, \aaa^2\bbb^2\ccc,\aaa^2\ccc\ddd^2, \aaa^2\ccc\ddd\eee, \aaa^2\ccc\eee^2, \\
		&\quad\,\,\, \aaa\bbb\ccc\eee^2, \ccc\ddd^3\eee, \ccc\ddd^2\eee^2, \ccc\ddd\eee^3\}.
	\end{align*}

	\subsection*{Case $H=\ddd_1\cdots$}
	Let the third of Figure \ref{a4b_pentagon} be the center tile $T_1$ in the partial neighborhood in Figure \ref{bde(d)}. Since $\bbb\ddd\eee$ is an $a^2b$-vertex and $H=\ddd_1\cdots$, we have the degree $3$ vertex $\eee_1\cdots=\bbb\ddd\eee$, which determines $T_5$ and $T_6$. The three pictures show three possible arrangements of $T_4$.
	\begin{figure}[htp]
		\centering
		\begin{tikzpicture}[>=latex,scale=0.85]

			\begin{scope}[xshift=-4cm]

				\fill (0,0.7) circle (0.1);

				\foreach \x in {1,2,3}
				\draw[rotate=72*\x]
				(18:0.7) -- (90:0.7);

				\foreach \x in {0,1,2}
				\draw[dotted,rotate=-72*\x]
				(18:0.7) -- (18:1.3) -- (-18:1.7) -- (-54:1.3) -- (-54:0.7);

				\draw
				(90:0.7) -- (18:0.7) -- (-54:0.7)
				(-54:0.7) -- (-54:1.3) -- (-90:1.7) -- (234:1.3) -- (234:0.7)
				(234:1.3) -- (198:1.7) -- (162:1.3) -- (162:0.7)
				(90:0.7) -- (120:1.5) -- (140:1.7) -- (162:1.3);

				\draw[dotted]
				(-54:1.3) -- (-18:1.7) -- (18:1.3) -- (18:0.7)
				(18:1.3) -- (40:1.7) -- (60:1.5) -- (90:0.7); 

				\draw[line width=1.5]
				(90:0.7) -- (162:0.7)
				(162:1.3) -- (198:1.7)
				(-54:1.3) -- (-90:1.7);

				\node at (90:0.4) {\small $\ddd$}; 
				\node at (162:0.45) {\small $\eee$};
				\node at (18:0.45) {\small $\bbb$};
				\node at (234:0.45) {\small $\ccc$};
				\node at (-54:0.45) {\small $\aaa$};

				\node at (-68:0.8) {\small $\bbb$}; 
				\node at (-112:0.8) {\small $\aaa$};
				\node at (-116:1.15) {\small $\ccc$};
				\node at (-90:1.45) {\small $\eee$};
				\node at (-62:1.15) {\small $\ddd$};

				\node at (176:0.8) {\small $\bbb$}; 
				\node at (220:0.8) {\small $\aaa$};
				\node at (224:1.15) {\small $\ccc$};
				\node at (198:1.45) {\small $\eee$};
				\node at (170:1.2) {\small $\ddd$};

				\node at (146:0.75) {\small $\ddd$}; 
				\node at (112:0.75) {\small $\eee$};
				\node at (152:1.2) {\small $\bbb$};
				\node at (140:1.5) {\small $\aaa$};
				\node at (125:1.35) {\small $\ccc$};

				\node[draw,shape=circle, inner sep=0.5] at (0,0) {\small $1$};
				\node[draw,shape=circle, inner sep=0.5] at (45:1.15) {\small $2$};\node[draw,shape=circle, inner sep=0.5] at (-18:1) {\small $3$};
				\node[draw,shape=circle, inner sep=0.5] at (-90:1.05) {\small $4$};
				\node[draw,shape=circle, inner sep=0.5] at (198:1) {\small $5$};
				\node[draw,shape=circle, inner sep=0.5] at (135:1.07) {\small $6$};
			\end{scope}

			\begin{scope}[xshift=0cm]

				\fill (0,0.7) circle (0.1);

				\foreach \x in {1,2,3}
				\draw[rotate=72*\x]
				(18:0.7) -- (90:0.7);

				\foreach \x in {0,1,2}
				\draw[dotted,rotate=-72*\x]
				(18:0.7) -- (18:1.3) -- (-18:1.7) -- (-54:1.3) -- (-54:0.7);

				\draw
				(90:0.7) -- (18:0.7) -- (-54:0.7)
				(-54:0.7) -- (-54:1.3) -- (-90:1.7) -- (234:1.3) -- (234:0.7)
				(234:1.3) -- (198:1.7) -- (162:1.3) -- (162:0.7)
				(90:0.7) -- (120:1.5) -- (140:1.7) -- (162:1.3);

				\draw[dotted]
				(-54:1.3) -- (-18:1.7) -- (18:1.3) -- (18:0.7)
				(18:1.3) -- (40:1.7) -- (60:1.5) -- (90:0.7); 

				\draw[line width=1.5]
				(90:0.7) -- (162:0.7)
				(162:1.3) -- (198:1.7)
				(-54:1.3) -- (-90:1.7);

				\node at (90:0.4) {\small $\ddd$}; 
				\node at (162:0.45) {\small $\eee$};
				\node at (18:0.45) {\small $\bbb$};
				\node at (234:0.45) {\small $\ccc$};
				\node at (-54:0.45) {\small $\aaa$};

				\node at (-68:0.8) {\small $\ccc$}; 
				\node at (-112:0.8) {\small $\aaa$};
				\node at (-116:1.15) {\small $\bbb$};
				\node at (-90:1.45) {\small $\ddd$};
				\node at (-62:1.15) {\small $\eee$};

				\node at (176:0.8) {\small $\bbb$}; 
				\node at (220:0.8) {\small $\aaa$};
				\node at (224:1.15) {\small $\ccc$};
				\node at (198:1.45) {\small $\eee$};
				\node at (170:1.2) {\small $\ddd$};

				\node at (146:0.75) {\small $\ddd$}; 
				\node at (112:0.75) {\small $\eee$};
				\node at (152:1.2) {\small $\bbb$};
				\node at (140:1.5) {\small $\aaa$};
				\node at (125:1.35) {\small $\ccc$};

				\node[draw,shape=circle, inner sep=0.5] at (0,0) {\small $1$};
				\node[draw,shape=circle, inner sep=0.5] at (45:1.15) {\small $2$};\node[draw,shape=circle, inner sep=0.5] at (-18:1) {\small $3$};
				\node[draw,shape=circle, inner sep=0.5] at (-90:1.05) {\small $4$};
				\node[draw,shape=circle, inner sep=0.5] at (198:1) {\small $5$};
				\node[draw,shape=circle, inner sep=0.5] at (135:1.07) {\small $6$};
			\end{scope}

			\begin{scope}[xshift=4cm]

				\fill (0,0.7) circle (0.1);

				\foreach \x in {1,2,3}
				\draw[rotate=72*\x]
				(18:0.7) -- (90:0.7);

				\foreach \x in {0,1,2}
				\draw[dotted,rotate=-72*\x]
				(18:0.7) -- (18:1.3) -- (-18:1.7) -- (-54:1.3) -- (-54:0.7);

				\draw
				(90:0.7) -- (18:0.7)
				(18:0.7) -- (-54:0.7)
				(-54:0.7) -- (-54:1.3) -- (-90:1.7) -- (234:1.3) -- (234:0.7)
				(234:1.3) -- (198:1.7) -- (162:1.3) -- (162:0.7)
				(90:0.7) -- (120:1.5) -- (140:1.7) -- (162:1.3);

				\draw[dotted]
				(18:1.3) -- (40:1.7) -- (60:1.5) -- (90:0.7) 
				(-54:1.3) -- (-18:1.7) -- (18:1.3) -- (18:0.7);

				\draw[line width=1.5]
				(90:0.7) -- (162:0.7)
				(162:1.3) -- (198:1.7)
				(234:1.3) -- (-90:1.7);

				\node at (90:0.4) {\small $\ddd$}; 
				\node at (162:0.45) {\small $\eee$};
				\node at (18:0.45) {\small $\bbb$};
				\node at (234:0.45) {\small $\ccc$};
				\node at (-54:0.45) {\small $\aaa$};

				\node at (-68:0.8) {\small $\aaa$}; 

				\node at (176:0.8) {\small $\bbb$}; 
				\node at (220:0.8) {\small $\aaa$};
				\node at (224:1.15) {\small $\ccc$};
				\node at (198:1.45) {\small $\eee$};
				\node at (170:1.2) {\small $\ddd$};

				\node at (146:0.75) {\small $\ddd$}; 
				\node at (112:0.75) {\small $\eee$};
				\node at (152:1.2) {\small $\bbb$};
				\node at (140:1.5) {\small $\aaa$};
				\node at (125:1.35) {\small $\ccc$};

				\node[draw,shape=circle, inner sep=0.5] at (0,0) {\small $1$};
				\node[draw,shape=circle, inner sep=0.5] at (45:1.15) {\small $2$};\node[draw,shape=circle, inner sep=0.5] at (-18:1) {\small $3$};
				\node[draw,shape=circle, inner sep=0.5] at (-90:1.05) {\small $4$};
				\node[draw,shape=circle, inner sep=0.5] at (198:1) {\small $5$};
				\node[draw,shape=circle, inner sep=0.5] at (135:1.07) {\small $6$};
			\end{scope}

		\end{tikzpicture}
		\caption{$H=\ddd_1\cdots$.}
		\label{bde(d)}
	\end{figure}

	In the first picture, we have $\aaa_1\bbb_4\cdots=\aaa\bbb^2,\aaa\bbb\ccc$.
	In the second picture, we have $\aaa_1\ccc_3\cdots=\aaa^2\ccc$, $\aaa\bbb\ccc$ or $\aaa\ccc^2$. If $\aaa_1\ccc_4\cdots=\aaa^2\ccc$, then $\bbb_1\cdots=\aaa\bbb^2,\aaa\bbb\ccc,\bbb^3,\bbb^2\ccc,\bbb\ccc^2$.
	In the third picture, we have $\aaa_5\ccc_1\cdots=\aaa\bbb\ccc$ or $\aaa\ccc^2$, $\aaa_1\aaa_4\cdots=\aaa^3$ or $\aaa^2\ccc$.

	\subsection*{Case $H=\eee_1\cdots$}
	Let the fourth of Figure \ref{a4b_pentagon} be the center tile $T_1$ in the partial neighborhood in Figure \ref{bde(e)}. Since $\bbb\ddd\eee$ is an $a^2b$-vertex and $H=\eee_1\cdots$, we have the degree $3$ vertex $\ddd_1\cdots=\bbb\ddd\eee$, which determines $T_2$. The four pictures show two possible arrangements of $T_3$. In the second, third and fourth pictures, we have the degree $3$ vertex $\bbb_1\ddd_3\cdots=\bbb\ddd\eee$, which determines $T_4$. In the second picture, we have the degree $3$ vertex $\aaa_1\ccc_4\cdots=\aaa_1\aaa_5\ccc_4$, which determines $\aaa_5$. In the third picture, we have the degree $3$ vertex $\aaa_1\ccc_4\cdots=\aaa_1\bbb_5\ccc_4$ , which determines $T_5$. In the fourth picture, we have the degree $3$ vertex $\aaa_1\ccc_4\cdots=\aaa_1\ccc_4\ccc_5$, which determines $T_5$.
	\begin{figure}[htp]
		\centering
		\begin{tikzpicture}[>=latex,scale=0.85]


			\begin{scope}[xshift=-5cm]

				\fill (0,0.7) circle (0.1);

				\foreach \x in {1,2,3}
				\draw[rotate=72*\x]
				(18:0.7) -- (90:0.7);

				\foreach \x in {0,1,2}
				\draw[dotted,rotate=-72*\x]
				(18:0.7) -- (18:1.3) -- (-18:1.7) -- (-54:1.3) -- (-54:0.7);

				\draw
				(18:0.7) -- (18:1.3) -- (-18:1.7)
				(90:0.7) -- (18:0.7)
				(18:0.7) -- (-54:0.7)
				(-54:0.7) -- (-54:1.3) -- (-90:1.7) -- (234:1.3) -- (234:0.7)
				(-18:1.7) -- (-54:1.3)
				(18:1.3) -- (40:1.7) -- (60:1.5) -- (90:0.7);

				\draw[dotted]
				(234:1.3) -- (198:1.7) -- (162:1.3) -- (162:0.7)
				(90:0.7) -- (120:1.5) -- (140:1.7) -- (162:1.3); 

				\draw[line width=1.5]
				(18:1.3) -- (-18:1.7)
				(90:0.7) -- (18:0.7)
				(-90:1.7) -- (-54:1.3);

				\node at (90:0.4) {\small $\eee$}; 
				\node at (162:0.45) {\small $\ccc$};
				\node at (18:0.45) {\small $\ddd$};
				\node at (234:0.45) {\small $\aaa$};
				\node at (-54:0.45) {\small $\bbb$};

				\node at (70:0.75) {\small $\ddd$}; 
				\node at (32:0.8) {\small $\eee$};
				\node at (26:1.2) {\small $\ccc$};
				\node at (39:1.5) {\small $\aaa$};
				\node at (55:1.3) {\small $\bbb$};

				\node at (4:0.8) {\small $\bbb$}; 
				\node at (-40:0.8) {\small $\aaa$};
				\node at (-44:1.15) {\small $\ccc$};
				\node at (-18:1.45) {\small $\eee$};
				\node at (10:1.15) {\small $\ddd$};

				\node at (-112:0.8) {\small $\aaa$};

				\node[draw,shape=circle, inner sep=0.5] at (0,0) {\small $1$};
				\node[draw,shape=circle, inner sep=0.5] at (42:1.15) {\small $2$};
				\node[draw,shape=circle, inner sep=0.5] at (-18:1) {\small $3$};
				\node[draw,shape=circle, inner sep=0.5] at (-90:1.05) {\small $4$};
				\node[draw,shape=circle, inner sep=0.5] at (198:1) {\small $5$};
				\node[draw,shape=circle, inner sep=0.5] at (135:1.07) {\small $6$};
			\end{scope}

			\begin{scope}[xshift=-1.5cm]

				\fill (0,0.7) circle (0.1);

				\foreach \x in {1,2,3}
				\draw[rotate=72*\x]
				(18:0.7) -- (90:0.7);

				\foreach \x in {0,1,2}
				\draw[dotted,rotate=-72*\x]
				(18:0.7) -- (18:1.3) -- (-18:1.7) -- (-54:1.3) -- (-54:0.7);

				\draw
				(18:0.7) -- (18:1.3) -- (-18:1.7)
				(90:0.7) -- (18:0.7)
				(18:0.7) -- (-54:0.7)
				(-54:0.7) -- (-54:1.3) -- (-90:1.7) -- (234:1.3) -- (234:0.7)
				(-18:1.7) -- (-54:1.3)
				(234:1.3) -- (198:1.7) -- (162:1.3) -- (162:0.7)
				(18:1.3) -- (40:1.7) -- (60:1.5) -- (90:0.7);

				\draw[dotted]
				(90:0.7) -- (120:1.5) -- (140:1.7) -- (162:1.3); 

				\draw[line width=1.5]
				(90:0.7) -- (18:0.7)
				(-54:0.7) -- (-54:1.3)
				(198:1.7) -- (162:1.3);

				\node at (90:0.4) {\small $\eee$}; 
				\node at (162:0.45) {\small $\ccc$};
				\node at (18:0.45) {\small $\ddd$};
				\node at (234:0.45) {\small $\aaa$};
				\node at (-54:0.45) {\small $\bbb$};

				\node at (70:0.75) {\small $\ddd$}; 
				\node at (32:0.8) {\small $\eee$};
				\node at (26:1.2) {\small $\ccc$};
				\node at (39:1.5) {\small $\aaa$};
				\node at (55:1.3) {\small $\bbb$};

				\node at (4:0.8) {\small $\bbb$}; 
				\node at (-40:0.8) {\small $\ddd$};
				\node at (-44:1.15) {\small $\eee$};
				\node at (-18:1.45) {\small $\ccc$};
				\node at (10:1.15) {\small $\aaa$};

				\node at (-68:0.8) {\small $\eee$}; 
				\node at (-112:0.8) {\small $\ccc$};
				\node at (-116:1.15) {\small $\aaa$};
				\node at (-90:1.45) {\small $\bbb$};
				\node at (-62:1.15) {\small $\ddd$};

				\node at (220:0.8) {\small $\aaa$}; 

				\node[draw,shape=circle, inner sep=0.5] at (0,0) {\small $1$};
				\node[draw,shape=circle, inner sep=0.5] at (42:1.15) {\small $2$};
				\node[draw,shape=circle, inner sep=0.5] at (-18:1) {\small $3$};
				\node[draw,shape=circle, inner sep=0.5] at (-90:1.05) {\small $4$};
				\node[draw,shape=circle, inner sep=0.5] at (198:1) {\small $5$};
				\node[draw,shape=circle, inner sep=0.5] at (135:1.07) {\small $6$};
			\end{scope}


			\begin{scope}[xshift=2cm]

				\fill (0,0.7) circle (0.1);

				\foreach \x in {1,2,3}
				\draw[rotate=72*\x]
				(18:0.7) -- (90:0.7);

				\foreach \x in {0,1,2}
				\draw[dotted,rotate=-72*\x]
				(18:0.7) -- (18:1.3) -- (-18:1.7) -- (-54:1.3) -- (-54:0.7);

				\draw
				(18:0.7) -- (18:1.3) -- (-18:1.7)
				(90:0.7) -- (18:0.7)
				(18:0.7) -- (-54:0.7)
				(-54:0.7) -- (-54:1.3) -- (-90:1.7) -- (234:1.3) -- (234:0.7)
				(-18:1.7) -- (-54:1.3)
				(234:1.3) -- (198:1.7) -- (162:1.3) -- (162:0.7)
				(18:1.3) -- (40:1.7) -- (60:1.5) -- (90:0.7);

				\draw[dotted]
				(90:0.7) -- (120:1.5) -- (140:1.7) -- (162:1.3); 

				\draw[line width=1.5]
				(90:0.7) -- (18:0.7)
				(-54:0.7) -- (-54:1.3)
				(198:1.7) -- (234:1.3);

				\node at (90:0.4) {\small $\eee$}; 
				\node at (162:0.45) {\small $\ccc$};
				\node at (18:0.45) {\small $\ddd$};
				\node at (234:0.45) {\small $\aaa$};
				\node at (-54:0.45) {\small $\bbb$};

				\node at (70:0.75) {\small $\ddd$}; 
				\node at (32:0.8) {\small $\eee$};
				\node at (26:1.2) {\small $\ccc$};
				\node at (39:1.5) {\small $\aaa$};
				\node at (55:1.3) {\small $\bbb$};

				\node at (4:0.8) {\small $\bbb$}; 
				\node at (-40:0.8) {\small $\ddd$};
				\node at (-44:1.15) {\small $\eee$};
				\node at (-18:1.45) {\small $\ccc$};
				\node at (10:1.15) {\small $\aaa$};

				\node at (-68:0.8) {\small $\eee$}; 
				\node at (-112:0.8) {\small $\ccc$};
				\node at (-116:1.15) {\small $\aaa$};
				\node at (-90:1.45) {\small $\bbb$};
				\node at (-62:1.15) {\small $\ddd$};

				\node at (176:0.8) {\small $\aaa$}; 
				\node at (220:0.8) {\small $\bbb$};
				\node at (224:1.15) {\small $\ddd$};
				\node at (198:1.45) {\small $\eee$};
				\node at (170:1.2) {\small $\ccc$};

				\node[draw,shape=circle, inner sep=0.5] at (0,0) {\small $1$};
				\node[draw,shape=circle, inner sep=0.5] at (42:1.15) {\small $2$};\node[draw,shape=circle, inner sep=0.5] at (-18:1) {\small $3$};
				\node[draw,shape=circle, inner sep=0.5] at (-90:1.05) {\small $4$};
				\node[draw,shape=circle, inner sep=0.5] at (198:1) {\small $5$};
				\node[draw,shape=circle, inner sep=0.5] at (135:1.07) {\small $6$};
			\end{scope}


			\begin{scope}[xshift=5.5cm]

				\fill (0,0.7) circle (0.1);

				\foreach \x in {1,2,3}
				\draw[rotate=72*\x]
				(18:0.7) -- (90:0.7);

				\foreach \x in {0,1,2}
				\draw[dotted,rotate=-72*\x]
				(18:0.7) -- (18:1.3) -- (-18:1.7) -- (-54:1.3) -- (-54:0.7);

				\draw
				(18:0.7) -- (18:1.3) -- (-18:1.7)
				(90:0.7) -- (18:0.7)
				(18:0.7) -- (-54:0.7)
				(-54:0.7) -- (-54:1.3) -- (-90:1.7) -- (234:1.3) -- (234:0.7)
				(-18:1.7) -- (-54:1.3)
				(234:1.3) -- (198:1.7) -- (162:1.3) -- (162:0.7)
				(18:1.3) -- (40:1.7) -- (60:1.5) -- (90:0.7);

				\draw[dotted]
				(90:0.7) -- (120:1.5) -- (140:1.7) -- (162:1.3); 

				\draw[line width=1.5]
				(90:0.7) -- (18:0.7)
				(-54:0.7) -- (-54:1.3)
				(198:1.7) -- (234:1.3);

				\node at (90:0.4) {\small $\eee$}; 
				\node at (162:0.45) {\small $\ccc$};
				\node at (18:0.45) {\small $\ddd$};
				\node at (234:0.45) {\small $\aaa$};
				\node at (-54:0.45) {\small $\bbb$};

				\node at (70:0.75) {\small $\ddd$}; 
				\node at (32:0.8) {\small $\eee$};
				\node at (26:1.2) {\small $\ccc$};
				\node at (39:1.5) {\small $\aaa$};
				\node at (55:1.3) {\small $\bbb$};

				\node at (4:0.8) {\small $\bbb$}; 
				\node at (-40:0.8) {\small $\ddd$};
				\node at (-44:1.15) {\small $\eee$};
				\node at (-18:1.45) {\small $\ccc$};
				\node at (10:1.15) {\small $\aaa$};

				\node at (-68:0.8) {\small $\eee$}; 
				\node at (-112:0.8) {\small $\ccc$};
				\node at (-116:1.15) {\small $\aaa$};
				\node at (-90:1.45) {\small $\bbb$};
				\node at (-62:1.15) {\small $\ddd$};

				\node at (176:0.8) {\small $\aaa$}; 
				\node at (220:0.8) {\small $\ccc$};
				\node at (224:1.15) {\small $\eee$};
				\node at (198:1.45) {\small $\ddd$};
				\node at (170:1.2) {\small $\bbb$};

				\node[draw,shape=circle, inner sep=0.5] at (0,0) {\small $1$};
				\node[draw,shape=circle, inner sep=0.5] at (42:1.15) {\small $2$};
				\node[draw,shape=circle, inner sep=0.5] at (-18:1) {\small $3$};
				\node[draw,shape=circle, inner sep=0.5] at (-90:1.05) {\small $4$};
				\node[draw,shape=circle, inner sep=0.5] at (198:1) {\small $5$};
				\node[draw,shape=circle, inner sep=0.5] at (135:1.07) {\small $6$};
			\end{scope}

		\end{tikzpicture}
		\caption{$H=\eee_1\cdots$.}
		\label{bde(e)}
	\end{figure}

	In the first picture, we have $\aaa_3\bbb_1\cdots=\aaa\bbb^2$ or $\aaa\bbb\ccc$, $\aaa_1\aaa_4\cdots=\aaa^3$ or $\aaa^2\ccc$.

	In the second picture, we have $\ccc_1\cdots=\aaa\bbb\ccc$, $\aaa\ccc^2$, $\bbb^2\ccc$, $\bbb\ccc^2$ or $\ccc^3$.

	In the third picture, we have $\aaa_5\ccc_1\cdots=\aaa^2\ccc$ or $\aaa\bbb\ccc$. If $\aaa_5\ccc_1\cdots=\aaa\bbb\ccc$, we have
	\[
	H=\{\aaa^2\ddd\eee, \ddd^3\eee, \aaa^3\ddd\eee, \aaa\ddd^3\eee,\aaa\ddd\eee^3, \ccc\ddd^3\eee\}.
	\]

	In the fourth picture, we have $\aaa_5\ccc_1\cdots=\aaa^2\ccc$, $\aaa\bbb\ccc$ or $\aaa\ccc^2$. If $\aaa_5\ccc_1\cdots=\aaa\ccc^2$, we have
	\[
	H=\{\aaa^2\ddd\eee, \aaa\ccc\ddd\eee,\ddd^2\eee^2, \ddd\eee^3,  \aaa^3\ddd\eee, \aaa^2\ccc\ddd\eee, \aaa\ccc^2\ddd\eee, \aaa\ddd^2\eee^2,\aaa\ddd^3\eee,\aaa\ddd\eee^3,  \ccc\ddd^2\eee^2, \ccc\ddd\eee^3\}.
	\]

	By the previous discussion, we get Table \ref{10}.

	\begin{table}[htp]
		\centering
		\caption{$\bbb\ddd\eee$ as the unique $a^2b$-vertex type.}\label{10}
		\resizebox{\linewidth}{!}{\begin{tabular}{|c|c|c|l|}
			\hline
			$a^2b$-Vertex& H &$a^3$-Vertex& \multicolumn{1}{c|}{Vertex} \\
			\hline\hline
			\multirow{20}{*}{$\bbb\ddd\eee$} & \multirow{5}{*}{$\aaa$} & \multirow{2}{*}{$\aaa^2\ccc$}& $\aaa\bbb^2,\aaa\bbb\ccc,  \aaa\ccc^2, \aaa\bbb^2\ccc, \aaa\bbb\ccc^2, \aaa\ccc^3, \aaa\bbb^3\ccc, \aaa\bbb^2\ccc^2, \aaa\bbb\ccc^3,\aaa\bbb\ccc\ddd^2,\aaa\bbb\ccc\eee^2,$ \\
			&&&$  \aaa\ccc^4, \aaa\ccc^2\ddd^2, \aaa\ccc^2\ddd\eee, \aaa\ccc^2\eee^2$  \\
			\cline{3-4}
			& &\multirow{2}{*}{$\aaa\bbb\ccc$} & $\aaa^3, \aaa^4, \aaa^3\bbb, \aaa^3\ccc, \aaa^2\ddd^2, \aaa^2\ddd\eee, \aaa^2\eee^2, \aaa\ccc\ddd^2, \aaa^5,\aaa^4\bbb,\aaa^4\ccc,\aaa^3\bbb^2,\aaa^3\ddd^2,$ \\
			&&&$\aaa^3\ddd\eee, \aaa^3\eee^2,\aaa^2\bbb\ddd^2,\aaa^2\bbb\eee^2, \aaa^2\ccc\ddd^2,\aaa^2\ccc\eee^2, \aaa\ccc^2\ddd^2,\aaa\ddd^3\eee,\aaa\ddd\eee^3$ \\
			\cline{3-4}
			& & \multirow{1}{*}{$\aaa\ccc^2$} & $ \aaa^3\ccc, \aaa^2\bbb\ccc, \aaa\ccc\ddd\eee, \aaa\ccc\eee^2, \aaa^4\ccc, \aaa^3\bbb\ccc, \aaa^2\bbb^2\ccc,\aaa^2\ccc\ddd^2, \aaa^2\ccc\ddd\eee,\aaa^2\ccc\eee^2$ \\
			\cline{2-4}
			& \multirow{4}{*}{$\bbb$} & \multirow{2}{*}{$\aaa^2\ccc$} & $\aaa\bbb^2,\aaa\bbb\ccc$, $\aaa\ccc^2, \aaa^2\bbb^2, \aaa\bbb^3, \aaa\bbb^2\ccc, \aaa\bbb\ccc^2, \bbb^2\ddd^2, \bbb\ccc\ddd^2,\aaa^2\bbb^3,\aaa\bbb^4,\aaa\bbb^3\ccc,$ \\
			&&&$  \aaa\bbb^2\ccc^2,\aaa\bbb\ccc^3, \aaa\bbb^2\ddd^2, \aaa\bbb^2\eee^2,\aaa\bbb\ccc\ddd^2,\aaa\bbb\ccc\eee^2, \bbb^3\ddd^2,\bbb^2\ccc\ddd^2,\bbb\ccc^2\ddd^2$ \\
			\cline{3-4}
			&& $\aaa\bbb\ccc$& $\aaa^3$ \\
			\cline{3-4}
			&& $\aaa\ccc^2$ & $\aaa^3$ \\
			\cline{2-4}
			& \multirow{6}{*}{$\ccc$} & $\aaa^3$ & $\aaa\bbb^2$ \\
			\cline{3-4}
			&&\multirow{2}{*}{$\aaa^2\ccc$} & $\aaa\bbb^2, \aaa\bbb\ccc, \aaa\ccc^2, \aaa\bbb^2\ccc, \aaa\bbb\ccc^2, \aaa\ccc^3, \bbb\ccc\ddd^2, \ccc^2\ddd^2, \ccc^2\ddd\eee,\aaa\bbb^3\ccc, \aaa\bbb^2\ccc^2,\aaa\bbb\ccc^3, $ \\
			&&&$ \aaa\bbb\ccc\ddd^2,\aaa\bbb\ccc\eee^2, \aaa\ccc^4, \aaa\ccc^2\ddd^2, \aaa\ccc^2\ddd\eee,\aaa\ccc^2\eee^2, \bbb^2\ccc\ddd^2,\bbb\ccc^2\ddd^2, \ccc^3\ddd^2, \ccc^3\ddd\eee$\\
			\cline{3-4}
			& & \multirow{1}{*}{$\aaa\bbb\ccc$} & $\aaa^3, \aaa^3\ccc, \aaa\ccc\ddd^2, \aaa^4\ccc,\aaa^2\ccc\ddd^2,\aaa^2\ccc\eee^2,\aaa\ccc^2\ddd^2, \ccc\ddd^4, \ccc\ddd^3\eee$ \\
			\cline{3-4}
			& & \multirow{2}{*}{$\aaa\ccc^2$} & $\aaa^3\ccc, \aaa^2\bbb\ccc, \aaa\ccc\ddd\eee, \aaa\ccc\ddd^2, \aaa\ccc\eee^2, \aaa^4\ccc, \aaa^3\bbb\ccc, \aaa^2\bbb^2\ccc, \aaa^2\ccc\ddd^2,\aaa^2\ccc\ddd\eee,\aaa^2\ccc\eee^2,$ \\
			&&&$  \aaa\bbb\ccc\eee^2,\ccc\ddd^3\eee, \ccc\ddd^2\eee^2, \ccc\ddd\eee^3$\\
			\cline{2-4}
			& \multirow{3}{*}{$\ddd$} & $\aaa^2\ccc$ & $\aaa\bbb^2,\aaa\bbb\ccc,\aaa\ccc^2, \bbb^3, \bbb^2\ccc,\bbb\ccc^2$ \\
			\cline{3-4}
			& & $\aaa\bbb\ccc$ & $\aaa^3$ \\
			\cline{3-4}
			& & $\aaa\ccc^2$ & $\aaa^3$ \\
			\cline{2-4}
			& \multirow{4}{*}{$\eee$} & $\aaa^3$ & $\aaa\bbb^2$\\
			\cline{3-4}
			&&$\aaa^2\ccc$ & $\aaa\bbb^2, \aaa\bbb\ccc,\aaa\ccc^2,\bbb^2\ccc,\bbb\ccc^2,\ccc^3$\\
			\cline{3-4}
			& & $\aaa\bbb\ccc$ & $\aaa^3, \aaa^2\ddd\eee, \ddd^3\eee, \aaa^3\ddd\eee, \aaa\ddd^3\eee,\aaa\ddd\eee^3, \ccc\ddd^3\eee$ \\
			\cline{3-4}
			& & $\aaa\ccc^2$ & $\aaa^2\ddd\eee, \aaa\ccc\ddd\eee, \ddd^2\eee^2, \ddd\eee^3, \aaa^3\ddd\eee, \aaa^2\ccc\ddd\eee, \aaa\ccc^2\ddd\eee, \aaa\ddd^2\eee^2,\aaa\ddd^3\eee, \aaa\ddd\eee^3, \ccc\ddd^2\eee^2, \ccc\ddd\eee^3$\\
			\hline
		\end{tabular} }

	\end{table}
	\end{proof}

	\begin{proposition}\label{bde,a2c}
		There is no non-symmetric $a^4b$-tiling with general angles, such that $\bbb\ddd\eee$ is the unique $a^2b$-vertex type and $\aaa^2\ccc$ is a vertex.
	\end{proposition}

	\begin{proof}
		By Lemma \ref{bde H}, one of the vertices in the last column of Table \ref{10}: $\{\aaa\bbb^2$, $\aaa\bbb\ccc$, $\aaa\ccc^2$, $\bbb^3$, $\bbb^2\ccc$, $\bbb\ccc^2$, $\ccc^3$, $\aaa^2\bbb^2$, $\aaa\bbb^3$, $\aaa\bbb^2\ccc$, $\aaa\bbb\ccc^2$, $\aaa\ccc^3$, $\bbb^2\ddd^2$, $\bbb\ccc\ddd^2$, $\ccc^2\ddd^2$, $\ccc^2\ddd\eee$, $\aaa^2\bbb^3$, $\aaa\bbb^4$, $\aaa\bbb^3\ccc$, $\aaa\bbb^2\ccc^2$, $\aaa\bbb^2\ddd^2$, $\aaa\bbb^2\eee^2$, $\aaa\bbb\ccc^3$, $\aaa\bbb\ccc\ddd^2$, $\aaa\bbb\ccc\eee^2$, $\aaa\ccc^4$, $\aaa\ccc^2\ddd^2$, $\aaa\ccc^2\ddd\eee$, $\aaa\ccc^2\eee^2$, $\bbb^3\ddd^2$, $\bbb^2\ccc\ddd^2$, $\bbb\ccc^2\ddd^2$, $\ccc^3\ddd^2$, $\ccc^3\ddd\eee\}$ must appear.

		For each combination of $3$ vertices, we will apply the irrational angle lemma and AAD to get Table \ref{bde+a2c}. The proof is completed by Table \ref{bde+a2c}, \ref{Tab-6.1}, \ref{Tab-6.3}, \ref{Tab-6.4}, \ref{Tab-6.5}, \ref{Tab-6.6} and \ref{Tab-6.7}.

		\begin{table}[htp]
			\centering
			\caption{$\bbb\ddd\eee$, $\aaa^2\ccc$ are vertices.}\label{bde+a2c}
			\resizebox{\linewidth}{!}{\begin{tabular}{|c|c|c|c|c|c|}
					\hline
					\multicolumn{2}{|c|}{Vertex}& Angles & Irrational Angle Lemma & AAD & Contradiction\\
					\hline\hline
					\multirow{34}{*}{$\bbb\ddd\eee,\aaa^2\ccc$} & $\aaa\bbb^2\eee^2$& $(1-\frac4f,\frac12-\eee+\frac2f,\frac8f,\frac32-\frac2f,\eee)$ & $n_2=n_5$& \multirow{3}{*}{$\thin^{\ccc}\eee^{\ddd}\thick^{\ddd}\eee^{\ccc}\thin\cdots$} & \multirow{2}{*}{No $\ddd^2\cdots$ by the} \\
					\cline{2-4}
					& $\aaa\bbb\ccc\eee^2$& $(1-\frac4f,3-2\ddd+\frac4f,\frac8f,\ddd,-1+\ddd-\frac4f)$ & $2n_2=n_4+n_5$& & \multirow{2}{*}{irrational angle lemma}\\
					\cline{2-4}
					&  $\aaa\ccc^2\eee^2$& $(1-\frac4f,\frac32-\ddd+\frac6f,\frac8f,\ddd,\frac12-\frac6f)$ &\multirow{4}{*}{$n_2=n_4$}& & \\
					\cline{2-3}
					\cline{5-6}
					& $\bbb^2\ddd^2$& $(1-\frac4f,1-\ddd,\frac8f,\ddd,1)$ & & \multirow{10}{*}{$\thin^{\bbb}\ddd^{\eee}\thick^{\eee}\ddd^{\bbb}\thin\cdots$} & \multirow{10}{*}{No $\eee^2\cdots$}\\
					\cline{2-3}
					& $\aaa\bbb^2\ddd^2$& $(1-\frac4f,\frac12-\ddd+\frac2f,\frac8f,\ddd,\frac32-\frac2f)$ & & & \\
					\cline{2-3}
					& $\bbb^2\ccc\ddd^2$& $(1-\frac4f,1-\ddd-\frac4f,\frac8f,\ddd,1+\frac4f)$ & &  & \\
					\cline{2-4}
					& $\bbb\ccc\ddd^2$& $(1-\frac4f,2-2\ddd-\frac8f,\frac8f,\ddd,\ddd+\frac8f)$ & \multirow{3}{*}{$2n_2=n_4+n_5$}& & \\
					\cline{2-3}
					& $\aaa\bbb\ccc\ddd^2$& $(1-\frac4f,1-2\ddd-\frac4f,\frac8f,\ddd,1+\ddd+\frac4f)$ & & & \\
					\cline{2-3}
					& $\bbb\ccc^2\ddd^2$& $(1-\frac4f,2-2\ddd-\frac{16}f,\frac8f,\ddd,\ddd+\frac{16}f)$ &  && \\
					\cline{2-4}
					& $\ccc^2\ddd^2$& $(1-\frac4f,1-\eee+\frac8f,\frac8f,1-\frac8f,\eee)$ &  \multirow{3}{*}{$n_2=n_5$}& & \\
					\cline{2-3}
					& $\aaa\ccc^2\ddd^2$& $(1-\frac4f,\frac32-\eee+\frac6f,\frac8f,\frac12-\frac6f,\eee)$ & & & \\
					\cline{2-3}
					& $\ccc^3\ddd^2$& $(1-\frac4f,1-\eee+\frac{12}f,\frac8f,1-\frac{12}f,\eee)$ & & & \\
					\cline{2-4}
					& $\bbb^3\ddd^2$& $(1-\frac4f,\frac23-\frac{2}3\ddd,\frac8f,\ddd,\frac43-\frac{1}3\ddd)$ & $2n_2+n_5=3n_4$& & \\
					\cline{2-6}
					& $\aaa^2\bbb^2$& $(1-\frac4f,\frac4f,\frac8f,\ddd,2-\ddd-\frac4f)$ &  \multirow{17}{*}{$n_4=n_5$} & &No degree $\ge4$ vertex $\ddd\eee\cdots$\\
					\cline{2-3}
					& $\aaa^2\bbb^3$& $(1-\frac4f,\frac8{3f},\frac8f,\ddd,2-\ddd-\frac8{3f})$ & & &(Similar to Proposition \ref{bde,a2b}) \\
					\cline{2-3}
					\cline{6-6}
					& $\aaa\bbb^2$& $(1-\frac4f,\frac12+\frac2f,\frac8f,\ddd,\frac32-\ddd-\frac2f)$ && &\multirow{4}{*}{See Table \ref{Tab-6.1}}\\
					\cline{2-3}
					&$\aaa\bbb\ccc$&\multirow{2}{*}{$(1-\frac{4}{f},1-\frac{4}{f},\frac{8}{f},\ddd,1-\ddd+\frac{4}{f})$} & & & \\
					\cline{2-2}
					& $\bbb^2\ccc$&  & & & \\
					\cline{2-3}
					& $\bbb^3$& $(1-\frac4f,\frac23,\frac8f,\ddd,\frac43-\ddd)$ & &  &\\
					\cline{2-3}
					\cline{6-6}
					& $\aaa\bbb^3$& $(1-\frac4f,\frac13+\frac4{3f},\frac8f,\ddd,\frac53-\ddd-\frac4{3f})$ & & &  \multirow{2}{*}{See Table \ref{Tab-6.3}}\\
					\cline{2-3}
					& $\aaa\bbb\ccc^2$& $(1-\frac4f,1-\frac{12}f,\frac8f,\ddd,1-\ddd+\frac{12}f)$ & & & \\
					\cline{2-3}
					\cline{6-6}
					& $\aaa\bbb^2\ccc$& $(1-\frac4f,\frac12-\frac2f,\frac8f,\ddd,\frac32-\ddd+\frac2f)$ & & & \multirow{2}{*}{See Table \ref{Tab-6.4}} \\
					\cline{2-3}
					& $\aaa\bbb^3\ccc$& $(1-\frac4f,\frac13-\frac4{3f},\frac8f,\ddd,\frac53-\ddd+\frac4{3f})$& & & \\
					\cline{2-3}
					\cline{6-6}
					& $\aaa\bbb^4$& $(1-\frac4f,\frac14+\frac1f,\frac8f,\ddd,\frac74-\ddd-\frac1f)$ & & & \multirow{2}{*}{See Table \ref{Tab-6.5}}\\
					\cline{2-3}
					& $\aaa\bbb^2\ccc^2$& $(1-\frac4f,\frac12-\frac6f,\frac8f,\ddd,\frac32-\ddd+\frac6f)$ & & & \\
					\cline{2-3}
					\cline{6-6}
					& $\aaa\bbb\ccc^3$& $(1-\frac4f,1-\frac{20}f,\frac8f,\ddd,1-\ddd+\frac{20}f)$ & & & \multirow{2}{*}{See Table \ref{Tab-6.6}} \\
					\cline{2-3}
					& $\aaa\ccc^2\ddd\eee$& $(1-\frac4f,1+\frac{12}f,\frac8f,\ddd,1-\ddd-\frac{12}f)$ & & & \\
					\cline{2-3}
					\cline{6-6}
					& $\ccc^2\ddd\eee$& $(1-\frac4f,\frac{16}f,\frac8f,\ddd,2-\ddd-\frac{16}f)$ & & & \multirow{2}{*}{See Table \ref{Tab-6.7}} \\
					\cline{2-3}
					& $\ccc^3\ddd\eee$& $(1-\frac4f,\frac{24}f,\frac8f,\ddd,2-\ddd-\frac{24}f)$ & & & \\
					\cline{2-3}
					\cline{6-6}
					& $\bbb\ccc^2$& $(1-\frac4f,2-\frac{16}f,\frac8f,\ddd,-\ddd+\frac{16}f)$ & & & See page \pageref{bde+a2c+bc2} \\
					\cline{2-6}
					& $\aaa\ccc^2$& \multirow{2}{*}{$(\frac23,2-\ddd-\eee,\frac23,\ddd,\eee),f=12$} &\multirow{4}{*}{}& & \multirow{4}{*}{See page \pageref{bde+a2c+ac3}} \\
					\cline{2-2}
					& $\ccc^3$& & & &   \\
					\cline{2-3}
					& $\aaa\ccc^3$& $(\frac45,2-\ddd-\eee,\frac25,\ddd,\eee),f=20$ &  & & \\
					\cline{2-3}
					& $\aaa\ccc^4$& $(\frac67,2-\ddd-\eee,\frac27,\ddd,\eee),f=28$ &  & &  \\
					\hline
			\end{tabular} }

		\end{table}

		\begin{table}[htp]
			\centering
			\caption{AVC for Cases $\{\bbb\ddd\eee,\aaa^2\ccc,\aaa\bbb^2\}$, $\{\bbb\ddd\eee,\aaa^2\ccc,\aaa\bbb\ccc(\bbb^2\ccc)\}$ and $\{\bbb\ddd\eee,\aaa^2\ccc,\bbb^3\}$.}\label{Tab-6.1}
			\begin{minipage}{1\textwidth}
				\resizebox{\linewidth}{!}{\begin{tabular}{c|c||c|c||c|c}
						\hline
						\multicolumn{2}{c||}{AVC for $\{\bbb\ddd\eee,\aaa^2\ccc,\aaa\bbb^2\}$}&\multicolumn{2}{c||}{AVC for $\{\bbb\ddd\eee,\aaa^2\ccc,\aaa\bbb\ccc(\bbb^2\ccc)\}$}&\multicolumn{2}{c}{AVC for $\{\bbb\ddd\eee,\aaa^2\ccc,\bbb^3\}$}\\
						\hline
						$f$&vertex&$f$&vertex&$f$&vertex\\
						\hline
						\hline
						all&$\bbb\ddd\eee,\aaa^2\ccc,\aaa\bbb^2$&all&$\bbb\ddd\eee,\aaa^2\ccc,\aaa\bbb\ccc,\bbb^2\ccc$&all&$\bbb\ddd\eee,\aaa^2\ccc,\bbb^3$\\
						$16s+12,s=1,2,\cdots$&$\aaa\bbb\ccc^{s+1},\aaa\ccc^{2s+2},\bbb^3\ccc^{s},$&$8s+4,s=2,3,\cdots$&$\aaa\ccc^{s+1},\bbb\ccc^{s+1},\ccc^s\ddd\eee,$&$24s,s=2,3,\cdots$&$\bbb^2\ccc^{2s},\bbb\ccc^{4s},\ccc^{6s},\ccc^{2s}\ddd\eee$\\
					    &$\bbb^2\ccc^{2s+1},\bbb\ccc^{3s+2},\ccc^{s+1}\ddd\eee,$&&$\ccc^{2s+1}$&$24s+12,s=2,3,\cdots$&$\aaa\bbb\ccc^{s+1},\aaa\ccc^{3s+2},\bbb^2\ccc^{2s+1},$\\
					    &$\ccc^{4s+3}$&&&&$\bbb\ccc^{4s+2},\ccc^{6s+3},\ccc^{2s+1}\ddd\eee$\\
						\hline
						\multicolumn{6}{c}{By the balance lemma, $\ccc^k\ddd\eee$ must appear. Its AAD gives $\aaa\thin\ddd$, $\aaa\thin\eee$, $\ddd\thin\ddd$, $\ddd\thin\eee$ or $\eee\thin\eee$, contradicting the AVC}\\
						\hline
				\end{tabular}}
			\end{minipage}
		\end{table}

			\begin{table}[htp]
				\centering
				\caption{AVC for Cases $\{\bbb\ddd\eee, \aaa^2\ccc,\aaa\bbb^3\}$ and $\{\bbb\ddd\eee,\aaa^2\ccc,\aaa\bbb\ccc^2\}$.}\label{Tab-6.3}
				\begin{minipage}{1\textwidth}
					\resizebox{\linewidth}{!}{\begin{tabular}{c|c||c|c}
						\hline
						\multicolumn{2}{c||}{AVC for $\{\bbb\ddd\eee,\aaa^2\ccc,\aaa\bbb^3\}$}&\multicolumn{2}{c}{AVC for $\{\bbb\ddd\eee, \aaa^2\ccc, \aaa\bbb\ccc^2\}$}\\
						\hline
						$f$&vertex&$f$&vertex\\
						\hline
						\hline
						all&$\bbb\ddd\eee,\aaa^2\ccc,\aaa\bbb^3$&all&$\bbb\ddd\eee, \aaa^2\ccc, \aaa\bbb\ccc^2,\bbb^2\ccc^3$ \\
						$24s+20,s=1,2,\cdots$&$ \aaa\bbb^2\ccc^{s+1},\aaa\bbb\ccc^{2s+2},\aaa\ccc^{3s+3},\bbb^5\ccc^s,\bbb^4\ccc^{2s+1}, $&$28$&$\aaa\bbb^2,\bbb^3\ccc,\ccc^2\ddd\eee,\aaa\ccc^4,\bbb\ccc^5,\ccc^7$\\
						&$\bbb^3\ccc^{3s+2},\bbb^2\ccc^{4s+3},\bbb\ccc^{5s+4}, \ccc^{6s+5},\ccc^{s+1}\ddd\eee $&$36$&$\bbb^3,\ccc^3\ddd\eee,\aaa\ccc^5,\bbb\ccc^6,\ccc^9$\\
						&&$8s+4,s=5,6,\cdots$&$\aaa\ccc^{s+1},\bbb\ccc^{s+2},\ccc^{2s+1},\ccc^{s-1}\ddd\eee$ \\
						\hline
						\multicolumn{4}{c}{By the balance lemma, $\ccc^k\ddd\eee$ must appear. Its AAD gives $\aaa\thin\ddd$, $\aaa\thin\eee$, $\ddd\thin\ddd$, $\ddd\thin\eee$ or $\eee\thin\eee$, contradicting the AVC} \\
						\hline
					\end{tabular}}
				\end{minipage}
			\end{table}

		\begin{table}[htp]
			\centering
			\caption{AVC for Cases $\{\bbb\ddd\eee, \aaa^2\ccc, \aaa\bbb^2\ccc\}$ and $\{\bbb\ddd\eee, \aaa^2\ccc, \aaa\bbb^3\ccc\}$.}\label{Tab-6.4}
			\begin{minipage}{1\textwidth}
				\resizebox{\linewidth}{!}{\begin{tabular}{c|c||c|c}
					\hline
					\multicolumn{2}{c||}{AVC for $\{\bbb\ddd\eee, \aaa^2\ccc, \aaa\bbb^2\ccc\}$} & \multicolumn{2}{c}{AVC for $\{\bbb\ddd\eee, \aaa^2\ccc, \aaa\bbb^3\ccc\}$}\\
					\hline
					$f$&vertex&$f$&vertex\\
					\hline
					\hline
					all&$\bbb\ddd\eee, \aaa^2\ccc,\aaa\bbb^2\ccc,\bbb^4\ccc$ &all&$\bbb\ddd\eee, \aaa^2\ccc, \aaa\bbb^3\ccc,\bbb^6\ccc$\\
					$16s+4,s=2,3,\cdots$&$\aaa\bbb\ccc^{s+1},\aaa\ccc^{2s+1},\bbb^3\ccc^{s+1},\bbb^2\ccc^{2s+1},$&$24s+4,s=2,3,\cdots$&$\aaa\bbb^2\ccc^{s+1},\aaa\bbb\ccc^{2s+1},\aaa\ccc^{3s+1},\bbb^5\ccc^{s+1},\bbb^4\ccc^{2s+1},$\\
					&$\bbb\ccc^{3s+1},\ccc^{4s+1},\ccc^s\ddd\eee$&&$\bbb^3\ccc^{3s+1},\bbb^2\ccc^{4s+1},\bbb\ccc^{5s+1},\ccc^{6s+1},\ccc^s\ddd\eee$\\
					\hline
					\multicolumn{4}{c}{By the balance lemma, $\ccc^k\ddd\eee$ must appear. Its AAD gives $\aaa\thin\ddd$, $\aaa\thin\eee$, $\ddd\thin\eee$, $\ddd\thin\ddd$ or $\eee\thin\eee$, contradicting the AVC} \\
					\hline
				\end{tabular}}
			\end{minipage}
		\end{table}

		\begin{table}[htp]
			\centering
			\caption{AVC for Cases $\{\bbb\ddd\eee, \aaa^2\ccc, \aaa\bbb^4\}$ and $\{\bbb\ddd\eee,\aaa^2\ccc,\aaa\bbb^2\ccc^2\}$.}\label{Tab-6.5}
			\begin{minipage}{1\textwidth}
				\resizebox{\linewidth}{!}{\begin{tabular}{c|c||c|c}
						\hline
						\multicolumn{2}{c||}{AVC for $\{\bbb\ddd\eee, \aaa^2\ccc, \aaa\bbb^4\}$} &\multicolumn{2}{c}{AVC for $\{\bbb\ddd\eee, \aaa^2\ccc, \aaa\bbb^2\ccc^2\}$} \\
						\hline
						$f$&vertex&$f$&vertex\\
						\hline
						\hline
						all&$\bbb\ddd\eee, \aaa^2\ccc, \aaa\bbb^4$&all&$\bbb\ddd\eee,\aaa^2\ccc,\aaa\bbb^2\ccc^2,\bbb^4\ccc^3$ \\
						$32s+28,s=1,2,\cdots$&$\aaa\bbb^3\ccc^{s+1},\aaa\bbb^2\ccc^{2s+2},\aaa\bbb\ccc^{3s+3},\aaa\ccc^{4s+4},$&$44$&$\aaa\bbb^{3},\ccc^{2}\ddd\eee,\aaa\bbb\ccc^{4},\bbb^5\ccc,\aaa\ccc^{6},\bbb^3\ccc^5,\bbb^2\ccc^{7},\bbb\ccc^{9},\ccc^{11}$\\
						&$\bbb^7\ccc^s,\bbb^6\ccc^{2s+1},\bbb^5\ccc^{3s+2},\bbb^4\ccc^{4s+3},$&$60$&$\bbb^5,\ccc^{3}\ddd\eee,\aaa\bbb\ccc^{5},\aaa\ccc^{8},\bbb^3\ccc^6,\bbb^2\ccc^{9},\bbb\ccc^{12},\ccc^{15}$\\
						&$\bbb^3\ccc^{5s+4},\bbb^2\ccc^{6s+5},\bbb\ccc^{7s+6},\ccc^{8s+7},\ccc^{s+1}\ddd\eee$&$16s+12,s=4,5,\cdots$&$\aaa\bbb\ccc^{s+2},\aaa\ccc^{2s+2},\bbb^3\ccc^{s+3},\bbb^2\ccc^{2s+3},\bbb\ccc^{3s+3},\ccc^{4s+3},\ccc^{s}\ddd\eee$ \\
						\hline
						\multicolumn{4}{c}{By the balance lemma, $\ccc^k\ddd\eee$ must appear. Its AAD gives $\aaa\thin\ddd$, $\aaa\thin\eee$, $\ddd\thin\eee$, $\ddd\thin\ddd$ or $\eee\thin\eee$, contradicting the AVC} \\
						\hline
				\end{tabular}}
			\end{minipage}
		\end{table}

		\begin{table}[htp]
			\centering
			\caption{AVC for Cases $\{\bbb\ddd\eee,\aaa^2\ccc,\aaa\bbb\ccc^3\}$ and $\{\bbb\ddd\eee,\aaa^2\ccc,\aaa\ccc^2\ddd\eee\}$.}\label{Tab-6.6}
			\begin{minipage}{1\textwidth}
				\resizebox{\linewidth}{!}{\begin{tabular}{c|c||c|c}
					\hline
					$f$&vertex&$f$&vertex\\
					\hline
					\hline
					all&$\bbb\ddd\eee,\aaa^2\ccc,\aaa\bbb\ccc^3,\bbb^2\ccc^5$&all&$\bbb\ddd\eee,\aaa^2\ccc,\aaa\ccc^2\ddd\eee,\ccc^3\ddd^2\eee^2$\\
					$36$&$\aaa\bbb^2\ccc,\ccc^2\ddd\eee,\aaa\ccc^5,\bbb^4\ccc,\bbb^3\ccc^3,\bbb\ccc^7,\ccc^9$&$28$&$\bbb\ccc^2,\aaa\ccc^4,\aaa\ddd^2\eee^2,\ccc^5\ddd\eee,\ccc\ddd^3\eee^3,\ccc^7$\\
					$44$&$\aaa\bbb^2,\bbb^3\ccc^2,\ccc^3\ddd\eee,\aaa\ccc^6,\bbb\ccc^8,\ccc^{11}$&$36$&$\bbb\ccc^3,\aaa\ccc^5,\ddd^3\eee^3,\ccc^6\ddd\eee,\ccc^9$\\
					$52$&$\bbb^3\ccc,\ccc^4\ddd\eee,\aaa\ccc^7,\bbb\ccc^9,\ccc^{13}$&$8s+4,s=5,6,\cdots$&$\aaa\ccc^{s+1},\bbb\ccc^{s-1},\ccc^{2s+1},\ccc^{s+2}\ddd\eee$\\
					$60$&$\bbb^3,\ccc^5\ddd\eee,\aaa\ccc^8,\bbb\ccc^{10},\ccc^{15}$&& \\
					$8s+4,s=8,9\cdots$&$\aaa\ccc^{s+1},\bbb\ccc^{s+3},\ccc^{2s+1},\ccc^{s-2}\ddd\eee$&&\\
					\hline
					\multicolumn{2}{c||}{By the balance lemma, $\ccc^k\ddd\eee$ must appear.} & \multicolumn{2}{c}{By the balance lemma, $\bbb\ccc^{k}$ must appear. }\\
					\multicolumn{2}{c||}{Its AAD gives $\aaa\thin\ddd$, $\aaa\thin\eee$, $\ddd\thin\ddd$, $\ddd\thin\eee$ or $\eee\thin\eee$, contradicting the AVC} & \multicolumn{2}{c}{Its AAD gives $\aaa\thin\bbb$, $\bbb\thin\bbb$ or $\ddd\thin\eee$, contradicting the AVC}\\
					\hline
				\end{tabular}}
			\end{minipage}
		\end{table}

		\begin{table}[htp]
			\centering
			\caption{AVC for Cases $\{\bbb\ddd\eee,\aaa^2\ccc,\ccc^2\ddd\eee\}$ and $\{\bbb\ddd\eee,\aaa^2\ccc,\ccc^3\ddd\eee\}$.}\label{Tab-6.7}
			\begin{minipage}{1\textwidth}
				\resizebox{\linewidth}{!}{\begin{tabular}{c|c||c|c}
						\hline
						\multicolumn{2}{c||}{AVC for $\{\bbb\ddd\eee,\aaa^2\ccc,\ccc^2\ddd\eee\}$}&\multicolumn{2}{c}{AVC for $\{\bbb\ddd\eee,\aaa^2\ccc,\ccc^3\ddd\eee\}$}\\
						\hline
						$f$&vertex&$f$&vertex\\
						\hline
						\hline
						all&$\bbb\ddd\eee,\aaa^2\ccc,\ccc^2\ddd\eee$&all&$\bbb\ddd\eee,\aaa^2\ccc,\ccc^3\ddd\eee$\\
						$16$&$\bbb\ccc^2,\ccc^4,\ddd^2\eee^2$&$20$&$\bbb\ccc^2,\aaa\ccc^3,\aaa\ccc\ddd\eee,\ccc^5,\ccc\ddd^2\eee^2$ \\
						$20$&$\aaa\bbb\ccc,\bbb^2\ccc,\aaa\ccc^3,\bbb\ccc^3,\ccc^{5}$&$24$&$\bbb\ccc^3,\ddd^2\eee^2,\ccc^6$ \\
						$8s,s=3,4,\cdots$&$\bbb^{x}\ccc^{2s-2x}(x=0,1,\cdots,s)$ &$28$&$\aaa\bbb\ccc,\bbb^2\ccc,\aaa\ccc^4,\bbb\ccc^4,\ccc^7$\\
						$8s+4,s=3,4,\cdots$&$\bbb^{x}\ccc^{2s+1-2x}(x=0,1,\cdots,s),$ &$8s,s=4,5,\cdots$&$\bbb^{x}\ccc^{2s-3x}(x=0,1,\cdots,\lfloor\frac{2s}3\rfloor)$\\
						&$\aaa\bbb^y\ccc^{s+1-2y}(y=0,1,\cdots,\lfloor\frac{s-1}2\rfloor)$& $8s+4,s=4,5,\cdots$&$\bbb^{x}\ccc^{2s+1-3x}(x=0,1,\cdots,\lfloor\frac{2s+1}3\rfloor),$\\
						&&&$\aaa\bbb^y\ccc^{s+1-3y}(y=0,1,\cdots,\lfloor\frac{s+1}3\rfloor)$\\
						\hline
						\multicolumn{4}{c}{By the balance lemma, $\bbb^k\ccc^l$ or $\aaa\bbb^k\ccc^l$ must appear. Its AAD gives $\aaa\thin\bbb$,  $\bbb\thin\bbb$, $\aaa\thin\ddd$, $\aaa\thin\eee$, $\ddd\thin\ddd$, $\ddd\thin\eee$ or $\eee\thin\eee$, contradicting the AVC} \\
						\hline
				\end{tabular}}
			\end{minipage}
		\end{table}

		\subsection*{Case. $\{\bbb\ddd\eee, \aaa^2\ccc,\bbb\ccc^2\}$} \label{bde+a2c+bc2}
		By the irrational angle lemma, we have $n_4=n_5$. The angle sums of $\bbb\ddd\eee, \aaa^2\ccc,\bbb\ccc^2$ and the angle sum for the pentagon imply
		\begin{align*}
			\aaa=1-\tfrac4f,\;
			\bbb=2-\tfrac{16}f,\;
			\ccc=\tfrac8f,\;
			\ddd+\eee=\tfrac{16}f.
		\end{align*}
		By the angle values, we get the
		\[
		\text{AVC}\sub\{\bbb\ddd\eee, \aaa^2\ccc,\aaa\bbb\ccc,\bbb\ccc^2,\aaa\ccc^{\frac{f+4}8},\aaa\ccc^{x}\ddd^{\frac{f-8x+4}{16}}\eee^{\frac{f-8x+4}{16}} ,\ccc^y\ddd^{\frac{f-4y}8}\eee^{\frac{f-4y}8} \},
		\]
		where $x=0,1,2,\cdots,\lfloor\tfrac{f+4}8\rfloor$, $y=0,1,2,\cdots,\lfloor\tfrac{f}4\rfloor$. If $f\neq4k(k\geq3)$, we know $\aaa\ccc^{x}\ddd^{\frac{f-8x+4}{16}}\eee^{\frac{f-8x+4}{16}},\ccc^y\ddd^{\frac{f-4y}8}\eee^{\frac{f-4y}8}$ are not vertices. Then $\ddd\cdots=\bbb\ddd\eee$. By the parity lemma, we have $\bbb\cdots=\bbb\ddd\eee$. Therefore $\bbb\ccc^2$ is not a vertex, a contradiction.

		If $f=12$, we get $\bbb=\ccc=\frac23$, a contradiction.
		If $f=16$, by the equation (\ref{coolsaet_eq1}), we have
		\begin{align*}
			\aaa=\tfrac34,\;
			\bbb=1,\;
			\ccc=\tfrac12,\;
			\tan(\ddd\pi)=\tfrac{(5\sqrt2+4)\sqrt{886-572\sqrt2}}{68}-\tfrac{\sqrt{2}}2+\tfrac12\approx1.223. \\
			\eee=1-\ddd\approx0.7182,\;
			\cos a=\tfrac{(2\sqrt2+5)\sqrt{443-286\sqrt2}}{68}+\tfrac{\sqrt2}4-\tfrac14\approx0.8182.
		\end{align*}
		Here $\delta$ and $\epsilon$ are normalized by $\pi$, while $a$ is in radians. Its $3D$ picture is shown in Figure \ref{f=16}. Therefore $\text{AVC}\sub\{\bbb\ddd\eee, \aaa^2\ccc, \aaa\bbb\ccc, \bbb\ccc^2, \ddd^2\eee^2,$ $\ccc^2\ddd\eee\}$. The AAD of $\bbb\ccc^2$ is $\thin^{\aaa}\bbb^{\ddd}\thin^{\aaa}\ccc^{\eee}\thin\ccc\thin$, $\thin^{\aaa}\bbb^{\ddd}\thin^{\eee}\ccc^{\aaa}\thin^{\aaa}\ccc^{\eee}\thin$ or $\thin^{\aaa}\bbb^{\ddd}\thin^{\eee}\ccc^{\aaa}\thin^{\eee}\ccc^{\aaa}\thin$. This implies $\aaa\ddd\cdots$ or $\aaa\eee\cdots$ is a vertex, contradicting the AVC.
		\begin{figure}[htp]
			\centering
			\includegraphics[scale=0.2]{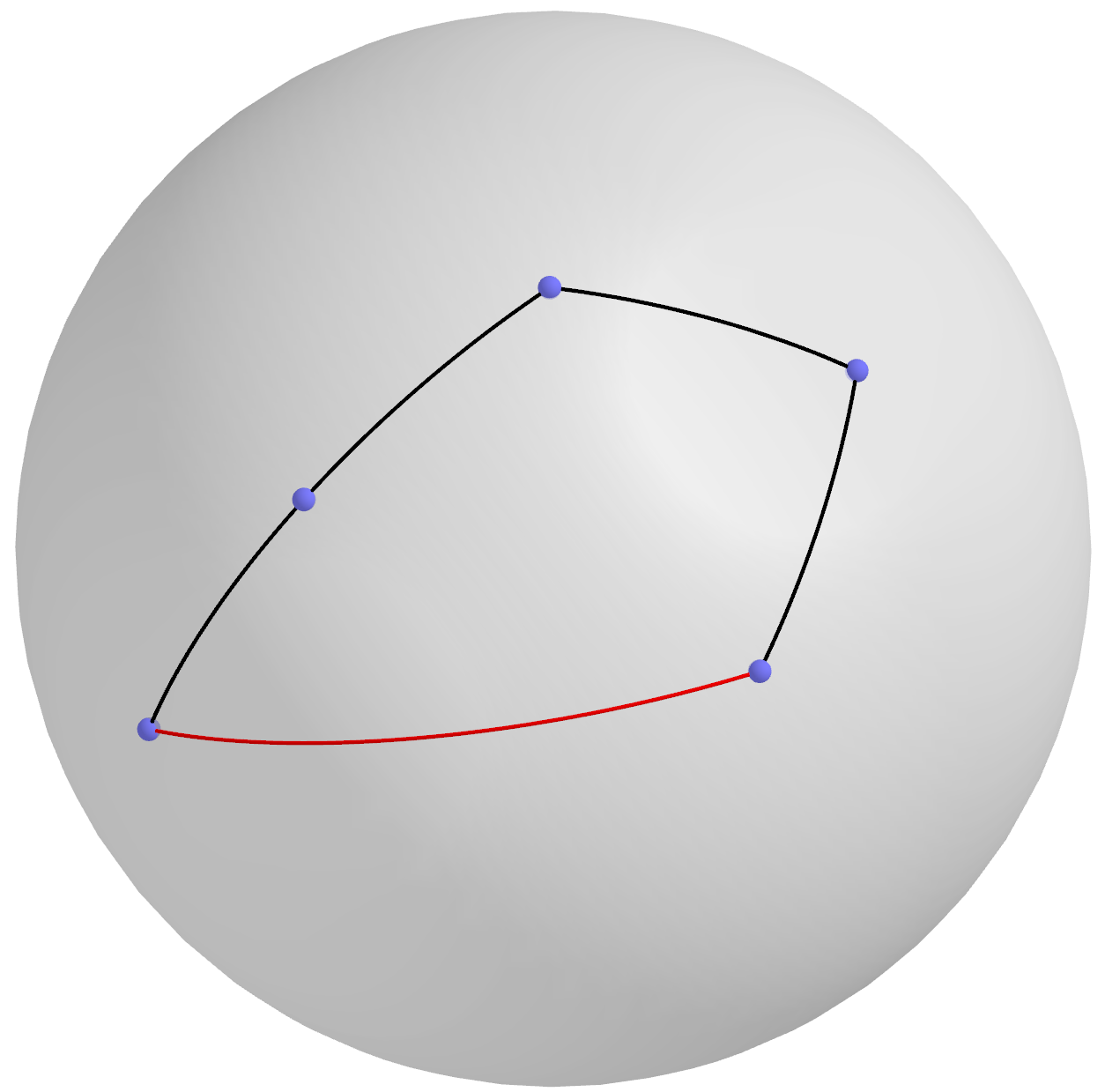}
			\caption{Case \{$\bbb\ddd\eee, \aaa^2\ccc,\bbb\ccc^2$\} for $f=16$.}
			\label{f=16}
		\end{figure}

		If $f\geq20$, we get four solutions of $\ddd$ to the equation (\ref{coolsaet_eq1}), but none of them lies in $(0,1)$. In other words, there is no such pentagon.

		\subsection*{Case. $\{\bbb\ddd\eee, \aaa^2\ccc,\aaa\ccc^3\}$} \label{bde+a2c+ac3}
			The angle sums of $\bbb\ddd\eee$, $\aaa^2\ccc$, $\aaa\ccc^3$ and the angle sum for the pentagon imply
			\begin{align*}
				\aaa=\tfrac45,\;
				\bbb+\ddd+\eee=2,\;
				\ccc=\tfrac25,\;
				f=20.
			\end{align*}

			\subsubsection*{Subcase. $\bbb<\ccc$}
			By Lemma \ref{geometry1}, we have $\ddd>\eee$. By $\bbb\ddd\eee$, we can obtain that $\ddd+\eee>\frac85$ and $\ddd>\frac45$. By $R(\ddd\eee\cdots)$, $R(\ddd^2\cdots)<$ all angles and the parity lemma, we get $\ddd\cdots=\bbb\ddd\eee$. Therefore $\text{AVC}\sub\{\bbb\ddd\eee, \aaa^2\ccc,\aaa\ccc^3,\ccc^5\}$.

			The AAD is $\aaa^2\ccc=\aaa^{\bbb}\thin^{\bbb}\aaa\cdots$, $\aaa^{\bbb}\thin^{\ccc}\aaa\cdots$ or $\aaa^{\bbb}\thin^{\aaa}\ccc\cdots$. This implies $\aaa\bbb\cdots$, $\bbb^2\cdots$ or $\bbb\ccc\cdots$ is a vertex, contradicting the AVC.

			\subsubsection*{Subcase. $1>\bbb>\ccc$}
			By Lemma \ref{geometry1}, we have $\ddd<\eee$. By $\bbb\ddd\eee$ and $1>\bbb>\ccc$, we know $\ddd+\eee>1$. By $\ddd<\eee$, we get $\eee>\frac12$. By $R(\eee^2\cdots)<\bbb$, $2\ddd$, $\ddd+\eee$, $2\eee$, $2\aaa$, $\aaa+\ccc$, $3\ccc$, we get $\eee^2\cdots=\ccc^2\eee^2$. If $\ccc^2\eee^2$ is a vertex, by the angle value, we have $\aaa=\tfrac45$, $\bbb+\ddd=\tfrac75$, $\ccc=\tfrac25$, $\eee=\tfrac35$. Applying the irrational angle lemma to this case, we have $n_2=n_4$. Then $\bbb\cdots=\bbb\ddd\cdots=\bbb\ddd\eee$. Therefore, $\ccc^2\eee^2$ is not a vertex. We get that $\eee^2\cdots$ is not a vertex. By $R(\aaa\ddd\eee\cdots)<\aaa$, $\bbb$, $\ccc$ and the balance lemma, we get no $\aaa\eee\cdots$. The AAD of $\aaa\ccc^3$ is $\thin^{\aaa}\ccc^{\eee}\thin^{\aaa}\ccc^{\eee}\thin\ccc\thin$ or $\thin^{\eee}\ccc^{\aaa}\thin^{\aaa}\ccc^{\eee}\thin\ccc\thin$. This implies $\aaa\eee\cdots$ or $\eee^2\cdots$ is a vertex, a contradiction.

			\subsubsection*{Subcase. $\bbb\geq1>\ccc$}
			By Lemma \ref{geometry1}, we have $\ddd<\eee$. By $\bbb\ddd\eee$ and $\bbb\geq1>\ccc$, we know $\ddd+\eee\leq1$. By $\ddd<\eee$, we get $\ddd<\frac12$. By the angle values, we get $\bbb\cdots=\bbb\ddd\eee$, $\bbb\ddd^k$, $\aaa\bbb\ddd^l$, $\bbb\ccc\ddd^m$ or $\bbb\ccc^2\ddd^n$.
			If $\bbb\ccc^2$ is a vertex, it actually belongs to Case. \{$\bbb\ddd\eee, \aaa^2\ccc,\bbb\ccc^2$\}. By the balance lemma, we get $\bbb\cdots=\bbb\ddd\eee$. Therefore, $\text{AVC}\sub\{\bbb\ddd\eee, \aaa^2\ccc,\aaa\ccc^3,\ccc^5\}$. We get the same contradiction as Subcase $\bbb<\ccc$.

		\vspace{9pt}

		\noindent Similar to Case \{$\bbb\ddd\eee, \aaa^2\ccc,\aaa\ccc^3$\}, three other cases \{$\bbb\ddd\eee, \aaa^2\ccc,\aaa\ccc^2$\}, \{$\bbb\ddd\eee, \aaa^2\ccc,$ $\ccc^3$\} and \{$\bbb\ddd\eee, \aaa^2\ccc,\aaa\ccc^4$\} admit no tilings.

	\end{proof}

	\begin{proposition}\label{bde,abc}
		There is no non-symmetric $a^4b$-tiling with general angles, such that $\bbb\ddd\eee$ is the unique $a^2b$-vertex type and $\aaa\bbb\ccc$ is a vertex.
	\end{proposition}

	\begin{proof}
		By Lemma \ref{bde H}, one of the vertices in the last column of Table \ref{10}: $\{\aaa^3$, $\aaa^4$, $\aaa^3\bbb$, $\aaa^3\ccc$, $\aaa^2\ddd^2$, $\aaa^2\ddd\eee$, $\aaa^2\eee^2$, $\aaa\ccc\ddd^2$, $\ddd^3\eee$, $\aaa^5$, $\aaa^4\bbb$, $\aaa^4\ccc$, $\aaa^3\bbb^2$, $\aaa^3\ddd^2$, $\aaa^3\ddd\eee$, $\aaa\ddd\eee^3$, $\aaa^3\eee^2$, $\aaa^2\bbb\ddd^2$, $\aaa^2\bbb\eee^2$, $\aaa^2\ccc\ddd^2$, $\aaa^2\ccc\eee^2$, $\aaa\ccc^2\ddd^2$, $\aaa\ddd^3\eee$, $\ccc\ddd^4$, $\ccc\ddd^3\eee\}$ must appear.

		For each combination of $3$ vertices, we will apply the irrational angle lemma and AAD to get Table \ref{bde+abc}. The proof is completed by the contradictions in the last column of Table \ref{bde+abc} and the last row of Table \ref{Tab-7.1}.

		\begin{table}[htp]
			\centering
			\caption{$\bbb\ddd\eee$, $\aaa\bbb\ccc$ are vertices.}\label{bde+abc}
			\resizebox{\linewidth}{!}{\begin{tabular}{|c|c|c|c|c|c|}
					\hline
					\multicolumn{2}{|c|}{Vertex}& Angles & Irrational Angle Lemma & AAD & Contradiction\\
					\hline\hline
					\multirow{23}{*}{$\bbb\ddd\eee,\aaa\bbb\ccc$} & $\aaa^2\eee^2$& $(\ddd-\frac4f,1-\frac4f,1-\ddd+\frac8f,\ddd,1-\ddd+\frac4f)$ & \multirow{2}{*}{$n_1+n_4=n_3+n_5$}& \multirow{4}{*}{$\thin^{\ccc}\eee^{\ddd}\thick^{\ddd}\eee^{\ccc}\thin\cdots$} & \multirow{4}{*}{No $\ddd^2\cdots$ by the} \\
					\cline{2-3}
					& $\aaa^2\bbb\eee^2$& $(-\frac12+\ddd-\frac{2}{f},1-\frac4f,\frac32-\ddd+\frac{6}{f},\ddd,1-\ddd+\frac4f)$ && & \multirow{4}{*}{irrational angle lemma}\\
					\cline{2-4}
					& $\aaa^3\eee^2$& $(\frac{2}3\ddd-\frac8{3f},1-\frac4f,1-\frac{2}3\ddd+\frac{20}{3f},\ddd,1-\ddd+\frac4f)$ &$2n_1+3n_4=2n_3+3n_5$& &\\
					\cline{2-4}
					&  $\aaa^2\ccc\eee^2$& \multirow{2}{*}{$(-1+2\ddd-\frac{12}f,1-\frac4f,2-2\ddd+\frac{16}f,\ddd,1-\ddd+\frac4f)$} & \multirow{2}{*}{$2n_1+n_4=2n_3+n_5$}& & \\
					\cline{2-2}
					&$\aaa\ddd\eee^3$&&&&\\
					\cline{2-6}
					& $\aaa^2\ddd^2$& $(1-\ddd,1-\frac4f,\ddd+\frac4f,\ddd,1-\ddd+\frac4f)$ & \multirow{2}{*}{$n_1+n_5=n_3+n_4$}& \multirow{9}{*}{$\thin^{\bbb}\ddd^{\eee}\thick^{\eee}\ddd^{\bbb}\thin\cdots$} & \multirow{9}{*}{No $\eee^2\cdots$}\\
					\cline{2-3}
					& $\aaa^2\bbb\ddd^2$& $(\frac12-\ddd+\frac{2}{f},1-\frac4f,\frac12+\ddd+\frac{2}{f},\ddd,1-\ddd+\frac4f)$ && & \\
					\cline{2-4}
					& $\aaa^3\ddd^2$& $(\frac23-\frac{2}3\ddd,1-\frac4f,\frac13+\frac{2}3\ddd+\frac4f,\ddd,1-\ddd+\frac4f)$ &$2n_1+3n_5=2n_3+3n_4$& & \\
					\cline{2-4}
					& $\aaa\ccc\ddd^2$& \multirow{2}{*}{$(1-\ccc+\frac4f,1-\frac4f,\ccc,\frac12-\frac2f,\frac12+\frac6f)$} & \multirow{2}{*}{$n_1=n_3$}& & \\
					\cline{2-2}
					& $\ddd^3\eee$& & & & \\
					\cline{2-4}
					& $\aaa^2\ccc\ddd^2$& \multirow{2}{*}{$(1-2\ddd-\frac4f,1-\frac4f,2\ddd+\frac8f,\ddd,1-\ddd+\frac4f)$} & \multirow{2}{*}{$2n_1+n_5=2n_3+n_4$}& & \\
					\cline{2-2}
					& $\aaa\ddd^3\eee$& & & & \\
					\cline{2-4}
					& $\aaa\ccc^2\ddd^2$& \multirow{2}{*}{$(2\ddd+\frac8f,1-\frac4f,1-2\ddd-\frac4f,\ddd,1-\ddd+\frac4f)$} & \multirow{2}{*}{$2n_1+n_4=2n_3+n_5$}& & \\
					\cline{2-2}
					& $\ccc\ddd^3\eee$&  & & & \\
					\cline{2-4}
					& $\ccc\ddd^4$& $(-1+4\ddd+\frac4f,1-\frac4f,2-4\ddd,\ddd,1-\ddd+\frac4f)$ & $4n_1+n_4=4n_3+n_5$ && \\
					\cline{2-6}
					& $\aaa^4$& $(\frac12,1-\frac4f,\frac12+\frac4f,\ddd,1-\ddd+\frac4f)$ & \multirow{11}{*}{$n_4=n_5$}&  & \multirow{4}{*}{No degree $\ge4$ vertex $\ddd\eee\cdots$} \\
					\cline{2-3}
					& $\aaa^3\bbb$& $(\frac13+\frac4{3f},1-\frac4f,\frac23+\frac8{3f},\ddd,1-\ddd+\frac4f)$ & &  &\multirow{4}{*}{(similar to Proposition \ref{bde,a2b})}\\
					\cline{2-3}
					& $\aaa^5$& $(\frac25,1-\frac4f,\frac35+\frac4f,\ddd,1-\ddd+\frac4f)$ & & & \\
					\cline{2-3}
					& $\aaa^4\bbb$& $(\frac14+\frac1f,1-\frac4f,\frac34+\frac3f,\ddd,1-\ddd+\frac4f)$ & & & \\
					\cline{2-3}
					\cline{6-6}
					& $\aaa^3$& $(\frac23,1-\frac4f,\frac13+\frac4f,\ddd,1-\ddd+\frac4f)$ & &  & \multirow{6}{*}{See Table \ref{Tab-7.1}} \\
					\cline{2-3}
					& $\aaa^3\ccc$& \multirow{2}{*}{$(\frac12-\frac2f,1-\frac4f,\frac12+\frac6f,\ddd,1-\ddd+\frac4f)$} & &  & \\
					\cline{2-2}
					& $\aaa^2\ddd\eee$& & & & \\
					\cline{2-3}
					& $\aaa^4\ccc$& \multirow{2}{*}{$(\frac13-\frac4{3f},1-\frac4f,\frac23+\frac{16}{3f},\ddd,1-\ddd+\frac4f)$} & & &  \\
					\cline{2-2}
					& $\aaa^3\ddd\eee$& & & &\\
					\cline{2-3}
					& $\aaa^3\bbb^2$& $(\frac8{3f},1-\frac4f,1+\frac4{3f},\ddd,1-\ddd+\frac4f)$ & & & \\
					\hline
			\end{tabular} }

		\end{table}

		\begin{table}[htp]
			\centering
			\caption{AVC for Cases $\{\bbb\ddd\eee,\aaa\bbb\ccc,\aaa^3\}$, $\{\bbb\ddd\eee,\aaa\bbb\ccc,\aaa^3\ccc(\aaa^2\ddd\eee)\}$, $\{\bbb\ddd\eee,\aaa\bbb\ccc,\aaa^4\ccc(\aaa^3\ddd\eee)\}$ and $\{\bbb\ddd\eee,\aaa\bbb\ccc,\aaa^3\bbb^2\}$.}\label{Tab-7.1}
			\begin{minipage}{1\textwidth}
				\resizebox{\linewidth}{!}{\begin{tabular}{c|c||c|c||c|c||c|c}
						\hline
						\multicolumn{2}{c||}{AVC for $\{\bbb\ddd\eee,\aaa\bbb\ccc,\aaa^3\}$}&\multicolumn{2}{c||}{AVC for $\{\bbb\ddd\eee,\aaa\bbb\ccc,\aaa^3\ccc(\aaa^2\ddd\eee)\}$}&\multicolumn{2}{c||}{AVC for $\{\bbb\ddd\eee,\aaa\bbb\ccc,\aaa^4\ccc(\aaa^3\ddd\eee)\}$}&\multicolumn{2}{c}{AVC for $\{\bbb\ddd\eee,\aaa\bbb\ccc,\aaa^3\bbb^2\}$}\\
						\hline
						$f$&vertex&$f$&vertex&$f$&vertex&$f$&vertex\\
						\hline
						\hline
						all&$\bbb\ddd\eee,\aaa\bbb\ccc,\aaa^3$&all&$\bbb\ddd\eee,\aaa\bbb\ccc,\aaa^3\ccc,\aaa^2\ddd\eee$&all&$\bbb\ddd\eee,\aaa\bbb\ccc,\aaa^4\ccc,\aaa^3\ddd\eee$&all&$\bbb\ddd\eee,\aaa\bbb\ccc,\aaa^3\bbb^2$\\
						$36$&$\ccc^2\ddd\eee,\aaa\ccc^3$&$36$&$\ccc^3$&$16$&$\aaa^2\bbb^2,\aaa^5\bbb,\aaa^8$&$8s,s=2,3,\cdots$&$\aaa^{6s}$\\
						&&&&&&$8s+4,s=2,3,\cdots$&$\aaa^{6s+3},\aaa^{3s+3}\bbb,\aaa^{3s+1}\ccc,\aaa^{3s}\ddd\eee$\\
						\hline
						\multicolumn{2}{c||}{By the balance lemma, $\ccc^2\ddd\eee$ must appear.} & \multicolumn{2}{c||}{By the balance lemma, $\aaa^2\ddd\eee$ must appear. }&\multicolumn{4}{c}{Can not balance.}\\
						\multicolumn{2}{c||}{Its AAD gives $\aaa\thin\eee$, $\bbb\thin\bbb$ or $\eee\thin\eee$,} & \multicolumn{2}{c||}{Its AAD gives $\bbb\thin\bbb$, $\ccc\thin\ddd$ or $\ccc\thin\eee$,}&\multicolumn{4}{c}{}\\
						\multicolumn{2}{c||}{contradicting the AVC} & \multicolumn{2}{c||}{contradicting the AVC}&\multicolumn{4}{c}{}\\
						\hline
				\end{tabular}}
			\end{minipage}
		\end{table}

	\end{proof}

	\begin{proposition}\label{bde,ac2}
		There is no non-symmetric $a^4b$-tiling with general angles, such that $\bbb\ddd\eee$ is the unique $a^2b$-vertex type and $\aaa\ccc^2$ is a vertex.
	\end{proposition}

	\begin{proof}
		By Lemma \ref{bde H}, one of the vertices in the last column of Table \ref{10}: $\{\aaa^3$, $\aaa^3\ccc$, $\aaa^2\bbb\ccc$, $\aaa^2\ddd\eee$, $\aaa\ccc\ddd^2$, $\aaa\ccc\ddd\eee$, $\aaa\ccc\eee^2$, $\ddd^2\eee^2$, $\ddd\eee^3$, $\aaa^4\ccc$, $\aaa^3\bbb\ccc$, $\aaa^3\ddd\eee$, $\aaa^2\bbb^2\ccc$, $\aaa^2\ccc\ddd^2$, $\aaa^2\ccc\ddd\eee$, $\aaa^2\ccc\eee^2$, $\aaa\bbb\ccc\eee^2$, $\aaa\ddd^3\eee$, $\aaa\ddd^2\eee^2$, $\aaa\ddd\eee^3$, $\ccc\ddd^3\eee$, $\ccc\ddd^2\eee^2$, $\ccc\ddd\eee^3\}$ must appear.

		For each combination of $3$ vertices, we will apply the irrational angle lemma and AAD to get Table \ref{bde+ac2}. The proof is completed by Table \ref{bde+ac2}, \ref{Tab-8.1}, \ref{Tab-8.2}, \ref{Tab-8.4} and \ref{Tab-8.5}, similar to Proposition \ref{ad2,ce2}.

		\begin{table}[htp]
			\centering
			\caption{$\bbb\ddd\eee,\aaa\ccc^2$ are vertices.}\label{bde+ac2}
			\resizebox{\linewidth}{!}{\begin{tabular}{|c|c|c|c|c|c|}
					\hline
					\multicolumn{2}{|c|}{Vertex}& Angles & Irrational Angle Lemma & AAD & Contradiction\\
					\hline\hline
					\multirow{23}{*}{$\bbb\ddd\eee,\aaa\ccc^2$} & $\aaa\ccc\eee^2$& $(\frac8f,\frac32-\ddd+\frac2f,1-\frac4f,\ddd,\frac12-\frac2f)$ & \multirow{2}{*}{$n_2=n_4$}& \multirow{6}{*}{$\thin^{\ccc}\eee^{\ddd}\thick^{\ddd}\eee^{\ccc}\thin\cdots$} & \multirow{6}{*}{No $\ddd^2\cdots$ by the irrational angle lemma} \\
					\cline{2-3}
					& $\aaa^2\ccc\eee^2$&$(\frac8f,\frac32-\ddd+\frac6f,1-\frac4f,\ddd,\frac12-\frac6f)$& &  & \\
					\cline{2-4}
					& $\ddd\eee^3$& $(\frac8f,\frac43-\frac{2}3\ddd,1-\frac4f,\ddd,\frac23-\frac{1}3\ddd)$ & \multirow{3}{*}{$2n_2+n_5=3n_4$}& & \\
					\cline{2-3}
					& $\aaa\ddd\eee^3$& $(\frac8f,\frac43-\frac{2}3\ddd+\frac8{3f},1-\frac4f,\ddd,\frac23-\frac{1}3\ddd-\frac8{3f})$ & & & \\
					\cline{2-3}
					& $\ccc\ddd\eee^3$& $(\frac8f,\frac53-\frac{2}3\ddd-\frac4{3f},1-\frac4f,\ddd,\frac13-\frac{1}3\ddd+\frac4{3f})$ & & & \\
					\cline{2-4}
					& $\aaa\bbb\ccc\eee^2$& $(\frac8f,3-2\ddd+\frac4{f},1-\frac4f,\ddd,-1+\ddd-\frac4{f})$ &$2n_2=n_4+n_5$ & & \\
					\cline{2-6}
					& $\aaa\ccc\ddd^2$& $(\frac8f,\frac32-\eee+\frac2f,1-\frac4f,\frac12-\frac2f,\eee)$ & \multirow{2}{*}{$n_2=n_5$}& \multirow{4}{*}{$\thin^{\bbb}\ddd^{\eee}\thick^{\eee}\ddd^{\bbb}\thin\cdots$} & \multirow{4}{*}{No $\eee^2\cdots$}\\
					\cline{2-3}
					& $\aaa^2\ccc\ddd^2$&$(\frac8f,\frac32-\eee+\frac6f,1-\frac4f,\frac12-\frac6f,\eee)$& & & \\
					\cline{2-4}
					& $\aaa\ddd^3\eee$& $(\frac8f,\frac43-\frac{2}3\eee+\frac8{3f},1-\frac4f,\frac23-\frac{1}3\eee-\frac8{3f},\eee)$ &\multirow{2}{*}{$2n_2+n_4=3n_5$} & & \\
					\cline{2-3}
					& $\ccc\ddd^3\eee$& $(\frac8f,1+2\ddd-\frac4f,1-\frac4f,\ddd,1-3\ddd+\frac4f)$ & & & \\
					\cline{2-6}
					& $\aaa\ccc\ddd\eee$& \multirow{2}{*}{$(\frac8f,1+\frac{4}f,1-\frac4f,\ddd,1-\ddd-\frac{4}f)$} & \multirow{10}{*}{$n_4=n_5$}&  & \multirow{4}{*}{See Table \ref{Tab-8.1}} \\
					\cline{2-2}
					& $\aaa\ddd^2\eee^2$& &&  &  \\
					\cline{2-3}
					& $\ddd^2\eee^2$& $(\frac8f,1,1-\frac4f,\ddd,1-\ddd)$ && &\\
					\cline{2-3}
					& $\aaa^2\bbb\ccc$& $(\frac8f,1-\frac{12}f,1-\frac4f,\ddd,1-\ddd+\frac{12}f)$ & & & \\
					\cline{6-6}
					\cline{2-3}
					& $\aaa^3\bbb\ccc$& $(\frac8f,1-\frac{20}f,1-\frac4f,\ddd,1-\ddd+\frac{20}f)$ & & & \multirow{2}{*}{See Table \ref{Tab-8.2}}\\
					\cline{2-3}
					& $\aaa^2\bbb^2\ccc$& $(\frac8f,\frac12-\frac6f,1-\frac4f,\ddd,\frac32-\ddd+\frac6f)$ & & &\\
					\cline{6-6}
					\cline{2-3}
					& $\aaa^2\ddd\eee$& $(\frac8f,\frac{16}f,1-\frac4f,\ddd,2-\ddd-\frac{16}f)$ & & & \multirow{2}{*}{See Table \ref{Tab-8.4}}\\
					\cline{2-3}
					& $\aaa^3\ddd\eee$& $(\frac8f,\frac{24}f,1-\frac4f,\ddd,2-\ddd-\frac{24}f)$ & & &\\
					\cline{6-6}
					\cline{2-3}
					& $\aaa^2\ccc\ddd\eee$& $(\frac8f,1+\frac{12}f,1-\frac4f,\ddd,1-\ddd-\frac{12}f)$ & & & \multirow{2}{*}{See Table \ref{Tab-8.5}}\\
					\cline{2-3}
					& $\ccc\ddd^2\eee^2$& $(\frac8f,\frac32-\frac2f,1-\frac4f,\ddd,\frac12-\ddd+\frac2f)$ & & &\\
					\cline{6-6}
					\cline{2-4}
					\cline{6-6}
					& $\aaa^3$& $(\frac23,2-\ddd-\eee,\frac23,\ddd,\eee),f=12$ & \multirow{3}{*}{} & & \multirow{3}{*}{Similar to Proposition \ref{bde as the unique $a^2b$-vertex} Case \{$\bbb\ddd\eee,\aaa^2\ccc,\aaa\ccc^3$\}}\\
					\cline{2-3}
					& $\aaa^3\ccc$& $(\frac25,2-\ddd-\eee,\frac45,\ddd,\eee),f=20$ &  & & \\
					\cline{2-3}
					& $\aaa^4\ccc$& $(\frac27,2-\ddd-\eee,\frac67,\ddd,\eee),f=28$ &  & & \\
					\hline
			\end{tabular} }

		\end{table}

		\begin{table}[htp]
			\centering
			\caption{AVC for Cases $\{\bbb\ddd\eee,\aaa\ccc^2,\aaa\ccc\ddd\eee(\aaa\ddd^2\eee^2)\}$, $\{\bbb\ddd\eee,\aaa\ccc^2,\ddd^2\eee^2\}$ and $\{\bbb\ddd\eee,\aaa\ccc^2,\aaa^2\bbb\ccc\}$.}\label{Tab-8.1}
			\begin{minipage}{1\textwidth}
				\resizebox{\linewidth}{!}{\begin{tabular}{c|c||c|c||c|c}
						\hline
						\multicolumn{2}{c||}{AVC for $\{\bbb\ddd\eee,\aaa\ccc^2,\aaa\ccc\ddd\eee(\aaa\ddd^2\eee^2)\}$}&\multicolumn{2}{c||}{AVC for $\{\bbb\ddd\eee,\aaa\ccc^2,\ddd^2\eee^2\}$}&\multicolumn{2}{c}{AVC for $\{\bbb\ddd\eee,\aaa\ccc^2,\aaa^2\bbb\ccc\}$}\\
						\hline
						$f$&vertex&$f$&vertex&$f$&vertex\\
						\hline
						\hline
						all&$\bbb\ddd\eee,\aaa\ccc^2,\aaa\ccc\ddd\eee,\aaa\ddd^2\eee^2$ &all&$\bbb\ddd\eee,\aaa\ccc^2,\ddd^2\eee^2$&	all&$\bbb\ddd\eee,\aaa\ccc^2,\aaa^2\bbb\ccc,\aaa^3\bbb^2$\\
						$8s+4,s=2,3,\cdots$&$\aaa^{2s+1},\aaa^{s}\bbb,\aaa^{s+1}\ccc,\aaa^{s+1}\ddd\eee$&$8s,s=2,3,\cdots$&$\aaa^{2s},\aaa^{s}\bbb,\aaa^{s}\ddd\eee$&$28$&$\bbb^2\ccc,\aaa^2\ddd\eee,\aaa\bbb^3,\aaa^4\ccc,\aaa^5\bbb,\aaa^7$\\
						&&&&$36$&$\bbb^3,\aaa^3\ddd\eee,\aaa^5\ccc,\aaa^6\bbb,\aaa^9$\\
						&&&&$8s+4,s=5,6,\cdots$&$\aaa^{2s+1},\aaa^{s+2}\bbb,\aaa^{s+1}\ccc,\aaa^{s-1}\ddd\eee$\\
						\hline
						\multicolumn{4}{c||}{By the balance lemma, $\aaa^{s}\bbb$ must appear. Its AAD gives $\bbb\thin\bbb$ or $\bbb\thin\ccc$, contradicting the AVC} &\multicolumn{2}{c}{By the balance lemma, $\aaa^k\ddd\eee$ must appear. }\\
						\multicolumn{4}{c||}{}&\multicolumn{2}{c}{Its AAD gives $\ccc\thin\ddd$ or $\ccc\thin\eee$, contradicting the AVC}\\
						\hline
				\end{tabular}}
			\end{minipage}
		\end{table}

		\begin{table}[htp]
			\centering
			\caption{AVC for Cases $\{\bbb\ddd\eee,\aaa\ccc^2,\aaa^3\bbb\ccc\}$ and $\{\bbb\ddd\eee,\aaa\ccc^2,\aaa^2\bbb^2\ccc\}$.}\label{Tab-8.2}
			\begin{minipage}{1\textwidth}
				\resizebox{\linewidth}{!}{\begin{tabular}{c|c||c|c}
						\hline
						\multicolumn{2}{c||}{AVC for $\{\bbb\ddd\eee,\aaa\ccc^2,\aaa^3\bbb\ccc\}$} &\multicolumn{2}{c}{AVC for $\{\bbb\ddd\eee,\aaa\ccc^2,\aaa^2\bbb^2\ccc\}$}\\
						\hline
						$f$&vertex&$f$&vertex\\
						\hline
						\hline
						all&$\bbb\ddd\eee,\aaa\ccc^2,\aaa^3\bbb\ccc,\aaa^5\bbb^2$&all&$\bbb\ddd\eee,\aaa\ccc^2,\aaa^2\bbb^2\ccc,\aaa^3\bbb^4$ \\
						$36$&$\aaa^2\ddd\eee,\aaa\bbb^2\ccc,\aaa\bbb^4,\aaa^5\ccc,\aaa^3\bbb^3,\aaa^7\bbb,\aaa^9$&$44$&$\aaa^2\ddd\eee,\bbb^3\ccc,\aaa^4\bbb\ccc,\aaa\bbb^5,\aaa^6\ccc,\aaa^5\bbb^3,\aaa^7\bbb^2,\aaa^9\bbb,\aaa^{11}$\\
						$44$&$\bbb^2\ccc,\aaa^3\ddd\eee,\aaa^2\bbb^3,\aaa^6\ccc,\aaa^8\bbb,\aaa^{11}$ &$60$&$\aaa^3\ddd\eee,\bbb^5,\aaa^5\bbb\ccc,\aaa^8\ccc,\aaa^6\bbb^3,\aaa^9\bbb^2,\aaa^{12}\bbb,\aaa^{15}$ \\
						$52$&$\aaa\bbb^3,\aaa^4\ddd\eee,\aaa^7\ccc,\aaa^9\bbb,\aaa^{13},$&$16s+12,s=4,5,\cdots$&$\aaa^{4s+3},\aaa^{3s+3}\bbb,\aaa^{2s+3}\bbb^2,\aaa^{s+3}\bbb^3,\aaa^{s+2}\bbb\ccc,$ \\
						$60$&$\bbb^3,\aaa^5\ddd\eee,\aaa^8\ccc,\aaa^{10}\bbb,\aaa^{15}$&&$\aaa^{2s+2}\ccc,\aaa^s\ddd\eee$\\
						$8s+4,s=8,9,\cdots$&$\aaa^{2s+1},\aaa^{s+3}\bbb,\aaa^{s+1}\ccc,\aaa^{s-2}\ddd\eee$& &\\
						\hline
						\multicolumn{4}{c}{By the balance lemma, $\aaa^k\ddd\eee$ must appear. Its AAD gives $\ccc\thin\ddd$ or $\ccc\thin\eee$, contradicting the AVC} \\
						\hline
				\end{tabular}}
			\end{minipage}
		\end{table}

		\begin{table}[htp]
			\centering
			\caption{AVC for Cases $\{\bbb\ddd\eee,\aaa\ccc^2,\aaa^2\ddd\eee\}$ and $\{\bbb\ddd\eee,\aaa\ccc^2,\aaa^3\ddd\eee\}$.}\label{Tab-8.4}
			\begin{minipage}{1\textwidth}
				\resizebox{\linewidth}{!}{\begin{tabular}{c|c||c|c}
						\hline
						\multicolumn{2}{c||}{AVC for $\{\bbb\ddd\eee,\aaa\ccc^2,\aaa^2\ddd\eee\}$}&\multicolumn{2}{c}{AVC for $\{\bbb\ddd\eee,\aaa\ccc^2,\aaa^3\ddd\eee\}$}\\
						\hline
						$f$&vertex&$f$&vertex\\
						\hline
						\hline
						all&$\bbb\ddd\eee,\aaa\ccc^2,\aaa^2\ddd\eee$&all&$\bbb\ddd\eee,\aaa\ccc^2,\aaa^3\ddd\eee$\\
						$16$&$\aaa^2\bbb,\aaa^4,\ddd^2\eee^2$ &$20$&$\aaa^2\bbb,\aaa^3\ccc,\aaa\ccc\ddd\eee,\aaa^5,\aaa\ddd^2\eee^2$\\
						$8s,s=3,4,\cdots$&$\aaa^{2s-2x}\bbb^{x}(x=0,1,\cdots,s)$ &$24$&$\aaa^3\bbb,\ddd^2\eee^2,\aaa^6$\\
						$8s+4,s=3,4,\cdots$&$\aaa^{2s+1-2x}\bbb^{x}(x=0,1,\cdots,s),$&$8s,s=4,5,\cdots$&$\aaa^{2s-3x}\bbb^{x}(x=0,1,\cdots,\lfloor\frac{2s}{3}\rfloor))$\\
						&$\aaa^{s+1-2y}\bbb^y\ccc(y=0,1,\cdots,\lfloor\frac{s-1}2\rfloor)$&$8s+4,s=3,4,\cdots$&$\aaa^{2s+1-3x}\bbb^{x}(x=0,1,\cdots,\lfloor\frac{2s+1}{3}\rfloor)),$\\
						&&&$\aaa^{s+1-3y}\bbb^y\ccc(y=0,1,\cdots,\lfloor\frac{s-1}{3}\rfloor)$\\
						\hline
						\multicolumn{4}{c}{By the balance lemma, $\aaa^k\bbb^l$ or $\aaa^k\bbb^l\ccc$ must appear. Its AAD gives $\bbb\thin\bbb$, $\bbb\thin\ccc$, $\ccc\thin\ddd$, $\ccc\thin\eee$, $\ddd\thin\ddd$, $\ddd\thin\eee$ or $\eee\thin\eee$, contradicting the AVC} \\
						\hline
				\end{tabular}}
			\end{minipage}
		\end{table}

		\begin{table}[htp]
			\centering
			\caption{AVC for Cases $\{\bbb\ddd\eee,\aaa\ccc^2,\aaa^2\ccc\ddd\eee\}$ and $\{\bbb\ddd\eee,\aaa\ccc^2,\ccc\ddd^2\eee^2\}$.}\label{Tab-8.5}
			\begin{minipage}{1\textwidth}
				\resizebox{\linewidth}{!}{\begin{tabular}{c|c||c|c}
						\hline
						\multicolumn{2}{c||}{AVC for $\{\bbb\ddd\eee,\aaa\ccc^2,\aaa^2\ccc\ddd\eee\}$}&\multicolumn{2}{c}{AVC for $\{\bbb\ddd\eee,\aaa\ccc^2,\ccc\ddd^2\eee^2\}$}\\
						\hline
						$f$&vertex&$f$&vertex\\
						\hline
						\hline
						all&$\bbb\ddd\eee,\aaa\ccc^2,\aaa^2\ccc\ddd\eee,\aaa^3\ddd^2\eee^2$&all&$\bbb\ddd\eee,\aaa\ccc^2,\ccc\ddd^2\eee^2$\\
						$28$&$\aaa^2\bbb,\aaa^4\ccc,\ccc\ddd^2\eee^2,\aaa^7,\aaa^5\ddd\eee,\aaa\ddd^3\eee^3$&$16s+12,s=1,2,\cdots$&$\aaa^{4s+3},\aaa^{s+1}\bbb,\aaa^{2s+2}\ccc,\aaa^{s+1}\ccc\ddd\eee,$ \\
						$36$&$\aaa^3\bbb,\aaa^5\ccc,\ddd^3\eee^3,\aaa^6\ddd\eee,\aaa^9$&&$\aaa^{3s+2}\ddd\eee,\aaa^{2s+1}\ddd^2\eee^2,\aaa^{s}\ddd^3\eee^3$ \\
						$8s+4,s=5,6,\cdots$&$\aaa^{2s+1},\aaa^{s-1}\bbb,\aaa^{s+1}\ccc,\aaa^{s+2}\ddd\eee$&&\\
						\hline
						\multicolumn{4}{c}{By the balance lemma, $\aaa^k\bbb$ must appear. Its AAD gives $\bbb\thin\bbb$ or $\bbb\thin\ccc$, contradicting the AVC} \\
						\hline
				\end{tabular}}
			\end{minipage}
		\end{table}
	\end{proof}

	\begin{proposition}\label{bde,a3,ab2}
		There is no non-symmetric $a^4b$-tiling with general angles, such that $\bbb\ddd\eee$ is the unique $a^2b$-vertex type and $\aaa^3$ is a vertex.
	\end{proposition}

	\begin{proof}
	 By Lemma \ref{bde H}, one of the vertices in the last column of Table \ref{10}: $\{\aaa\bbb^2,\aaa\bbb\ccc,\aaa\ccc^2\}$ must appear. By Proposition \ref{bde,abc}, \ref{bde,ac2},  we may assume that $\aaa\bbb\ccc$, $\aaa\ccc^2$ are not vertices. So $\aaa\bbb^2$ is a vertex.	By the irrational angle lemma, we have $n_4=n_5$. The angle sums of $\bbb\ddd\eee, \aaa^3,\aaa\bbb^2$ and the angle sum for the pentagon imply
		\begin{align*}
			\aaa=\tfrac23,\;
			\bbb=\tfrac23,\;
			\ccc=\tfrac13+\tfrac4f,\;
			\ddd+\eee=\tfrac43.
		\end{align*}
		By $\bbb\neq\ccc$, we get $f\ge14$. By $\ccc<R(\ddd\eee\cdots)=\aaa=\bbb<2\ccc,\ddd+\eee$ and $n_4=n_5$, we have $\ddd\cdots=\ddd\eee\cdots=\bbb\ddd\eee$ since $\bbb\ddd\eee$ is the unique $a^2b$-vertex type. Then $\#\bbb\ge\#\bbb\ddd\eee+\#\aaa\bbb^2=f+\#\aaa\bbb^2>f$, a contradiction.
	\end{proof}

	\section{Tilings with a unique $a^2b$-vertex type $\aaa\ddd\eee$}
	\label{ade as the unique $a^2b$-vertex}

	This section is devoted to the classification of edge-to-edge tilings of the sphere by congruent pentagons with the edge combination $a^4b$, where $\aaa\ddd\eee$ is the unique $a^2b$-vertex type. By Lemma \ref{l2}, we know that there exists an $a^3$-vertex around the special tile. This means one of $\{\aaa^3$, $\bbb^3$, $\ccc^3$, $\aaa^2\bbb$, $\aaa^2\ccc$, $\aaa\bbb^2$, $\aaa\ccc^2$, $\bbb^2\ccc$, $\bbb\ccc^2$, $\aaa\bbb\ccc\}$ must appear. Up to the symmetry of interchanging $\bbb\leftrightarrow\ccc$ and $\ddd\leftrightarrow\eee$, it is enough to consider the following cases $\{\aaa^3,\aaa^2\ccc,\aaa\ccc^2,\aaa\bbb\ccc,\bbb^{2}\ccc,\ccc^3\}$.

	\begin{proposition}\label{ade,2bc}
		In a non-symmetric $a^4b$-tiling with general angles, if $\aaa\ddd\eee$ is the unique $a^2b$-vertex type and $\bbb^{2}\ccc$ is a vertex, then the tiling is the non-symmetric $3$-layer earth map tiling $T(4m\,\aaa\ddd\eee,2m\,\bbb^2\ccc,2\ccc^m)$ with $4m$ tiles for some $m\ge4$,  or its standard flip modification $T((8k+4)\aaa\ddd\eee,4k\,\bbb^2\ccc,4\bbb\ccc^{k+1})$ for odd $m=2k+1$.
	\end{proposition}

	\begin{proof}

    The angle sums of $\aaa\ddd\eee$, $\bbb^2\ccc$ and the angle sum for the pentagon imply
    \begin{align*}
    	\aaa+\ddd+\eee=2,\;
    	\bbb=1-\tfrac{4}{f},\;
    	\ccc=\tfrac8f.
    \end{align*}
	We have $1>\beta>\frac23>\gamma$. By Lemma \ref{geometry1}, this implies $\delta<\epsilon$. We divide the discussion into three cases.

	\subsection*{Case. $\aaa\ge1$}
	Then $\ddd,\eee<1$. By Lemma \ref{geometry2}, we get $\bbb+2\ddd>1$ and $\ccc+2\eee>1$. By $R(\aaa\cdots)=\ddd+\eee<\bbb+\ccc,\bbb+2\ddd,\ccc+2\eee$ and the parity lemma, we get $\aaa\cdots=\aaa\ddd\eee,\aaa\ccc^x(x\ge2)$ or $\aaa\ccc^y\ddd^z(z\ge2)$.

	If $\aaa\ccc^x$ is a vertex, by the irrational angle lemma and $\aaa\ddd\eee$, we have $n_4=n_5$ and $\aaa\ccc^y\ddd^z$ cannot be a vertex.  Then we have the unique AAD $\aaa\ccc\cdots=\aaa\ccc^{x}=\thin^{\ccc}\aaa^{\bbb}\thin^{\eee}\ccc^{\aaa}\thin\cdots\thin^{\eee}\ccc^{\aaa}\thin$ since $\aaa^2\cdots$ and $\aaa\bbb\cdots$ are not vertices. This gives a vertex $\bbb^{\aaa}\thin^{\ccc}\eee\cdots$. By $\bbb+2\ddd>1$ and $\ccc+2\eee>1$, we get $R(\bbb\ddd\eee\cdots)<\ccc+\ddd+\eee$. By $R(\bbb^2\cdots)=\ccc<x\ccc=\ddd+\eee$, we get $\bbb^2\cdots=\bbb^2\ccc$.   Then $\bbb\eee\cdots=\bbb\ddd\eee\cdots=\bbb\ccc^y\ddd\eee$ or 	$\bbb\ddd^2\eee^2$. Both possibilities will lead to some contradiction in the following.

	By $\aaa\thin\ccc\cdots=\aaa^{\bbb}\thin^{\eee}\ccc\cdots$ or $\aaa^{\ccc}\thin^{\aaa}\ccc\cdots$, we have the AAD  $\bbb\thin\eee\cdots=\bbb\ddd^2\eee^2=\thin^{\ddd}\bbb^{\aaa}\thin^{\ccc}\eee^{\ddd}\thick^{\eee}\ddd^{\bbb}\thin^{\bbb}\ddd^{\eee}\thick^{\ddd}\eee^{\ccc}\thin$, $\thin^{\ddd}\bbb^{\aaa}\thin^{\ccc}\eee^{\ddd}\thick^{\eee}\ddd^{\bbb}\thin^{\ccc}\eee^{\ddd}\thick^{\eee}\ddd^{\bbb}\thin$ or $\thin^{\ddd}\bbb^{\aaa}\thin^{\ccc}\eee^{\ddd}\thick^{\ddd}\eee^{\ccc}\thin^{\bbb}\ddd^{\eee}\thick^{\eee}\ddd^{\bbb}\thin$. This gives a vertex $\bbb^{\ddd}\thin^{\ddd}\bbb\cdots$ or $\thin^{\aaa}\bbb^{\ddd}\thin^{\eee}\ccc^{\aaa}\thin\cdots$.
	Then $\bbb^{\ddd}\thin^{\ddd}\bbb\cdots=\bbb^2\ccc=\thin^{\aaa}\bbb^{\ddd}\thin^{\ddd}\bbb^{\aaa}\thin^{\aaa}\ccc^{\eee}\thin$. This gives a vertex $\aaa^2\cdots$, a contradiction. Then $\thin^{\aaa}\bbb^{\ddd}\thin^{\eee}\ccc^{\aaa}\thin\cdots=\bbb^2\ccc$, $\bbb\ccc^{\frac{f+4}{8}}$ or $\bbb\ccc^y\ddd\eee$. Either one gives a vertex $\aaa^2\cdots$ or  $\aaa\bbb\cdots$, a contradiction. Therefore, $\bbb\ddd^2\eee^2$ is not a vertex.

	We have the unique AAD
	$\bbb\thin\eee\cdots=\bbb\ccc^y\ddd\eee=\thin^{\ddd}\bbb^{\aaa}\thin^{\ccc}\eee^{\ddd}\thick^{\eee}\ddd^{\bbb}\thin^{\eee}\ccc^{\aaa}\thin\cdots\thin^{\eee}\ccc^{\aaa}\thin$. This gives a new $\thin^{\aaa}\bbb^{\ddd}\thin^{\ccc}\eee^{\ddd}\thick\cdots=\bbb\ccc^y\ddd\eee=\thin^{\aaa}\bbb^{\ddd}\thin^{\ccc}\eee^{\ddd}\thick^{\eee}\ddd^{\bbb}\thin\ccc\thin\cdots$. This gives a vertex $\aaa^2\cdots$ or  $\aaa\bbb\cdots$, a contradiction.  Therefore, $\bbb\ccc^y\ddd\eee$ is not a vertex.

	In summary, $\aaa\ccc^x$ is not a vertex and $\aaa\cdots=\aaa\ddd\eee$ or $\aaa\ccc^y\ddd^z(z\ge2)$. To balance $\aaa$ and $\ddd$, we have $\aaa\cdots=\aaa\ddd\eee$. By the angle values, we have $\text{AVC}\sub\{\aaa\ddd\eee,\bbb^2\ccc,\bbb\ccc^{\frac{f+4}{8}},\ccc^{\frac{f}4}\}$.

    \subsection*{Case. $\aaa<1$ and $\eee\ge1$}

    By $R(\eee\cdots)=\aaa+\ddd<1<2\bbb,\bbb+\ccc$, the parity lemma and the balance lemma, we get $\eee\cdots=\ddd\cdots=\aaa\ddd\eee$ or $\ccc^x\ddd\eee(x\ge2)$. If $\ccc^x\ddd\eee$ is a vertex, we have the unique AAD $\ccc^{x}\ddd\eee=\thin^{\ccc}\eee^{\ddd}\thick^{\eee}\ddd^{\bbb}\thin^{\aaa}\ccc^{\eee}\thin\cdots\thin^{\aaa}\ccc^{\eee}\thin$ since $\bbb\eee\cdots$ and $\eee^2\cdots$ are not vertices. This gives a vertex $\thin^{\bbb}\aaa^{\ccc}\thin^{\ddd}\bbb^{\aaa}\thin\cdots$.
    By $\aaa=x\ccc\ge2\ccc$ and  $\ccc\thin\eee\cdots=\ccc^x\ddd\eee=\ccc^{\eee}\thin^{\ccc}\eee\cdots$, we have the AAD $\thin^{\bbb}\aaa^{\ccc}\thin^{\ddd}\bbb^{\aaa}\thin\cdots=\aaa^{y}\bbb\ccc^{z}(z\ge0)=\thin\aaa\thin\cdots\thin^{\bbb}\aaa^{\ccc}\thin^{\ddd}\bbb^{\aaa}\thin^{\eee}\ccc^{\aaa}\thin\cdots\thin^{\eee}\ccc^{\aaa}\thin$.

    By $R(\bbb^2\cdots)=\ccc$ and $\aaa=x\ccc\ge2\ccc$, we get $\bbb^2\cdots=\bbb^2\ccc$. Then we get the unique AAD $\bbb^2\ccc=\thin^{\aaa}\bbb^{\ddd}\thin^{\aaa}\bbb^{\ddd}\thin^{\aaa}\ccc^{\eee}\thin$ since $\ddd^2\cdots$ and $\ddd\thin\eee\cdots$ are not vertices.  Therefore there is no $\aaa^{\bbb}\thin^{\bbb}\aaa\cdots$ and $\bbb\ccc\cdots=\bbb^2\ccc$ or $\aaa^y\bbb\ccc^z(y\ge0)$.

    Then AAD $\aaa^{\bbb}\thin^{\ccc}\aaa\cdots$ implies a vertex $\thin^{\ddd}\bbb^{\aaa}\thin^{\aaa}\ccc^{\eee}\thin\cdots=\aaa^{y}\bbb\ccc^{z}$. This gives a vertex  $\bbb\eee\cdots$ or $\ccc^{\eee}\thin^{\ccc}\aaa\cdots$, or $\eee\thin\eee\cdots$, or $\ddd\thin\eee\cdots$, a contradiction. Therefore there is no $\aaa^{\bbb}\thin^{\ccc}\aaa\cdots$. This implies that there is no $\aaa\thin\aaa\thin\aaa\cdots$ and $\aaa\thin\aaa\cdots=\aaa^{\ccc}\thin^{\ccc}\aaa\cdots$.

    Then we have the unique AAD $\thin^{\bbb}\aaa^{\ccc}\thin^{\ddd}\bbb^{\aaa}\thin\cdots=\aaa\bbb\ccc^{z}=\thin^{\bbb}\aaa^{\ccc}\thin^{\ddd}\bbb^{\aaa}\thin^{\eee}\ccc^{\aaa}\thin\cdots\thin^{\eee}\ccc^{\aaa}\thin$. Similarly, we have the unique AAD $\aaa\thin\bbb\cdots=\aaa^y\bbb\ccc^z=\aaa\thin^{\ddd}\bbb\cdots=\thin^{\bbb}\aaa^{\ccc}\thin^{\ddd}\bbb^{\aaa}\thin^{\eee}\ccc^{\aaa}\thin\cdots\thin^{\eee}\ccc^{\aaa}\thin$. However, the AAD $\thin^{\bbb}\aaa^{\ccc}\thin^{\ddd}\bbb^{\aaa}\thin^{\eee}\ccc^{\aaa}\thin\cdots\thin^{\eee}\ccc^{\aaa}\thin$ implies a new $\thin^{\bbb}\aaa^{\ccc}\thin^{\aaa}\bbb^{\ddd}\thin\cdots$. The contradiction implies $\aaa^{\ccc}\thin^{\ddd}\bbb\cdots$ is not a vertex.

    In summary, $\ccc^x\ddd\eee$ is not a vertex. Then $\eee\cdots=\aaa\ddd\eee$. By the angle values, we have $\text{AVC}\sub\{\aaa\ddd\eee,\bbb^2\ccc,\bbb\ccc^{\frac{f+4}{8}},\ccc^{\frac{f}4}\}$.

    \subsection*{Case. $\aaa<1$ and $\eee<1$}
    Then it is convex and $\aaa+\ddd>1$, $\ddd+\eee>1$. By Lemma \ref{geometry2}, we get $\aaa+2\ccc>1$, $2\aaa+\ccc>1$, $2\aaa+\bbb>1$ and $\bbb+2\ddd>1$. By $R(\eee^2\cdots)<\aaa,\bbb+\ccc,2\ddd$ and the parity lemma, we get $\eee^2\cdots=\ccc^x\eee^2(x\ge2)$. If $\ccc^x\eee^2$ is a vertex, by the irrational angle lemma, we have $n_1=n_4$. Then AAD $\ccc^x\eee^2=\eee^{\ddd}\thick^{\ddd}\eee\cdots$ implies that $\ddd^2\cdots$ is a vertex. By $n_1=n_4$, we have $\ddd^2\cdots=\aaa^2\ddd^2\cdots$, contradicting $2\aaa+2\ddd>2$. Therefore $\ccc^x\eee^2$ is not a vertex. By $R(\ddd\eee\cdots)=\aaa<1<2\bbb,\bbb+\ccc$ and the parity lemma, we get $\eee\cdots=\ddd\eee\cdots=\aaa\ddd\eee$ or $\ccc^y\ddd\eee(y\ge2)$. Then the same proof as Case $\aaa<1$ and $\eee\ge1$ gives  $\text{AVC}\sub\{\aaa\ddd\eee,\bbb^2\ccc,\bbb\ccc^{\frac{f+4}{8}},\ccc^{\frac{f}4}\}$.

    \vspace{5pt}
    In all three cases above, we get $\text{AVC}\sub\{\aaa\ddd\eee,\bbb^2\ccc,\bbb\ccc^{\frac{f+4}{8}},\ccc^{\frac{f}4}\}$. Then we get the non-symmetric $3$-layer earth map tiling and its standard flip modification by exactly the same proof as the ``Subcase. $\ccc^{\frac{f}4}$ is a vertex'' and  ``Subcase. $\bbb\ccc^{\frac{f+4}{8}}$ is a vertex'' on page \pageref{c 1}.

	\end{proof}

	By Proposition \ref{ade,2bc},
	we can assume that $\bbb^2\ccc$ and $\bbb\ccc^2$ are not vertices anymore in the rest of this section.

	\begin{proposition}\label{ade,c3}
		In a non-symmetric $a^4b$-tiling with general angles, if $\aaa\ddd\eee$ is the unique $a^2b$-vertex type and $\ccc^{3}$ is a vertex, then all tilings are
		\begin{enumerate}
			\item Pentagonal subdivisions of the octahedron $T(24\aaa\ddd\eee,8\ccc^3,6\bbb^4)$, with $24$ tiles.
			\item Pentagonal subdivisions of the icosahedron $T(60\aaa\ddd\eee,20\ccc^3,12\bbb^5)$, with $60$ tiles.
		\end{enumerate}
	\end{proposition}

	\begin{proof}
		Note that the pentagonal subdivision of the tetrahedron is not included, because it is symmetric. We will first apply Lemma \ref{hdeg} to get a new vertex $\bbb\cdots$ (excluding $\bbb^2\ccc$, $\bbb\ccc^2$), then apply the irrational angle lemma to control the AVC, and finally try to get all tilings by the AAD. We summarize such computations and deductions in Table \ref{table 28}, excluding the symmetric ($\bbb=\ccc$) case $\{\aaa\ddd\eee,\ccc^3,\bbb^3\}$. The proof is completed by Table \ref{table 28} and \ref{Tab-9.1}, similar to Proposition \ref{ad2,ce2}.


		\begin{table*}[htp]
			\centering
			\caption{$\aaa\ddd\eee,\ccc^3$ are vertices.}\label{table 28}
			~\\
			\resizebox{\linewidth}{!}{\begin{tabular}{|c|c|c|c|c|}
					\hline
					\multicolumn{2}{|c|}{Vertex} & Angles & Irrational Angle Lemma  & Contradiction \\
					\hline\hline
					\multirow{7}{*}{$\aaa\ddd\eee, \ccc^{3}$}&$\aaa^2\bbb$&$(\frac56-\frac2f,\frac13+\frac4f,\frac23,\ddd,\frac76-\ddd+\frac2f)$&\multirow{4}{*}{$n_4=n_5$}&\multirow{4}{*}{See Table \ref{Tab-9.1}}\\
					\cline{2-3}
					&$\aaa\bbb^{2}$&$(\frac43-\frac8f,\frac13+\frac4f,\frac23,\ddd,\frac23-\ddd+\frac8f)$&&\\
					\cline{2-3}
					&$\aaa\bbb\ccc$&$(1-\frac4f,\frac13+\frac4f,\frac23,\ddd,1-\ddd+\frac4f)$&&\\
					\cline{2-3}
					&$\aaa\bbb^{3}$&$(1-\frac{12}f,\frac13+\frac4f,\frac23,\ddd,1-\ddd+\frac{12}f)$
					&&\\
					\cline{2-5}
					&$\bbb^3\ccc$&$(2-\ddd-\eee,\frac49,\frac23,\ddd,\eee),f=36$&\multirow{3}{*}{} &  \multirow{1}{*}{Similar to Proposition \ref{bde as the unique $a^2b$-vertex} Case \{$\bbb\ddd\eee,\aaa^2\ccc,\aaa\ccc^3$\}}\\
					\cline{2-3}
					\cline{5-5}
					&$\bbb^4$&$(2-\ddd-\eee,\frac12,\frac23,\ddd,\eee),f=24$ & & \multirow{2}{*}{Admit the unique pentagonal subdivision tiling}\\
					\cline{2-3}
					&$\bbb^5$&$(2-\ddd-\eee,\frac25,\frac23,\ddd,\eee),f=60$&&\\
					\hline
			\end{tabular}}
		\end{table*}

		\begin{table}[htp]
			\centering
			\caption{AVC for Cases $\{\aaa\ddd\eee,\ccc^3,\aaa^2\bbb\}$, $\{\aaa\ddd\eee,\ccc^3,\aaa\bbb^2\}$, $\{\aaa\ddd\eee,\ccc^3,\aaa\bbb\ccc\}$,  and $\{\aaa\ddd\eee,\ccc^3,\aaa\bbb^3\}$.}\label{Tab-9.1}
			\begin{minipage}{1\textwidth}
				\resizebox{\linewidth}{!}{\begin{tabular}{c|c||c|c||c|c||c|c}
						\hline
						\multicolumn{2}{c||}{AVC for $\{\aaa\ddd\eee,\ccc^3,\aaa^2\bbb\}$}&\multicolumn{2}{c||}{AVC for $\{\aaa\ddd\eee,\ccc^3,\aaa\bbb^2\}$}&\multicolumn{2}{c||}{AVC for $\{\aaa\ddd\eee,\ccc^3,\aaa\bbb\ccc\}$}&\multicolumn{2}{c}{AVC for $\{\aaa\ddd\eee,\ccc^3,\aaa\bbb^3\}$}\\
						\hline
						$f$&vertex&$f$&vertex&$f$&vertex&$f$&vertex\\
						\hline
						\hline
						all&$\aaa\ddd\eee,\ccc^3,\aaa^2\bbb$&all&$\aaa\ddd\eee,\ccc^3,\aaa\bbb^2$&all&$\aaa\ddd\eee,\ccc^3,\aaa\bbb\ccc$&all&$\aaa\ddd\eee,\ccc^3,\aaa\bbb^3$\\
						$60$&$\aaa\bbb^3,\bbb^2\ddd\eee,\bbb^5$&$24$&$\bbb^4,\bbb^2\ddd\eee,\ddd^2\eee^2$&$36$&$\bbb^3\ccc,\bbb^2\ddd\eee$&$60$&$\aaa^2\bbb,\bbb^2\ddd\eee,\bbb^5$ \\
						&&$36$&$\bbb^3\ccc,\bbb\ccc\ddd\eee$&&&&\\
						&&$60$&$\bbb^5,\bbb^3\ddd\eee,\bbb\ddd^2\eee^2$&&&&\\
						\hline
						\multicolumn{8}{c}{By the balance lemma $\bbb^k\ddd\eee,\bbb\ccc\ddd\eee$ or $\bbb^k\ddd^2\eee^2$ must appear. Its AAD gives $\aaa\thin\aaa$, $\aaa\thin\ccc$, $\bbb\thin\ccc$, $\ccc\thin\ddd$ or  $\ddd\thin\eee$, contradicting the AVC}\\
				\end{tabular}}
			\end{minipage}
		\end{table}

	\subsection*{Case. $\{\aaa\ddd\eee, \ccc^3, \bbb^4\}$}
		The angle sums of $\aaa\ddd\eee, \ccc^3, \bbb^4$ and the angle sum for the pentagon imply $\bbb=\tfrac12<\ccc=\tfrac23$. By Lemma \ref{geometry1}, we have $\ddd>\eee$.

		\subsubsection*{Subcase. $\ddd\ge1$}
		By $R(\ddd\eee\cdots)=\aaa<2\bbb$, $\bbb+\ccc$, $2\ccc$ and the parity lemma, we get  $\ddd\cdots=\ddd\eee\cdots=\aaa\ddd\eee$. Therefore $\text{AVC}=\{\aaa\ddd\eee, \ccc^3, \bbb^4\}$.

		\subsubsection*{Subcase. $\ddd<1$, $\aaa\ge1$}
		By the angle values, we get $\aaa\cdots=\aaa\ddd\eee$, $\aaa\bbb^2$ or $\aaa\eee^{x}(x\ge4)$. By the previous discussion, we get that $\aaa\bbb^2$ is not a vertex. Then we know $\aaa\cdots=\aaa\ddd\eee$ or $\aaa\eee^{x}(x\ge4)$. By the balance lemma, we have $\aaa\cdots=\aaa\ddd\eee$. Therefore $\text{AVC}=\{\aaa\ddd\eee, \ccc^3, \bbb^4\}$.

		\subsubsection*{Subcase. $\ddd<1$, $\aaa<1$}
		By $\aaa<1$, we get $\ddd+\eee>1$. Then by  $R(\ddd^2\cdots)<\aaa,2\bbb,\bbb+\ccc,2\ccc,2\ddd,\ddd+\eee,2\eee$ and the parity lemma, we get that $\ddd^2\cdots$ is not a vertex. By $R(\ddd\eee\cdots)=\aaa<2\bbb$, $\bbb+\ccc$, $2\ccc$ and the parity lemma, we get $\ddd\cdots=\ddd\eee\cdots=\aaa\ddd\eee$. Therefore $\text{AVC}=\{\aaa\ddd\eee, \ccc^3, \bbb^4\}$.

		In conclusion, we always have $\text{AVC}=\{\aaa\ddd\eee, \ccc^3, \bbb^4\}$.

    \vspace{6pt}

	Similarly, for Case \{$\aaa\ddd\eee, \ccc^3, \bbb^5$\}, we get $\text{AVC}=\{\aaa\ddd\eee, \ccc^3, \bbb^5\}$. Therefore, we
	can easily obtain pentagonal subdivisions of the octahedron and icosahedron (see Figure \ref{poufen} for 3D pictures, see Table \ref{tab-1} for geometric data). Three pentagonal subdivisions and their moduli spaces were described in \cite{lwy}.

	\end{proof}

	\begin{proposition}\label{ade,a2c}
		There is no non-symmetric $a^4b$-tiling with general angles, such that $\aaa\ddd\eee$ is the unique $a^2b$-vertex type and $\aaa^{3}$, $\aaa^{2}\ccc$ or $\aaa\ccc^{2}$ is a vertex.
	\end{proposition}

	\begin{proof}
		For each case, we will first apply Lemma \ref{hdeg} to get a new vertex $\bbb\cdots$ (excluding $\bbb^2\ccc$, $\bbb\ccc^2$), then apply the irrational angle lemma to control the AVC. We summarize such computations and deductions in Table \ref{table ade,a2c}, excluding the symmetric ($\bbb=\ccc$) ones: Case $\{\aaa\ddd\eee,\aaa^2\ccc,\aaa^2\bbb\}$, Case $\{\aaa\ddd\eee,\aaa\ccc^2,\aaa\bbb^2\}$ and $\{\aaa\ddd\eee,\aaa\ccc^2,\aaa\bbb\ccc\}$. The proof is completed by Table \ref{table ade,a2c} and \ref{Tab-10.2}.

		\begin{table*}[htp]
			\centering
			\caption{All possible cases with $\aaa^{3},\aaa^{2}\ccc$ or $\aaa\ccc^{2}$.}\label{table ade,a2c}
			~\\
			\resizebox{\linewidth}{!}{\begin{tabular}{|c|c|c|c|c|}
					\hline
					\multicolumn{2}{|c|}{Vertex} & Angles & Irrational Angle Lemma & Contradiction \\
					\hline\hline
					\multirow{9}{*}{$\aaa\ddd\eee,\aaa^{3}$}& \multirow{2}{*}{$\aaa\bbb\ccc$} &$(\frac23,\bbb,\frac43-\bbb,\ddd,\frac43-\ddd)$, $f=12$
					&&\multirow{13}{*}{No degree $\ge4$ vertex $\ddd\eee\cdots$}\\
					&&and all vertices have degree $3$&&\\
					\cline{2-4}
					&$\aaa^2\bbb$&\multirow{3}{*}{$(\frac23,\frac23,\frac13+\frac4f,\ddd,\frac43-\ddd)$}&\multirow{20}{*}{$n_4=n_5$}& \multirow{13}{*}{(Similar to Proposition \ref{bde,a2b})}\\
					\cline{2-2}
					& $\aaa\bbb^2$&
					& &\\
					\cline{2-2}
					& $\bbb^3$&
					& &\\
					\cline{2-3}
					& $\aaa\bbb^{3}$&$(\frac23,\frac49,\frac59+\frac4f,\ddd,\frac43-\ddd)$
					& &\\
					\cline{2-3}
					&$\bbb^{3}\ccc$&$(\frac23,\frac12-\frac2f,\frac12+\frac6f,\ddd,\frac43-\ddd)$
					& &\\
					\cline{2-3}
					&$\bbb^{4}$&$(\frac23,\frac12,\frac12+\frac4f,\ddd,\frac43-\ddd)$
					& &\\
					\cline{2-3}
					& $\bbb^{5}$&$(\frac23,\frac25,\frac35+\frac4f,\ddd,\frac43-\ddd)$
					&&\\
					\cline{1-3}
					\multirow{7}{*}{$\aaa\ddd\eee,\aaa^{2}\ccc$}&$\aaa\bbb^{2}$&$(\frac45-\frac8{5f},\frac35+\frac4{5f},\frac25+\frac{16}{5f},\ddd,\frac65-\ddd+\frac8{5f})$
					& &\\
					\cline{2-3}
					&$\aaa\bbb^{3}$&$(\frac57-\frac{12}{7f},\frac37+\frac4{7f},\frac47+\frac{24}{7f},\ddd,\frac97-\ddd+\frac{12}{7f})$
					& &\\
					\cline{2-3}
					&$\bbb^{3}\ccc$&$(\frac34-\frac3f,\frac12-\frac2f,\frac12+\frac6f,\ddd,\frac54-\ddd+\frac3f)$
					& &\\
					\cline{2-3}
					& $\bbb^{4}$&$(\frac34-\frac2f,\frac12,\frac12+\frac4f,\ddd,\frac54-\ddd+\frac2f)$
					& &\\
					\cline{2-3}
					&$\bbb^{5}$&$(\frac7{10}-\frac2f,\frac25,\frac35+\frac4f,\ddd,\frac{13}{10}-\ddd+\frac2f)$
					& &\\
					\cline{2-3}
					\cline{5-5}
					&$\aaa\bbb\ccc$&$(1-\frac4f,1-\frac4f,\frac8f,\ddd,1-\ddd+\frac4f)$
					& &\multirow{3}{*}{See Table \ref{Tab-10.2}}\\
					\cline{2-3}
					&$\bbb^3$&$(\frac56-\frac2f,\frac23,\frac13+\frac4f,\ddd,\frac76-\ddd+\frac2f)$
					& &\\
					\cline{1-3}
					\multirow{6}{*}{$\aaa\ddd\eee,\aaa\ccc^{2}$}&$\bbb^3$&$(\frac43-\frac8f,\frac23,\frac13+\frac4f,\ddd,\frac23-\ddd+\frac8f)$
					& &\\
					\cline{2-3}
					\cline{5-5}
					&$\aaa^2\bbb$&$(\frac45-\frac8{5f},\frac25+\frac{16}{5f},\frac35+\frac4{5f},\ddd,\frac65-\ddd+\frac8{5f})$
					& &\multirow{6}{*}{No degree $\ge4$ vertex $\ddd\eee\cdots$}\\
					\cline{2-3}
					&$\aaa\bbb^{3}$&$(\frac45-\frac{24}{5f},\frac25+\frac8{5f},\frac35+\frac{12}{5f},\ddd,\frac65-\ddd+\frac{24}{5f})$
					& &\\
					\cline{2-3}
					&$\bbb^{3}\ccc$&$(1-\frac{12}f,\frac12-\frac2f,\frac12+\frac6f,\ddd,1-\ddd+\frac{12}f)$
					& &\\
					\cline{2-3}
					&$\bbb^{4}$&$(1-\frac8f,\frac12,\frac12+\frac4f,\ddd,1-\ddd+\frac8f)$
					& &\\
					\cline{2-3}
					&$\bbb^{5}$&$(\frac45-\frac8f,\frac25,\frac35+\frac4f,\ddd,\frac65-\ddd+\frac8f)$& &\\
					\hline
			\end{tabular}}
		\end{table*}

		\begin{table}[htp]
			\centering
			\caption{AVC for Cases $\{\aaa\ddd\eee,\aaa^2\ccc, \aaa\bbb\ccc\}$, $\{\aaa\ddd\eee,\aaa^2\ccc, \bbb^3\}$ and $\{\aaa\ddd\eee,\aaa\ccc^2,\bbb^3\}$.}\label{Tab-10.2}
			\begin{minipage}{1\textwidth}
				\resizebox{\linewidth}{!}{\begin{tabular}{c|c||c|c||c|c}
					\hline
					\multicolumn{2}{c||}{AVC for $\{\aaa\ddd\eee,\aaa^2\ccc, \aaa\bbb\ccc\}$}&\multicolumn{2}{c||}{AVC for $\{\aaa\ddd\eee,\aaa^2\ccc,\bbb^3\}$}&\multicolumn{2}{c}{AVC for $\{\aaa\ddd\eee,\aaa\ccc^2,\bbb^3\}$}\\
					\hline
					$f$&vertex&$f$&vertex&$f$&vertex\\
					\hline
					\hline
					all&$\aaa\ddd\eee,\aaa^2\ccc,\aaa\bbb\ccc$&all&$\aaa\ddd\eee,\aaa^2\ccc,\bbb^3$&all&$\aaa\ddd\eee,\aaa\ccc^2,\bbb^3$\\
					$8s+4,s=2,3,\cdots$&$\aaa\ccc^{s+1},\bbb\ccc^{s+1},\ccc^{2s+1},\ccc^{s}\ddd\eee$&$60$&$\aaa\ccc^3,\ccc^2\ddd\eee,\ccc^5$&$24$&$\ccc^4,\ccc^2\ddd\eee,\ddd^2\eee^2$\\
					&&&&$36$&$\bbb\ccc^3,\bbb\ccc\ddd\eee$\\
					&&&&$60$&$\ccc^5,\ccc^3\ddd\eee,\ccc\ddd^2\eee^2$\\
					\hline
					\multicolumn{4}{c||}{By the balance lemma, $\ccc^k\ddd\eee$ must appear.}&\multicolumn{2}{c}{By the balance lemma, $\ccc^k\ddd\eee, \ccc^k\ddd^2\eee^2$ or $\bbb\ccc\ddd\eee$ must appear.}\\
					\multicolumn{4}{c||}{Its AAD gives $\bbb\thin\ddd$, $\bbb\thin\eee$, $\ddd\thin\eee$ or $\eee\thin\eee$, contradicting the AVC}&\multicolumn{2}{c}{Its AAD gives $\aaa\thin\aaa$, $\aaa\thin\bbb$,  $\bbb\thin\ccc$, $\bbb\thin\ddd$, $\bbb\thin\eee$ or $\ddd\thin\eee$, contradicting the AVC}\\
					\hline
				\end{tabular}}
			\end{minipage}
		\end{table}

	\end{proof}

	\begin{proposition}\label{ade,abc}
		There is no non-symmetric $a^4b$-tiling with general angles, such that $\aaa\ddd\eee$ is the unique $a^2b$-vertex type and $\aaa\bbb\ccc$ is a vertex.
	\end{proposition}

	\begin{proof}
		By Proposition \ref{ade,2bc}, \ref{ade,c3}, \ref{ade,a2c}, we may assume  that $\aaa^3$, $\aaa^2\ccc$ $(\aaa^2\bbb)$, $\aaa\ccc^2$ $(\aaa\bbb^2)$, $\bbb^2\ccc$ $(\bbb\ccc^2)$, $\ccc^3$ $(\bbb^3)$ are not vertices. Then $\aaa\ddd\eee$, $\aaa\bbb\ccc$ are the only degree $3$ vertices. By \eqref{vcountf} and \eqref{vcountv}, $\#\aaa\ge v_3=20+\sum_{k\ge 4}(3k-10)v_k>f=12+2\sum_{k\ge 4}(k-3)v_k$, a contradiction.
	\end{proof}

	\section{Some non-edge-to-edge quadrilateral tilings}\label{Conclusion}

    All pentagonal prototiles in Theorem \ref{theorem1} and \ref{theorem2} have fixed $\bbb$, $\ccc$ values, and may degenerate to quadrilaterals when $\aaa$, $\ddd$, or $\eee$  equals $\pi$ as shown in Figure \ref{degradation}.

	\begin{figure}[htp]
		\centering
		\begin{tikzpicture}[>=latex,scale=0.9]
			\begin{scope}[xshift=-3 cm,scale=0.5*0.75]
				\draw (-4,0)--(4,0)--(2.5,4)--(-5.5,4)--(-4,0);
				\draw[line width=1.5] (-4,0)--(4,0);
				\node at (-3.5,0.5){\small $\ddd$};\node at (3,0.5){\small $\eee$};
				\node at (-1.5,3.6){\small $\aaa$};\node at (-4.5,3.5){\small $\bbb$};\node at (2,3.5){\small $\ccc$};
				\node at (0.3,4.5){$a$};  \node at (-3.5,4.5){$a$}; \node at (-5.5,2.4){$a$}; \node at (3.8,2.4){$a$}; \node at (0,-0.5){$b$};
				\fill (-1.5,4) circle (0.1);
			\end{scope}

			\begin{scope}[xshift=0.5 cm,scale=0.4*0.75]
				\draw (0,0)--(8,0)--(4,5)--(-3.5,5)--(-4.5,0)--(0,0);
				\draw[line width=1.5] (0,0)--(8,0);
				\node at (0.3,0.58){\small $\ddd$};\node at (-3.8,0.65){\small $\bbb$};
				\node at (3.7,4.3){\small $\ccc$};\node at (-3,4.5){\small $\aaa$};\node at (6.3,0.65){\small $\eee$};
				\node at (3.5,-0.8){$b$};  \node at (-2,-0.8){$a$}; \node at (-5,3){$a$}; \node at (6.5,3){$a$}; \node at (0.5,5.6){$a$};
			\end{scope}

		    \begin{scope}[xshift=4.9 cm,scale=0.4*0.75]
		    	\draw (0,0)--(8,0)--(4,5)--(-3.5,5)--(-4.5,0)--(0,0);
		    	\draw[line width=1.5] (0,0)--(8,0);
		    	\node at (0.3,0.58){\small $\eee$};\node at (-3.8,0.65){\small $\ccc$};
		    	\node at (3.7,4.3){\small $\bbb$};\node at (-3,4.5){\small $\aaa$};\node at (6.3,0.65){\small $\ddd$};
		    	\node at (3.5,-0.8){$b$};  \node at (-2,-0.8){$a$}; \node at (-5,3){$a$}; \node at (6.5,3){$a$}; \node at (0.5,5.6){$a$};
		    \end{scope}
		\end{tikzpicture}
		\caption{Some quadrilaterals as degenerate pentagons.}
		\label{degradation}
	\end{figure}

\begin{figure}[htp]
		\centering
		\begin{overpic}[scale=0.125]{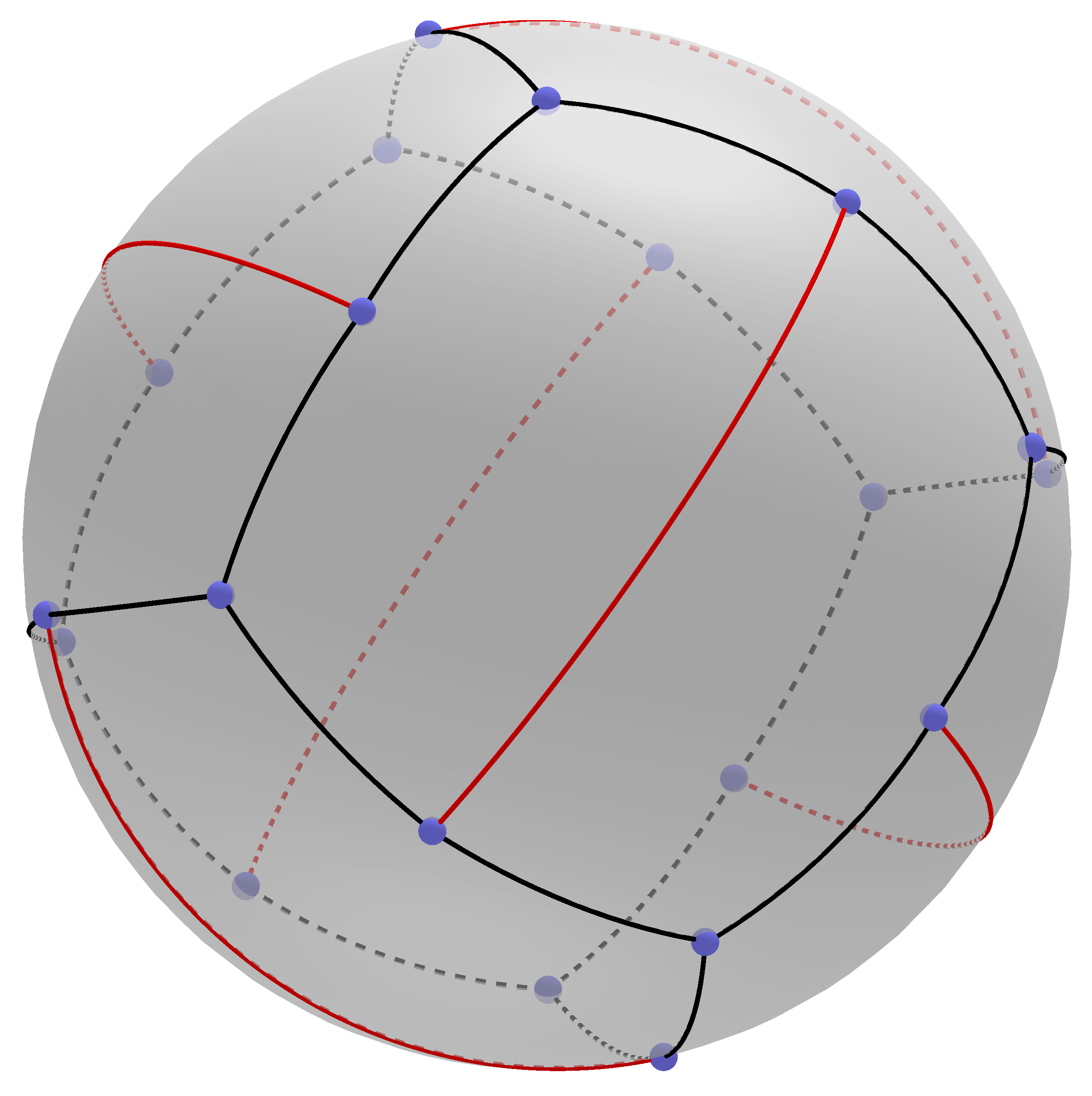}
			\put(85,5){\large $f=12$}
		\end{overpic} \hspace{30pt}
	    \begin{overpic}[scale=0.105]{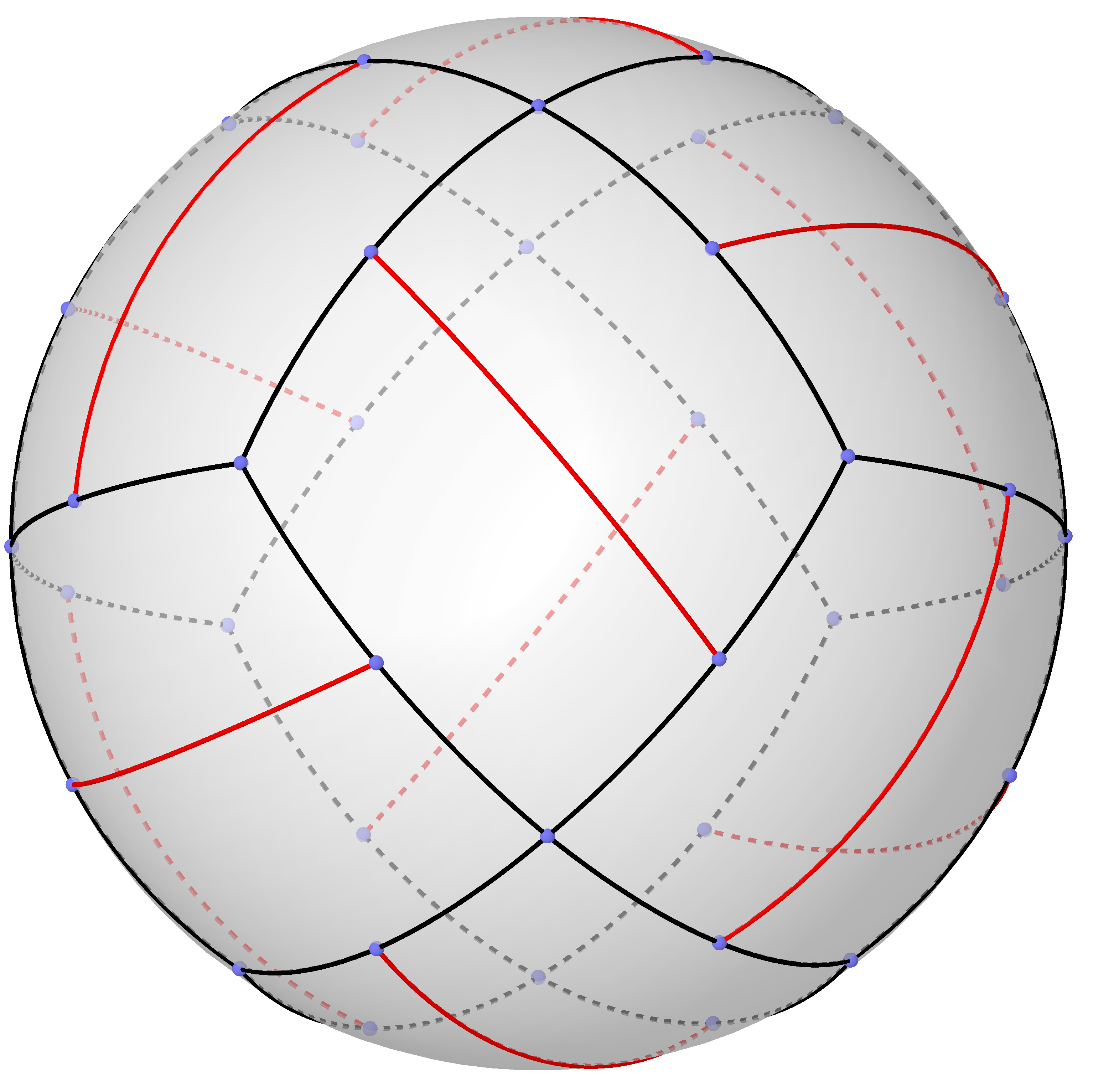}
	    	\put(85,5){\large $f=24$}
	    \end{overpic}  \hspace{30pt}
        \begin{overpic}[scale=0.097]{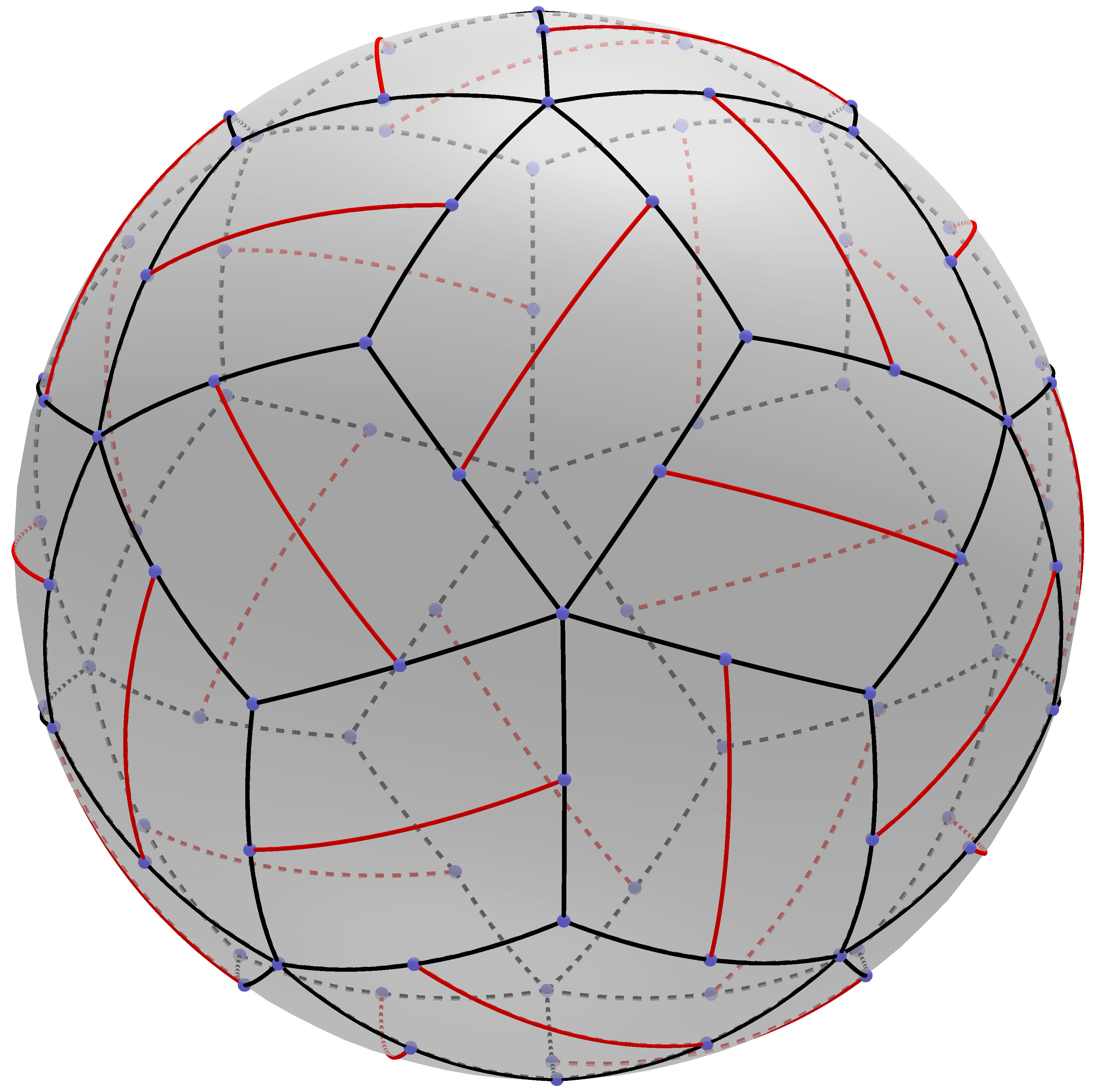}
        	\put(85,5){\large $f=60$}
        \end{overpic}

		\caption{The pentagonal subdivisions of the Platonic solids with $\aaa=\pi$.}
		\label{quadrilateral tilings1}
	\end{figure}

\begin{figure}[htp]
		\centering
		\begin{overpic}[scale=0.13]{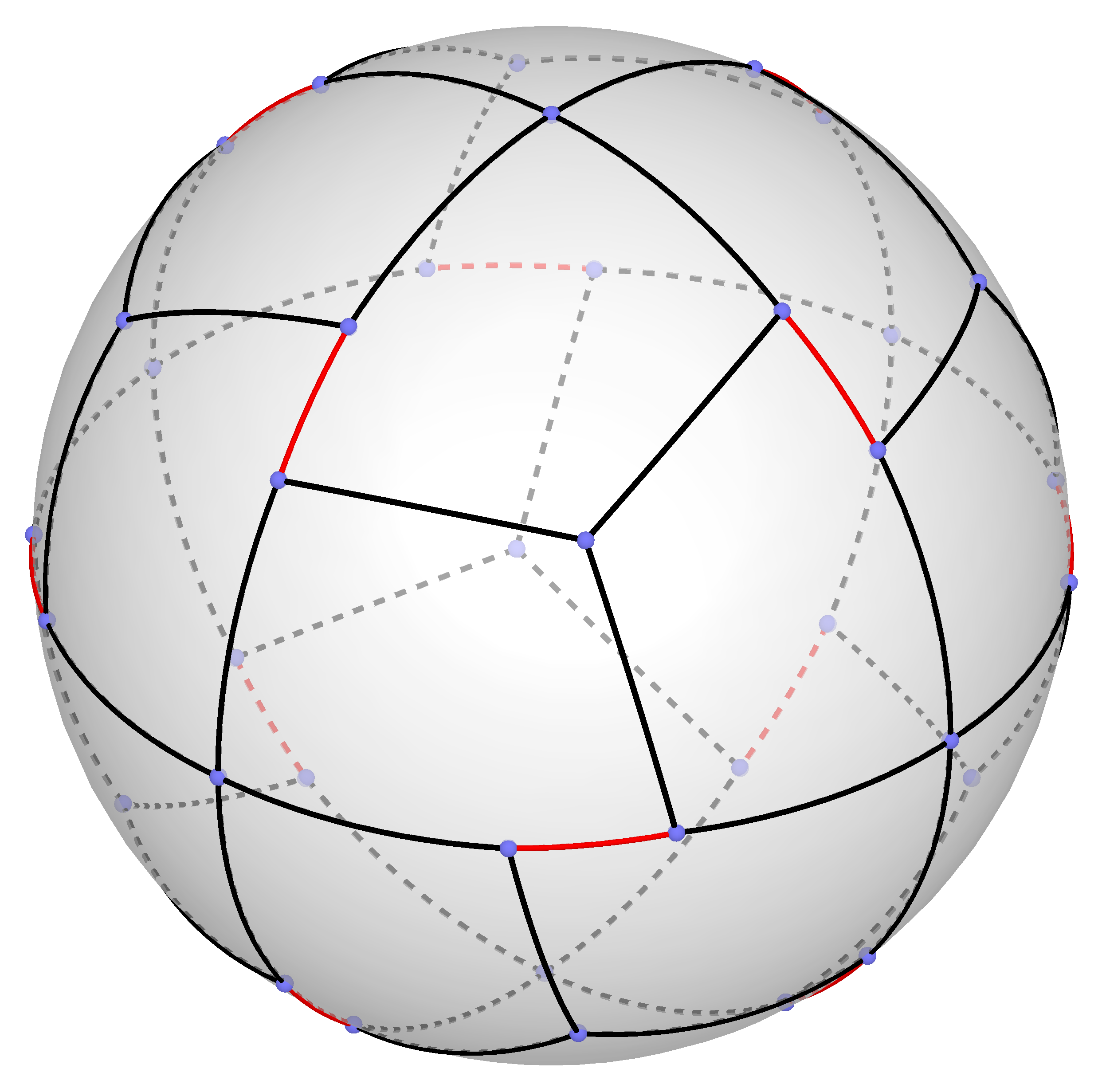}
			\put(85,5){\large $f=24$}
		\end{overpic} \hspace{60pt}
		\begin{overpic}[scale=0.098]{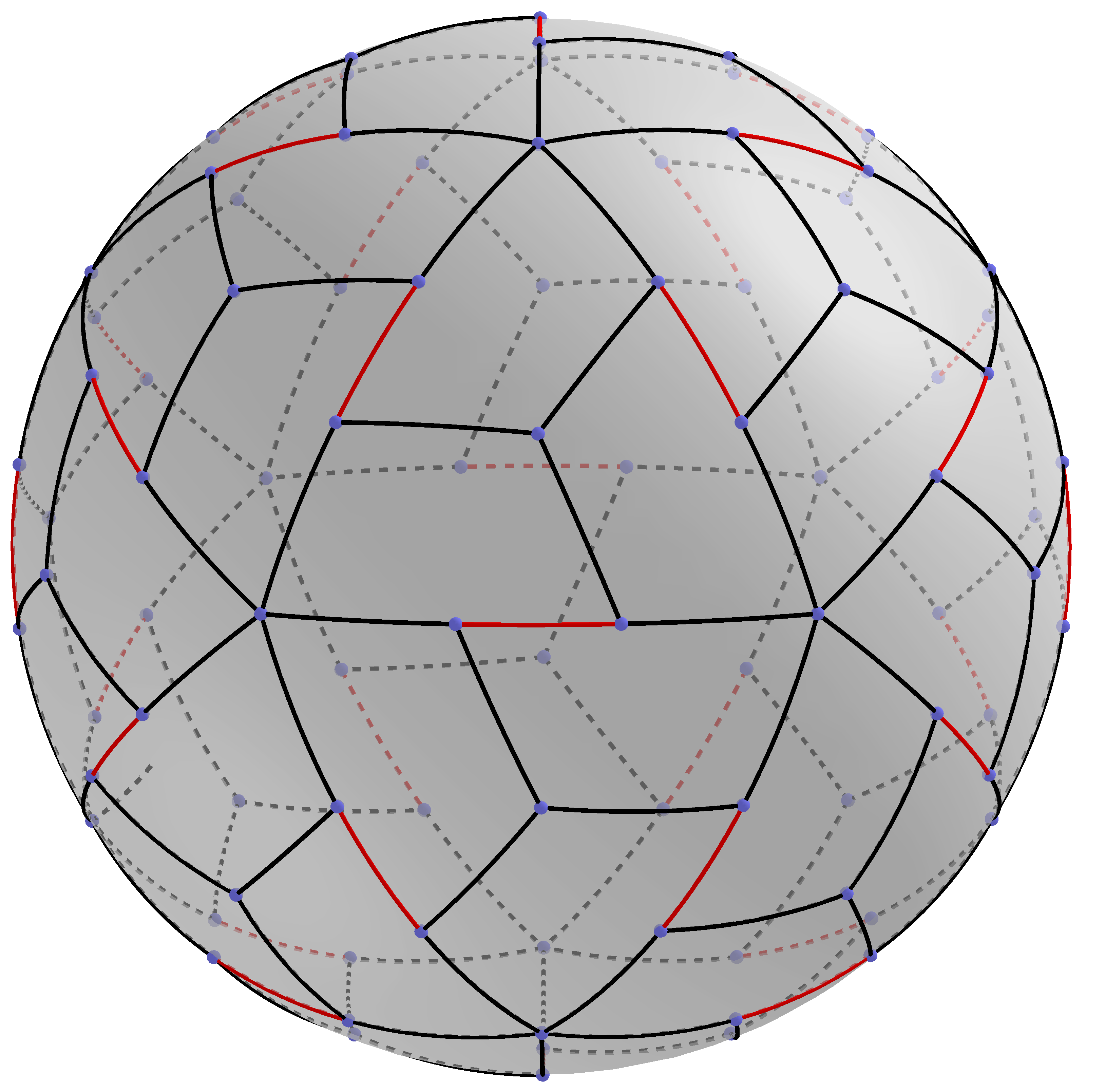}
			\put(85,5){\large $f=60$}
		\end{overpic}
		\caption{The pentagonal subdivisions of the Platonic solids with $\ddd=\pi$.}
		\label{quadrilateral tilings2}
	\end{figure}

\begin{figure}[htp]
		\centering
		\includegraphics[scale=0.0456]{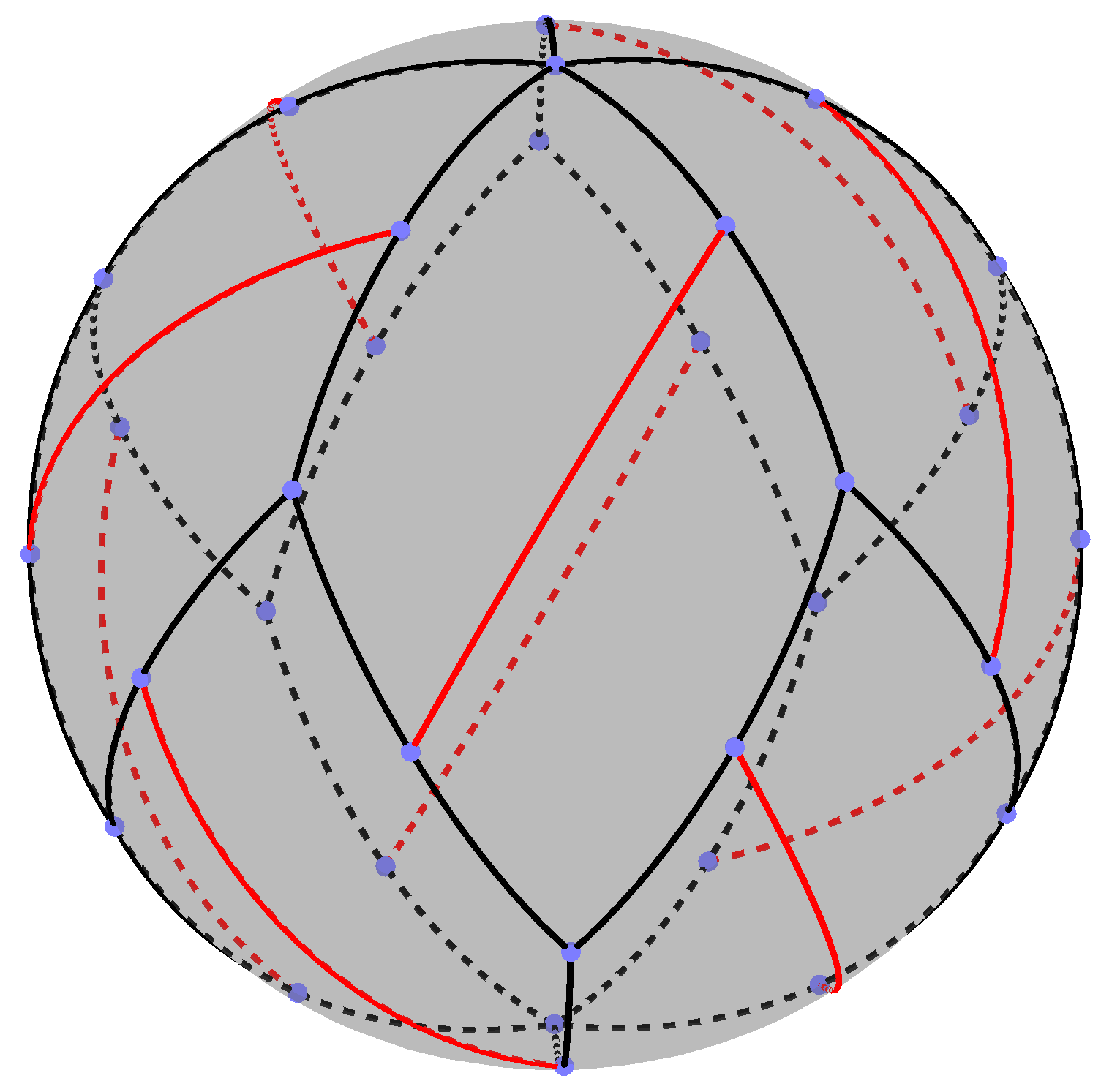}\hspace{10pt}
		\begin{overpic}[scale=0.0303]{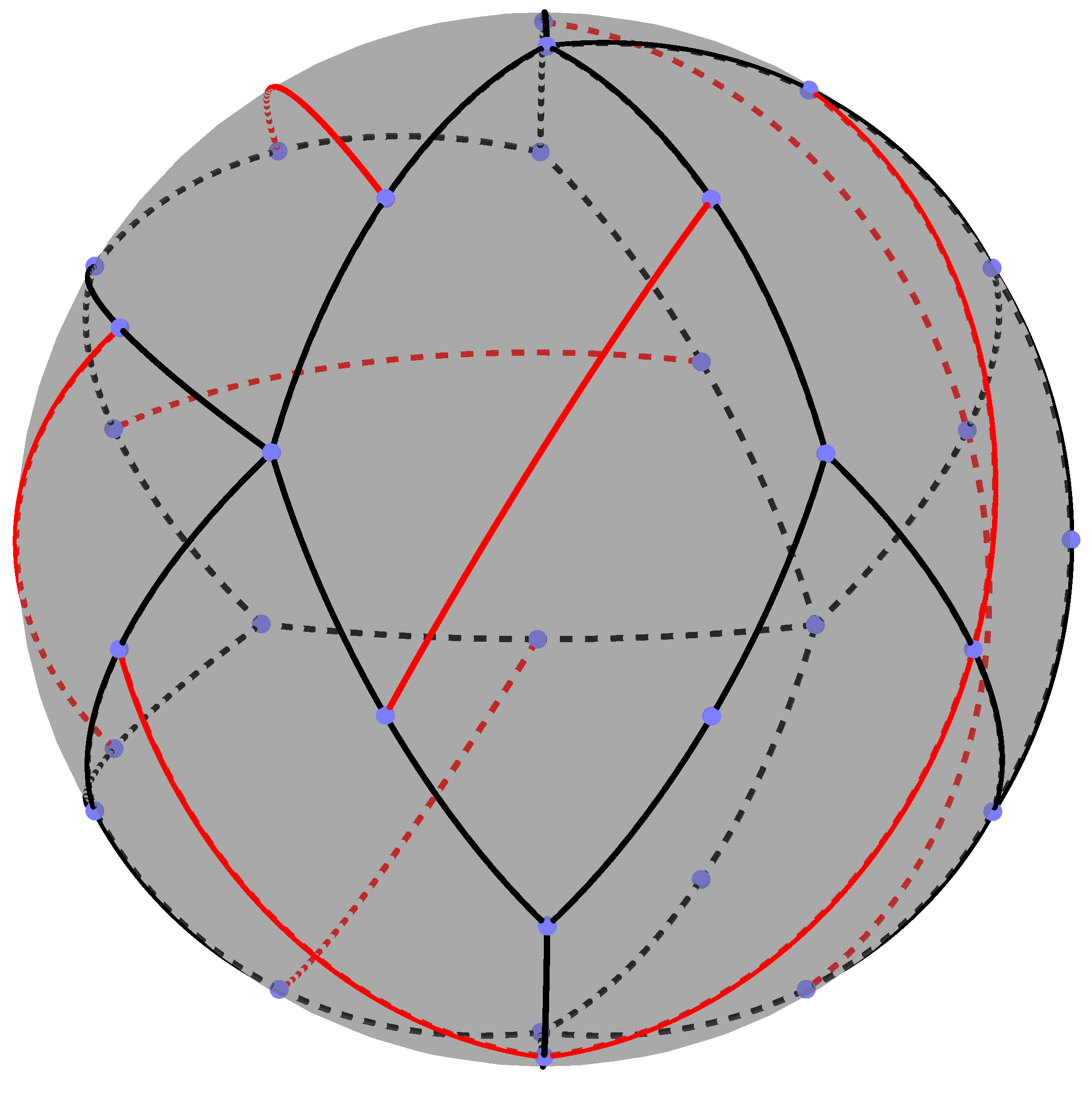}
			\put(-30,0){\large $\aaa=\pi$}
		\end{overpic} \hspace{10pt}
	    \includegraphics[scale=0.0408]{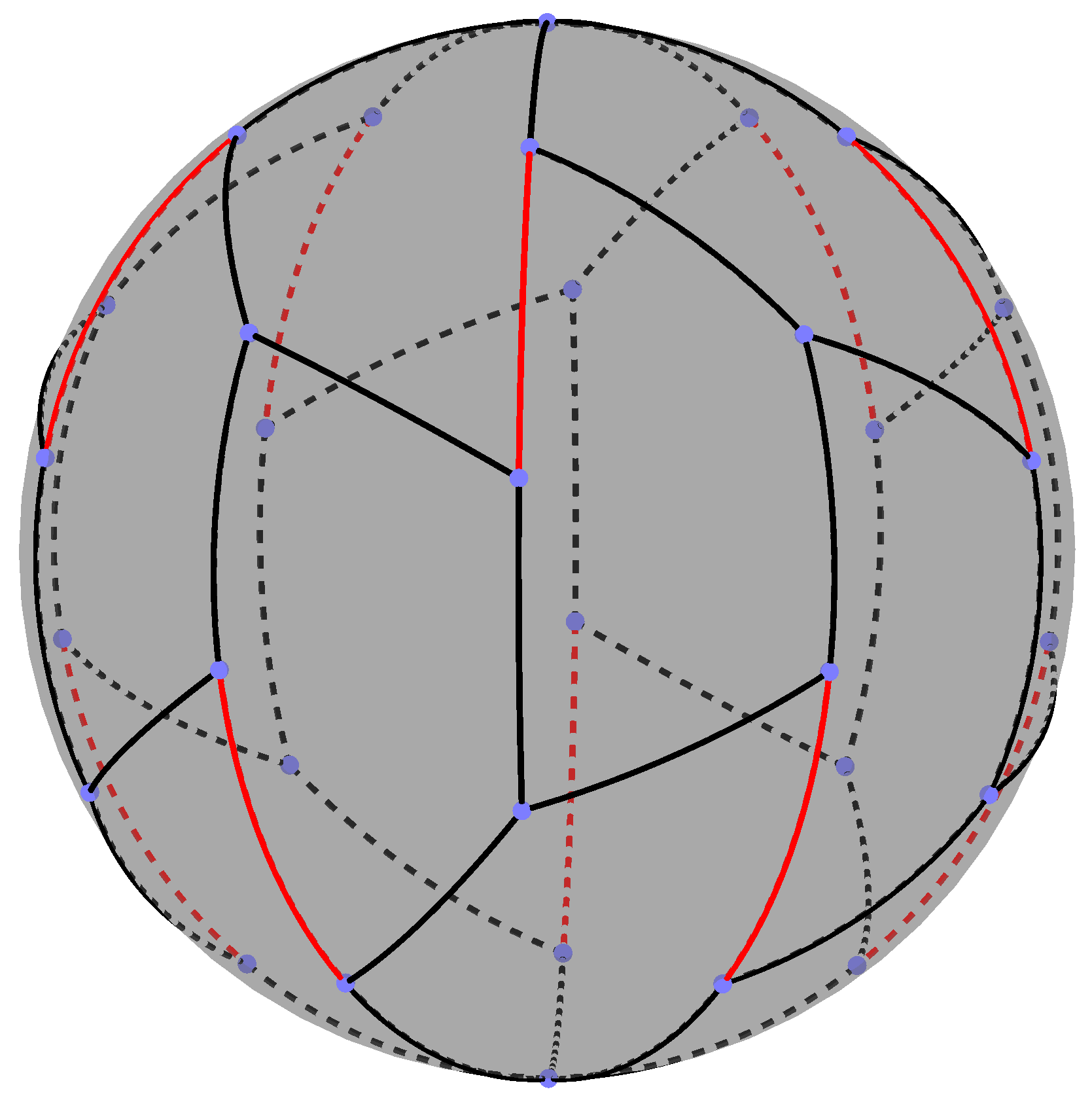}\hspace{10pt}
		\begin{overpic}[scale=0.0428]{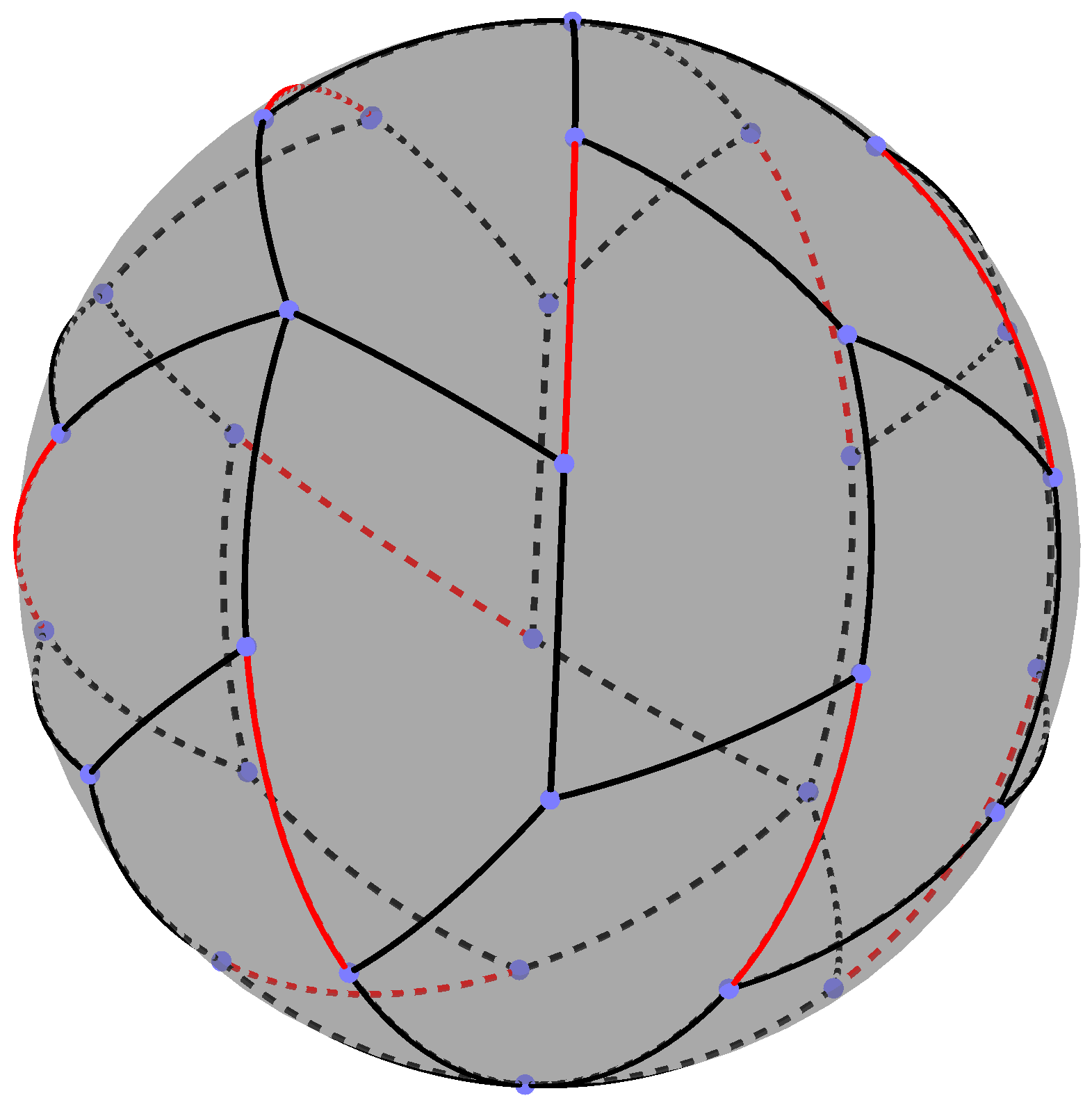}
			\put(-30,0){\large $\eee=\pi$}
		\end{overpic}

		\caption{The degenerate non-symmetric $3$-layer earth map tilings with $20$ tiles and their standard flips.}
		\label{quadrilateral tilings4}
	\end{figure}

	Then it induces some non-edge-to-edge quadrilateral tilings as follows:
	\begin{enumerate}
		\item Three pentagonal subdivisions of the Platonic solids with $12$, $24$ and $60$ tiles and $\aaa=\pi$ in Figure \ref{quadrilateral tilings1}, and two pentagonal subdivisions of the Platonic solids with $24$ and $60$ tiles and $\ddd=\pi$ in Figure \ref{quadrilateral tilings2};
		\item A sequence of unique pentagons with $\aaa$ or $\eee=\pi$ admitting non-symmetric $3$-layer earth map tilings $T(4m\,\aaa\ddd\eee,2m\,\bbb^2\ccc,2\ccc^m)$ with $4m$ tiles for any $m\ge4$ (see the first and third pictures of Figure \ref{quadrilateral tilings4}), among which each odd $m=2k+1$ case always admits a standard flip modification: $T((8k+4)\aaa\ddd\eee,4k\,\bbb^2\ccc,4\bbb\ccc^{k+1})$ (see the second and fourth pictures of Figure \ref{quadrilateral tilings4}).
	\end{enumerate}

\end{document}